%% file: main.tex
\documentclass[10pt]{amsart}

\pdfoutput=1
\usepackage[top=1.4in,bottom=1.275in,left=1.2in,right=1.2in]{geometry}
\usepackage{amsmath}
\usepackage{amssymb}
\usepackage{amsthm}
\usepackage{amscd}
\usepackage{amsaddr}
\usepackage[only,llbracket,rrbracket]{stmaryrd}
\usepackage{calligra}
\usepackage{booktabs}
\usepackage{array}
\newcolumntype{?}{!{\vrule width 1pt}}
\usepackage{lmodern}
\usepackage[toc,page]{appendix}

\usepackage{bm}
\usepackage{times}
\usepackage{microtype}
\usepackage[font=small,figurewithin=section]{caption}
\usepackage[labelformat=simple]{subcaption}

\usepackage[hypertexnames=false]{hyperref}
\usepackage[noabbrev,nameinlink,capitalize]{cleveref}

\usepackage{tikz}
\usetikzlibrary{matrix,arrows,decorations.pathmorphing}
\usepackage{enumitem}
\usepackage{algorithm}
\usepackage[noend]{algpseudocode}

\algnewcommand\algorithmicto{\textbf{to}}
\algnewcommand\RETURN{\State \textbf{return} }
\numberwithin{equation}{section}
\makeatletter
\g@addto@macro\bfseries{\boldmath}
\makeatother
\hypersetup{colorlinks=true,citecolor=black,filecolor=black,linkcolor=black,urlcolor=black}
\setlist{nolistsep}

\SetSymbolFont{stmry}{bold}{U}{stmry}{m}{n}

\allowdisplaybreaks

\input{preamble}

\title[Socratis Petrides and Leszek Demkowicz]{An adaptive multigrid solver for DPG methods with applications in linear acoustics and electromagnetics}
\author[Socratis Petrides and Leszek Demkowicz]{Socratis Petrides$^1$ and Leszek Demkowicz$^2$}
\address{$^1$Center for Applied Scientific Computing, Lawrence Livermore National Laboratory}
\address{$^2$Oden Institute for Computational Engineering and Sciences, The University of Texas at Austin}
\email{petrides1@llnl.gov}
\email{leszek@oden.utexas.edu}

\begin{document}

	\input{abstract}

	\maketitle

	\input{intro}

	\input{uwacoustics}

	\input{analysis}

	\input{mg}

	\input{results}

	\input{conclusion}

	\section*{Acknowledgments}

\thanks{
This work has been supported with grants by the Air Force Office of Scientific Research (AFOSR grant FA9550-12-1-0484 and FA9550-17-1-0090).
The authors would also like to express their gratitude to Dr. Jay Gopalakrishnan for his suggestions and his insightful comments throughout this work. 

This work was performed under the auspices of the U.S. Department of Energy by Lawrence Livermore National Laboratory under Contract DE-AC52-07NA27344, LLNL-JRNL-815524.
This document was prepared as an account of work sponsored by an agency of the United States government. 
Neither the United States government nor Lawrence Livermore National Security, LLC, nor any of their employees makes any warranty, expressed or implied, or assumes any legal liability or responsibility for the accuracy, completeness, or usefulness of any information, apparatus, product, or process disclosed, or represents that its use
would not infringe privately owned rights. Reference herein to any specific commercial product, process, or service
by trade name, trademark, manufacturer, or otherwise does not necessarily constitute or imply its endorsement,
recommendation, or favoring by the United States government or Lawrence Livermore National Security, LLC.
The views and opinions of authors expressed herein do not necessarily state or reflect those of the United States
government or Lawrence Livermore National Security, LLC, and shall not be used for advertising or product
endorsement purposes.
}
	\appendix
	\input{interpolant}
	\bibliographystyle{siam}
	\bibliography{shortref,ref}
\end{document}

%% file: preamble.tex
\newcommand{\half}{\frac{1}{2}}
\let\div\relax 
\DeclareMathOperator{\div}{div}
\DeclareMathOperator{\curl}{curl}
\DeclareMathOperator*{\argmin}{arg\,min}
\newcommand{\bsL}{\boldsymbol{L}}
\newcommand{\Aop}{\mathrm{A}} 
\newcommand{\Bop}{\mathrm{B}} 
\newcommand{\Rop}{\mathrm{R}} 
\newcommand{\Top}{\mathrm{T}} 
\newcommand{\sfA}{\mathsf{A}}
\newcommand{\sfB}{\mathsf{B}}
\newcommand{\sfP}{\mathsf{P}}
\newcommand{\sfG}{\mathsf{G}}
\newcommand{\sfQ}{\mathsf{Q}}
\newcommand{\sfu}{\mathsf{u}}
\newcommand{\sfv}{\mathsf{v}}
\newcommand{\sfr}{\mathsf{r}}
\newcommand{\sfz}{\mathsf{z}}
\newcommand{\sfellt}[2]{\raisebox{0.0\height}{$#1{}\text{\textsf{\fontfamily{PTSans-TLF}\selectfont l}}$}}
\newcommand{\sfell}{\mathpalette\sfellt\relax}
\newcommand{\C}{\mathbb{C}} 

\DeclareMathOperator*{\supp}{supp}

\newcommand{\mthu}{\mathrm{u}} 
\newcommand{\mthv}{\mathrm{v}} 
\newtheoremstyle{boldremark}
    {\dimexpr\topsep/2\relax} 
    {\dimexpr\topsep/2\relax} 
    {}          
    {}          
    {\bfseries} 
    {.}         
    {.5em}      
    {}          
\theoremstyle{plain}
\newtheorem*{definition*}{Definition}

\theoremstyle{plain}
\newtheorem{theorem}{Theorem}
\newtheorem{lemma}{Lemma}

\theoremstyle{boldremark}


%% file: abstract.tex

\begin{abstract}
We propose an adaptive multigrid preconditioning technology for solving linear systems arising from Discontinuous Petrov-Galerkin (DPG) discretizations. Unlike standard multigrid techniques, this preconditioner involves only trace spaces defined on the mesh skeleton, and it is suitable for adaptive hp-meshes. The key point of the construction is the integration of the iterative solver with a fully automatic and reliable mesh refinement process provided by the DPG technology. The efficacy of the solution technique is showcased with numerous examples of linear acoustics and electromagnetic simulations, including simulations in the high-frequency regime, problems which otherwise would be intractable. Finally, we analyze the one-level preconditioner (smoother) for uniform meshes and we demonstrate that theoretical estimates of the condition number of the preconditioned linear system can be derived based on well established theory for self--adjoint positive definite operators. 
\end{abstract}

%% file: intro.tex

\section{Introduction} \label{sec:sec1}
\subsection{Background}

Numerical simulation of wave propagation phenomena has been a very challenging problem in the computational science community for many years. The main difficulty is the following. Most commonly used numerical methods for PDEs, such as the Finite Element or the Finite Difference methods lose their stability, unless the wave frequency is ``resolved''. That is, the underlying discretization mesh has to be ``fine enough'' in order to satisfy the \emph{Nyquist criterion} and overcome the so called \emph{pollution effect}. The Nyquist criterion requires a minimum number of finite elements or discretization points per wavelength for the method to remain stable. Moreover, as the frequency increases, even if the Nyquist criterion is met, the solution gets ``polluted'', i.e., becomes out of phase or very diffusive. 
Pollution can be minimized by using high order discretizations, but it can not be eliminated. To the best of our knowledge, to this day there is no \emph{pollution-free} method for simulations in two or three space dimensions \cite{Babuska1}. 

Fine meshes and high order discretizations result in very large linear systems and therefore extremely hard to solve. State of the art direct solvers based on LU decompositions for sparse matrices fail to scale efficiently and they usually run out of memory. While for two dimensional problems multi-frontal direct solvers for sparse matrices can be optimal, i.e., for a linear system of size $N$ the complexity and memory estimates  are $\mathcal{O}(N^{3/2})$ and $\mathcal{O}(N log N)$ respectively, in three space dimensions this is not the case. For a 3D problem on a uniform mesh, the computational cost is estimated to be $\mathcal{O}(N^2)$ and the required memory is $\mathcal{O}(N^{4/3})$ \cite{mumps1}.

An alternative approach is to employ iterative solution techniques. This kind of solvers do not need to perform any kind of factorization. On the contrary, they are based on recursive matrix--vector multiplications, operations that have both complexity and memory requirements of $\mathcal{O}(N)$. Usually, the optimal strategy is to construct a ``good'' approximation of the inverse of the system matrix and use it as a preconditioner to a Krylov subspace method such as the Generalized Minimum Residual (GMRES) or the Conjugate Gradient (CG) method. Unfortunately, constructing an optimal preconditioner for time--harmonic wave simulations is still a notoriously hard problem. The time--harmonic wave operator (acoustics or Maxwell) is highly indefinite and cutting edge preconditioners for elliptic problems (such as multilevel, multigrid, DDM) fail to work. A common strategy for preconditioning indefinite operators involve the construction of a preconditioner for GMRES which shifts the spectrum to the positive half plane. The resulting solver is not always reliable, especially when the grid is very coarse. In addition, multigrid techniques lose their efficiency, since in order to converge, they usually require expensive coarse grid solves because the coarse mesh has to be fine enough to maintain discrete stability.

Several recent approaches show more promising results \cite{gander2019class}. In particular well known attempts include but 
not limited to, shifted Laplacian techniques \cite{ErnstGander,LahayeD,Gander2}, domain decomposition methods \cite{KimSeungil,JaySchwarz,DBLP,Stolk}, multigrid methods \cite{Stolk14,Jaymultigrid,XavierV2} and stabilized FEM techniques based on artificial absorption \cite{absorption_2,absorption_1}. A lot of attention has been paid to the sweeping preconditioner \cite{EngquistY1,VionA} and some of its variants such us the polarized traces \cite{Zepeda}, the L--sweeps method \cite{taus2020sweeps} and the Source transfer method \cite{chen2013source,leng2020diagonal}. These are domain decomposition methods which introduce absorbing boundary conditions in the interior of the domain (often a PML) to transfer the information among the subdomains. Although it is hard to exploit parallelization, they are provably very optimal for problems in unbounded domains. However, their efficiency can significantly deteriorate for quasi--resonant problems.

\subsection{Preconditioning the DPG system}

There are several works on preconditioning DPG systems for elliptic problems. To the best of our knowledge the first one chronologically, was done by Barker et al. in \cite{Barker2014}. In this work the authors analyzed the one level additive Schwarz preconditioner for the Poisson problem and showcased, both theoretically and numerically, convergence of the preconditioned CG solver, independent of the mesh size $h$. A few years later, Barker et al. in \cite{barker} constructed a scalable preconditioner for the primal DPG method again for the Poisson problem. The key point of the implementation was the norm equivalence of the DPG bilinear form with the standard $H^1$ and $H(\div)$ norms. Naturally, existing $H^1$ and $H(\div)$ algebraic multigrid (AMG) technologies could be utilized to precondition the DPG system. Following a similar approach in \cite{LiXu}, the authors extended theoretical analysis for the Poisson problem to the two level setting. Their work was based on well known results on preconditioning the $H^1$, $H(\div)$ and $H(\curl)$ spaces \cite{arnold1997preconditioning,dryja1994domain,hiptmair1999multilevel}. Finally, a geometric multigrid preconditioner for the Poisson and Stokes problem is presented in \cite{ROBERTS20172018}.

Designing a preconditioner for wave problems is much more challenging. While, it is relatively straight forward to construct and analyze preconditioners for elliptic problems using the idea of norm equivalence, for wave problems the equivalence constants are frequency dependent. As far as we know there are two works so far attempting to construct robust and efficient preconditioners for wave problems but none of them provide any theoretical results. In \cite{LiXu} the authors present a numerical study on a one level additive Schwarz preconditioner when applied to the Helmholtz problem. Their results indicate uniform convergence of the solver with respect to the frequency $\omega$ and the mesh size $h$. The convergence is however sensitive to the overlap and the size of the Schwarz patches. The second attempt on preconditioning the Helmholtz problem is described in \cite{Gopalakrishnan2015}. This work introduces a multiplicative Schwarz (overlapping block Gauss--Seidel) preconditioner for the primal DPG formulation. The preconditioner is shown to converge at a rate independent of the polynomial order $p$ and the frequency $\omega$, but no dependence study on the mesh size $h$ is provided.

For our construction we take a slightly different path than most of the works described above. The preconditioner is constructed directly from the underlying bilinear form of the wave problem. Using the properties of the ultraweak DPG method we invoke a numerical experiment to examine the dependence of the condition number of the preconditioned system with respect to the polynomial order of approximation, the frequency, the Schwarz patch size and the mesh size. Even though the analysis is done only for a one level additive smoother, it gives us useful insights for the construction of a multilevel preconditioner. We note that our results are consistent with the results of related work described above.

\subsection{Outline} 
The paper is organized as follows. First, we give a brief summary of the DPG methods, using as a model problem the linear acoustics equations. Moving on, in \cref{sec:sec3}, we provide some theoretical results on the one-level preconditioner (smoother) for these equations. The analysis is complemented with a numerical experiment that is described in detail in \cref{appendix:interpolant}. In \cref{sec:sec4}, we present the multigrid algorithm and in \cref{sec:sec5} we showcase with numerous numerical results the efficacy of this method for both acoustics and electromagnetic simulations. Lastly, in \cref{sec:sec6}, we summarize our conclusions and outline future work.

%% file: uwacoustics.tex

\section{The DPG method and linear acoustics} \label{sec:sec2}
The DPG method for the time-harmonic linear acoustics equations was analyzed in dept in \cite{Zitelli,demk2}. In this section we briefly mention the DPG setup and provide some background for the analysis of the one-level additive Schwarz smoother for preconditioning the ultraweak DPG formation of the time-harmonic acoustics equations.

\subsection{The ideal Petrov--Galerkin method: a brief outline}

We consider the following abstract setting. Given a continuous bilinear (sesquilinear) form $b(\cdot,\cdot)$ defined on the product $U\times V$ of Hilbert spaces $U$ (trial space) and $V$ (test space), and a continuous linear (anti--linear) form $l(\cdot)$ defined on $V$, we want to find the solution to the problem: 
\begin{equation*}
   \left\{
      \begin{aligned}
         &u \in U \\
         &b(u,v) = l(v), \quad v \in V. 
      \end{aligned}
   \right.
\end{equation*}
We assume that the above continuous problem is well posed, i.e, the continuous inf--sup condition \cite{babuska} for the bilinear form $b(\cdot,\cdot)$ holds. 
Given now a finite--dimensional trial space $U_h \subset U$, we can recast the above linear problem as a minimization problem in the dual norm  $\|\cdot\|_V^{'}$, i.e., we now seek the solution to 
\begin{equation*}
u_h = \argmin_{w_h \in U_h} \frac{1}{2}\|l-\Bop w_h\|^2_{V^{'}},
\end{equation*}
where $\Bop : U \rightarrow V^{'}$ is defined as 
\begin{equation*}
\langle \Bop u,v\rangle_{V^{'}\times V} = b(u,v) 
\end{equation*} 
and $\langle \cdot , \cdot \rangle_{V^{'}\times V}$ denotes the duality pairing on ${V^{'}\times V}$. By introducing the \emph{Riesz operator} for the test space, the dual norm $|\cdot\|_{V^{'}}$ can be replaced with $\|\Rop_V^{-1}(\cdot)\|_{V}$ and arrive at a simpler (to compute with) minimization problem: 
\begin{equation}
\label{eq:minres2}
u_h = \argmin_{w_h \in U_h} \frac{1}{2}\|\Rop_V^{-1}(l-\Bop w_h)\|^2_{V}.
\end{equation} 
The Riesz operator $\Rop_V:V \rightarrow V^{'}$ is an isometric isomorphism defined by
\begin{equation*}
(\Rop_Vy)(v) = \langle \Rop_v y , v \rangle_{V^{'}\times V} = (y,v)_V, \quad y,v \in V.
\end{equation*}
\cref{eq:minres2} leads to two additional characterizations for the DPG method. Taking the G$\hat{\text{a}}$teaux derivative, leads to the following linear problem:
\begin{equation}\label{eq:eq5}
   \left\{
      \begin{aligned}
         &u_h \in U_h\subset U\\
         &(\Rop_V^{-1}(l-\Bop u_h),\Rop_V^{-1}\Bop\delta u_h) = 0, \quad \delta u_h \in U_h
      \end{aligned}
   \right.
\end{equation}
We introduce the \emph{trial--to--test operator} $\Top:U_h \rightarrow V$, defined by 
\begin{equation}\label{eq:eq6}
(\Top \delta u_h,\delta v)_V = b(\delta u_h,\delta v) \,\,\,\, \delta u_h \in U_h, \, \delta v \in V, \quad \text{ i.e., }  \Top = \Rop_V^{-1} \Bop 
\end{equation}
and we define the \emph{optimal test space} as $V_h^{\text{opt}}: = \Top(U_h)$, i.e, the image of the trial to test operator T acting on the trial space $U_h$. Therefore, problem \eqref{eq:eq5} leads to the following Petrov--Galerkin scheme with \emph{optimal test functions}.
\begin{equation*}
   \left\{
      \begin{aligned}
         &u_h \in U_h\subset U\\
         &b(u_h,v_h) = l(v_h),  \quad v_h \in V_h^{\text{opt}} \\
      \end{aligned}
   \right.
\end{equation*}
Clearly, the above construction of the test space always delivers a symmetric positive definite system and it is \emph{optimal} in then sense that it guarantees uniform and unconditional discrete stability \cite{dpg_opt}. 
Indeed, for $v_h \in V_h^{\text{opt}}$ we have:
\begin{equation*}
\frac{b(u_h,v_h)}{\|v_h\|_V} = \frac{b(u_h,\Top u_h)}{\|\Top u_h\|_V} = \|\Top u_h\|_V = \sup_{v\in V} \frac{|(\Top u_h,v)_V|}{\|v\|_V} = \sup_{v\in V} \frac{|b(u_h,v)|}{\|v\|_V}. 
\end{equation*}
Therefore the \emph{discrete inf--sup} stability condition \cite{babuska}, is trivially satisfied. 

\noindent
Additionally, if we define $\psi := \Rop_V^{-1}(l-\Bop u_h)$ to be the \emph{error representation function}, \eqref{eq:eq5} leads to the following mixed problem:
\begin{equation*}
\left\{
\begin{alignedat}{3}
    &u_h \in U_h, \, \psi \in V,\\
    &(\psi,v)_V+b(u_h,v)&&=l(v), &&\quad v\in V,\\
    &b(w_h,\psi)&&=0, &&\quad w_h\in U_h
\end{alignedat}
\right.
\end{equation*}
where now we solve simultaneously for the original unknown $u_h$ and the error representation function $\psi$. Note that the norm of $\psi$ offers a built--in a--posteriori error indicator. Indeed, if we define the \emph{energy} norm as  
\begin{equation*}
\|u\|_E := \|\Bop u\|_{V^{'}} = \sup_{v\in V} \frac{|b(u,v)|}{\|v\|_V}
\end{equation*}
then we can compute the Galerkin error measured in the \emph{energy} norm, by computing the norm of the error representation function (residual). 
\begin{equation*}
\|u_h - u\|_E = \|\Bop(u_h - u)\|_{V^{'}} = \|\Bop u_h - l\|_{V^{'}} = \|\Rop_V^{-1}(\Bop u_h - l)\|_V = \|\psi\|_V. 
\end{equation*}

\subsection{The ultraweak formulation -- linear acoustics}
Consider the \emph{operator form} of the time--harmonic linear acoustics equations with impedance boundary conditions in a simply connected domain $\Omega$:
\begin{equation*}
\left\{
 \begin{alignedat}{3}
  &i \omega p + \div u    &&= f_1,      && \quad \text{ in }  \Omega\\
  &i \omega u + \nabla p  &&= f_2,      && \quad \text{ in } \Omega\\
  &p - u \cdot n          &&= 0,        && \quad \text{ on } \partial \Omega
   \end{alignedat}
\right .   
\end{equation*}
Here, $p$ and $u$ denote the pressure and the velocity of the wave respectively. Define the group variable $\mathfrak{u} := (p,u)$, the group load $\mathfrak{f} := (f_1,f_2) \in L^2(\Omega) \times (L^2(\Omega))^d =: \boldsymbol{L^2}(\Omega)$ and the acoustics operator 
$\Aop((p,u)) = (i \omega p + \div u, i \omega u + \nabla p)$. Then the problem can be compactly written as:
\begin{equation}
\label{eq:first_order_system}
\left\{
\begin{aligned}
&\mathfrak{u} \in D(\Aop) \\
&\Aop \mathfrak{u} = \mathfrak{f}
\end{aligned}
\right.
\end{equation}
where $D(\Aop)$ denotes the the domain of $\Aop$:
\begin{equation*}
D(\Aop) := \{ \mathfrak{u} \in \boldsymbol{L^2}(\Omega)\, : \Aop\mathfrak{u} \in \boldsymbol{L^2}(\Omega), \,  p - u \cdot n \,\,= 0  \text{ on } \partial \Omega \}.
\end{equation*}
By testing (\ref{eq:first_order_system}) with the group variable $\mathfrak{v}:=(q,v)$, and using integration by parts, the regularity is passed from the trial to the test space. The resulting variational formulation is called the \emph{ultraweak formulation}:
\begin{equation*}
\left\{
\begin{aligned}
&\mathfrak{u} \in \boldsymbol{L^2}(\Omega) \\
&(\mathfrak{u}, \Aop^\ast \mathfrak{v}) = (\mathfrak{f},\mathfrak{v}), \quad \mathfrak{v} \in D(\Aop^\ast),
\end{aligned}
\right.
\end{equation*}
where $D(\Aop^\ast)$ is equipped with the adjoint graph norm,
\begin{equation*}
\Vert \mathfrak{v} \Vert_V^2 := \Vert \Aop^\ast \mathfrak{v} \Vert^2 + \alpha^2 \Vert \mathfrak{v} \Vert^2 \, ,
\end{equation*}
with a scaling parameter $\alpha = \mathcal{O}(1)$. 
Observe that for this specific problem, the operator is (formally) {\em skew--adjoint}, $\Aop^\ast = - \Aop$, and 
\begin{equation*}
 D(\Aop^{*})= \{\mathfrak{v} \in \bsL^2(\Omega)\, : \Aop^\ast\mathfrak{v} \in \bsL^2(\Omega) :  
q + v \cdot n \,\,= 0  \text{ on } \partial \Omega \}.
\end{equation*}
Note that with the above choice of norm on the test space and $\alpha = 0$ the ultraweak formulation delivers $L^2$--projection. In practice, when the space is broken, $\alpha$ has to be non--zero in order for the norm to be \emph{localizable}. A common choice is $\alpha = 1$. It easy to see that the two norms are robustly equivalent. This in turn gives a pollution free method for one--space dimension problems (see \cite{sp,demk2,Zitelli,dpg_e} for additional details). Additionally, we emphasize that other equivalent (in terms of stability) formulations can be explored \cite{sp2,demk_varform,fuentes2017coupled,MR3543022}, however, especially for wave operators, the ultraweak formulation is provably superior because of its approximability properties \cite{Zitelli,demk2,petrides2019adaptive,nagaraj20193d,henneking2019dpg,henneking2020numerical}.

The boundary conditions on the test functions (built into the definition of $D(\Aop^\ast)$), can be eliminated at the expense of introducing extra unknowns $\hat{\mathfrak{u}}:=(\hat{p}, \hat{u}_n)$. We arrive then at the linear problem
\begin{equation*}
\left\{
\begin{aligned}
&\mathfrak{u} \in \bsL^2(\Omega),\, \hat{\mathfrak{u}} \in \hat{U} \\
&(\mathfrak{u}, \Aop^\ast \mathfrak{v}) + \langle \hat{\mathfrak{u}}, \mathfrak{v} \rangle
= (\mathfrak{f},\mathfrak{v}), \quad \mathfrak{v} \in H_{\Aop^\ast}(\Omega) \, ,
\end{aligned}
\right.
\end{equation*}
where
\begin{equation*}
H_{\Aop^\ast}(\Omega) := \{ \mathfrak{v} \in \bsL^2(\Omega) \, : \, \Aop^\ast \mathfrak{v} \in \bsL^2(\Omega) \}
\end{equation*}
is equipped with the same graph norm as above. The flux $\hat{\mathfrak{u}}$ is defined as the trace of a function that lives in the original
energy graph space,
\begin{equation*}
\hat{U} := \mathrm{tr}\,  D(\Aop) \,. 
\end{equation*}
More precisely, the trace space is defined by:
\begin{equation}\label{eq:trace_space}
\begin{aligned}
\hat{U} &= \{ (\hat{p},\hat{u}_n) \,  : \, \text{exists } (p,u) \in D(\Aop) \text{ such that } (\hat{p},\hat{u}_n) = \mathrm{tr} (p,u) \} \\
& =  \{ (\hat{p},\hat{u}_n) \in H^{\half}(\partial \Omega) \times H^{-\half}(\partial \Omega) \, : \, \hat{p} - \hat{u}_n = 0 \text{ on } \partial \Omega \} \, .
\end{aligned}
\end{equation}
and it is equipped with the minimum energy extension norm,
\begin{equation*}
\Vert \hat{\mathfrak{u}} \Vert_{\hat{U}} := \inf_{\begin{array}{l}\mathfrak{u} \in D(\Aop)\\ \mathrm{tr} \,\mathfrak{u}
= \hat{\mathfrak{u}} \end{array}} \Vert \mathfrak{u} \Vert_{H_\Aop(\Omega)}
\end{equation*}
where 
\begin{equation}\label{eq:ha}
\Vert \mathfrak{u} \Vert_{H_\Aop(\Omega)}^2 := \Vert \Aop \mathfrak{u} \Vert^2 + \Vert \mathfrak{u} \Vert^2 \, .
\end{equation}
\subsection{The Discontinuous Petrov--Galerkin method: ``breaking'' the test space} \label{sec:DPG_UW}
Following the same procedure as above, but with broken (discontinuous) test functions instead, we arrive at the \emph{ultaweak variational formulation with broken test spaces}:
\begin{equation*}
\left\{
\begin{aligned}
&\mathfrak{u} \in \bsL^2(\Omega),\, \hat{\mathfrak{u}} \in \hat{U}_h \\
&\underbrace{(\mathfrak{u}, \Aop^\ast \mathfrak{v}) + \langle \hat{\mathfrak{u}}, \mathfrak{v} \rangle}_{=:b((\mathfrak{u},
\hat{\mathfrak{u}}),\mathfrak{v})}
= (\mathfrak{f},\mathfrak{v}), \quad \mathfrak{v} \in V(\Omega_h) \, .
\end{aligned}
\right.
\end{equation*}
The broken test space is defined as:
\begin{equation*}
V(\Omega_h) := \{ v \in \bsL^2(\Omega) \, : \Aop^\ast_h \mathfrak{v} \in \bsL^2(\Omega) \} 
\end{equation*}
and it is equipped with the adjoint graph norm,
\begin{equation*}
\Vert \mathfrak{v} \Vert^2_{V(\Omega_h)} := \sum_K  
\underbrace{(\Vert \Aop^\ast_h \mathfrak{v} \Vert^2 + \Vert \mathfrak{v} \Vert^2)}_{=: \Vert \mathfrak{v} \Vert^2_{V(K)}}  \, .
\end{equation*}
Here, the symbol $h$ indicates that the operator is understood {\em element--wise}.
Fluxes and traces are defined now on the whole mesh skeleton $\Gamma_h$ and equipped with the minimum energy extension norm.
The trace space is given by
\begin{equation*}
\hat{U}_h  =  \{ (\hat{p},\hat{u}_n) \in H^{\half}(\Gamma_h) \times H^{-\half}(\Gamma_h) \, : \, \hat{p} - \hat{u}_n = 0 \text{ on } \partial \Omega \} \, .
\end{equation*}
Additionally, we define the {\em energy norm} by
\begin{equation}\label{eq:energy_norm}
\Vert (\mathfrak{u},\hat{\mathfrak{u}})\Vert_E^2 
:= \left( \sup_{\mathfrak{v} \in V(\Omega_h)} \frac{b((\mathfrak{u},\hat{\mathfrak{u}}),\mathfrak{v})}
{\Vert \mathfrak{v} \Vert_{V(\Omega_h)}} \right)^2
= \sum_K 
\underbrace{ \left(\sup_{\mathfrak{v} \in V(K)} \frac{\vert( \mathfrak{u},\Aop^\ast_h \mathfrak{v}) + \langle \hat{\mathfrak{u}},
\mathfrak{v} \rangle \vert}
{\Vert \mathfrak{v} \Vert_{V(K)}}\right)^2 }_{=:\Vert (\mathfrak{u},\hat{\mathfrak{u}})\Vert_{E,K}^2  }  
\end{equation}
The above energy norm satisfies the same property as standard Sobolev energy norms, i.e., the global norm squared ($\Vert (\mathfrak{u},\hat{\mathfrak{u}})\Vert_E^2$) equals the sum of element norms squared ($\Vert (\mathfrak{u},\hat{\mathfrak{u}})\Vert_{E,K}^2$) \cite{Carstensen}. This result is very crucial for the analysis of the preconditioner since it allows us to consider a single element, when proving norm equivalence. 

\subsection{The practical DPG method} We note that for practical computations one has to derive the so called \emph{practical DPG method} \cite{Gopalakrishnan1}. That is, the inversion of the Riesz operator can only be approximated, since computing it exactly would require the solution of an infinite dimensional boundary value problem. The procedure is then to consider a truncated finite dimensional test space $V^r \subset V$ with $dim(V^r) >> dim(U_h)$, the so called \emph{enriched test space}. Consequently, both the trial--to--test operator $T$ and the error representation function $\psi$ are also approximated. Various studies exist that examine the effect of this approximation 
to the discrete stability with the overall conclusion that stability is not significantly affected. We refer the reader to \cite{srir_soc,Carstensen,Gopalakrishnan1,CH17} for further reading.
%

%% file: analysis.tex
%
\section{Analysis of the one-level preconditioner}\label{sec:sec3}
We are now ready to provide a theoretical convergence analysis for the one level additive Schwarz smoother. The analysis is based on the \emph{subspace correction theory} \cite{Xu,Xu2,xu1} for self--adjoint positive definite operators. We start by stating the subspace correction theorem and we proceed by verifying its assumptions for the DPG ultraweak formulation for the time--harmonic acoustics problem presented in the previous section. 
\subsection{Additive Schwarz preconditioner and the subspace correction theory}
Let $\sfB$ be self--adjoint and positive definite defined on the vector space $U$. Suppose that $U_i$, $i=1,\ldots J$, are closed subspaces of the Hilbert space $(U, (\cdot,\cdot))$. Additionally, let $\sfB_i$ be self--adjoint and positive definite operator defined on $U_i$ by 
\begin{equation*}
(\sfB_i \sfu_i, \sfv_i) = (\sfB\sfu_i,\sfv_i), \quad \forall \sfu_i, \sfv_i \in U_i,
\end{equation*}
and let $\sfQ_i:U\rightarrow U_i$ denote the $(\cdot,\cdot)$-orthogonal projection onto $U_i$. Then, the operator
$$\sfA = \sum_{i=1}^J \sfB^{-1}_i \sfQ_i$$
is called the \emph{additive} preconditioner based on subspaces ${U_i}$ and operators $\sfB_i$. The additive Schwarz algorithm is given by:
\begin{algorithm}[H]
\caption{Additive Schwarz/parallel subspace correction}\label{alg:PSC}
  \begin{algorithmic}[1]
    \Procedure{PSC}{$\sfu_k,\sfu_{k+1}$}\Comment{Given $\sfu_k$ compute $\sfu_{k+1}$}

    \State $\sfr= \sfell - \sfB \sfu_k$
    \Comment{Compute initial residual}

    \State $\sfr_i = \sfQ_i \sfr$
    \Comment{Project the residual on to $U_i$}

    \State $\sfB_i \sfz_i = \sfr_i$
    \Comment{Solve on the subspace the local problem}

    \State $\sfu_{k+1}= \sfu_k + \theta \sum_{i=i}^J \sfz_i$
    \Comment{Correct $\sfu_k$ on each subspace, $\theta \in (0,1]$}
    \EndProcedure
  \end{algorithmic}
\end{algorithm}
\begin{theorem}[\it{Subspace correction} \cite{Xu2,Xu,xu1}]\label{subspace_thm}
  For the above setting, assume the following two statements hold: 
  \begin{itemize}
    \item (strengthened Cauchy--Schwarz inequality) 
    there exists $\beta>0$ such that for all $\sfu_i, \sfv_i \in U_i$ 
    \begin{equation*}
      \sum^J_{i=1}\sum^J_{j=1} |(\sfu_i, \sfv_j)_\sfB| \le \beta \Big(\sum^J_{i=1} \|\sfu_i\|_\sfB^2 \Big)^{1/2} \Big(\sum^J_{j=1}\|\sfv_j\|_\sfB^2 \Big)^{1/2}
    \end{equation*}
    \item (stable decomposition) there exists $\alpha > 0$ such that $\forall \sfu \in U$, there exists a decomposition  $\sfu = \sum^J_{i=1} \sfu_i, \text{ with } \sfu_i \in U_i,$ that satisfies
    \begin{equation*}
      \sum^J_{i=1} \|\sfu_i\|^2_\sfB \le \alpha^{-1} \|\sfu\|^2_\sfB.
    \end{equation*} 
  \end{itemize}
  Then, the following equivalence relation is true:
  \begin{equation}\label{eq:equival2}
    \alpha (\sfu,\sfu)_\sfB \le (\sfP \sfu, \sfu)_\sfB \le \beta (\sfu,\sfu)_\sfB
  \end{equation}
  where $\sfP = \sum_{i=1}^J\sfP_i$ and $\sfP_i: U \rightarrow U_i$ is the $(\cdot,\cdot)_\sfB$-orthogonal projector, i.e.,
  \begin{equation*}
    (\sfP_i\sfu, \sfv_i)_\sfB = (\sfu, \sfv_i)_\sfB, \quad \forall \sfv_i \in U_i.
  \end{equation*}
  \em  
  Note that $\sfB_i \sfP_i = \sfQ_i \sfB$ and therefore $\sfP = \sfA \sfB$ is the preconditioned matrix. Consequently, \cref{eq:equival2} gives an estimate of the condition number $\kappa(\sfA \sfB) = \beta/\alpha$. A detailed proof of this theorem can be found in \cite{xu1}. 
\end{theorem}
%
\subsection{Verification of the theorem assumptions in the DPG setting}
\subsubsection{\bf Setup}
Let $\mthu = (\mathfrak{u},\mathfrak{\hat{u}})$ denote the group variable including the field and trace variables of the ultraweak formulation. We want to precondition the energy norm $\|\cdot\|_E$ given by \cref{eq:energy_norm}. Let $b_E(\cdot, \cdot$) be the Hermitian and coercive form corresponding to the energy norm, i.e.,
\begin{equation*}
  b_E(\mthu, \mthv) = b(\mthu,\Top\mthv) = (\Top\mthu,\Top\mthv)_V
\end{equation*}
where $\Top$ is the trial--to--test operator for the broken formulation defined according to \cref{eq:eq6}. We introduce $\{\Omega_i\}_{i=1}^{J}$ a finite cover of $\Omega$ such that each $\bar{\Omega}_i$ is the support of a vertex shape function. We denote the size of a vertex patch by $\delta$. In addition, we assume that the cover satisfies the \emph{finite overlap property}, i.e., there exists an integer $r$ such that each point of $\Omega$ is contained in at most $r$ of the sets $\Omega_i$. Equivalently, let $\chi_i$ be the characteristic function of $\Omega_i$. Then
\begin{equation}\label{eq:overlap}
  \sum_{i=1}^J \|\chi_i\|^2_{L^\infty} = \sum_{i=1}^J \|\chi_i\|_{L^\infty} \le r 
\end{equation}
Finally, the local energy subspaces corresponding to the partition are given by:
\begin{equation*}
  U_i = \bsL^2(\Omega_i) \times \hat{U}_i,
\end{equation*}
where $\hat{U}_i:=\{\hat{\mathfrak{u}} \in \hat{U}_h : \hat{\mathfrak{u}}=0 \text{ on } \Gamma_h - \Omega_i\}$, and $\Gamma_h$ is the mesh skeleton.
\subsubsection{\bf Strengthened Cauchy--Schwarz inequality}
Verifying the first assumption, of \cref{subspace_thm} is straight forward. Let $\sfu_i \in U_i$ and $\sfv_j \in U_j$. Then, for any inner product $(\cdot,\cdot)_B$ the inequality 
\begin{equation*}
  |(\sfu_i, \sfv_j)_\sfB| \le \varepsilon_{ij} \|\sfu_i\|_\sfB \|\sfv_j\|_\sfB
\end{equation*}
is true for $\varepsilon_{ij} \le 1$. Indeed, from Cauchy--Schwarz inequality, $\varepsilon_{ij} = 1$ when $(\sfu_i, \sfv_j)_\sfB \ne 0$, but it can be chosen to be $0$ when $(\sfu_i, \sfv_j)_\sfB = 0$. Taking the sum over $i$ and $j$ we have:
\begin{equation*}
   \begin{aligned}
      \sum^J_{i=1}\sum^J_{j=1} |(\sfu_i, \sfv_j)_\sfB|
      & \le \sum^J_{i=1} \sum^J_{j=1} \varepsilon_{ij} \|\sfu_i\|_\sfB \|\sfv_j\|_\sfB \\
      & \le \rho(\mathcal{E}) \Big(\sum^J_{i=1} \|\sfu_i\|_\sfB^2 \Big)^{1/2} \Big(\sum^J_{j=1}\|\sfv_j\|_\sfB^2 \Big)^{1/2}
   \end{aligned}
\end{equation*}
where $\rho(\mathcal{E}) = \|\mathcal{E}\|_2$ is the spectral radius of the matrix $\varepsilon_{ij}$.
Note that the entries of the matrix $\varepsilon_{ij}$ are zero for non overlapping subdomains. Indeed, suppose that $\supp\{\mthu_i\} \subseteq \Omega_i$. Then $\supp\{\Top\mthu_i\} \subseteq \Omega_i$. This is a direct consequence of the use of broken test spaces, i.e., $\Top$, the DPG \emph{trial-to-test operator} is local and therefore $\Top\mthu_i$ is discontinuous. Therefore
\begin{equation*}
   b_E(\mthu_i, \mthv_j) = b(\mthu_i, \Top \mthv_j) = (\Top \mthu_i, \Top \mthv_j)_V = 0, 
\end{equation*}
for all $\,\mthv_j \in \Omega_j\,$ such that $\,\supp\{\mthv_j\} \cap \Omega_i = \emptyset$,
and consequently
\begin{equation*}
   b_E(\mthu_i, \mthv_j) = 0, 
\end{equation*}
for $\Omega_i, \Omega_j$ disjoint and $\supp\{\mthu_i\} \subseteq \Omega_i$, $\supp\{\mthv_j\} \subseteq \Omega_j$. Finally, by the finite overlap assumption, we obtain an upper bound for the spectral radius: 
\begin{equation*}
   \rho(\mathcal{E}) \le r
\end{equation*}
and the result reads:
\begin{equation*}
  \sum^J_{i=1}\sum^J_{j=1} |b_E(\sfu_i, \sfv_j)| \le r \Big(\sum^J_{i=1} \|\sfu_i\|_E^2 \Big)^{1/2} \Big(\sum^J_{j=1}\|\sfv_j\|_E^2 \Big)^{1/2}
\end{equation*}
\subsubsection{\bf Stable Decomposition}
Proving the second assumption of \cref{subspace_thm} is a bit more involved. This 
is exactly the well-known {\em stable decomposition} assumption, which in simple terms it says that if $\sfu$ is decomposed into patch contributions, then the sum of energies stored in the patches must be controlled by the energy of $\sfu$.

\noindent
We start with a DPG global stability result \cite{Carstensen,demk2},
\begin{equation*}
\Vert \mathfrak{u} \Vert^2 + \Vert \hat{\mathfrak{u}} \Vert_{\hat{U}}^2
\leq \left[ \frac{1}{\gamma^2} + (1 + \frac{1}{\gamma})^2 \right] \, \Vert (\mathfrak{u}, \hat{\mathfrak{u}}) \Vert_E^2  \, .
\end{equation*}
Above $\gamma$ denotes the global boundedness below constant for operator $A$. For the acoustics operator and the case of impedance BC, constant $\gamma$ is independent of frequency $\omega$ \cite{demk2}. Combining the stability estimate with continuity we derive the equivalence relation:
\begin{equation*}
\gamma_1^2 (\|\mathfrak{u}\|^2 + \|\mathfrak{\hat{u}}\|^2_{\hat{U}}) \le 
\|(\mathfrak{u}, \hat{\mathfrak{u}})\|^2_E \le M_1^2 (\|\mathfrak{u}\|^2 + 
\|\mathfrak{\hat{u}}\|^2_{\hat{U}})\, ,
\end{equation*}
where $\gamma_1^{-2} = \frac{1}{\gamma^2} + (1 + \frac{1}{\gamma})^2$ 
and $M_1$ is a continuity constant.
The above relation allows us to construct a stable decomposition for the original trial norm instead of the energy norm. We will prove the result by considering separate cases for the $L^2$ space and the trace space $\hat{U}$.
For simplicity we consider a the two--dimensional case for the rest of the proof. The three--dimensional case is fully analogous. 
\medskip
\paragraph[Stable decomposition for the $L^2$ space]{\bf Stable decomposition for the $L^2$ space}
Let $Y_p = \mathcal{Q}^{(p,q)} = \mathcal{P}^p\otimes \mathcal{P}^q$ denote the space of polynomials of order less or equal $p,q$ with respect to $x,y$ respectively. These spaces are described in detail in \cite[Ch.2]{petrides2019adaptive}. Additionally, let $\{\Phi_j\}_{j=1}^J$ be a partition of unity subordinate to the covering $\Omega_j$ of $\Omega$, so that
\begin{equation*}
   \sum_{j=1}^J \Phi_j = 1, \qquad 0\le\Phi_j\le 1, \qquad \supp(\Phi_j) \subset \Omega_j
\end{equation*}
Given $\mathfrak{u} \in Y_p$ we define the following decomposition
\begin{equation*}
   \mathfrak{u}_j = \sfP \Phi_j \mathfrak{u}
\end{equation*}
where $\sfP:\mathcal{Q}^{(p+1,q+1)} \rightarrow \mathcal{Q}^{(p,q)}$ is the $L^2$--orthogonal projection.
Obviously
\begin{equation*}
\sum_{i=1}^J \mathfrak{u}_j = \sum_{j=1}^J \sfP \Phi_j \mathfrak{u} = \sfP \underbrace{\sum_{i=1}^J \Phi_j}_{=1} \mathfrak{u}_j = \sfP \mathfrak{u} = \mathfrak{u}
\end{equation*}
It easy to see that the decomposition is stable. Indeed
\begin{equation*}
   \|\mathfrak{u}_j\| = \|\sfP \Phi_j \mathfrak{u}\| = \|\sfP \Phi_j \mathfrak{u}\|_{L^2(\Omega_j)} \le \|\sfP\|\|\Phi_j\|_{L^\infty} \|\mathfrak{u}\|_{L^2(\Omega_j)} \le \|\mathfrak{u}\|_{L^2(\Omega_j)}
\end{equation*}
Summing up for all subdomains, and by the finite overlap property \eqref{eq:overlap} we have
\begin{equation*}
\sum_{j=1}^J \|\mathfrak{u}_j\|^2 \le \sum_{j=1}^J \|\mathfrak{u}\|^2_{L^2(\Omega_j)} \le  r \|\mathfrak{u}\|^2
\end{equation*}
\paragraph[Stable decomposition for the trace space $\hat{U}$ space]{\bf Stable decomposition for the trace space $\hat{U}$}

We shall consider a single element $K \in \Omega_j$. We denote by $H_A(K)$ the graph energy space of functions defined on element $K$. Recall that $\bar{\Omega}_j$ is defined by the support of a vertex shape function $\Phi_j$ of size $\delta$. Then, $\{\Phi_j\}_{i=1}^J$ is a partition of unity,
\begin{equation*}
  \sum_{j=1}^J \Phi_j = 1, \qquad 0\le\Phi_j\le 1, \qquad \supp(\Phi_j) \subset \Omega_j
  \quad \text{and} \quad 
  \|\nabla \Phi_j\|_{L^\infty} \le C \delta^{-1} 
\end{equation*}
Let $H_A(K)$ be the graph energy space of functions defined on element $K$, we denote by
$H_A^j(K)$ the subspace of functions vanishing on $K - \supp(\Phi_j)$. 

\medskip
\paragraph{\bf Operator dependent projection based interpolant.} Let $\mathfrak{\hat{u}}$ be a sufficiently regular function defined on the skeleton. Additionally, let's assume the existence of a projection-based interpolant $\hat{\Pi}\mathfrak{\hat{u}}$ (which will be defined explicitly later on) taking values on the mesh skeleton such that
\begin{equation}\label{eq:proj_prop}
\|\hat{\Pi}\Phi_j\mathfrak{\hat{u}}\|_{\hat{U}} \le \underbrace{\|\hat{\Pi}\|_{\hat{U}}}_{=:C_{\hat{\Pi}}} \|\Phi_j\mathfrak{\hat{u}}\|_{\hat{U}}
\end{equation}
Construction of such an interpolant is operator-dependent, and we expect constant $C_{\hat{\Pi}}$
to grow mildly with the frequency $\omega$
for the acoustics operator,  and be independent of the polynomial order $p$ and the fine mesh size $h$. We investigate the dependence on $\omega$, $h$, and $p$ numerically by constructing such an interpolation operator. The construction and the numerical experiments are described in the next section.
Consider now a function on the skeleton $\mathfrak{\hat{u}}$. We define the corresponding decomposition by
\begin{equation*}
  \mathfrak{\hat{u}}_j = \hat{\Pi} \Phi_j \mathfrak{\hat{u}} \, .
\end{equation*}
Clearly,
\begin{equation*}
\sum_{j=1}^J \mathfrak{\hat{u}}_j  = \sum_{j=1}^J \hat{\Pi} \Phi_j \mathfrak{\hat{u}}
=  \hat{\Pi} (\underbrace{\sum_{j=1}^J \Phi_j}_{=1} \mathfrak{\hat{u}}) = \hat{\Pi} \mathfrak{\hat{u}} 
= \mathfrak{\hat{u}} 
\end{equation*}
and by the postulated property \eqref{eq:proj_prop}
\begin{equation*}
\|\mathfrak{\hat{u}}_j\|_{\hat{U}} = \|\hat{\Pi}\Phi_j\mathfrak{\hat{u}}\|_{\hat{U}} \le C_{\hat{\Pi}} \|\Phi_j\mathfrak{\hat{u}}\|_{\hat{U}} \le C_{\hat{\Pi}} \|\Phi_j\mathfrak{u}\|_{H_A}
\end{equation*}
where $\mathfrak{u}$ is an extension of $\mathfrak{\hat{u}}$ to the element such that $tr\mathfrak{u} = \mathfrak{\hat{u}}$. 
The last term bounds nicely for the acoustics operator,
\begin{equation*}
\|\Phi_j\mathfrak{u}\|^2_{H_\Aop} = \|\Aop(\Phi_j \mathfrak{u}))\|^2+\|(\Phi_j \mathfrak{u})\|^2
\end{equation*}
and
\begin{equation*}
\begin{aligned}
\|\Aop(\Phi_j(u,p))\|^2 &=\|i \omega \Phi_j p +\text{div}(\Phi_j u) \|^2+\|i \omega u+\nabla(\Phi_j p)\|^2 \\
                &\le 2\Big(\|\Phi_j (i \omega p + \text{div }u)\|^2 + \|\Phi_j(i \omega u+\nabla p)\|^2 + \|\nabla \Phi_j \cdot u \|^2 + \| \nabla \Phi_j p\|^2 \Big) \\
                &\le 2\Big( \|\Aop (u,p) \|^2 + \delta^{-2} (\| u \|^2 + \|  p \|^2)\Big) \, 
\end{aligned}
\end{equation*}
Therefore 
\begin{equation*}
\|\mathfrak{\hat{u}}_j\|_{\hat{U}} \le C_{\hat{\Pi}} C_{\delta} \|\mathfrak{u}\|_{H_A}, \quad \text{where } \quad C_{\delta} = \mathcal{O}(\delta^{-1})
\end{equation*}
Note that the constant $C_\delta$ is independent of frequency. 
The above inequality is true for an arbitrary extension $\mathfrak{u}$ of $\mathfrak{\hat{u}}$. Taking the infimum on the right hand side with respect to extensions and summing over all subdomains we obtain the stable decomposition for the trace space. 
Combining everything together we have the following lemma.
\begin{lemma}
For every function $\sfu = (\mathfrak{u},\mathfrak{\hat{u}}) \in \bsL^2(\Omega) \times \hat{U}$ there exists a stable decomposition
\begin{equation*}
\sfu = \sum_{j=1}^J \sfu_j,\quad \sfu_j \in \bsL^2(\Omega_j) \times \hat{U}_j,
\end{equation*}
such that
\begin{equation*}
\sum_{j=1}^J \| \sfu_j \|^2_{E} \le C \| \sfu \|^2_{E} \,.
\end{equation*}
where $C = \mathcal{O}(r C^2_{\hat{\Pi}} C^2_\delta)$. Here $r$ comes from the finite overlap property and $C_\delta = \mathcal{O}(\delta^{-1})$. 
This in turn gives an upper bound for the condition number of the preconditioned DPG system
\begin{equation*}
\kappa(\sfA\sfB) \le C \frac{r^2}{\delta^2}, \quad \text{where} \quad C = \mathcal{O}(C^2_{\hat{\Pi}}).
\end{equation*}
\end{lemma}
It remains to study the dependence of the interpolation norm $C_{\hat{\Pi}}$ on the discretization size $h$, polynomial order $p$ and the frequency $\omega$. For this, we design a numerical experiment to compute the interpolation norm. We present the results below. Further details describing the design of the experiment are given in \cref{appendix:interpolant}.
\subsection{Numerical study of \texorpdfstring{$\|\hat{\Pi}\|_{\hat{U}}$}{}}\label{one_level_res}
We examine the dependence of the norm of the interpolation operator described above, on the polynomial order $p$, the discretization size $h$, and the frequency $\omega$. We present the results in the tables below for polynomial orders $p=2,4,6$. For each case of $p$ we also present convergence of the preconditioned CG solver. The subdomains are defined to be the support of a vertex function defined on a quadrilateral uniform mesh of size $h=1/2$ (i.e, 9 subdomains with a fixed overlap size $\delta = 1/2$. 
\begin{table}[H]
\begin{subtable}{0.36\textwidth}\centering
\begin{tabular}{c??ccccc}\toprule\toprule
  h \textbackslash $\omega$ & $\pi$ & $2\pi$ & $4\pi$ & $8\pi$ & $16\pi$\\
 \bottomrule
 \toprule
 1/2  &   15  &  16  &   17   &  {\color{red} 11}  &  {\color{red} 8} \\
 1/4  &   16  &  18  &   20   &  {\color{red} 14}  &  {\color{red} 9} \\
 1/8  &   16  &  18  &   22   &  23  &  {\color{red} 10} \\
 1/16 &   16  &  18  &   23   &  24  &  25 \\
 1/32 &   16  &  19  &   23   &  25  &  25 \\
 \bottomrule 
 \bottomrule 
 \end{tabular}
 \caption{Iteration count}
\end{subtable}
\hspace{50pt}
\begin{subtable}{0.49\textwidth}
   \centering
\begin{tabular}{c??ccccc}\toprule\toprule
  h \textbackslash $\omega$ & $\pi$ & $2\pi$ & $4\pi$ & $8\pi$ & $16\pi$\\
 \bottomrule
 \toprule
         1/2  &  2.553 & 4.327  & 6.088  & {\color{red}3.768} &  {\color{red}2.268} \\
         1/4  &  2.540 & 4.427  & 7.598  & {\color{red}6.320} &  {\color{red}4.328} \\
         1/8  &  2.510 & 4.415  & 6.972  & 8.247 &  {\color{red}6.061} \\
         1/16 &  2.500 & 4.362  & 6.841  & 8.379 &  9.555 \\
         1/32 &  2.497 & 4.343  & 6.785  & 8.651 &  9.621 \\
 \bottomrule
 \bottomrule
\end{tabular}
\caption{Value of the constant $\|\hat{\Pi}\|_{\hat{U}}$}
\end{subtable}
\caption{Polynomial order $p=2$. Left: iteration count for CG preconditioned with additive Schwarz smoother with fixed $\delta = 1/2$. Right: The value of the interpolation norm.}
\end{table}
\raggedbottom
\vspace{-20pt}
\begin{table}[H]
\begin{subtable}[t]{0.36\textwidth}\centering
\begin{tabular}{c??ccccc}\toprule\toprule
  h \textbackslash $\omega$ & $\pi$ & $2\pi$ & $4\pi$ & $8\pi$ & $16\pi$\\
 \bottomrule
 \toprule
   1/2  &   16  &  18  &   22   &  {\color{red} 17 }  & {\color{red} 9} \\
   1/4  &   16  &  18  &   22   &  22  &  {\color{red} 11 } \\
   1/8  &   16  &  18  &   22   &  25  &  23 \\
   1/16 &   17  &  19  &   23   &  25  &  25 \\
   1/32 &   17  &  19  &   23   &  25  &  26 \\
 \bottomrule
 \bottomrule
\end{tabular}
 \caption{Iteration count}
\end{subtable}
\hspace{50pt}
\begin{subtable}[t]{0.5\textwidth}
   \centering
\begin{tabular}{c??ccccc}\toprule\toprule
  h \textbackslash $\omega$ & $\pi$ & $2\pi$ & $4\pi$ & $8\pi$ & $16\pi$\\
 \bottomrule
 \toprule
   1/2  &  1.516 & 1.976 & 3.222 & {\color{red}3.168} & {\color{red}1.871} \\
   1/4  &  1.497 & 2.021 & 3.328 & 5.262 & {\color{red}2.582} \\
   1/8  &  1.491 & 2.046 & 3.432 & 5.122 & 5.891 \\
   1/16 &  1.490 & 2.055 & 3.492 & 5.260 & 6.478 \\
   1/32 &  1.489 & 2.057 & 3.512 & 5.352 & 6.556 \\
 \bottomrule
 \bottomrule
\end{tabular}
\caption{Value of the constant $\|\hat{\Pi}\|_{\hat{U}}$}
\end{subtable}
\caption{Polynomial order $p=4$. Left: iteration count for CG preconditioned with additive Schwarz smoother with fixed $\delta = 1/2$. Right: The value of the interpolation norm.}
\end{table}
\raggedbottom
\vspace{-20pt}
\begin{table}[H]
\begin{subtable}{0.36\textwidth}\centering
\begin{tabular}{c??ccccc}\toprule\toprule
  h \textbackslash $\omega$ & $\pi$ & $2\pi$ & $4\pi$ & $8\pi$ & $16\pi$\\
 \bottomrule
 \toprule
   1/2  &  17   &   19  &   23  &   24  &   {\color{red} 12}  \\
   1/4  &  16   &   18  &   22  &   23  &   22  \\
   1/8  &  17   &   18  &   22  &   24  &   24  \\
   1/16 &  17   &   18  &   22  &   25  &   25  \\
   1/32 &  17   &   18  &   23  &   25  &   26  \\
 \bottomrule
 \bottomrule
\end{tabular}
 \caption{Iteration count}
\end{subtable}
\hspace{50pt}
\begin{subtable}{0.5\textwidth}
   \centering
\begin{tabular}{c??ccccc}\toprule\toprule
  h \textbackslash $\omega$ & $\pi$ & $2\pi$ & $4\pi$ & $8\pi$ & $16\pi$\\
 \bottomrule
 \toprule
   1/2  & 1.271 & 1.522 & 1.874 & 3.514 & {\color{red}2.391} \\
   1/4  & 1.268 & 1.512 & 1.753 & 3.093 & 3.026 \\
   1/8  & 1.267 & 1.510 & 1.734 & 2.750 & 3.583 \\
   1/16 & 1.267 & 1.510 & 1.730 & 2.685 & 3.724 \\
   1/32 & 1.266 & 1.510 & 1.730 & 2.671 & 3.962 \\
 \bottomrule
 \bottomrule
\end{tabular}
\caption{Value of the constant $\|\hat{\Pi}\|_{\hat{U}}$}
\end{subtable}
\caption{Polynomial order $p=6$. Left: iteration count for CG preconditioned with additive Schwarz smoother with fixed $\delta = 1/2$. Right: The value of the interpolation norm.}
\end{table}
\raggedbottom
For these experiments we simulate a plane wave in the square domain $\Omega = (0,1)^2$, using impedance boundary conditions on the entire boundary of the domain. We run our simulations for five different frequencies and we perform successive uniform $h$-refinements, starting with a mesh of size $h=1/2$. The CG solver is terminated when the $l_2$-norm of the residual drops below $10^{-6}$. The numbers in red color denote that the error is above 90\%, i.e, the mesh is not fine enough to resolve the wave.

The following interesting observations can be derived from the results. First, the solver shows convergence, independent of the polynomial order. This, can also be verified in the value of the constant $\hat{\Pi}_{\hat{U}}$. Gopalakrishnan and Sch{\"o}berl observed the same behavior in their construction of a multiplicative Schwarz preconditioner in \cite{Gopalakrishnan2015}. The second observation involves the dependence on the mesh size $h$. As we can clearly see, both the constant and the CG convergence are independent of $h$. 

Lastly, the number of iterations of the solver grows mildly with respect to the frequency, an observation that is also reflected in the value of the constant. However, in cases where the mesh is too coarse to resolve a high frequency wave, the solver is more efficient with respect to the number of iterations (see numbers in red). This behavior, might not seem to be of great value for uniform meshes, but in case of adaptive refinements, this could be very beneficial. As we have seen in the previous chapter, the DPG method being unconditionally stable, allows for adaptive refinements starting from very coarse meshes. Therefore, integrating an iterative solver within the adaptive process, even in the pre-asymptotic region, seems to be a promising direction to follow.  

Overall, the computed interpolation norm gives a very useful insight on the effectiveness of the preconditioner. We emphasize that, there is no reason to believe that the bound provided by the interpolation norm is sharp. However, the overall trend of the value of the norm, is consistent with the convergence behavior of the solver. 

Note that in practical computations the interior degrees of freedom ($L^2$ field variables) can be eliminated in an element--wise fashion. The size of the final system is then significantly reduced, since the condensed system involves only the trace unknowns. A well known Schur complement result guarantees that the upper bound of the condition number of the condensed and original operators are the same \cite[Appendix C]{petrides2019adaptive}. This allows us to analyze the original system even if in actual computations we use the statically condensed system.

%% file: mg.tex
%

\section{Extension to multigrid}
\label{sec:sec4}
We would like to extend our construction to the two-- (multi--) level setting, i.e, accelerate the preconditioner by coupling (in the multiplicative way) the additive Schwarz smoother with a coarse grid solve. The main reason is the unavoidable dependence of the one--level preconditioner on the number of subdomains and the size of overlap. In order to, keep the number of iterations of the solver under control, both the number of subdomains and the size of the overlap have to remain constant and unfortunately this adds significant work on the local solves within each subdomain as $h$ decreases. On the other hand, if the number of subdomains increases, the cost of each local solve remains constant, but the overall number of iterations grows. A preliminary implementation for a two--level setting for acoustics for two dimensional problems was introduced in \cite{sp}. Here, we extend that work to three dimensional simulations and to the multi--level setting.
\medskip
\subsection{Discussion on implementation}
\begin{figure}[H]
 \captionsetup[subfigure]{labelformat=empty}
 \centering
\begin{subfigure}{0.12\textwidth}
    \centering
    \includegraphics[width=0.9\linewidth]{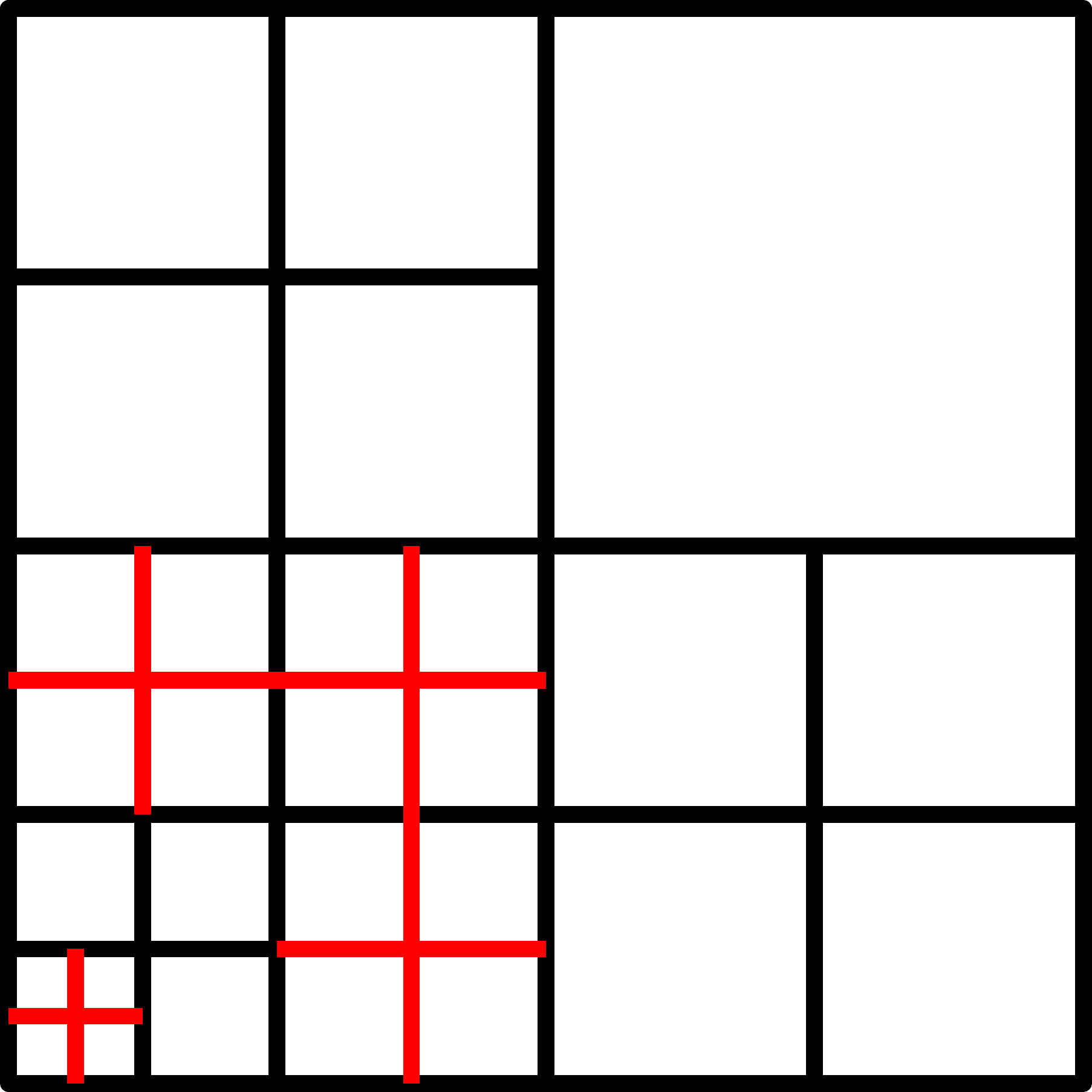}
    \caption{{Fine grid}}
\end{subfigure}
\begin{subfigure}{0.03\textwidth}
    \centering
    \includegraphics[width=1.0\linewidth]{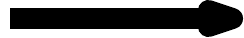}
  \end{subfigure}
  \begin{subfigure}{0.12\textwidth}
    \centering
    \includegraphics[width=0.9\linewidth]{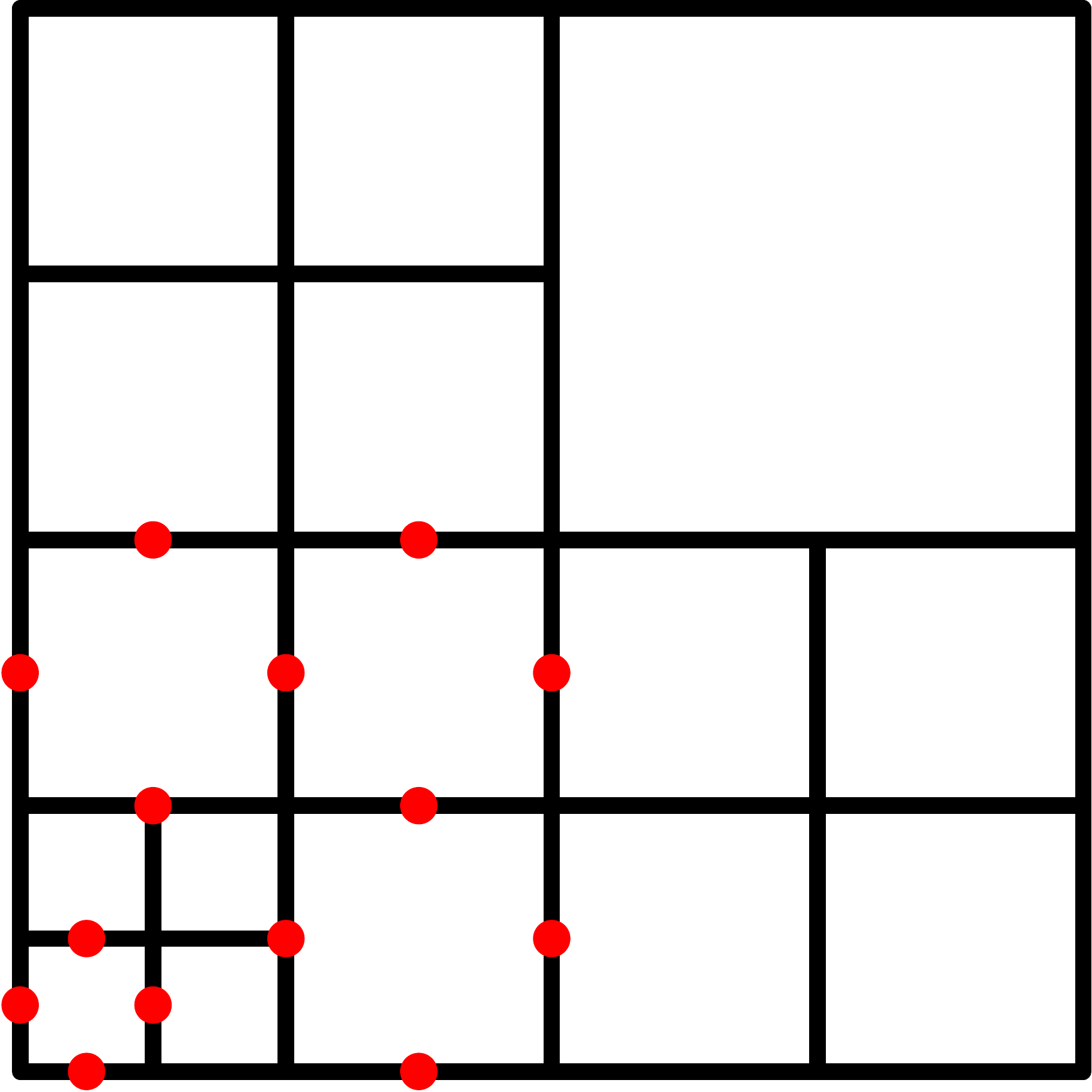}
\caption{{Macro grid}}
\end{subfigure}
\hspace{160pt}
  \begin{subfigure}{0.12\textwidth}
    \centering
    \includegraphics[width=0.9\linewidth]{figures/mg_schematic/macrogrid2.pdf}
\caption{{Macro grid}}
\end{subfigure}
\begin{subfigure}{0.03\textwidth}
    \centering
    \includegraphics[width=1.0\linewidth]{figures/mg_schematic/arrow_right.pdf}
  \end{subfigure}
\begin{subfigure}{0.12\textwidth}
    \centering
    \includegraphics[width=0.9\linewidth]{figures/mg_schematic/finegrid2.pdf}
\caption{{Fine grid}}
 \end{subfigure}

\begin{subfigure}{0.05\textwidth}
    \centering
    \includegraphics[width=0.3\linewidth]{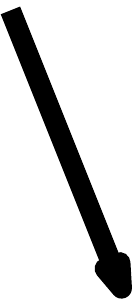}
  \end{subfigure}
  \hspace{170pt}
\begin{subfigure}{0.05\textwidth}
    \centering
    \includegraphics[width=0.3\linewidth]{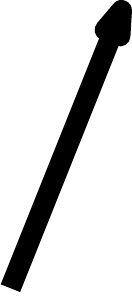}
  \end{subfigure}

\begin{subfigure}{0.12\textwidth}
    \centering
    \includegraphics[width=0.9\linewidth]{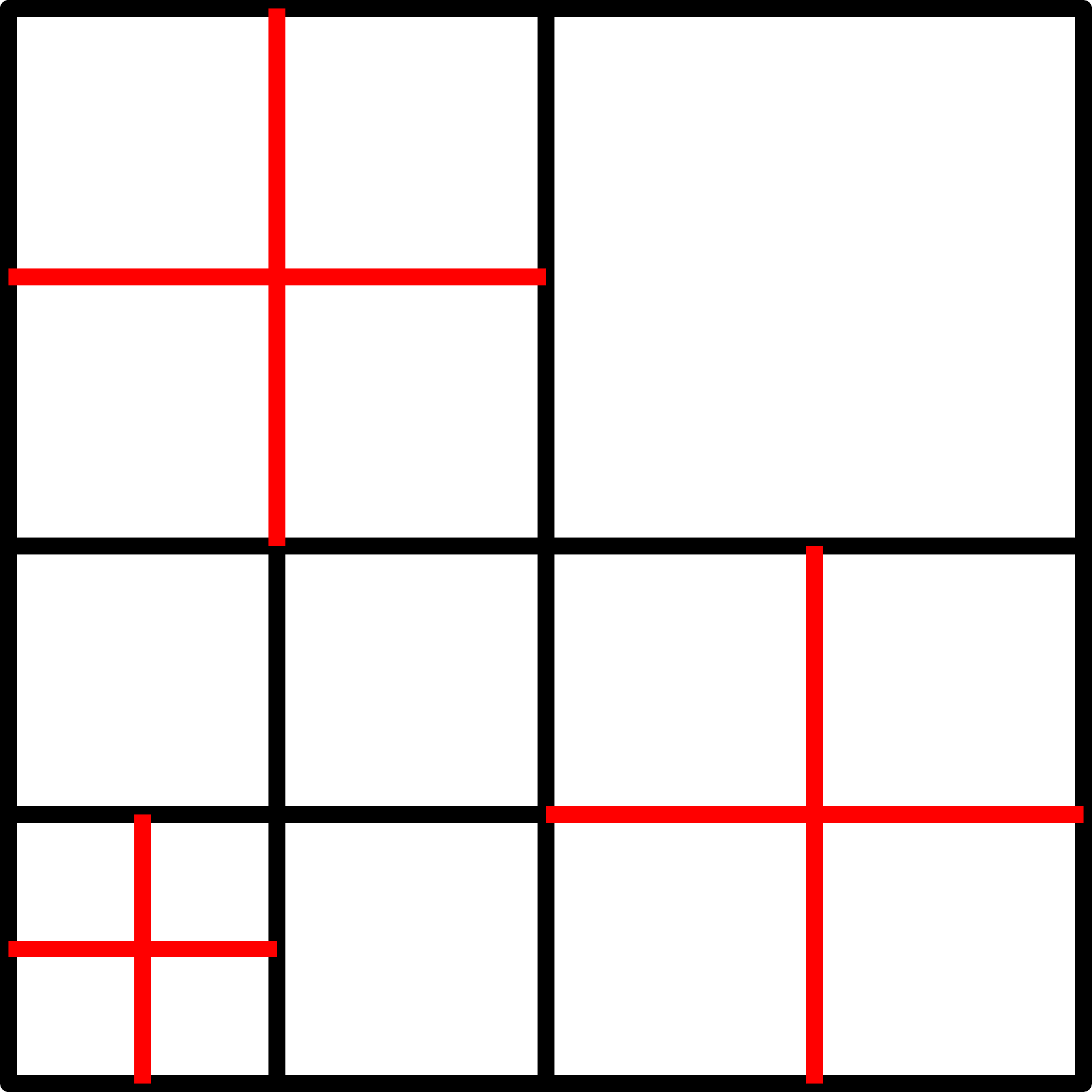}
\caption{Interm. grid}
  \end{subfigure}
\begin{subfigure}{0.03\textwidth}
    \centering
    \includegraphics[width=1.0\linewidth]{figures/mg_schematic/arrow_right.pdf}
  \end{subfigure}
  \begin{subfigure}{0.12\textwidth}
    \centering
    \includegraphics[width=0.9\linewidth]{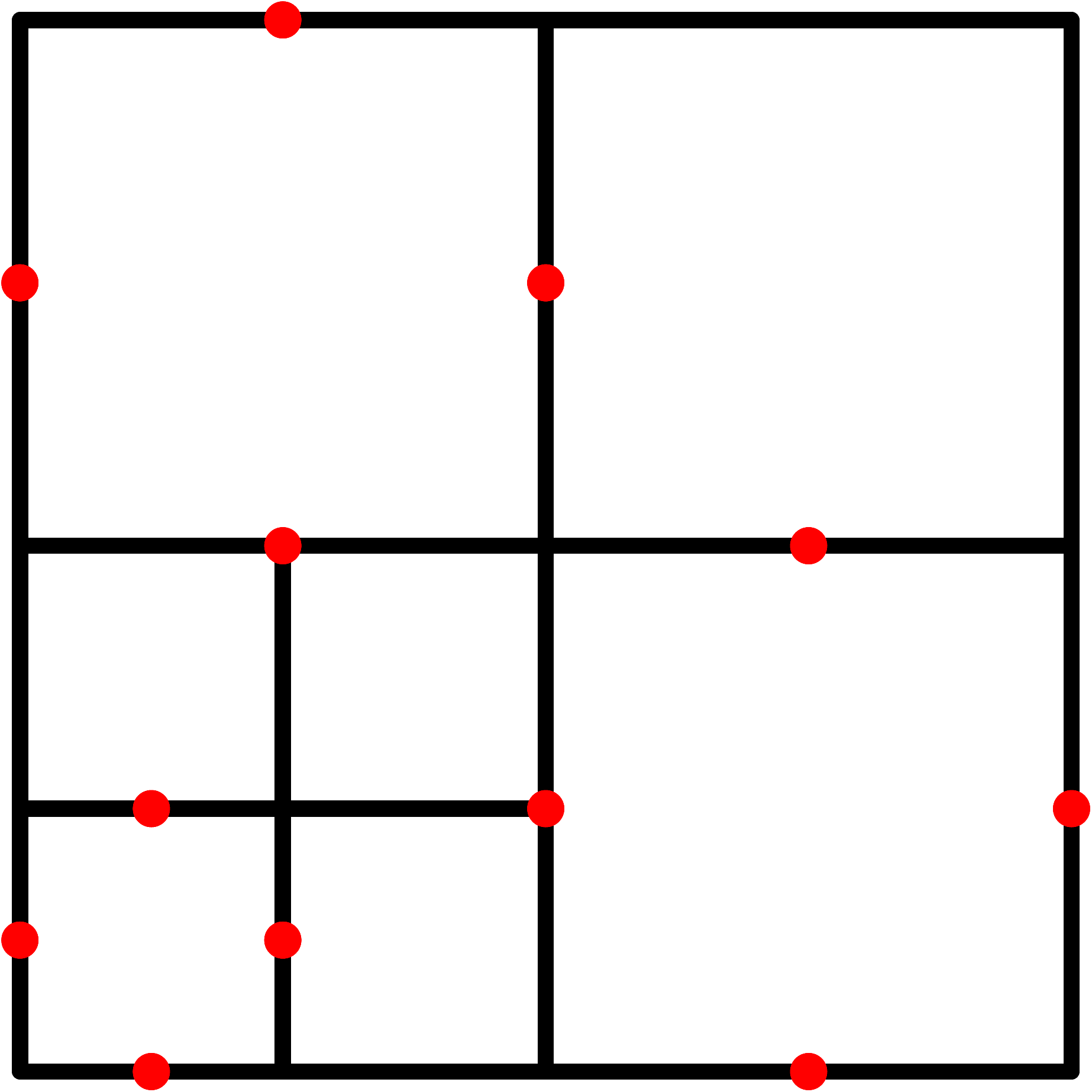}
\caption{{Macro grid}}
\end{subfigure}
\hspace{15pt}
  \begin{subfigure}{0.12\textwidth}
    \centering
    \includegraphics[width=0.9\linewidth]{figures/mg_schematic/macrogrid1.pdf}
\caption{{Macro grid}}
\end{subfigure}
\begin{subfigure}{0.03\textwidth}
    \centering
    \includegraphics[width=1.0\linewidth]{figures/mg_schematic/arrow_right.pdf}
  \end{subfigure}
\begin{subfigure}{0.12\textwidth}
    \centering
    \includegraphics[width=0.9\linewidth]{figures/mg_schematic/finegrid1.pdf}
\caption{{Interm. grid}}
  \end{subfigure}

\begin{subfigure}{0.04\textwidth}
    \centering
    \includegraphics[width=0.4\linewidth]{figures/mg_schematic/arrow_down2.pdf}
  \end{subfigure}
  \hspace{10pt}
\begin{subfigure}{0.04\textwidth}
    \centering
    \includegraphics[width=0.4\linewidth]{figures/mg_schematic/arrow_up2.pdf}
  \end{subfigure}

\begin{subfigure}{0.12\textwidth}
    \centering
    \includegraphics[width=1.0\linewidth]{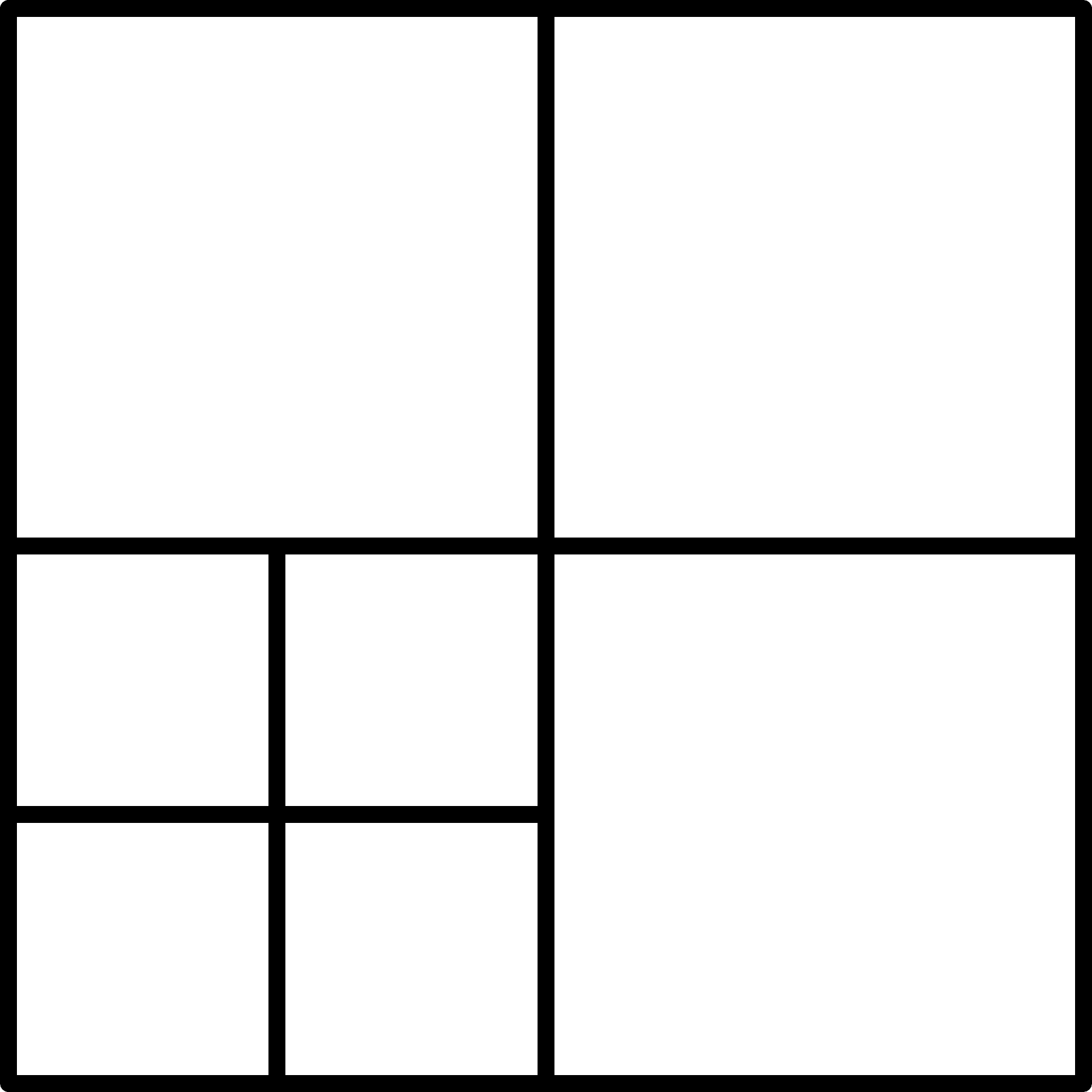}
\caption{{Coarse grid}}
  \end{subfigure}
  \caption{Multigrid v-cycle schematic. For demonstration purposes the schematic is in 2D. In 3D it is fully analogous}
\label{fig:mg_constr}
\end{figure}

The major components of the multigrid preconditioner are: a) inter-grid transfer operators, b) smoother and c) coarse grid solver. We describe the construction of each component for the ultraweak formulation for the acoustics problem presented in Section \ref{sec:DPG_UW}. 
\medskip
\paragraph{\bf Inter-grid transfer operators}
Recall that the ultraweak formulation involves two sets of variables, i.e., the field variables $(u,p) \in U=\boldsymbol{L^2}(\Omega)$ and the \emph{trace} variables 
$(\hat{p},\hat{u}_n) \in \hat{U}$ given by \cref{eq:trace_space}. Since the field variables are only $L^2$--conforming, they can be element--wise eliminated from the global system. The resulting linear system therefore involves only the trace unknowns. Constructing, a prolongation operator from a coarse trace space to an $h$--refined fine trace space is not so straight forward. The reason is simply because new edges (faces in 3D) which are created after an $h$--refinement hold new interface variables which have no ancestors. Therefore the usual prolongation operator (natural inclusion) based on constrained approximation \cite{Demkbook2} is not well defined. Following our work in \cite{sp}, we overcome this difficulty by considering a two--step procedure for the prolongation. First, we construct the so called \emph{macro--grid}. This is the mesh that is constructed by eliminating all the trace variables that do not lie on the coarse grid skeleton using a Schur complement restriction operator. The final inter--grid transfer operator is then defined to be the composition of a natural inclusion operator from the coarse grid to the macro grid and a Schur--complement extension operator from the macro--grid to the fine grid. This construction is demonstrated in Figure \ref{fig:mg_constr}. For the $p$--refinement case with hierarchical shape functions, the new degrees of freedom are simply set to zero. This is the standard inclusion operator for the $p$--multigrid algorithm. Lastly, the restriction operator is defined to be the transpose of the prolongation operator. 
\begin{figure}[H]
\begin{subfigure}{0.49\textwidth}
\centering
   \includegraphics[width=0.5\linewidth]{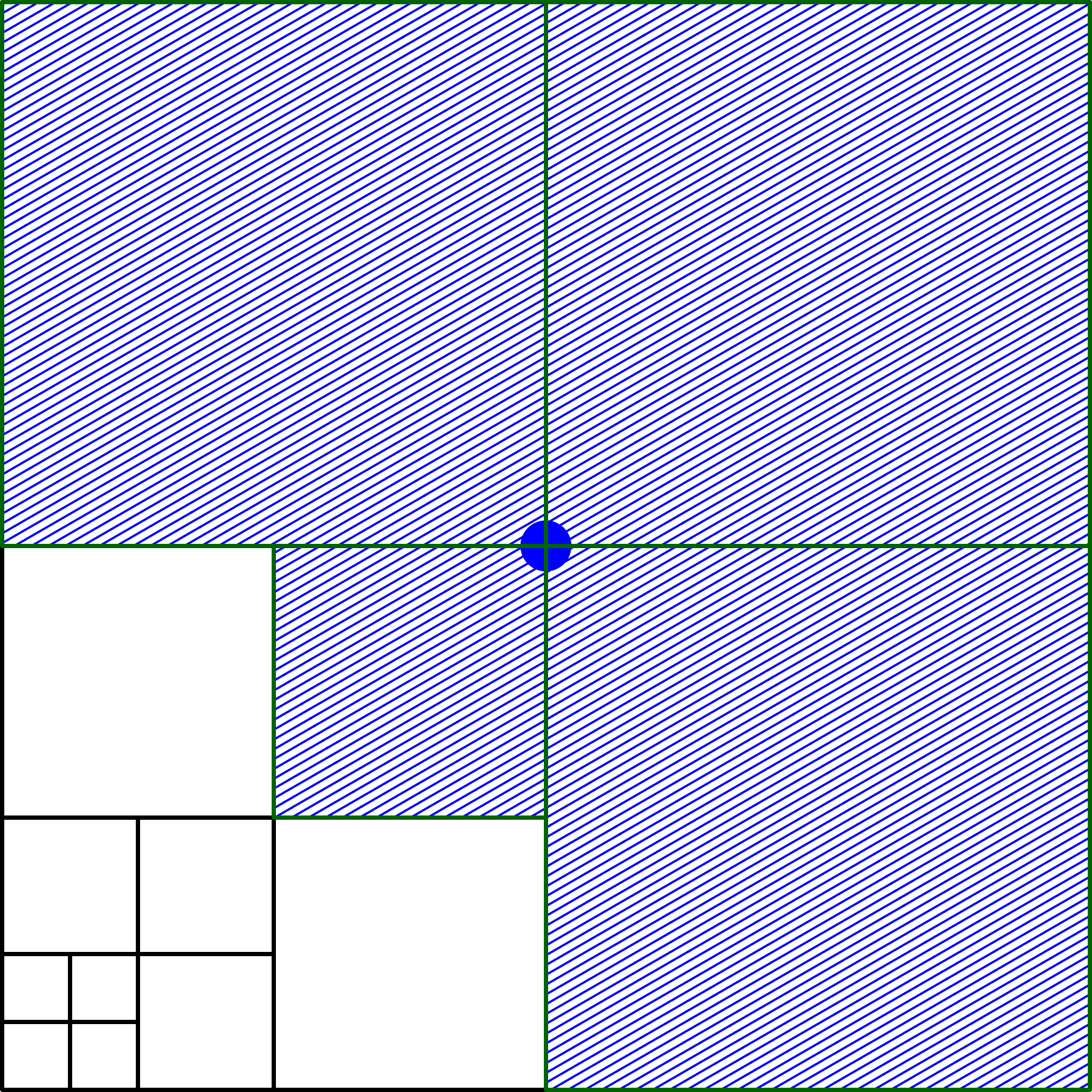}
   \caption{Vertex patch - coarse grid}
\end{subfigure}
\begin{subfigure}{0.49\textwidth}
   \centering
   \includegraphics[width=0.5\linewidth]{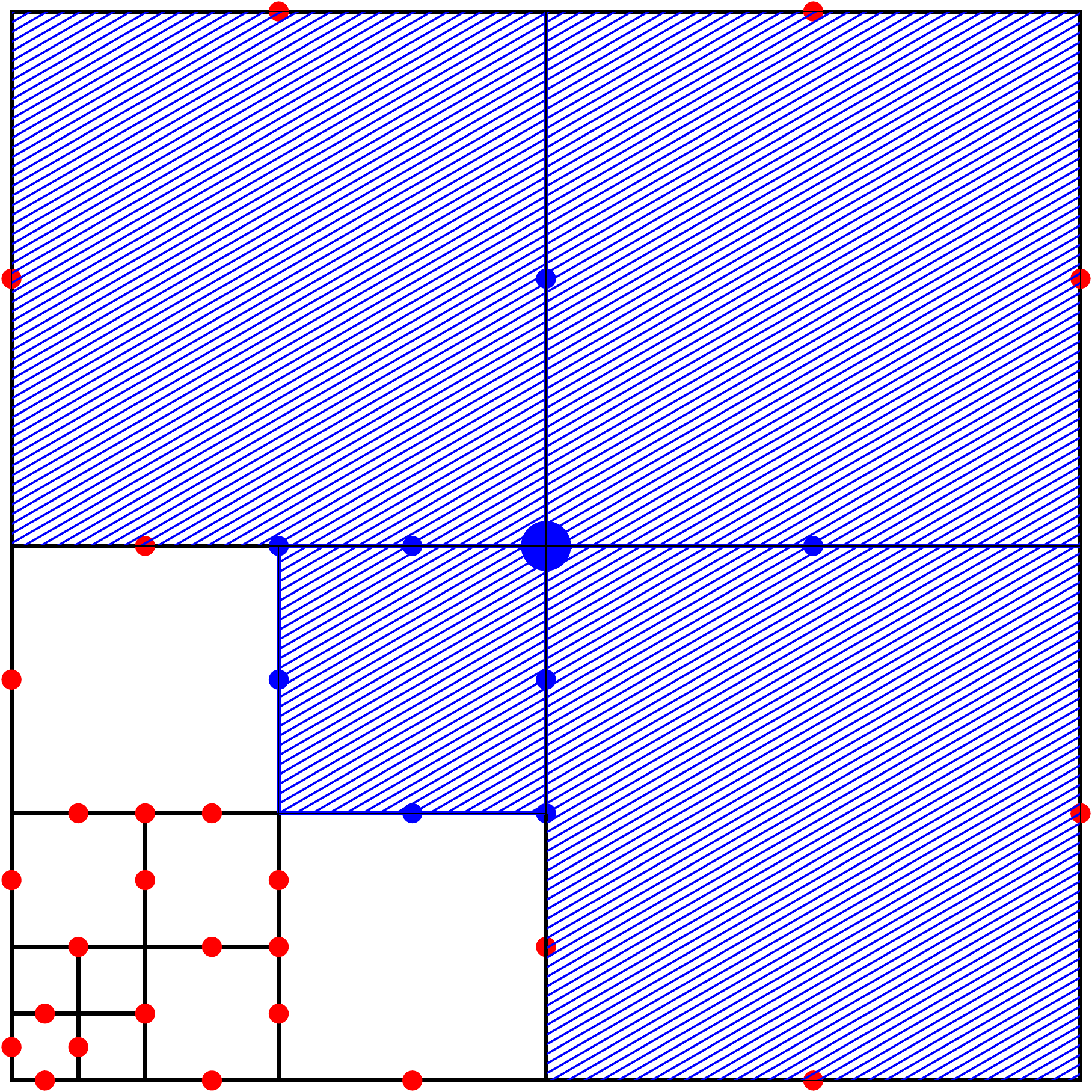}
   \caption{Vertex patch - macro grid}
\end{subfigure}
\caption{Construction of a smoother patch. A smoother patch is defined by the support of a coarse grid vertex basis function.}
\label{fig:patch}
\end{figure}
\paragraph{\bf Smoother}
As a smoother, we use the additive Schwarz procedure described in the previous section. The Schwarz patches are defined to be the support of the vertex basis functions of a coarser grid which is 1 level down the multigrid hierarchy, and the smoothing solves are performed on the macro-grids, see Figure \ref{fig:patch} for illustration. 
\medskip
\paragraph{\bf Coarse-grid solve} For the coarse grid correction we use a direct solver, specifically the Multi--frontal Sparse Direct Solver MUMPS \cite{mumps1,MUMPS2}. 
Note that the coarse grid has the field variables eliminated, and its size is usually significantly smaller than the size of the fine--grid system.

%% file: results.tex

\section{Numerical results}
\label{sec:sec5}
This chapter is devoted to numerical results for relatively large simulations, using the Conjugate Gradient (CG) solver preconditioned with the multigrid technology. Our interest is in computationally challenging acoustics and electromagnetic problems. We present various numerical experiments in acoustics and Maxwell equations in both uniform and adaptive refinement setting. All the simulations were performed on a single node with 24 cores and 256G of RAM memory, using the Finite Element Library \emph{hp3D} \cite{Demkbook1,Demkbook2,keith}.
%
\input{results_maxwell}
\input{results_acoustics}

%

%% file: results_maxwell.tex

\subsection{Time harmonic Maxwell equations}

We first consider the time harmonic form of the Maxwell equations given by 
\begin{equation*}
 \left\{
 \begin{aligned}
 i \omega \mu H + \nabla \times E & = 0   \\
-i \omega \epsilon E + \nabla \times H & = J  
   \end{aligned}
 \right.
\end{equation*}
where $E$ is the electric field, $H$ is the magnetic field, $\omega$ is the angular frequency, $\epsilon$ is the permittivity and $\mu$ and permeability of the material. 

\subsubsection{Comparison with standard multigrid methods}\label{sec:comparison_gmres}
For our first Maxwell experiment we compare our DPG multigrid technology with the multigrid preconditioner for the standard Galerkin method described in \cite{Jaymultigrid}. Following the experiment in \cite{Jaymultigrid}, we consider the computational domain $\Omega = (0,1)^3$ and a \emph{perfect electric conductor} (PEC) material, i.e, the boundary condition is $n \times E = 0 $ on $\partial \Omega$. For the discretization we use a uniform hexahedral mesh of the lowest order. Define the space $H_0^{-\frac{1}{2}}(\curl, \Gamma_h)$ on the mesh skeleton as 
\begin{equation*}
H_0^{-\frac{1}{2}}(\curl, \Gamma_h) := \{E \in H^{-\frac{1}{2}}(\curl, \Gamma_h): n\times E = 0 \text{ on } \partial \Omega\}
\end{equation*}
where $H^{-\frac{1}{2}}(\curl,\Gamma_h)$ is defined as \cite[Ch. 2]{petrides2019adaptive}.
Then the ultraweak DPG formulation is
\begin{equation*}
\left\{
\begin{alignedat}{3}
& E, H \in (L^2(\Omega))^3, \\
& \hat{E} \in H_0^{-\frac{1}{2}}(\text{curl},\Gamma_h), \,\, \hat{H} \in H^{-\frac{1}{2}}(\text{curl},\Gamma_h) \\
 &i \omega \mu (H,F)+(E,\nabla_h \times F) + \langle n \times \hat{E}, F \rangle_{\Gamma_h} &&= 0, && \quad F \in H(\text{curl}, \Omega_h), \\
 -&i \omega \epsilon (E,G)+(H,\nabla_h \times G) + \langle n \times \hat{H} , G \rangle_{\Gamma_h}  &&= (J,G),&& \quad G \in H(\text{curl}, \Omega_h)
\end{alignedat}
\right.
\end{equation*}
where the letter $h$ denotes element--wise operations. The simulation is driven by the right hand side which is chosen such that the exact solution for the electric field is given by the finite element lift of the function 
$$E_{\text{ex}}(x,y,z) = [y(1-y)z(1-z),yx(1-x)z(1-z),x(1-x)y(1-y)].$$

We choose $\epsilon = \mu = 1$. The initial guess for the CG solver is set to zero and the iterations are terminated when the norm of the residual is reduced by a factor of $10^{-6}$. 
In \cref{tab:DPG_MG_iter1,tab:DPG_MG_iter10} we present the iteration count of the CG solver when preconditioned with our multigrid technology for $\omega = 1$ and $10$. Likewise, \cref{tab:FEM_MG_iter1,tab:FEM_MG_iter10}, retrieved from \cite{Jaymultigrid}, show the iteration count of the GMRES solver preconditioned with the standard multigrid technology. We note that the smoother used in \cite{Jaymultigrid}, is of multiplicative type (block Gauss-Seidel). In the tables, $h$ and $H$ denote the fine and the coarse grid discretization size respectively. For the DPG multigrid implementation, On each level we perform one pre- and one post- smoothing step, so that we keep the symmetry. At the coarsest level the problem is solved exactly with a direct solver. 
\begin{table}[H]\centering
\begin{subtable}[t]{0.33\textwidth}\centering
\begin{tabular}{c??cccc}\toprule\toprule
  h \textbackslash H & 1/2 & 1/4 & 1/8 &1/16\\
 \bottomrule
 \toprule
 1/4  & 7 &    &     &    \\
 1/8  & 7 &  7 &     &    \\
 1/16 & 7 &  7 &  7  &    \\
 1/32 & 6 &  7 &  6  &  7 \\
 1/64 & 6 &  6 &  6  &  6 \\
 \bottomrule
 \bottomrule
\end{tabular}
\caption{CG preconditioned with MG for DPG}
   \label{tab:DPG_MG_iter1}
\end{subtable}
\hspace{50pt}
\begin{subtable}[t]{0.33\textwidth}
   \centering
\begin{tabular}{c??cccc}\toprule\toprule
  h \textbackslash H & 1/2 & 1/4 & 1/8 &1/16\\
 \bottomrule
 \toprule
 1/4  & 6  &      &     &    \\
 1/8  & 7  &  7   &     &    \\
 1/16 & 9  &  10  &  8  &    \\
 1/32 & 10 &  10  &  9  &  7 \\
 1/64 & 11 &  11  &  9  &  8 \\
 \bottomrule
 \bottomrule
\end{tabular}
\caption{GMRES preconditioned with MG for FEM}
   \label{tab:FEM_MG_iter1}
\end{subtable}
\caption{Iteration count for $\omega = 1$. Observe the uniform convergence with respect to $h$ and $H$.}
\end{table}

\begin{table}[H]\centering
\begin{subtable}[t]{0.33\textwidth}\centering
\begin{tabular}{c??cccc}\toprule\toprule
  h \textbackslash H & 1/2 & 1/4 & 1/8 &1/16\\
 \bottomrule
 \toprule
 1/4  & 9 &    &        &       \\
 1/8  & 11 &  10 &      &      \\
 1/16 & 14 &  12 &  9   &      \\
 1/32 & 14 &  14 &  12  &  10 \\
 1/64 & 15 &  15 &  12  &  12 \\
 \bottomrule
 \bottomrule
\end{tabular}
\caption{CG preconditioned with MG for DPG}
   \label{tab:DPG_MG_iter10}
\end{subtable}
\hspace{50pt}
\begin{subtable}[t]{0.33\textwidth}
   \centering
\begin{tabular}{c??cccc}\toprule\toprule
  h \textbackslash H & 1/2 & 1/4 & 1/8 &1/16\\
 \bottomrule
 \toprule
 1/4  & 3$^*$  &      &     &    \\
 1/8  & 2$^*$  &  37   &     &    \\
 1/16 & 3$^*$  &  48  &  18  &    \\
 1/32 & 2$^*$ &  78$^*$  & 22  &  16 \\
 1/64 & 2$^*$ &  78$^*$  & 21  &  17 \\
 \bottomrule
 \bottomrule
\end{tabular}
\caption{GMRES preconditioned with MG for FEM}
   \label{tab:FEM_MG_iter10}
\end{subtable}
\caption{Iteration count for $\omega = 10$. The number of iterations for the DPG method grows mildly with the frequency but always converges to the true solution. Uniform convergence is achieved when a fine enough coarse grid is used. On the contrary the GMRES method fails to deliver reliable solutions when the coarse grid is in the pre--asymptotic region.}
\end{table}
The purpose of this comparison is not to compare number to number the 
iteration count for each multigrid technology, but rather to observe 
the general convergence trend of our preconditioner compared to the 
preconditioner for the standard Galerkin method. 
The entries $n^*$ in the \cref{tab:FEM_MG_iter10} denote that even 
though the GMRES solver did converge, the final iterate differed from 
the true solution by more than $10^{-3}$ (measured in the appropriate norm). 
This is a well known flaw of the GMRES solver, i.e, the residual of the GMRES algorithm might be small enough and the stopping criterion is met, but the output solution is far from the true solution. This happens when the coarse grid is not fine enough for the corresponding frequency and stability is lost. On the contrary, this undesirable convergence behavior is not happening for the Conjugate Gradient solver. The discrete pre--asymptotic stability of the DPG method, along with the theory of self--adjoint preconditioners, ensure that the CG solver always converge to the true solution. However, the convergence does depended on how ``good'' the coarse grid is with respect to the frequency. A similar dependence on the frequency is observed in our 2D simulations \cite{sp}. For both preconditioners, uniform convergence with respect to the frequency is recovered when the coarse grid is fine enough. 

\medskip
\subsubsection{Fichera ``oven'' problem}
\begin{figure}[H]
 \centering
\includegraphics[width=0.5\linewidth]{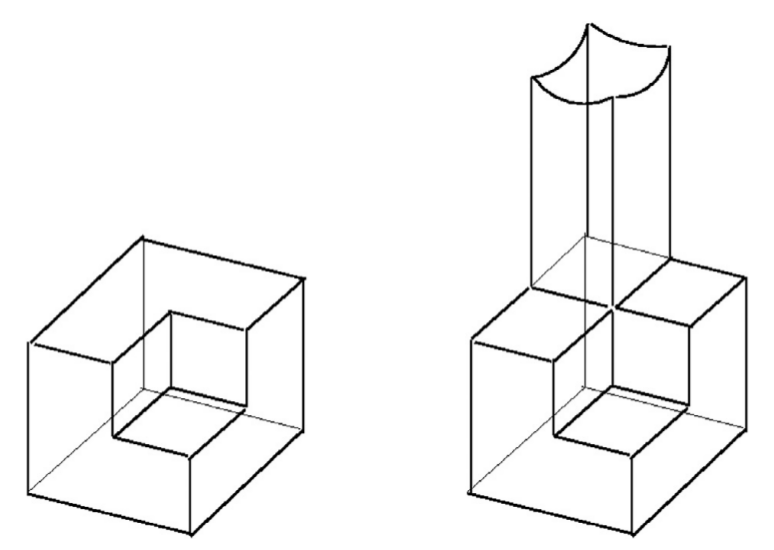}
\caption{Fichera corner with a truncated infinite waveguide attached at the top (retrieved, September 10, 2020 from \cite{Carstensen})}\label{fig:fichera_domain}
\end{figure}

For our second Maxwell experiment we solve the Fichera ``oven'' problem which was first presented in \cite{Carstensen}. The adaptive nature of the DPG method makes it suitable for this problem because the solution is expected to be singular. The set up is as follows. For the construction of the domain we start with the cube $(0,2)^3$ which is uniformly refined into eight cubes and then one is removed creating the Fichera corner. Then, an infinite waveguide is attached at the top and it's truncated at a unit distance from the Fichera corner (see \cref{fig:fichera_domain}). We choose $\epsilon = \mu =1$ and $\omega = 5$. The simulation is driven by a non--homogeneous electric boundary condition on the waveguide and a homogeneous electric boundary condition elsewhere. That is,
\begin{equation*}
n \times E = 
\left\{
\begin{aligned}
& n \times E_d &&\text{ across the waveguide section}\\
& 0  &&\text{ elsewhere}
\end{aligned}
\right.
\end{equation*}
where $E_d = (\sin \pi x_2, 0, 0) $ is the first propagating mode.
%

We start the simulation with a mesh of eight cubes with a uniform order of approximation $p=3$ and perform successive $h$--adaptive refinements. The adaptive refinements are driven by the built--in DPG error indicator (the norm of the error representation function $\psi$), and terminated when the norm of the residual decreases by one order of magnitude. Note that an exact solution for this problem is not known. 
%
%
\input{fichera_figs}

\begin{figure}[H]
\begin{subfigure}[t]{0.49\textwidth}
\centering
  \includegraphics[width=1.0\linewidth]{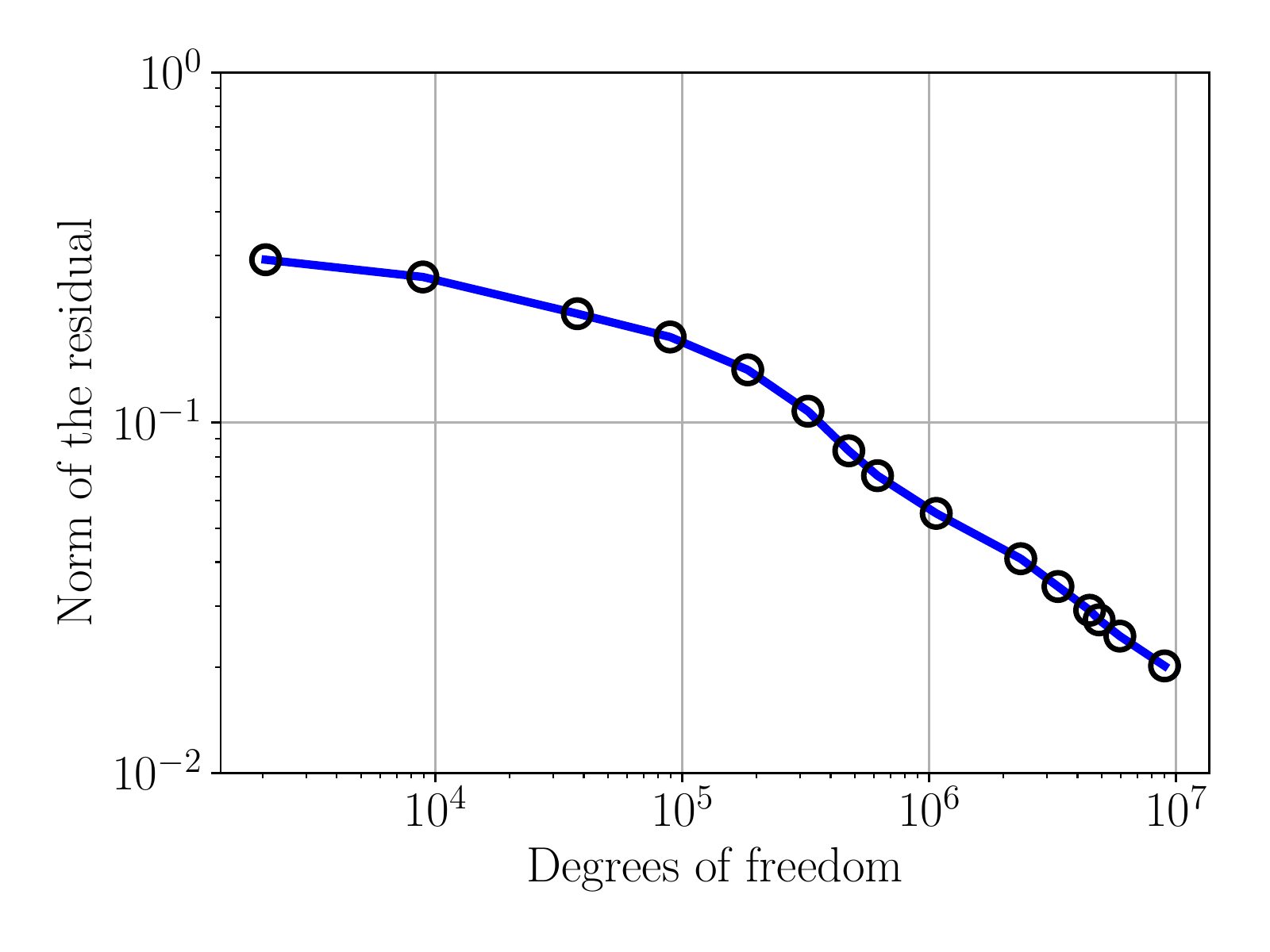}
  \caption{Residual convergence}\label{fig:fichera_err_convergence}
\end{subfigure}
\begin{subfigure}[t]{0.49\textwidth}
\centering
  \includegraphics[width=1.0\linewidth]{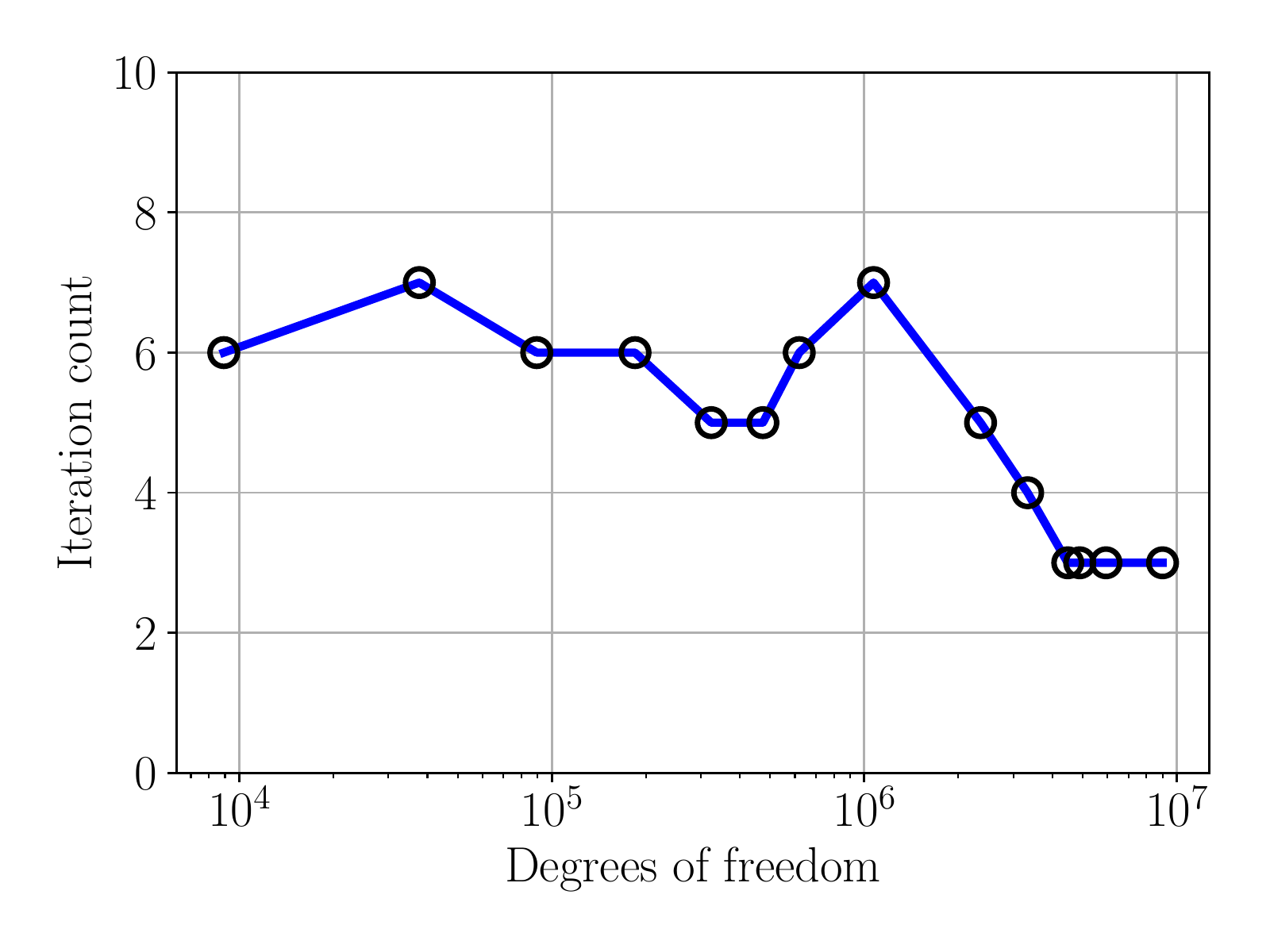}
  \caption{Iteration count for the preconditioned CG solver}\label{fig:fichera_cg_convergence}
\end{subfigure}
\caption{Fichera problem: convergence of the residual (left) and the CG solver(right)}\label{fig:fichera_convergence}
\end{figure}
%
The multigrid preconditioner setting is as follows. Starting from the finest mesh, several coarser adaptive meshes, which belong to the history of refinements, are selected in the multigrid cycle. An exact solve is performed at a coarse grid, where the size of the system is significantly smaller than the fine grid system and small enough to be solve efficiently with a direct solver. The selection is made at run time and depends on the computer architecture and the current memory available.    

In \cref{fig:fichera_prob} we show the evolution of the mesh along with the numerical solution for the real part of the x-component of the electric field. As we can see, the ultraweak DPG method, captures the singularities at the reentrant corner and edges very fast. Observe how the solution changes qualitatively and quantitatively with the resolution of the singularities. The convergence of the residual is displayed in \cref{fig:fichera_err_convergence}. Since there is no known exact solution to this problem, the only way to quantify the convergence is through the DPG residual. Clearly the norm of the residual converges to zero as we proceed with refinements. In \cref{fig:fichera_cg_convergence} we illustrate the convergence behavior of the CG solver preconditioned with multigrid. For this simulation, at most five levels were used in the multigrid cycles. At all intermediate multigrid levels, we perform 10 smoothing steps with a relaxation parameter selected according to the finite overlap property (here $\theta = 0.2$). The CG iterations are initiated with a zero initial guess and they are terminated when the norm of the discrete residual is reduced below $10^{-6}$. As we can observe, the solver shows uniform convergence with respect to the discretization size throughout the adaptive refinement process. We emphasize that the problem was too large to be solved with a direct solver. 

\medskip
\subsubsection{Scattering of a Gaussian beam from a cube}
Our last Maxwell example involves the simulation of a high--frequency Gaussian beam scattering from a rigid cube. We consider the domain $(0,1)^3$ where a cube of side length $a = 1/7$ centered at $(0.5,0.5,0.5)$ is removed ($\Omega = (0,1)^3\backslash (\frac{3}{7},\frac{4}{7})^3$). The simulation is driven by an impedance boundary condition on the outer cube and homogeneous electric boundary condition on the inner cube. The impedance data around the origin correspond to a high--frequency Gaussian beam propagating inside the domain at a specific angle. Away from the source the impedance data smoothly decay to zero in order to simulate absorbing boundary conditions.
\noindent
The DPG ultraweak formulation for this setting is given by:
\begin{equation*}
\left\{
\begin{alignedat}{3}
& E, H \in (L^2(\Omega))^3, \\
&\hat{E} \in H_{\Gamma_2}^{-\frac{1}{2}}(\text{curl},\Gamma_h), \,\, \hat{H} \in H^{-\frac{1}{2}}(\text{curl},\Gamma_h) \\
&n \times (n \times \hat{E}) - n \times \hat{H} = \hat{g}, \,\, \text{ on } \Gamma_1,  \\
 &i \omega \mu (H,F)+(E,\nabla_h \times F) + \langle n \times \hat{E} , F \rangle_{\Gamma_h} &&= 0, \quad F \in H(\text{curl}, \Omega_h), \\
 -&i \omega \epsilon (E,G)+(H,\nabla_h \times G) + \langle n \times \hat{H} , G \rangle_{\Gamma_h} && = 0, \quad G \in H(\text{curl}, \Omega_h)
\end{alignedat}
\right.
\end{equation*} 
where, $\Gamma_1$, $\Gamma_2$ are the boundaries of the outer and the inner cube respectively and 
\begin{equation*}
H_{\Gamma_2}^{-\frac{1}{2}}(\text{curl},\Gamma_h) := \{\hat{E} 
\in H^{-\frac{1}{2}}(\text{curl},\Gamma_h): n \times \hat{E}=0 \text{ on } \Gamma_2\}.
\end{equation*}

We start the simulation with a uniform mesh of $342$ cubes of order $p=3$ and we initiate adaptive $hp$--refinements. We follow an ad--hoc refinement strategy; an element marked for refinement is h-refined unless its size is less than half a wavelength, in which case it's p-refined (two elements per wavelength). In order to resolve the anticipated singularities on the scatterer, the elements adjacent to its corners and edges are marked and forced to be h-refined when needed. 

Finally, when an adaptive $p$--refinement is performed, we follow the following rule (minimum rule). When an element is marked for a $p$--refinement we first find the neighboring elements with respect to its faces. The order of a neighboring element is then increased by one if it is less than the intended order of the marked element. The last step is to assign orders to edges and faces. The order of a face is defined to be the minimum of the orders of its neighboring elements. After all faces are assigned their order, the edge orders take the value of the minimum of the orders of the faces they belong to. With this rather complicated rule, we ensure that a $p$--unrefinement is not possible (the sequence of meshes remains nested), and the constrained approximation technology for handling hanging nodes \cite{Demkbook2} is well defined.  We mention that a maximum rule was also tested successfully. However, in order to maintain stability the enriched order of the test space had to be increased, and this often led to significantly increased element computation times.

\begin{figure}[H]
\begin{subfigure}[b]{0.48\textwidth}
      \centering
      \includegraphics[width=1.0\linewidth]{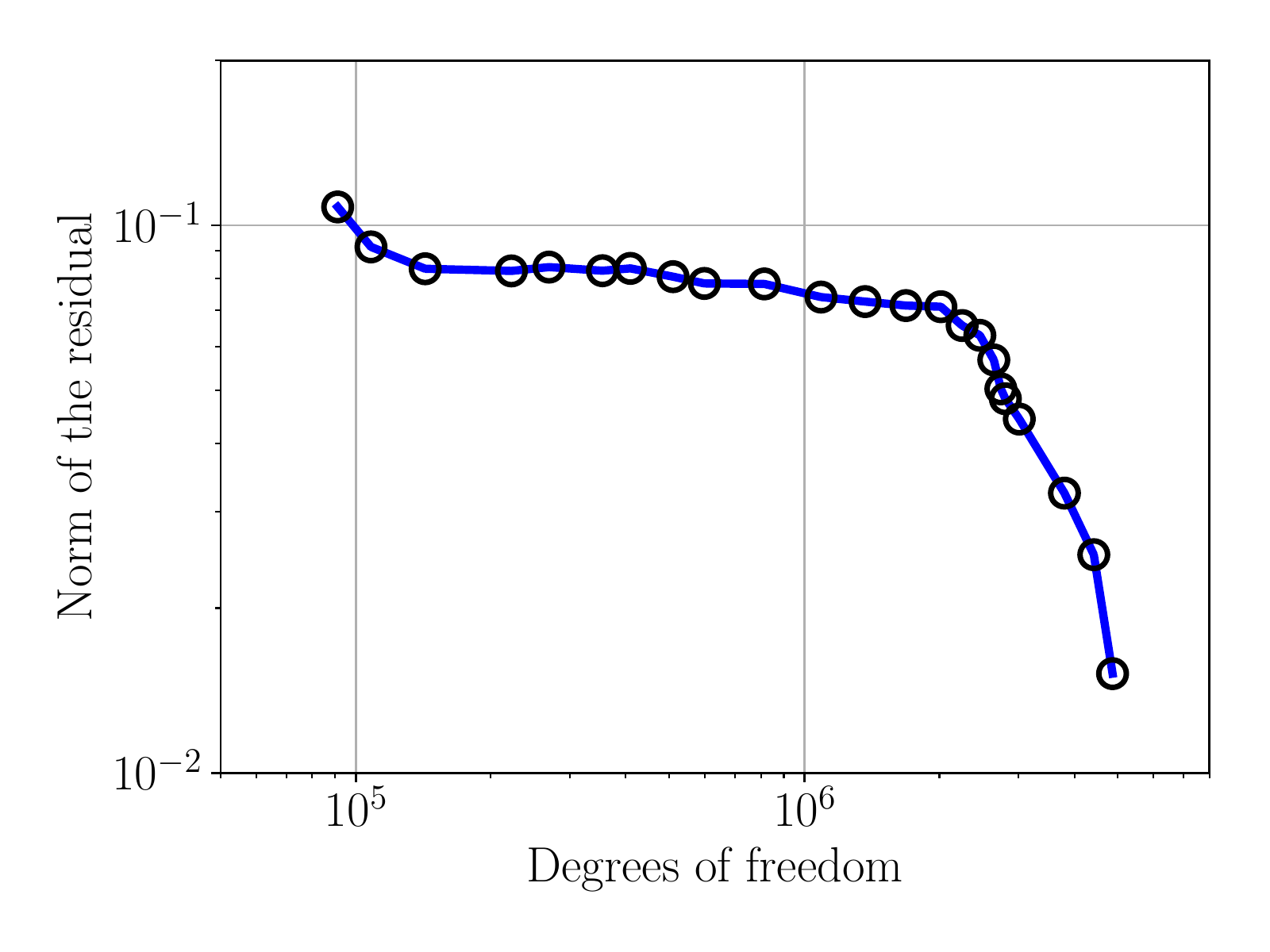}
      \caption{Residual convergence}\label{fig:cube_res_conv}
\end{subfigure}
\begin{subfigure}[b]{0.48\textwidth}
      \centering
      \includegraphics[width=1.0\linewidth]{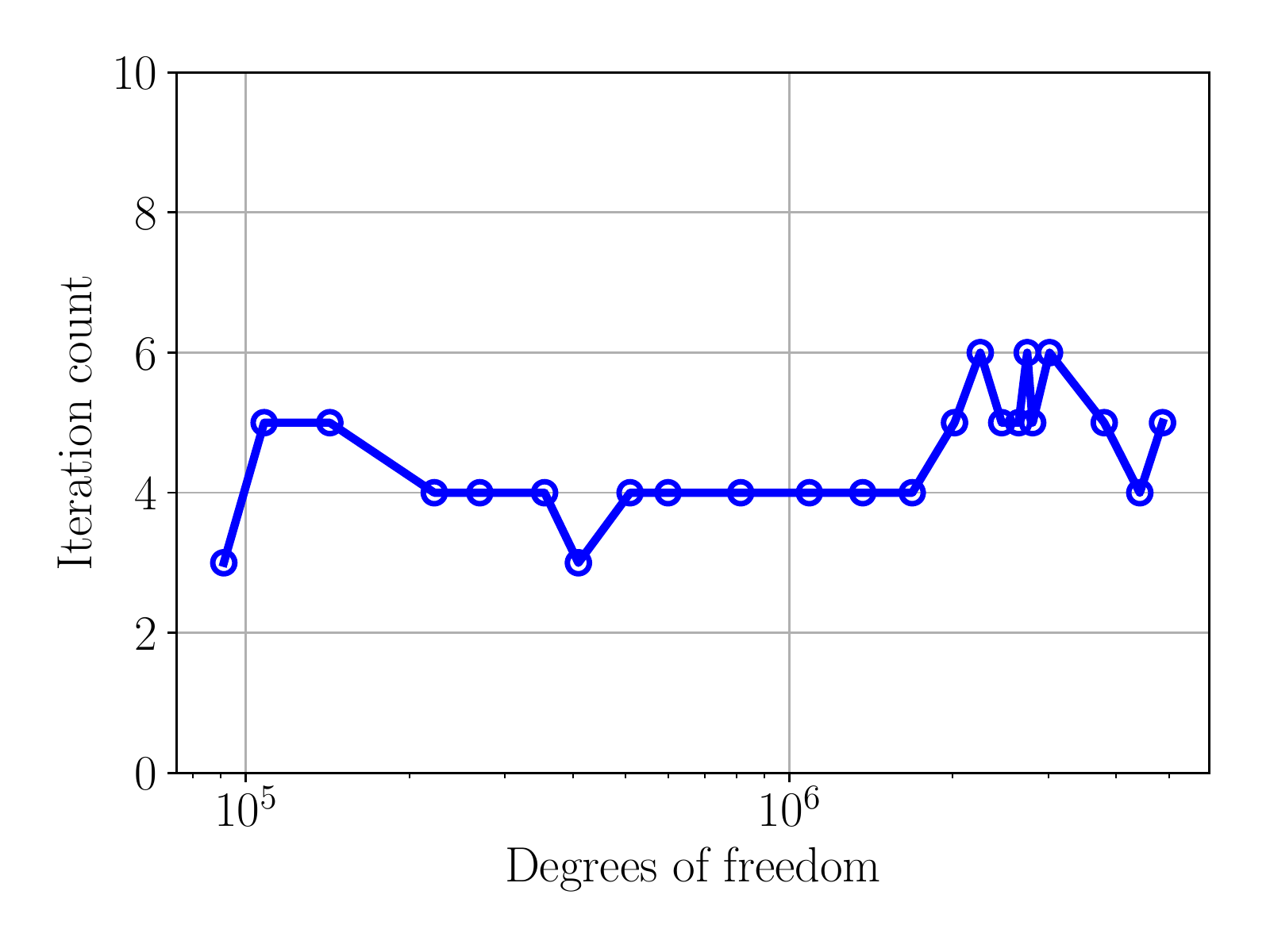}
      \caption{Preconditioned CG iteration count}\label{fig:cube_cg_iter}
\end{subfigure}
    \caption{Convergence of the residual and the preconditioned CG solver}
\end{figure}

We run the simulation for $\epsilon = \mu = 1$ and $\omega = 50\pi$. This frequency corresponds to approximately $40$ wavelengths inside the computational domain. The multigrid setting is the same an in the previous example. For this experiment, at most eight multigrid levels were used. Adaptivity is guided by the norm of the residual which is shown to be driven to zero in \cref{fig:cube_res_conv}.  The iteration count of the CG solver preconditioned with our multigrid technology is presented in \cref{fig:cube_cg_iter}. Note that the number of iterations of the CG solver remains under control throughout the adaptive refinement process. Lastly, the evolution of the mesh and the numerical solution of the real part of $E_x$ are shown in \cref{fig:max_cube2,fig:max_efield} respectively. Note that in these figures we show only the part of the domain below the plane defined by the point $(0.5,0.5,0.5)$ and the normal vector $(-0.5,-0.5,1)$.
%
%
%
\input{beam_20_figs}

%% file: fichera_figs.tex
\begin{figure}[H]
\hspace{-5pt}
\begin{subfigure}[b]{0.22\textwidth}
\includegraphics[trim=250 0 250 0,clip,width=1.0\linewidth]{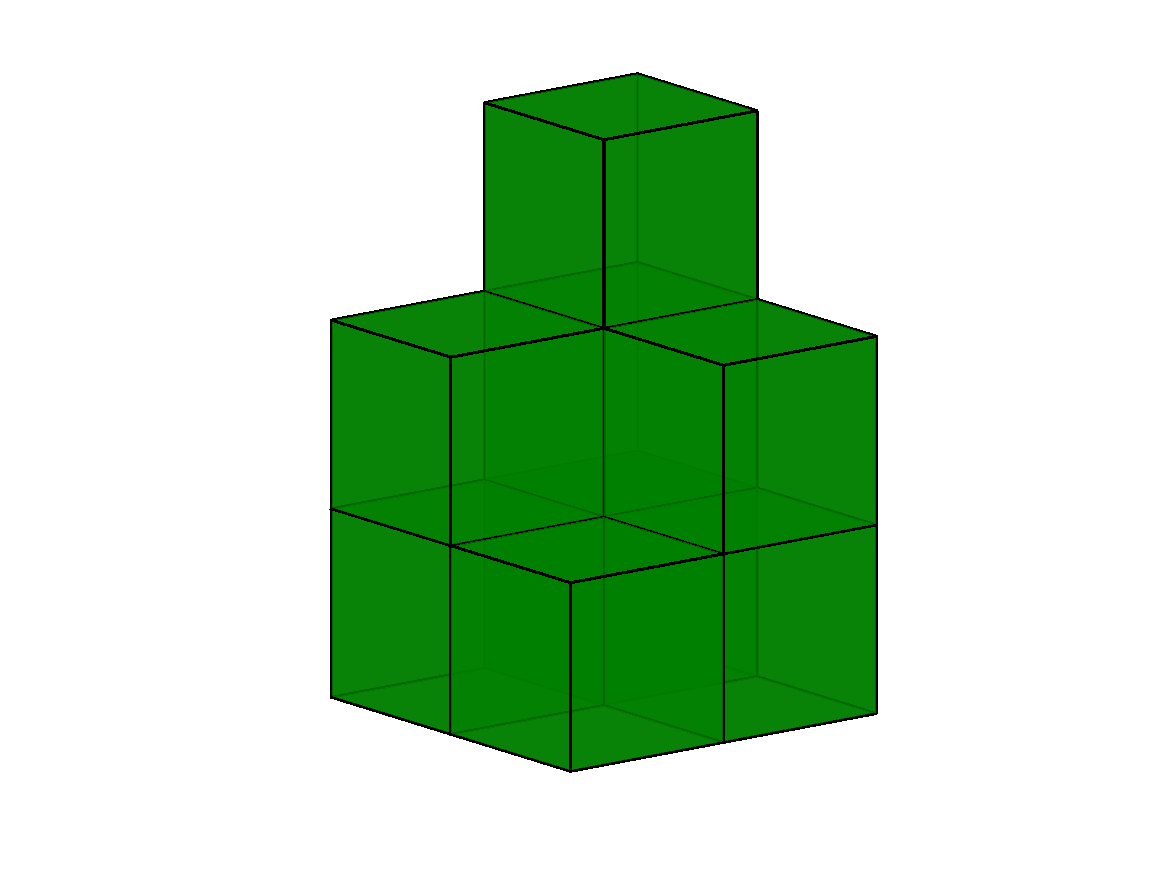}
\end{subfigure}
\begin{subfigure}[b]{0.22\textwidth}
\includegraphics[trim=210 0 210 0,clip,width=1.0\linewidth]{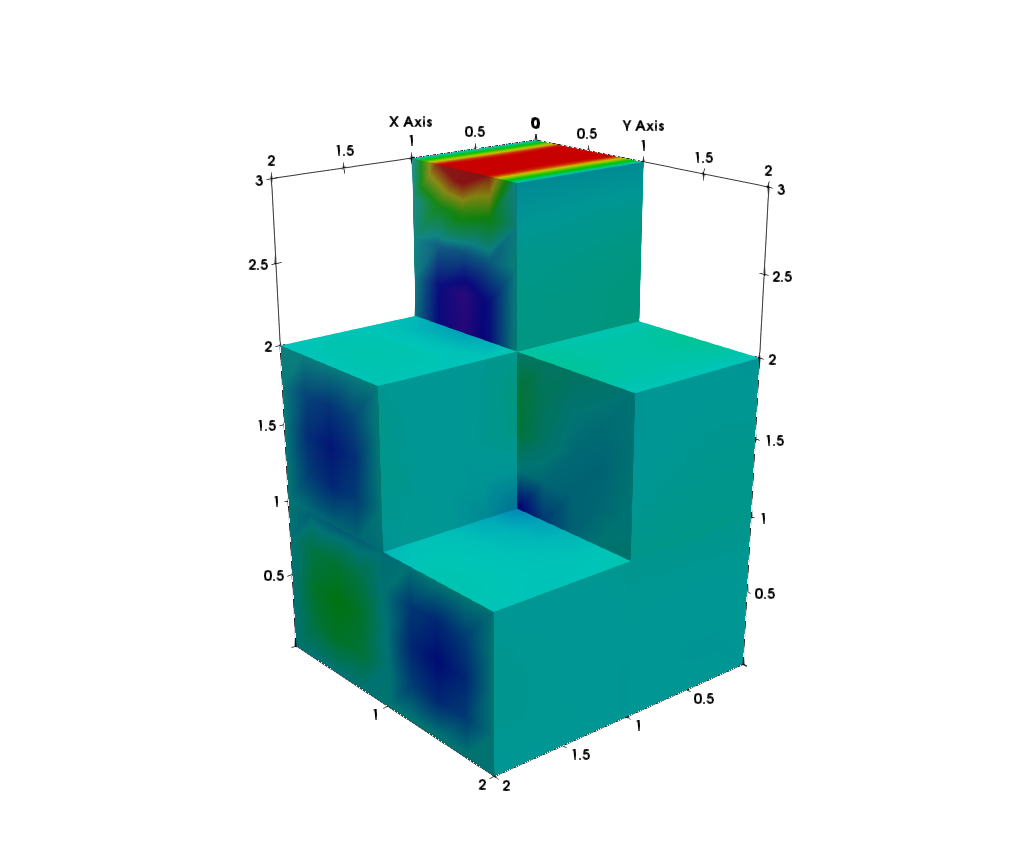}
\end{subfigure}
\begin{subfigure}[b]{0.22\textwidth}
\includegraphics[trim=210 0 210 0,clip,width=1.0\linewidth]{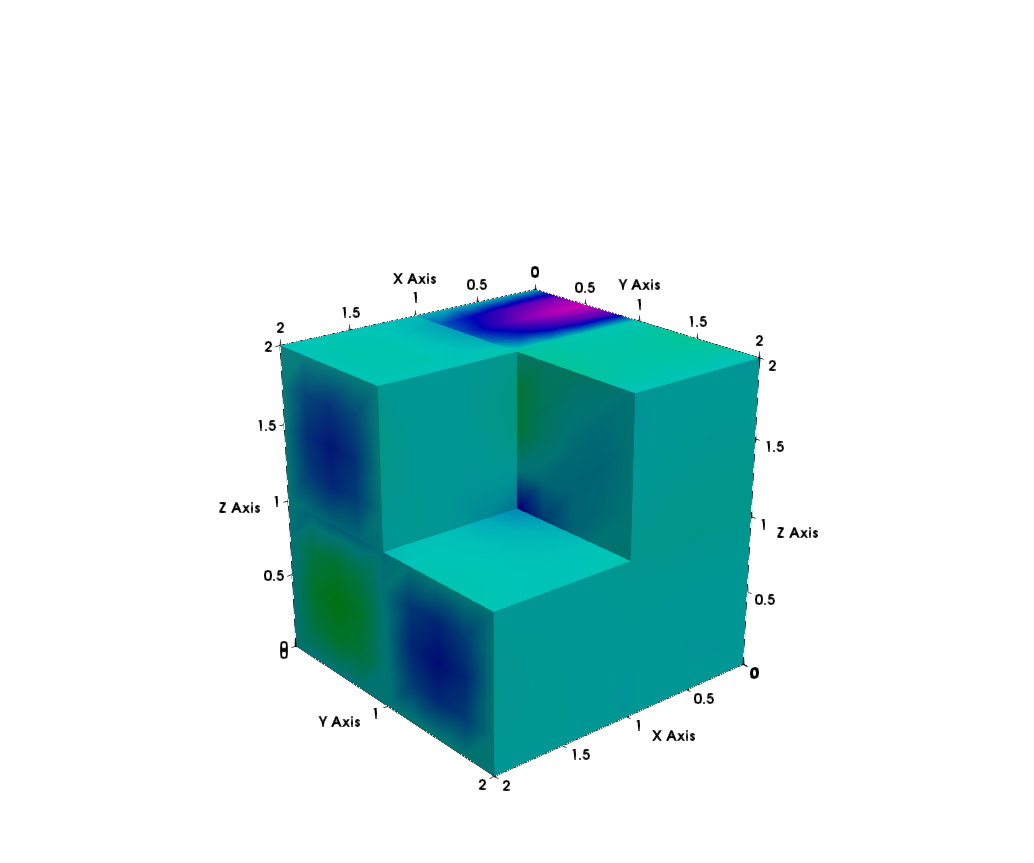}
\end{subfigure}
\begin{subfigure}[b]{0.22\textwidth}
\includegraphics[trim=210 0 210 0,clip,width=1.0\linewidth]{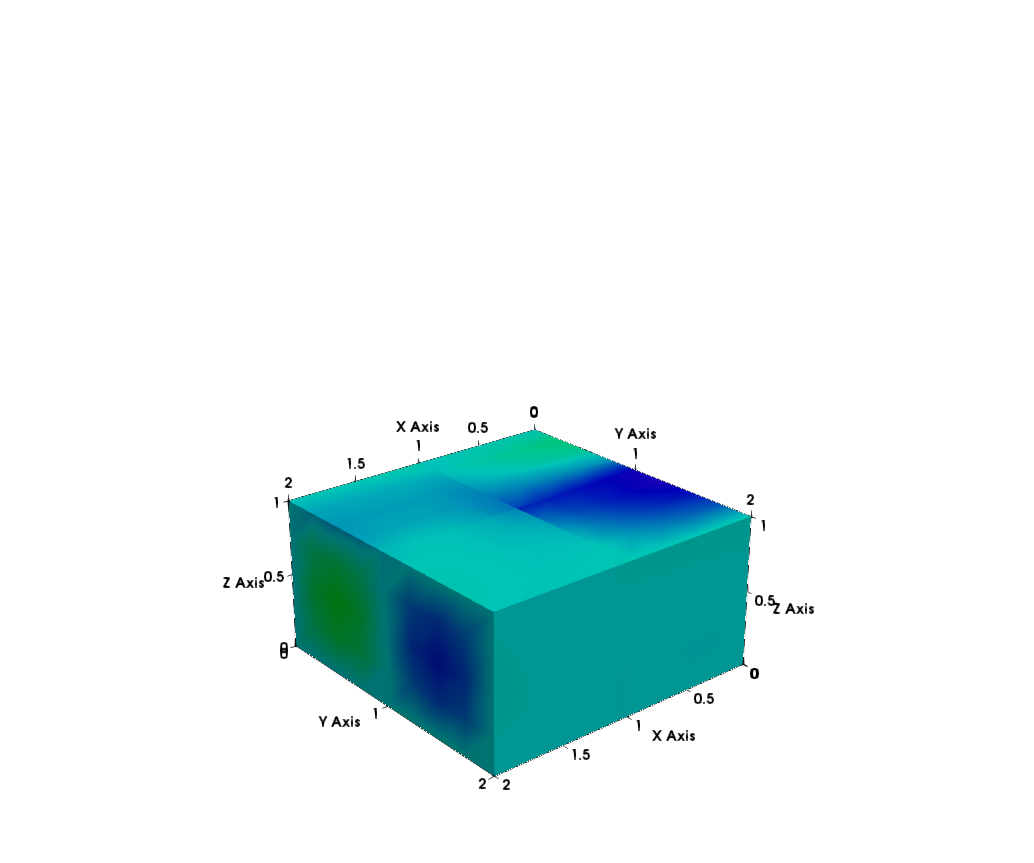}
\end{subfigure}
\begin{subfigure}[b]{0.10\textwidth}
\includegraphics[width=1.0\linewidth]{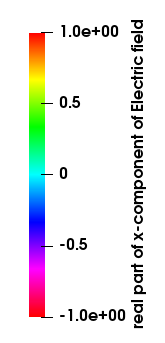}
\end{subfigure}
\end{figure}
\vspace{-25pt}
\begin{figure}[H]
\hspace{-5pt}
\begin{subfigure}[b]{0.22\textwidth}
\includegraphics[trim=250 0 250 0,clip,width=1.0\linewidth]{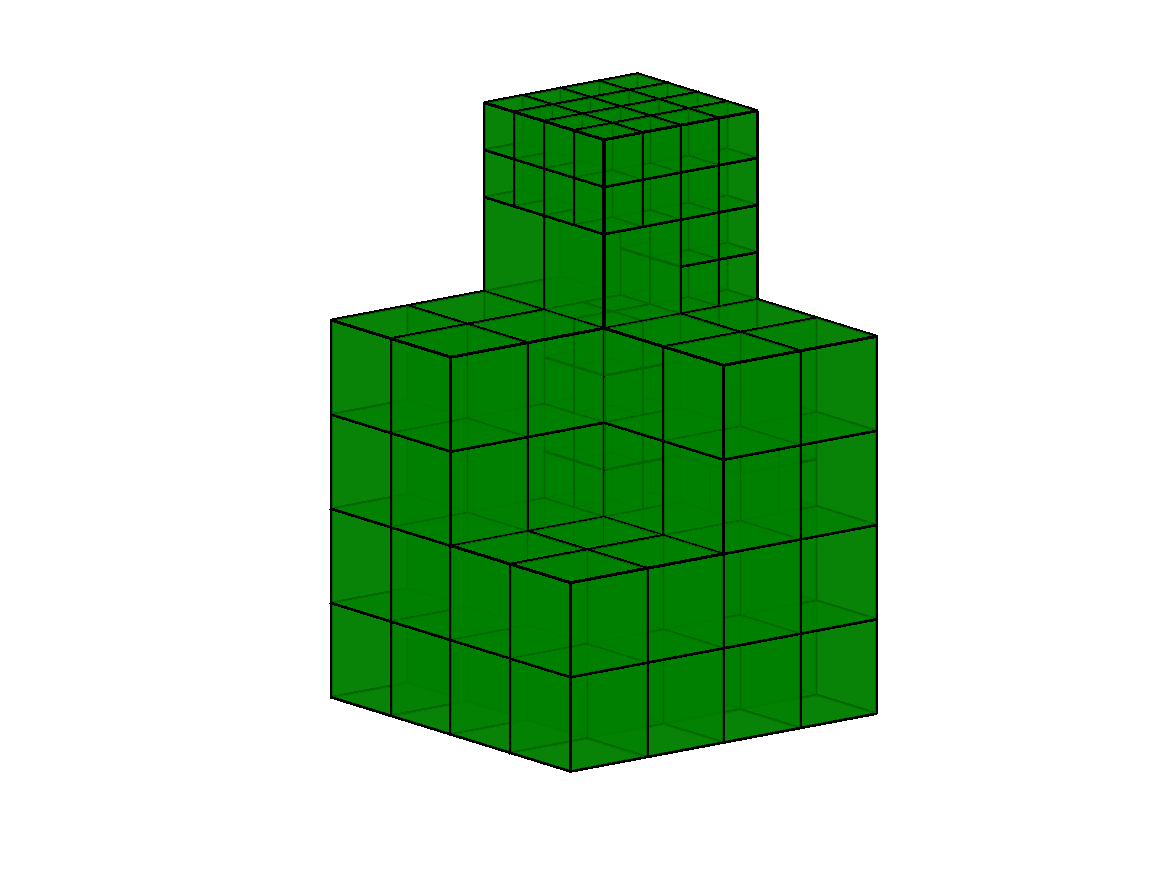}
\end{subfigure}
\begin{subfigure}[b]{0.22\textwidth}
\includegraphics[trim=210 0 210 0,clip,width=1.0\linewidth]{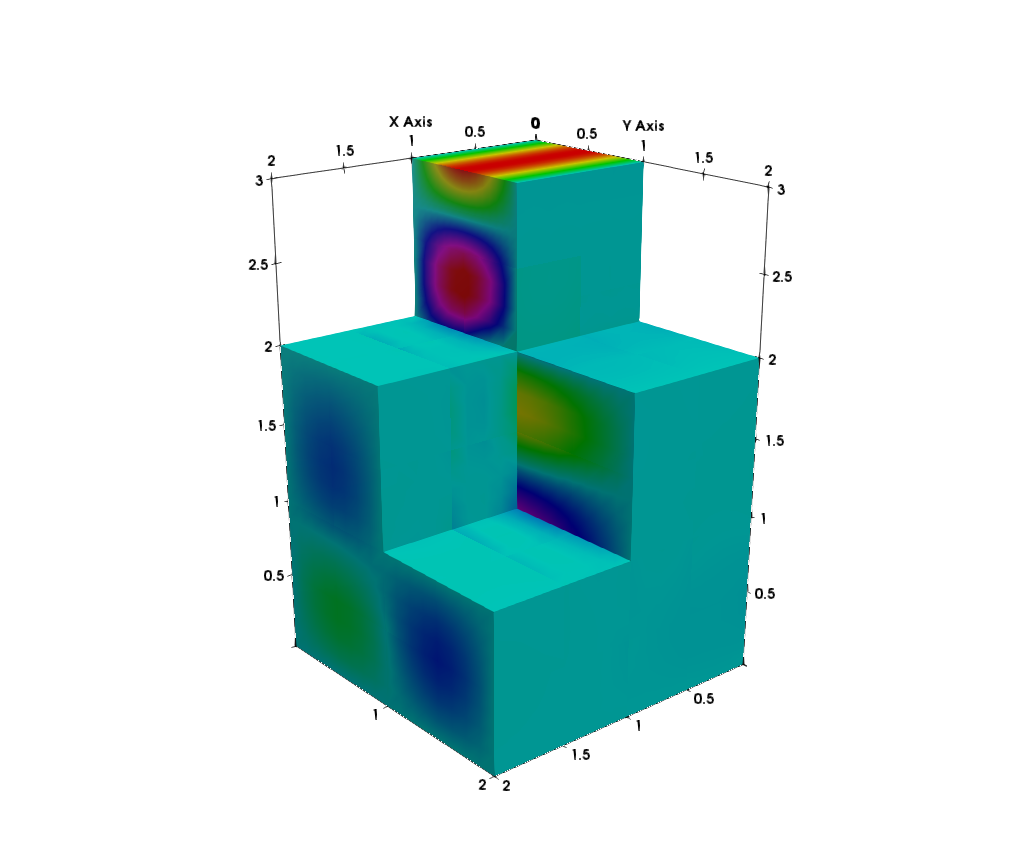}
\end{subfigure}
\begin{subfigure}[b]{0.22\textwidth}
\includegraphics[trim=210 0 210 0,clip,width=1.0\linewidth]{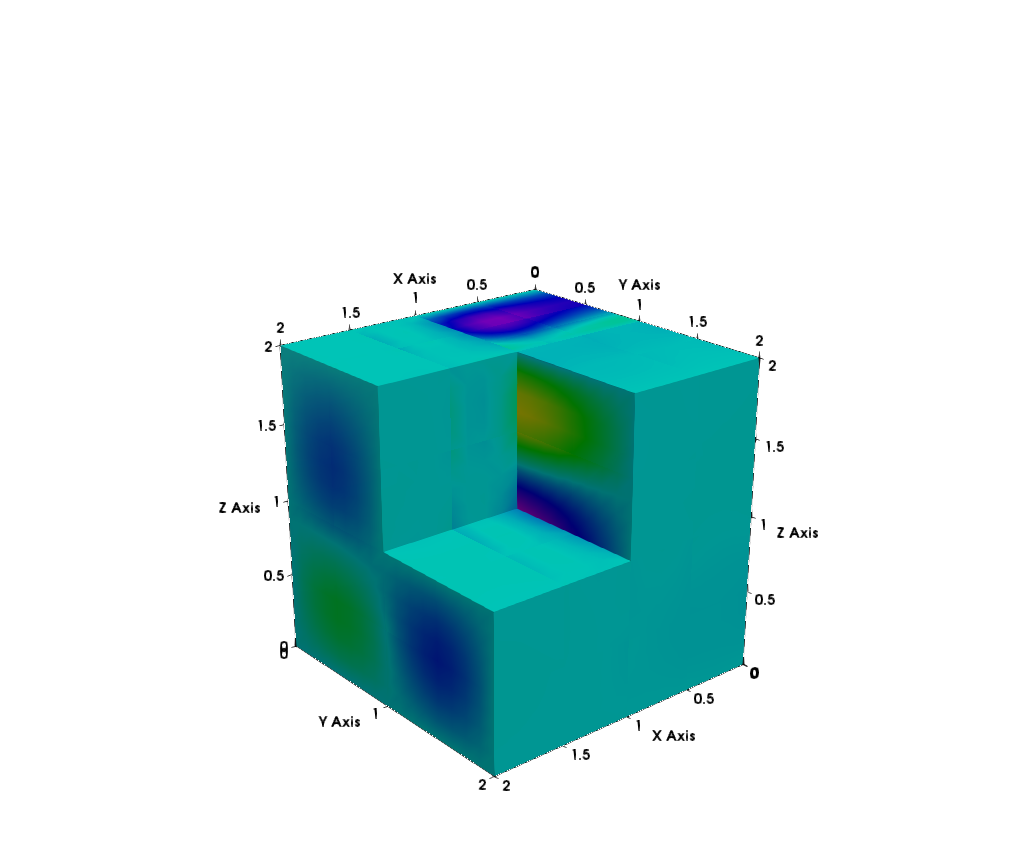}
\end{subfigure}
\begin{subfigure}[b]{0.22\textwidth}
\includegraphics[trim=210 0 210 0,clip,width=1.0\linewidth]{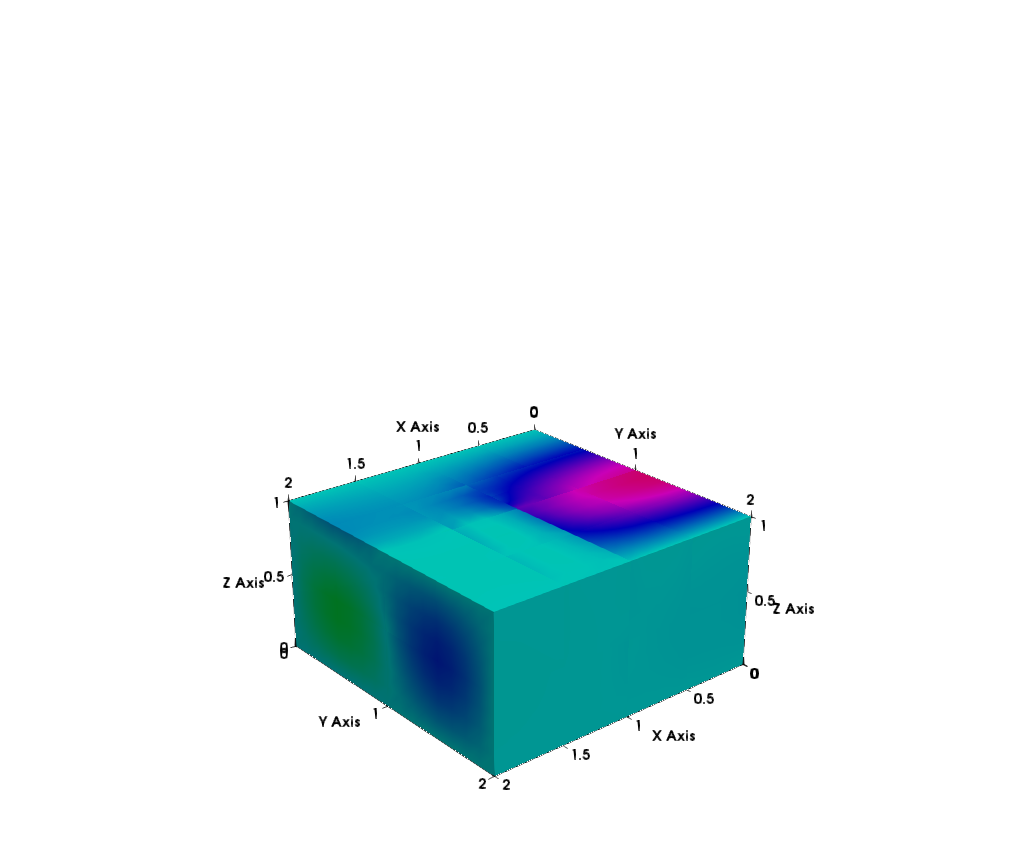}
\end{subfigure}
\begin{subfigure}[b]{0.10\textwidth}
\includegraphics[width=1.0\linewidth]{figures/results/em/fichera/cmap.png}
\end{subfigure}
\end{figure}
\vspace{-25pt}
\begin{figure}[H]
\hspace{-5pt}
\begin{subfigure}[b]{0.22\textwidth}
\includegraphics[trim=250 0 250 0,clip,width=1.0\linewidth]{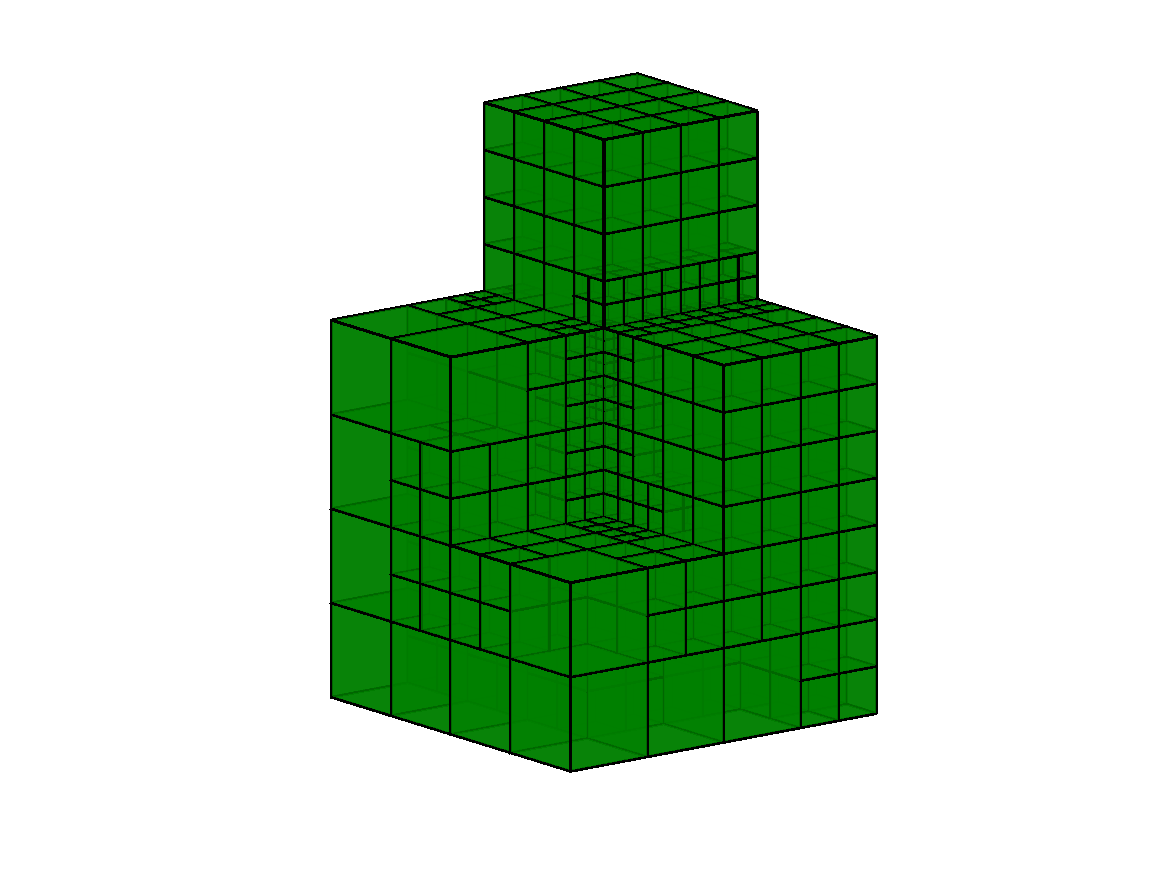}
\end{subfigure}
\begin{subfigure}[b]{0.22\textwidth}
\includegraphics[trim=210 0 210 0,clip,width=1.0\linewidth]{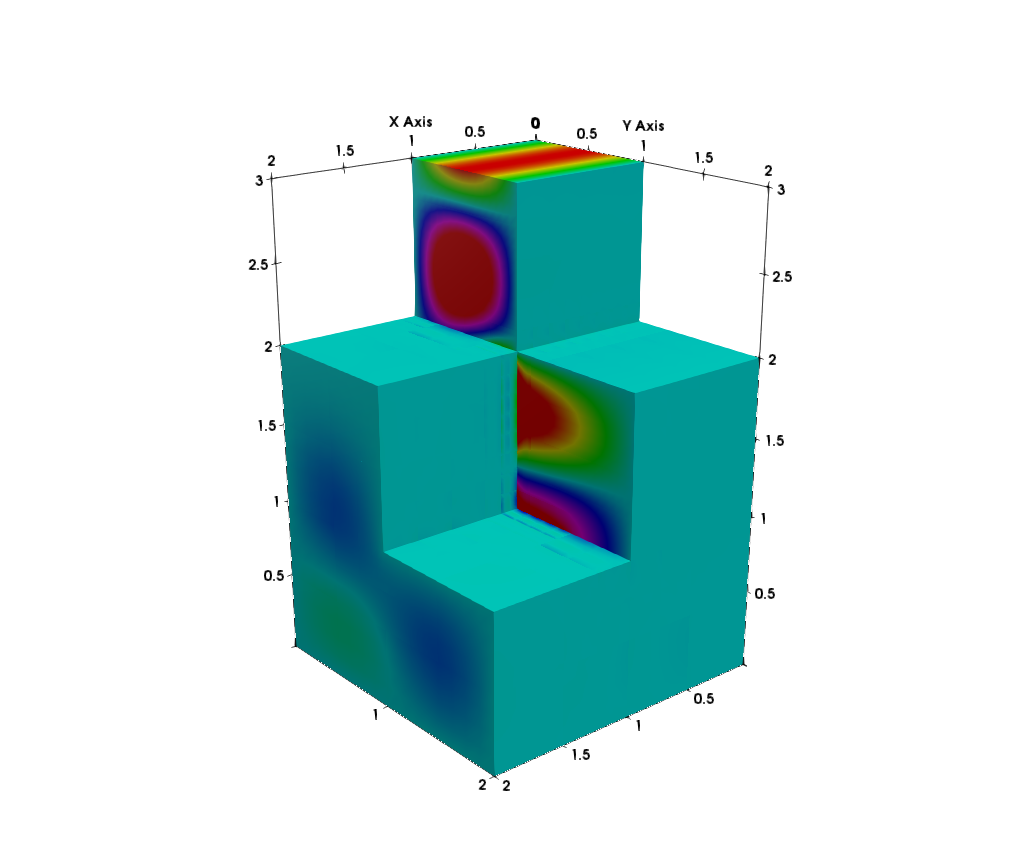}
\end{subfigure}
\begin{subfigure}[b]{0.22\textwidth}
\includegraphics[trim=210 0 210 0,clip,width=1.0\linewidth]{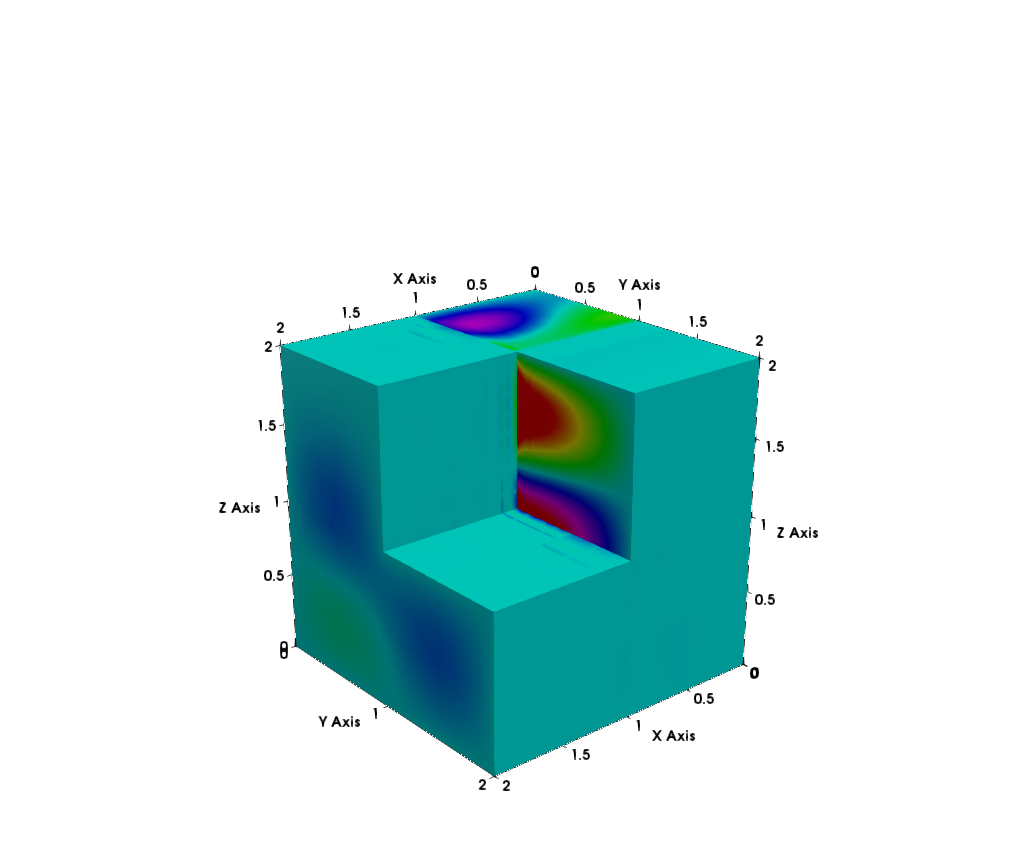}
\end{subfigure}
\begin{subfigure}[b]{0.22\textwidth}
\includegraphics[trim=210 0 210 0,clip,width=1.0\linewidth]{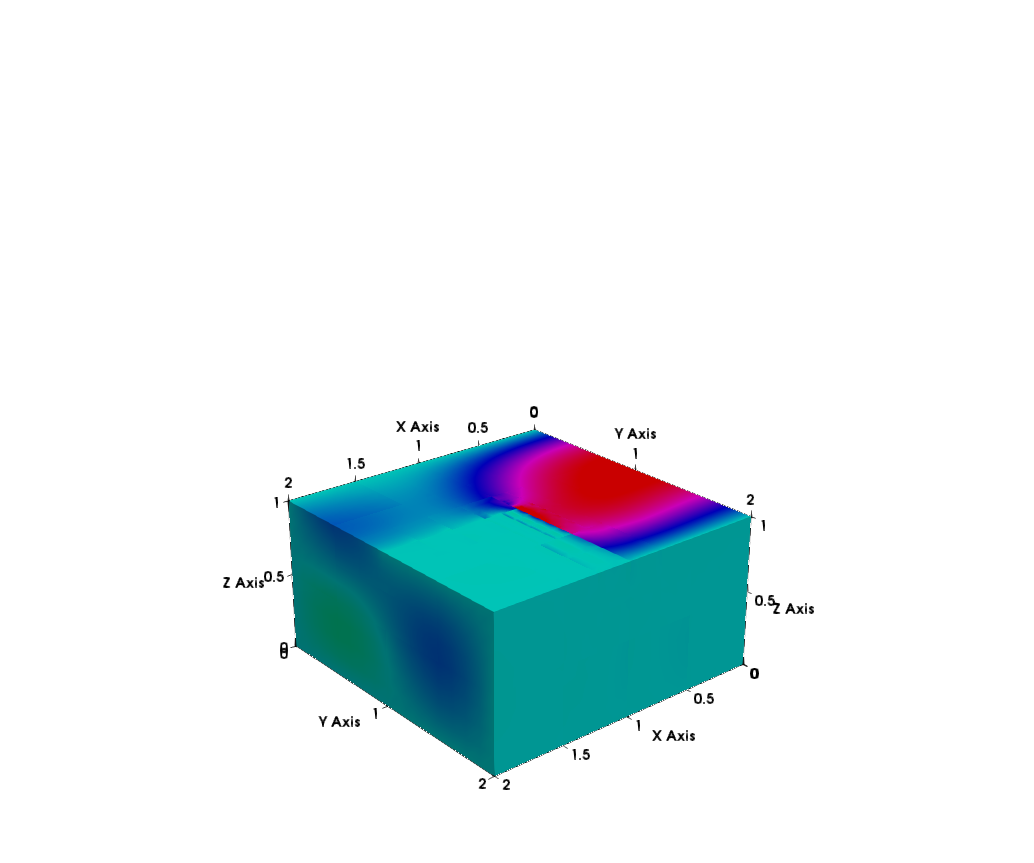}
\end{subfigure}
\begin{subfigure}[b]{0.10\textwidth}
\includegraphics[width=1.0\linewidth]{figures/results/em/fichera/cmap.png}
\end{subfigure}
\end{figure}
\vspace{-25pt}
\begin{figure}[H]
\hspace{-5pt}
\begin{subfigure}[b]{0.22\textwidth}
\includegraphics[trim=250 0 250 0,clip,width=1.0\linewidth]{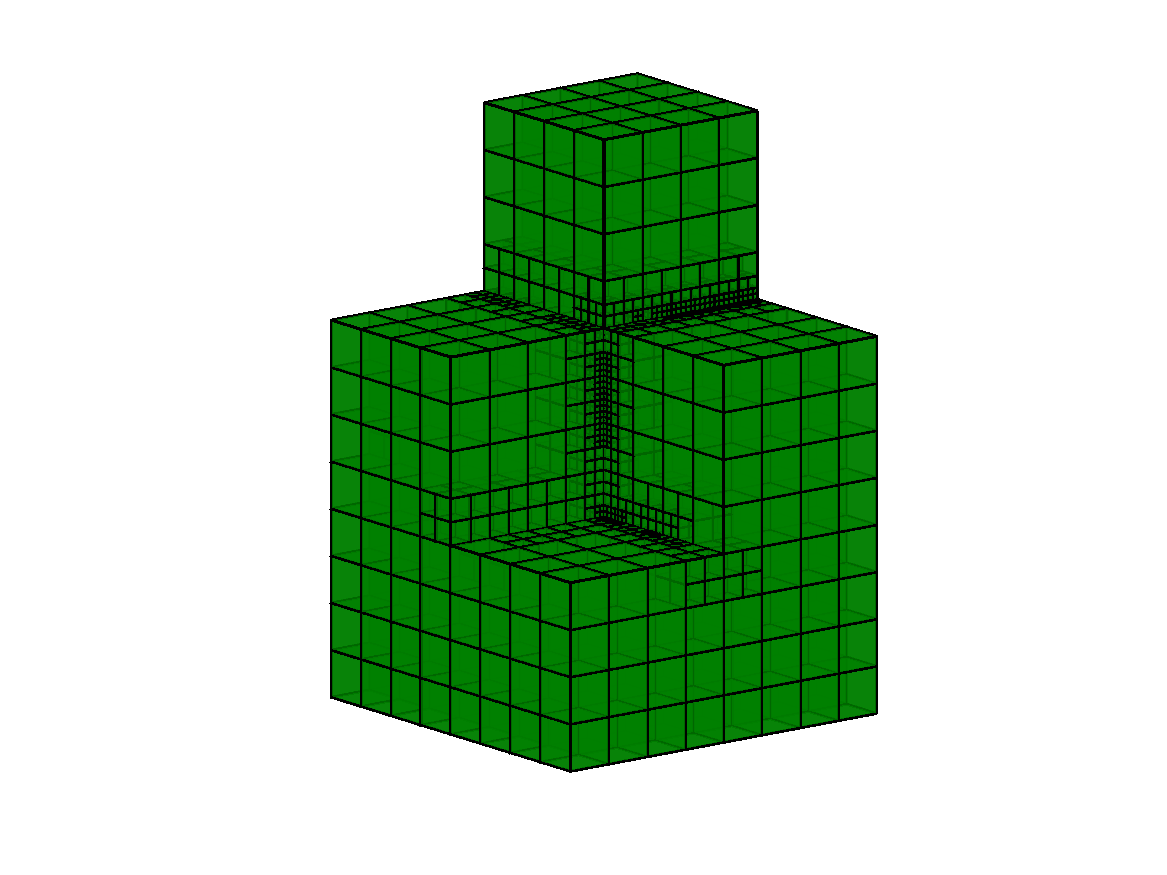}
\end{subfigure}
\begin{subfigure}[b]{0.22\textwidth}
\includegraphics[trim=210 0 210 0,clip,width=1.0\linewidth]{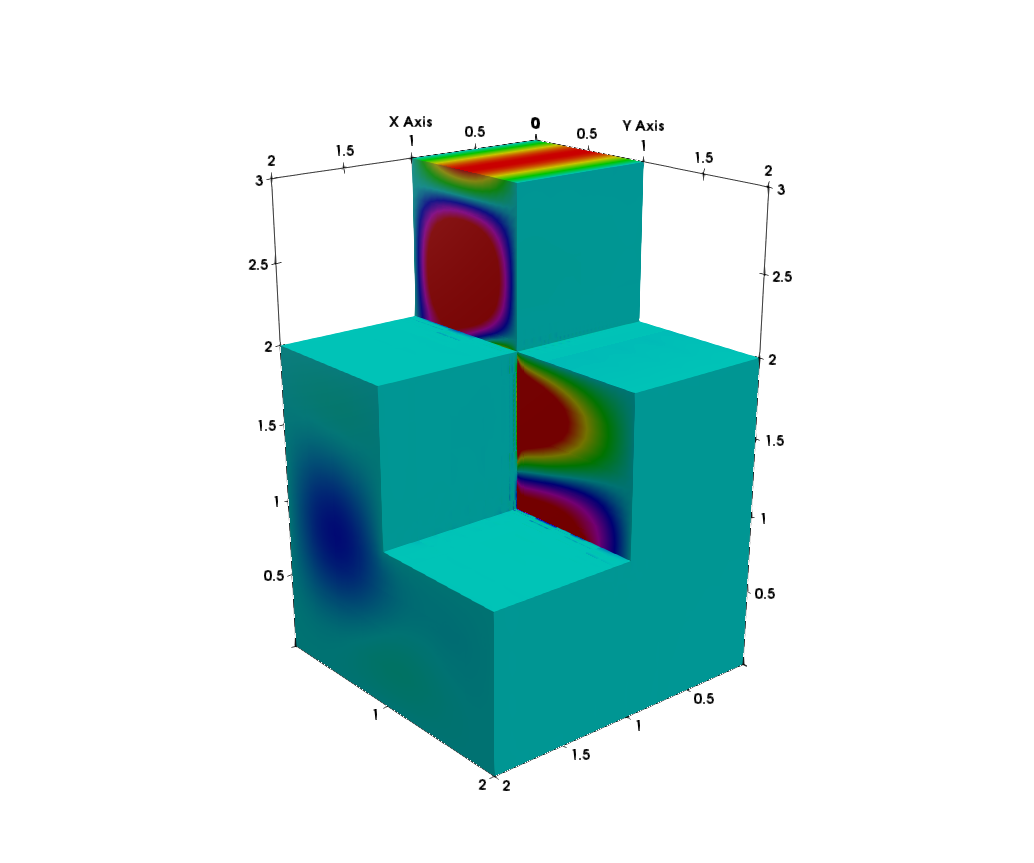}
\end{subfigure}
\begin{subfigure}[b]{0.22\textwidth}
\includegraphics[trim=210 0 210 0,clip,width=1.0\linewidth]{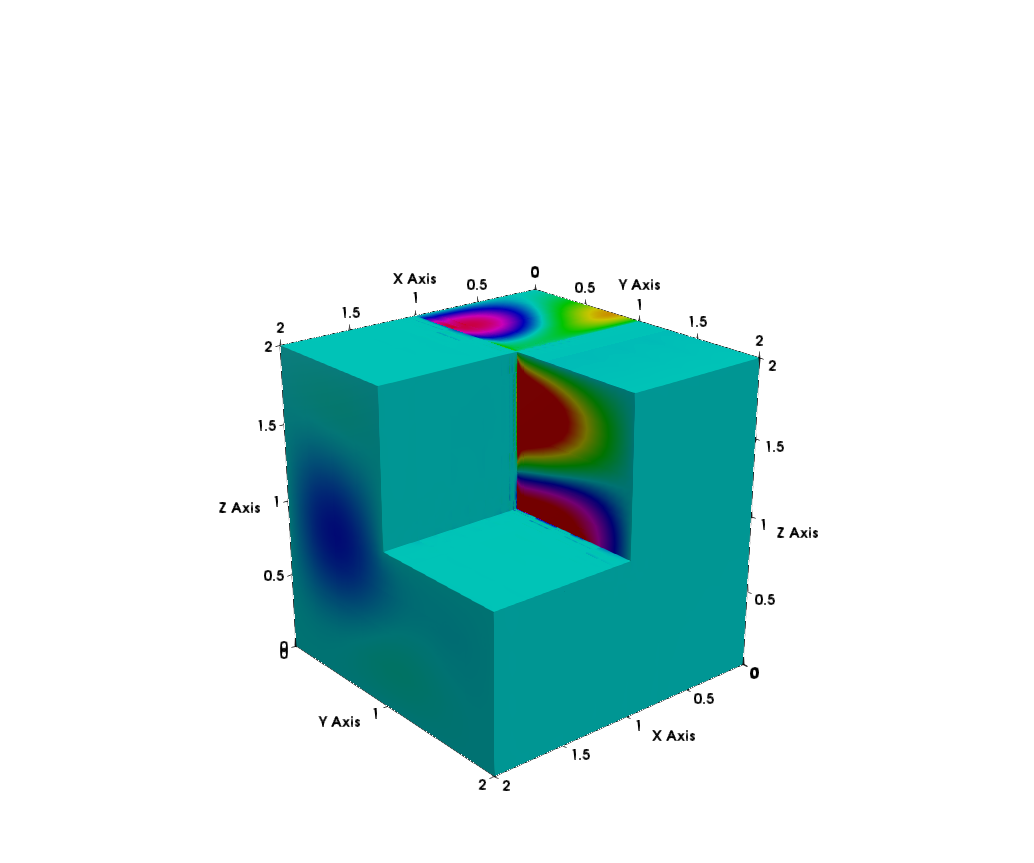}
\end{subfigure}
\begin{subfigure}[b]{0.22\textwidth}
\includegraphics[trim=210 0 210 0,clip,width=1.0\linewidth]{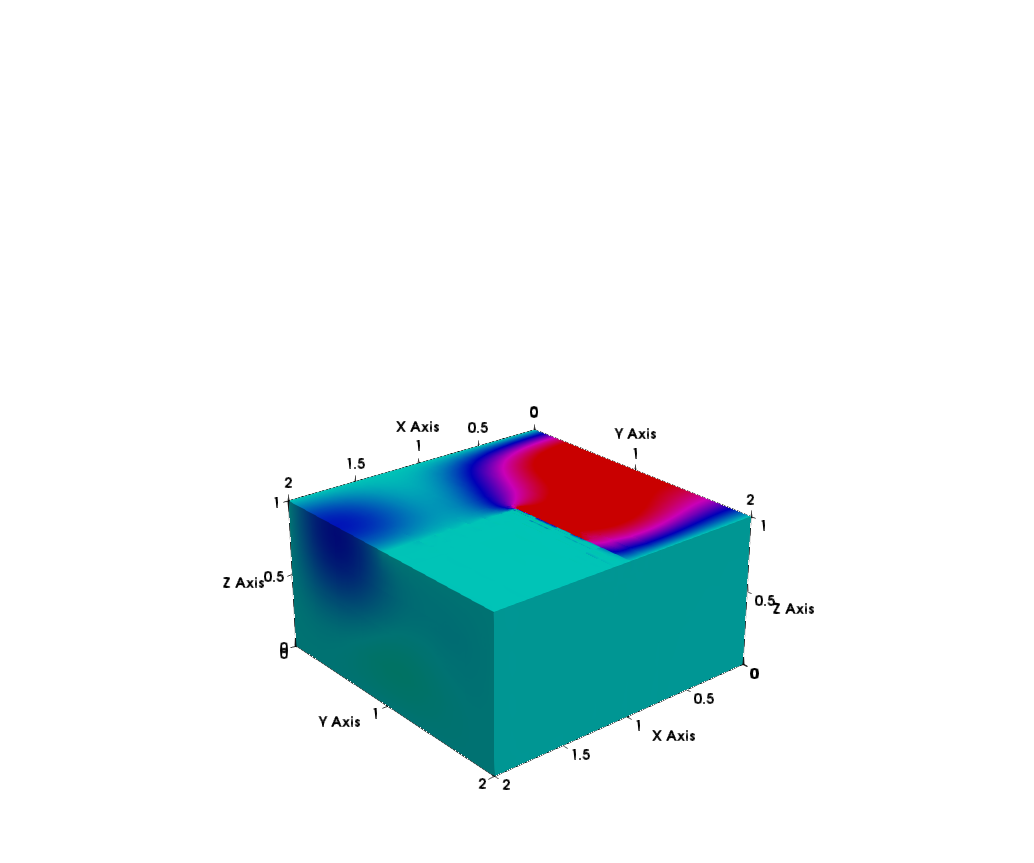}
\end{subfigure}
\begin{subfigure}[b]{0.10\textwidth}
\includegraphics[width=1.0\linewidth]{figures/results/em/fichera/cmap.png}
\end{subfigure}


\hspace{-5pt}
\begin{subfigure}[b]{0.22\textwidth}
\includegraphics[trim=250 0 250 0,clip,width=1.0\linewidth]{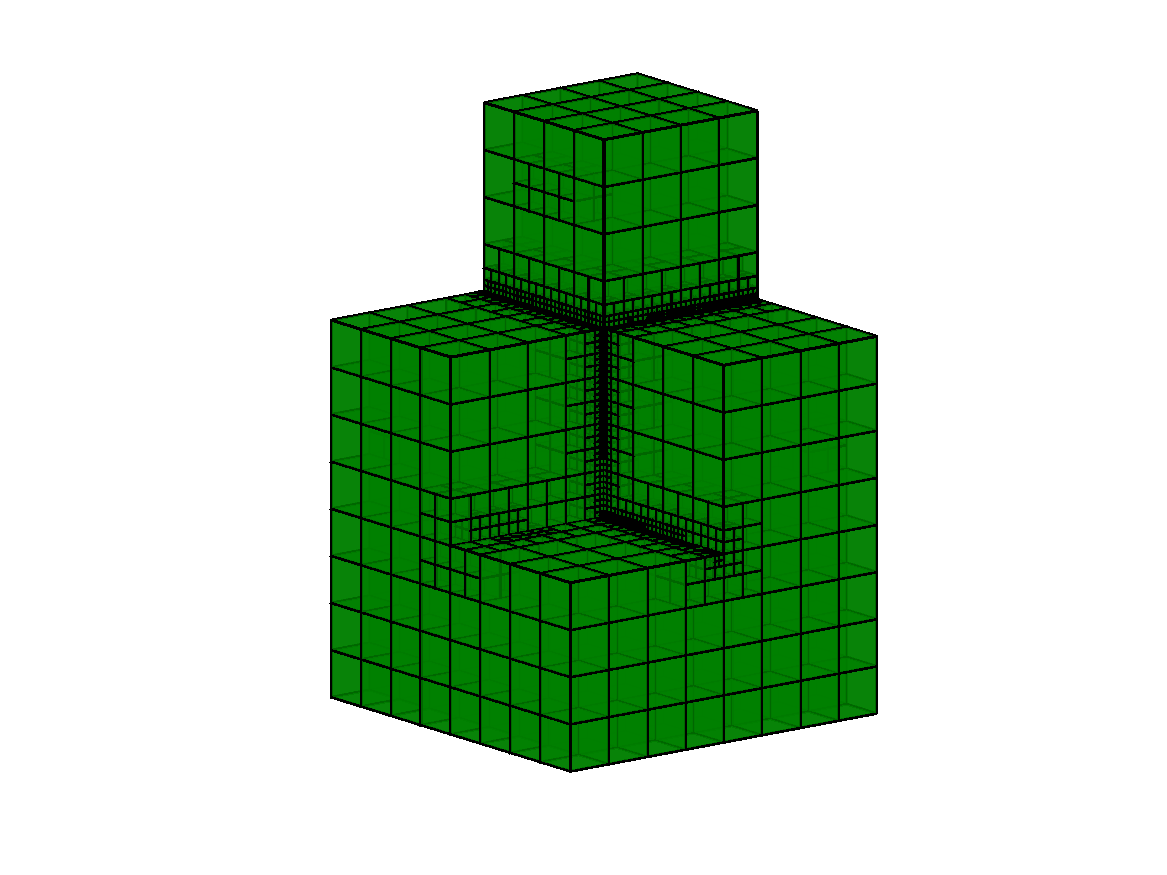}
\end{subfigure}
\begin{subfigure}[b]{0.22\textwidth}
\includegraphics[trim=210 0 210 0,clip,width=1.0\linewidth]{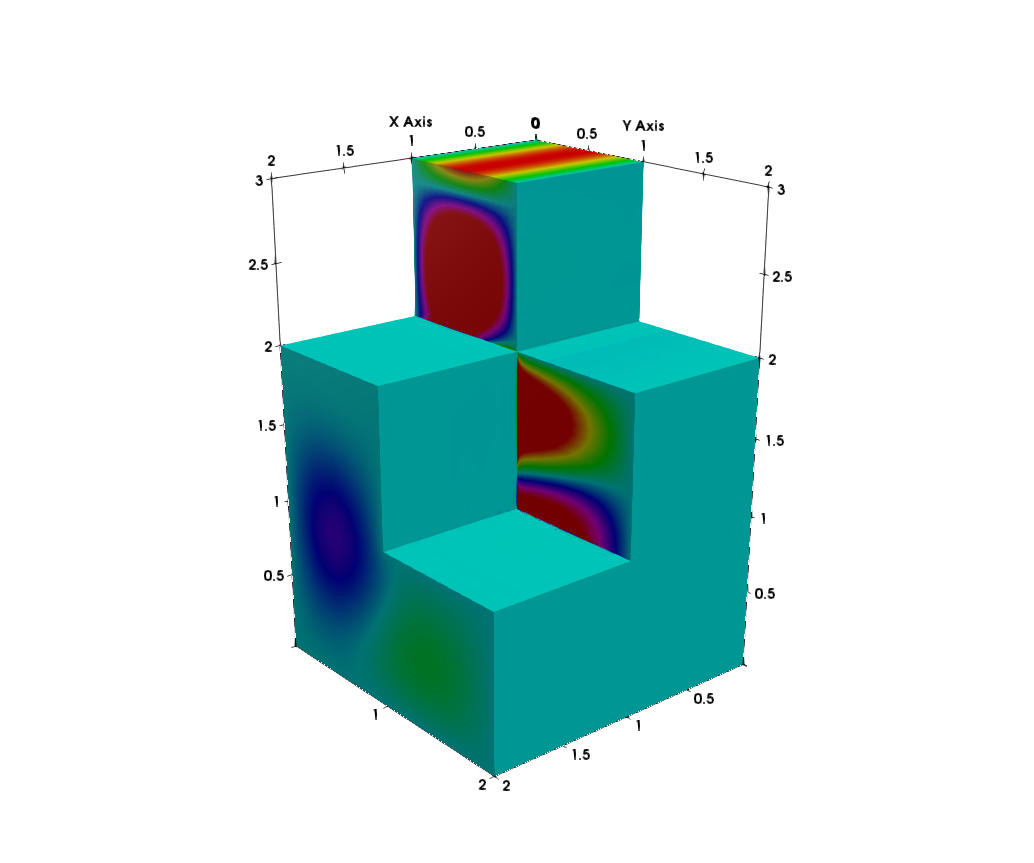}
\end{subfigure}
\begin{subfigure}[b]{0.22\textwidth}
\includegraphics[trim=210 0 210 0,clip,width=1.0\linewidth]{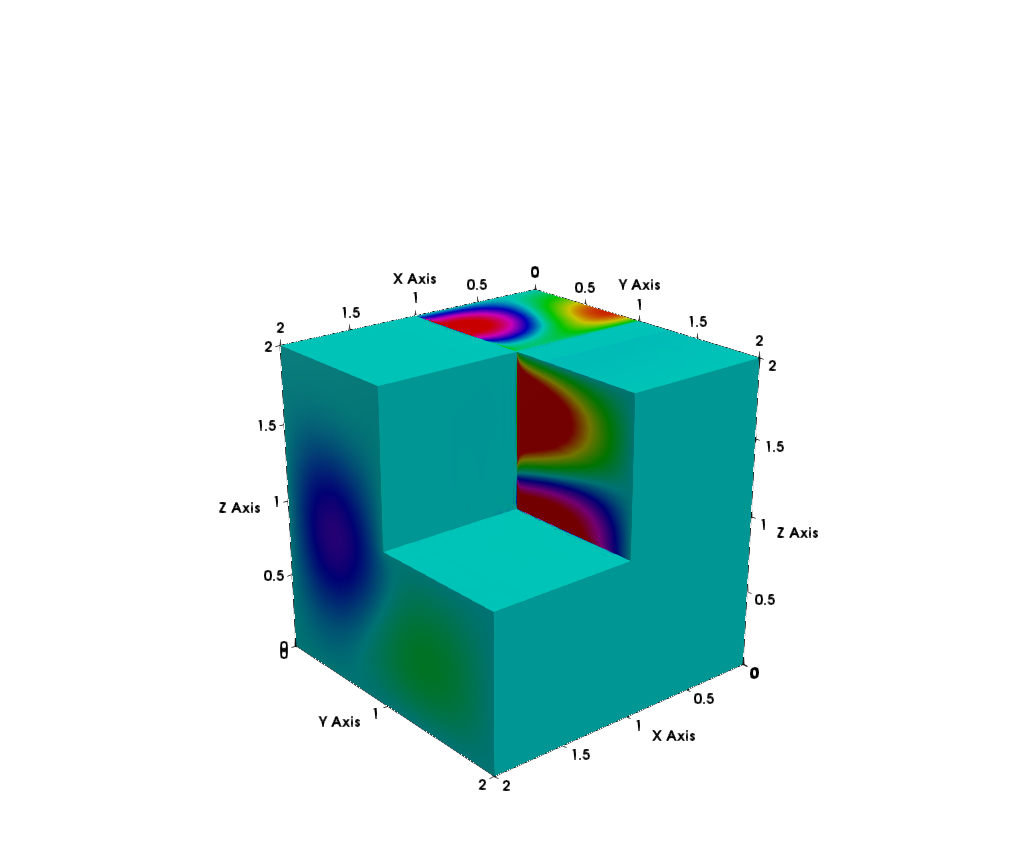}
\end{subfigure}
\begin{subfigure}[b]{0.22\textwidth}
\includegraphics[trim=210 0 210 0,clip,width=1.0\linewidth]{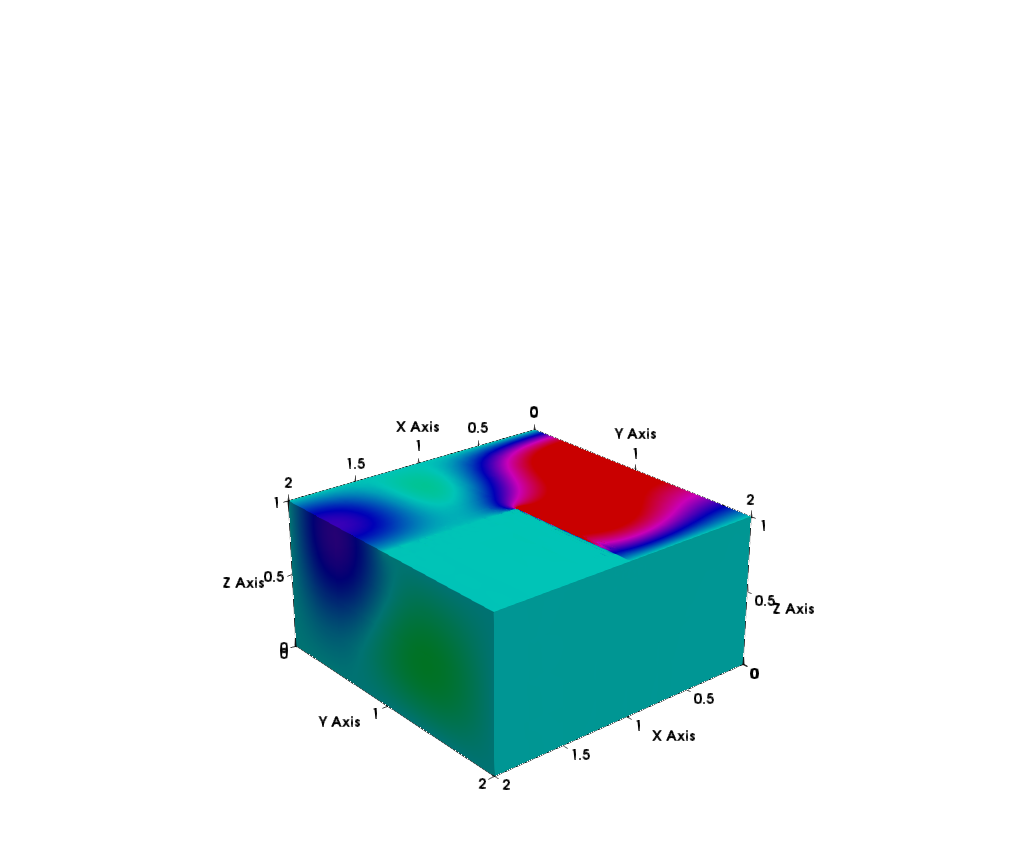}
\end{subfigure}
\begin{subfigure}[b]{0.10\textwidth}
\includegraphics[width=1.0\linewidth]{figures/results/em/fichera/cmap.png}
\end{subfigure}
\end{figure}
\vspace{-25pt}
\begin{figure}[H]
\hspace{-5pt}
\begin{subfigure}[b]{0.22\textwidth}
\includegraphics[trim=250 0 250 0,clip,width=1.0\linewidth]{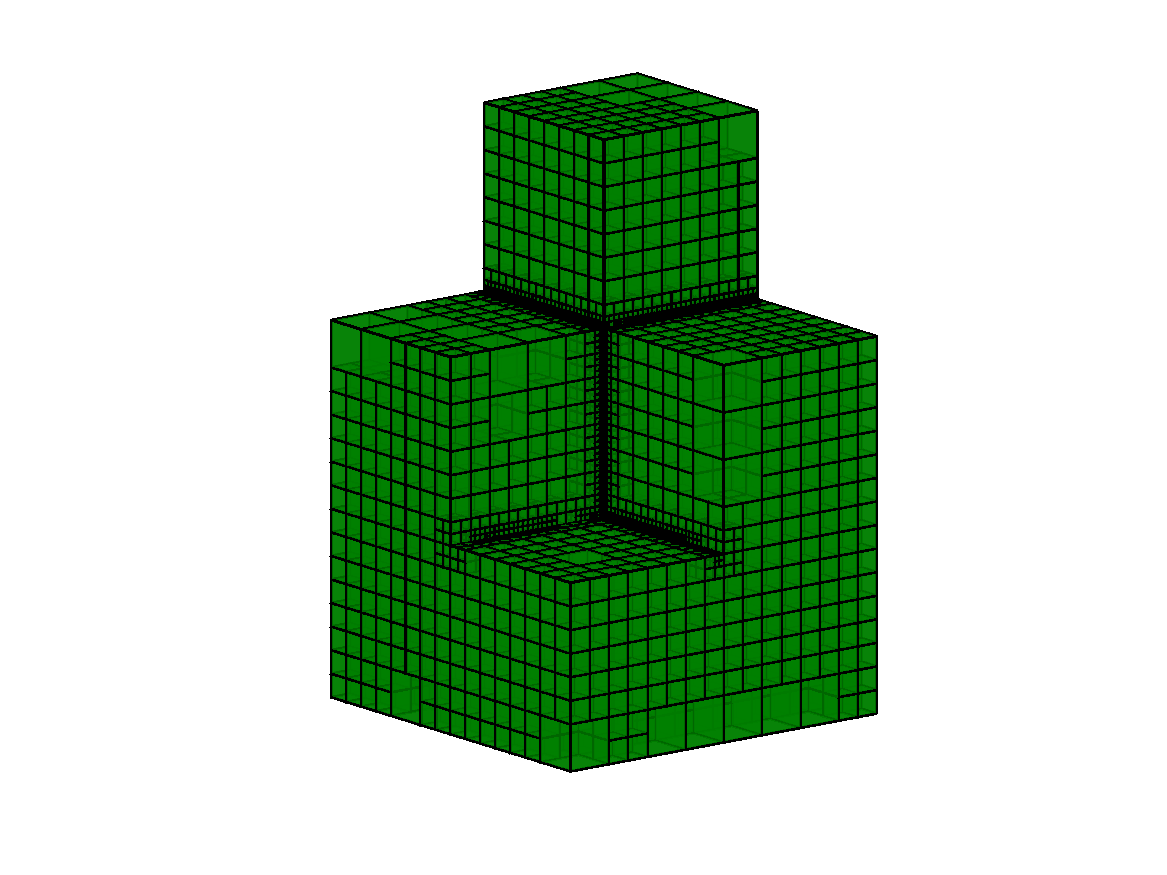}
\end{subfigure}
\begin{subfigure}[b]{0.22\textwidth}
\includegraphics[trim=210 0 210 0,clip,width=1.0\linewidth]{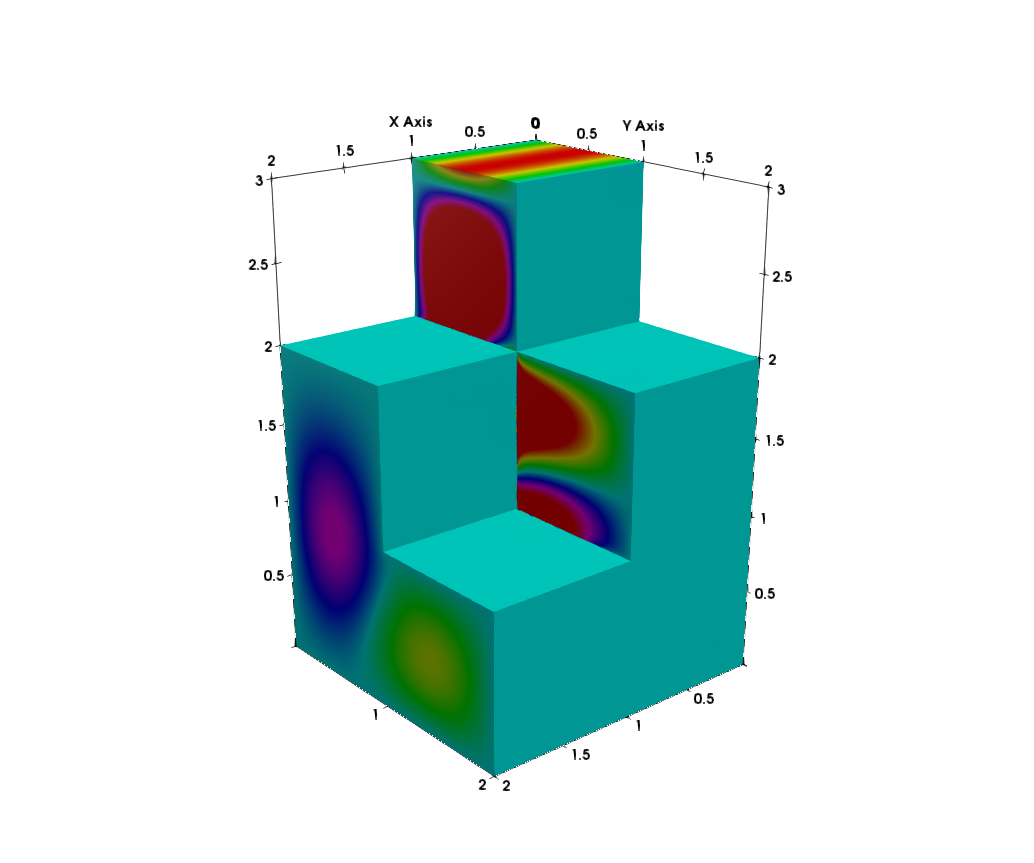}
\end{subfigure}
\begin{subfigure}[b]{0.22\textwidth}
\includegraphics[trim=210 0 210 0,clip,width=1.0\linewidth]{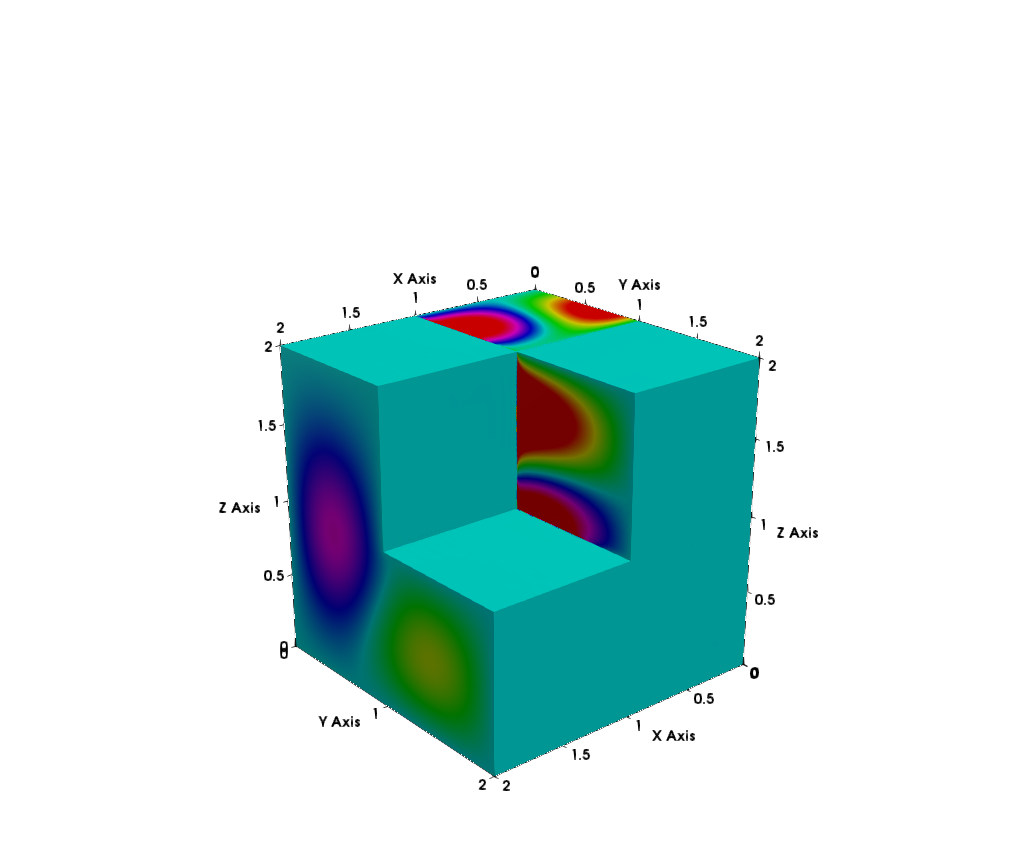}
\end{subfigure}
\begin{subfigure}[b]{0.22\textwidth}
\includegraphics[trim=210 0 210 0,clip,width=1.0\linewidth]{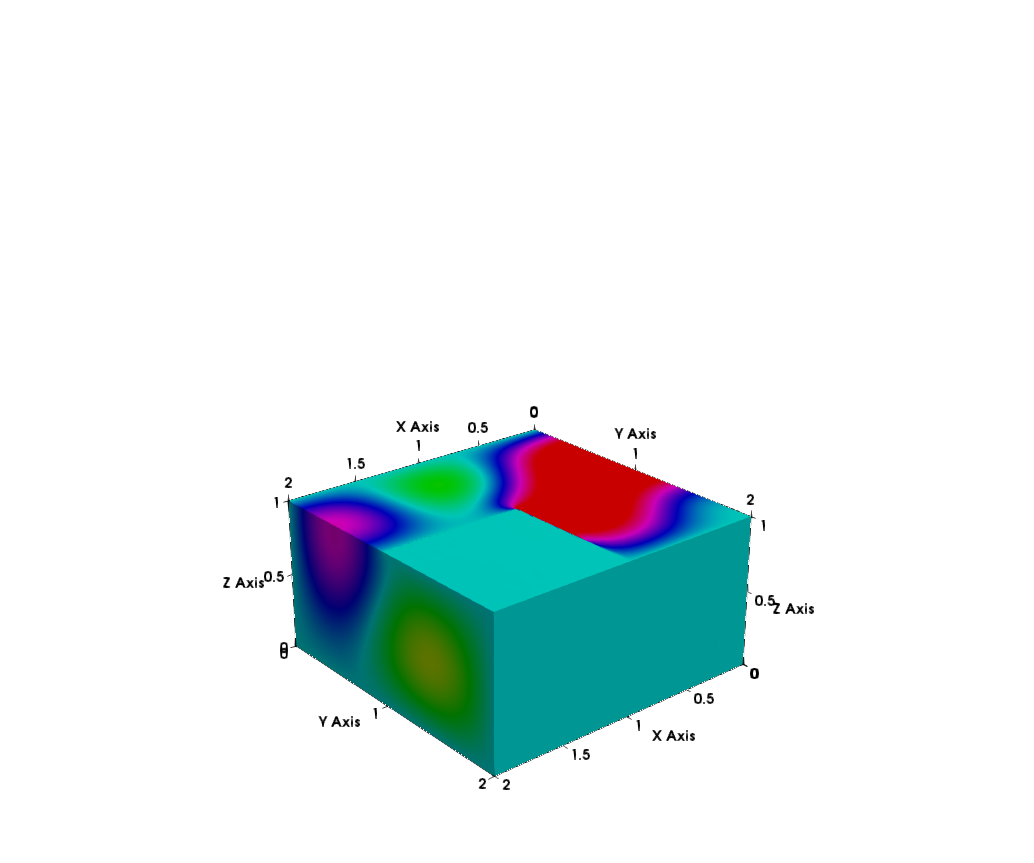}
\end{subfigure}
\begin{subfigure}[b]{0.10\textwidth}
\includegraphics[width=1.0\linewidth]{figures/results/em/fichera/cmap.png}
\end{subfigure}
\end{figure}
\vspace{-25pt}
\begin{figure}[H]
\hspace{-5pt}
\begin{subfigure}[b]{0.22\textwidth}
\includegraphics[trim=250 0 250 0,clip,width=1.0\linewidth]{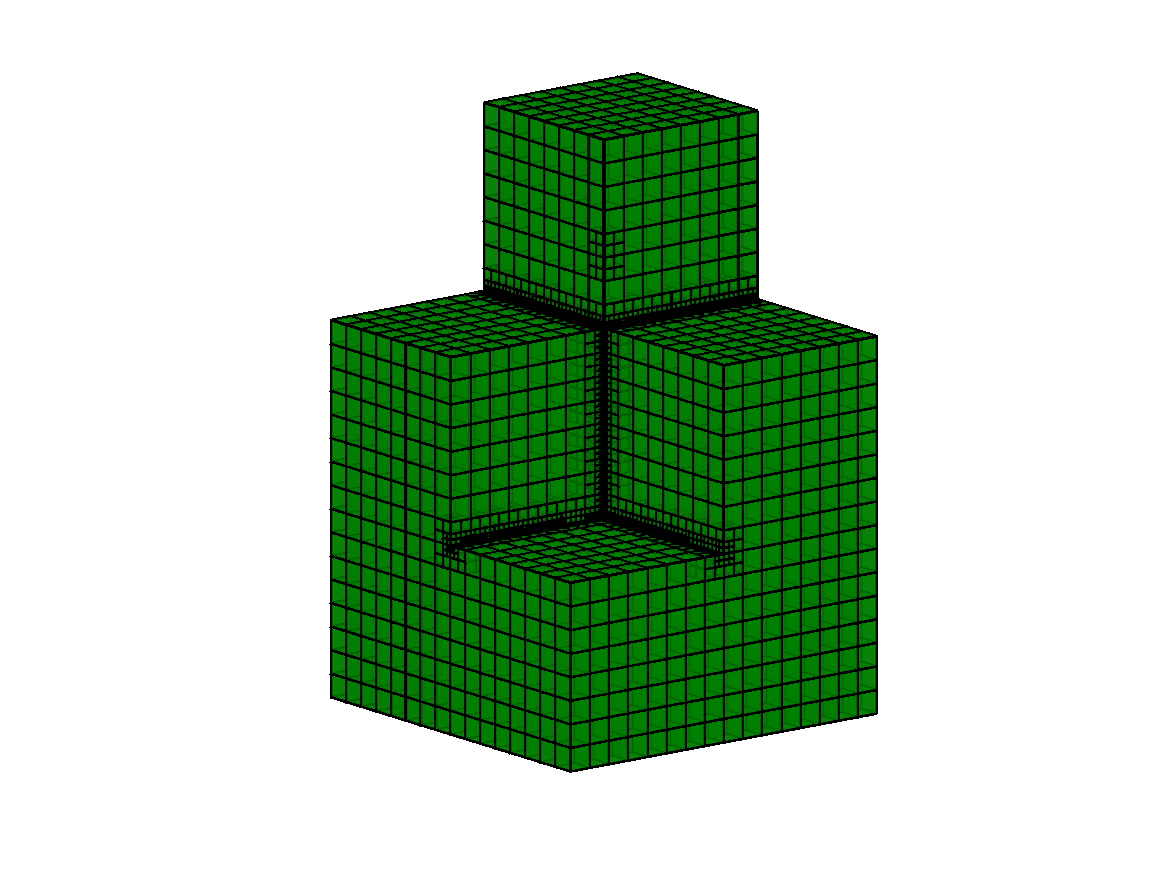}
\end{subfigure}
\begin{subfigure}[b]{0.22\textwidth}
\includegraphics[trim=210 0 210 0,clip,width=1.0\linewidth]{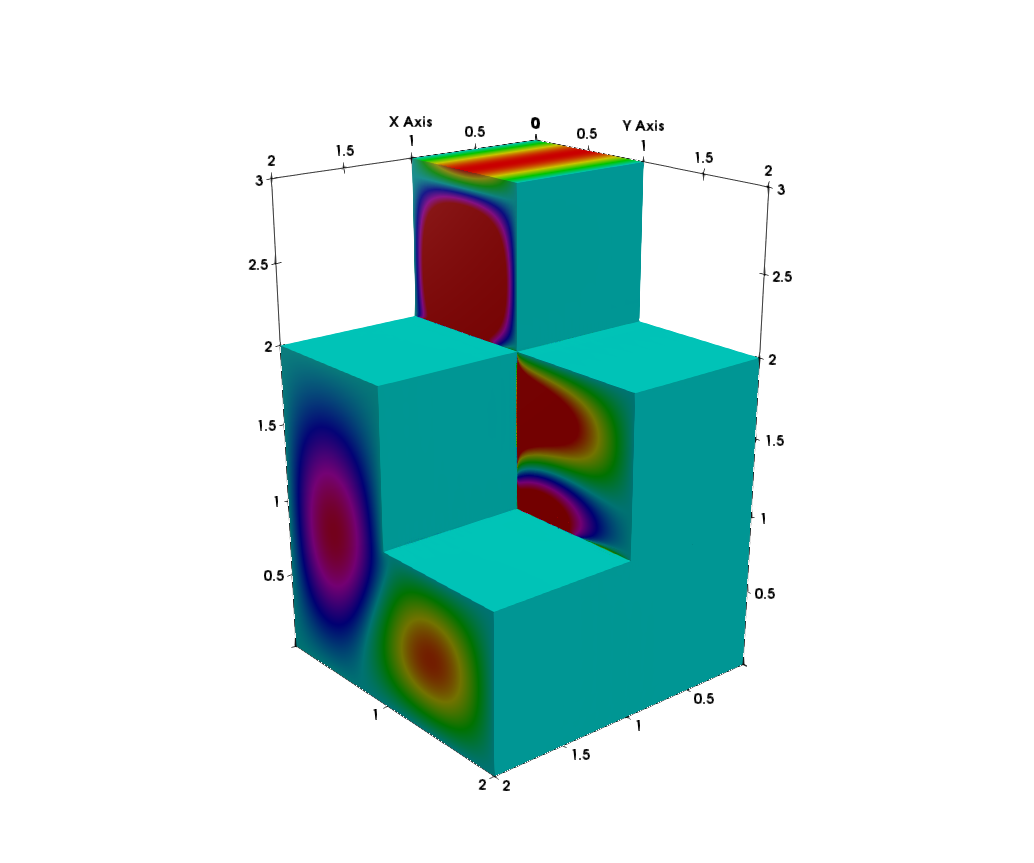}
\end{subfigure}
\begin{subfigure}[b]{0.22\textwidth}
\includegraphics[trim=210 0 210 0,clip,width=1.0\linewidth]{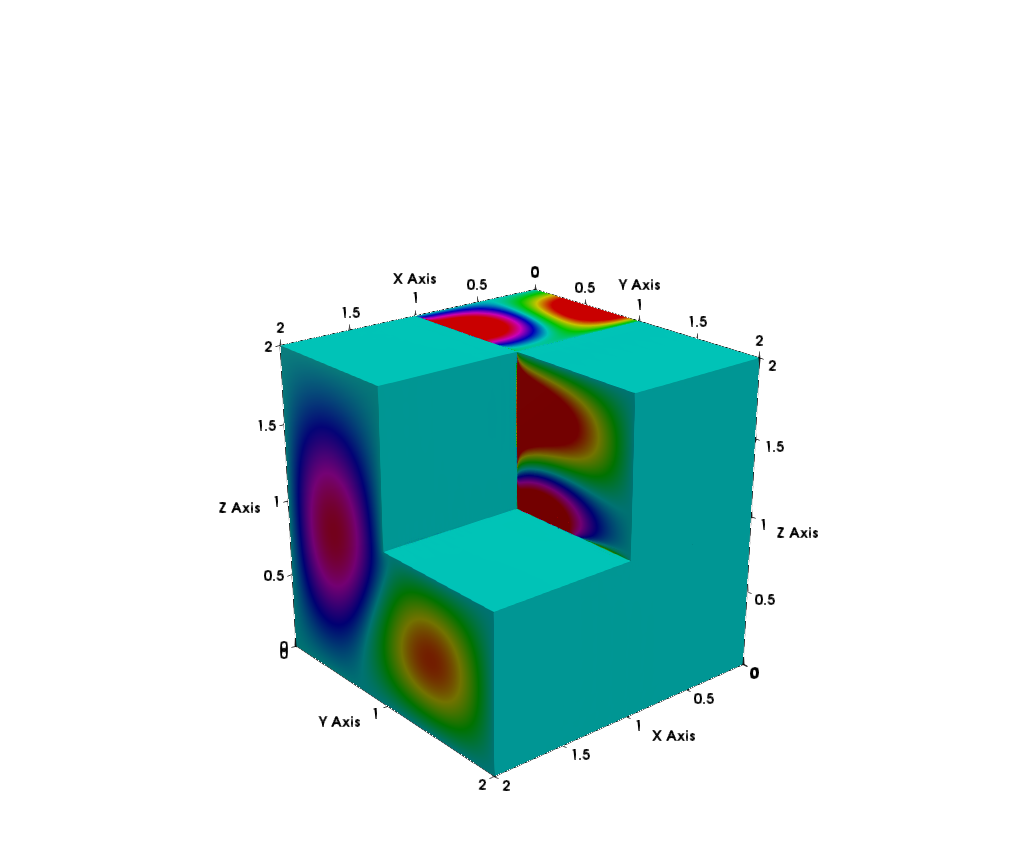}
\end{subfigure}
\begin{subfigure}[b]{0.22\textwidth}
\includegraphics[trim=210 0 210 0,clip,width=1.0\linewidth]{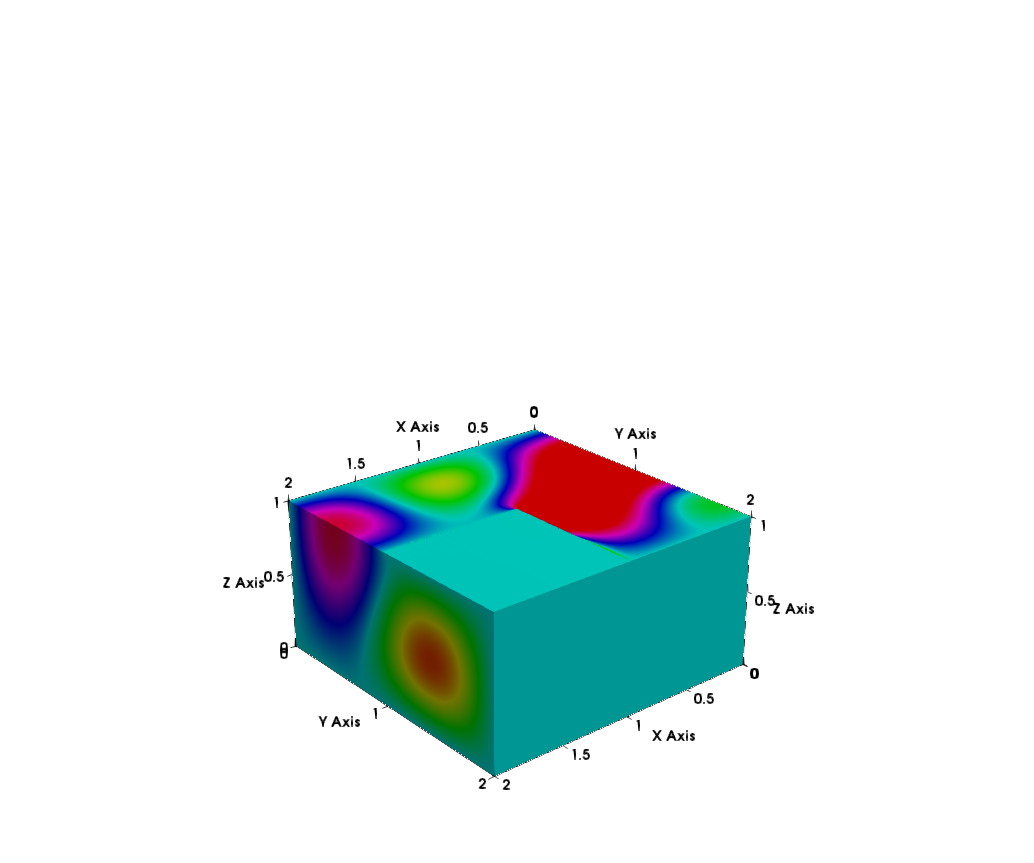}
\end{subfigure}
\begin{subfigure}[b]{0.10\textwidth}
\includegraphics[width=1.0\linewidth]{figures/results/em/fichera/cmap.png}
\end{subfigure}
\end{figure}
\vspace{-25pt}
\begin{figure}[H]
\hspace{-5pt}
\begin{subfigure}[b]{0.22\textwidth}
\includegraphics[trim=250 0 250 0,clip,width=1.0\linewidth]{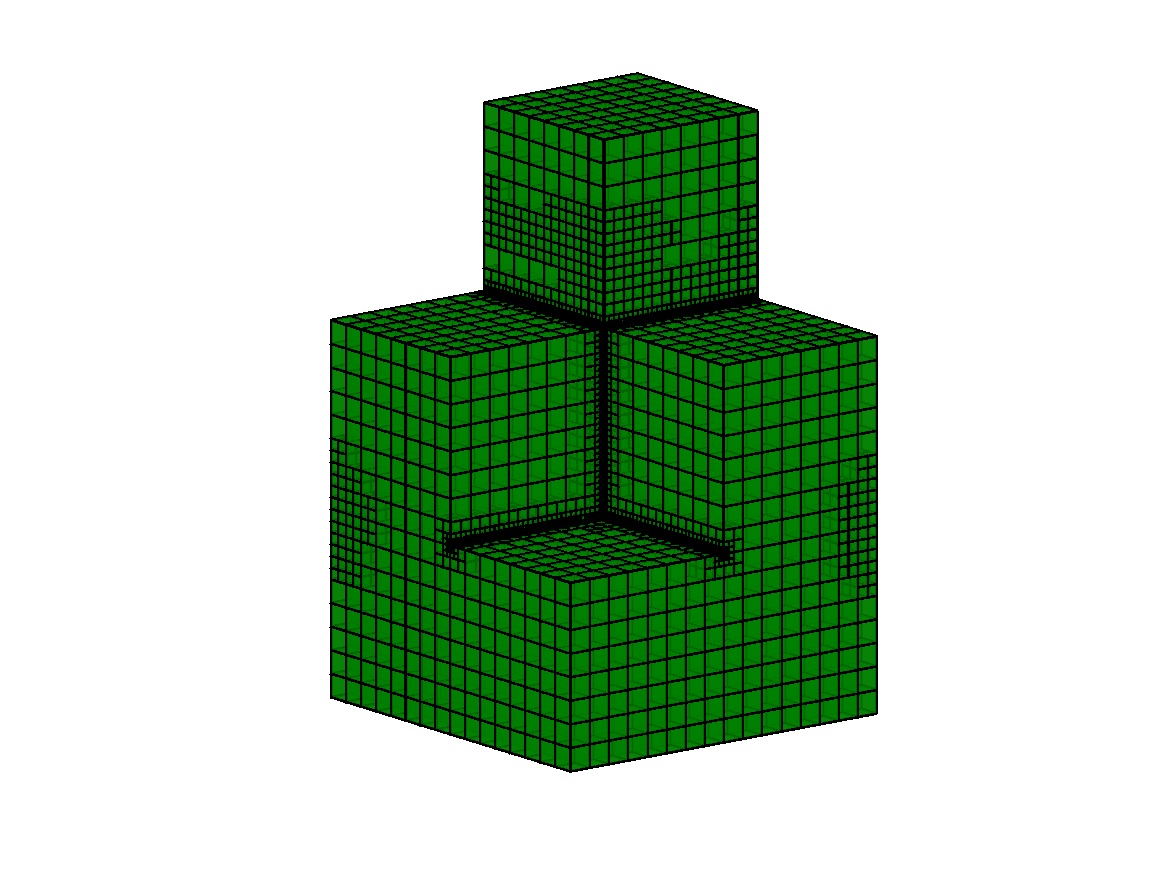}
\end{subfigure}
\begin{subfigure}[b]{0.22\textwidth}
\includegraphics[trim=210 0 210 0,clip,width=1.0\linewidth]{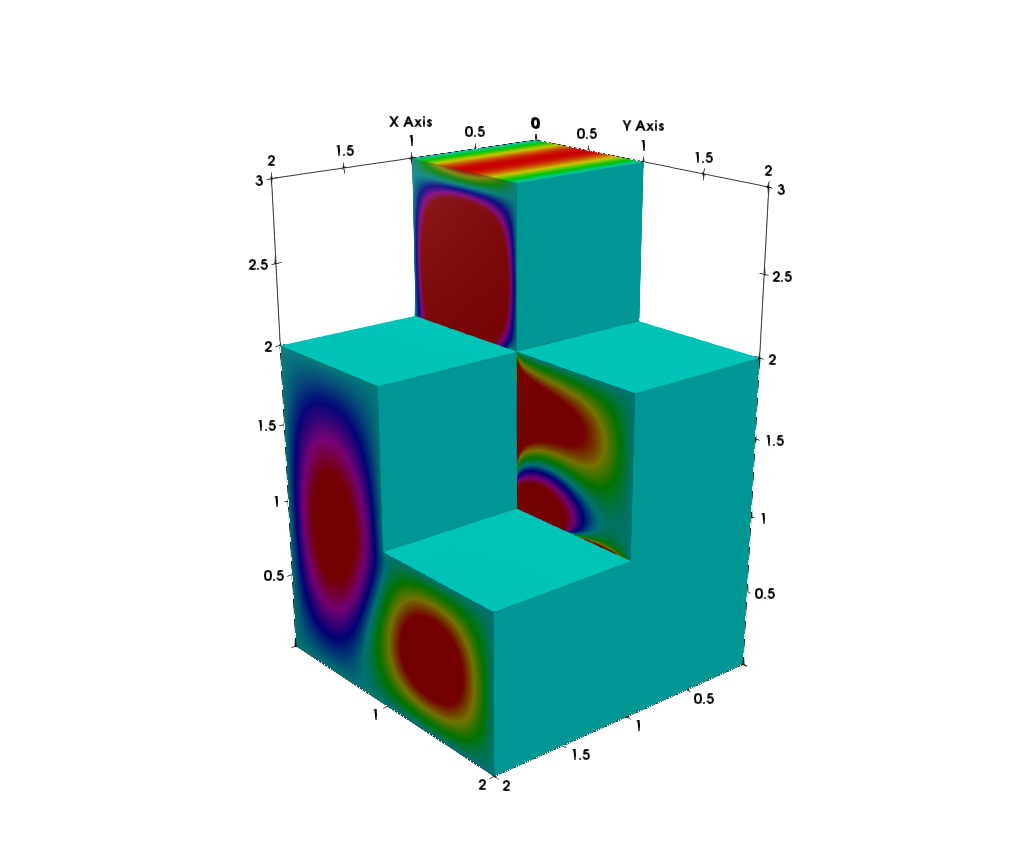}
\end{subfigure}
\begin{subfigure}[b]{0.22\textwidth}
\includegraphics[trim=210 0 210 0,clip,width=1.0\linewidth]{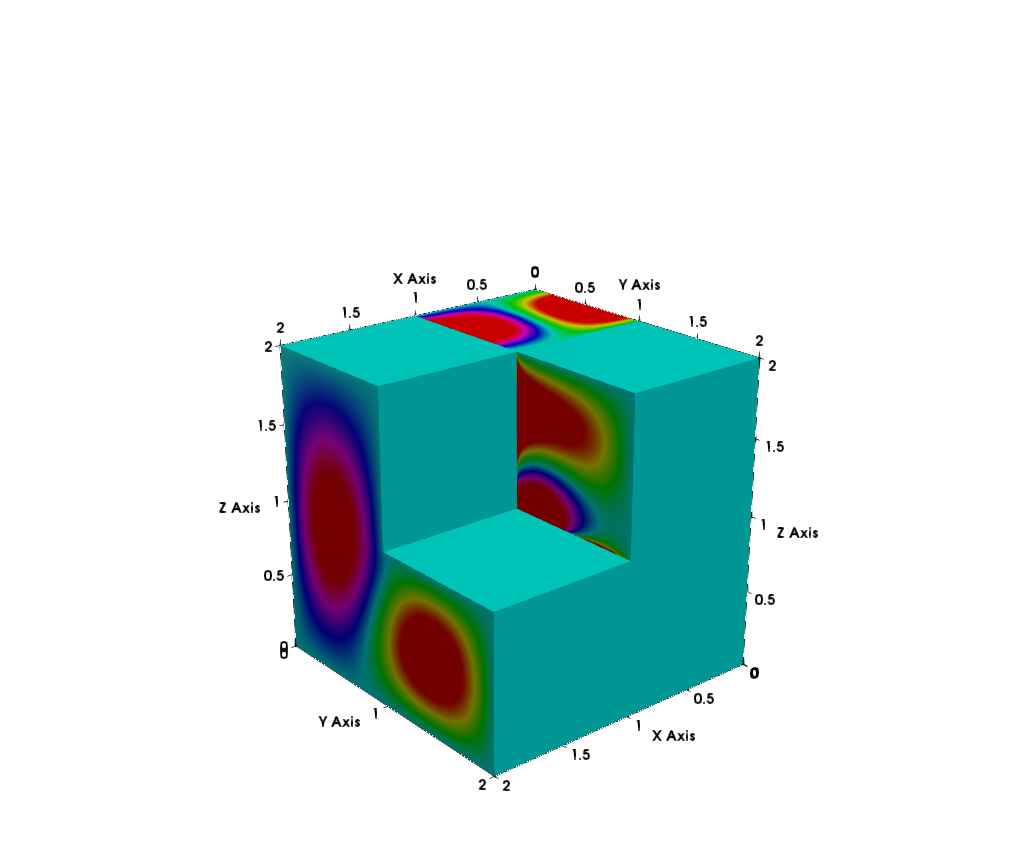}
\end{subfigure}
\begin{subfigure}[b]{0.22\textwidth}
\includegraphics[trim=210 0 210 0,clip,width=1.0\linewidth]{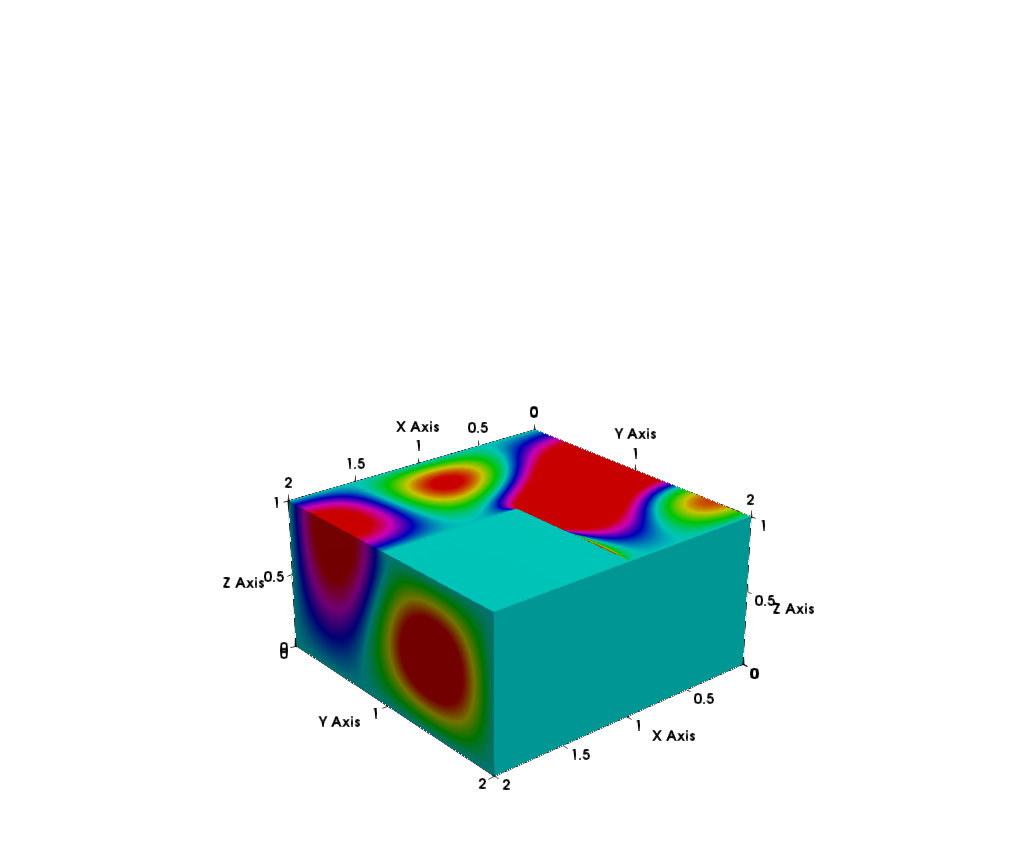}
\end{subfigure}
\begin{subfigure}[b]{0.10\textwidth}
\includegraphics[width=1.0\linewidth]{figures/results/em/fichera/cmap.png}
\end{subfigure}
\setcounter{figure}{2}
\caption{Evolution of the mesh and the numerical solution for the real part of the x-component of the electric field. These results are for the meshes 1,3,5,7,9,11,13 and 15.}\label{fig:fichera_prob}
\end{figure}

%% file: beam_20_figs.tex

\begin{figure}[H]
      \centering
\captionsetup[subfigure]{labelformat=empty}
\begin{subfigure}[b]{0.24\textwidth}
      \centering
      \includegraphics[trim=220 0 220 0,clip,width=1.0\linewidth]{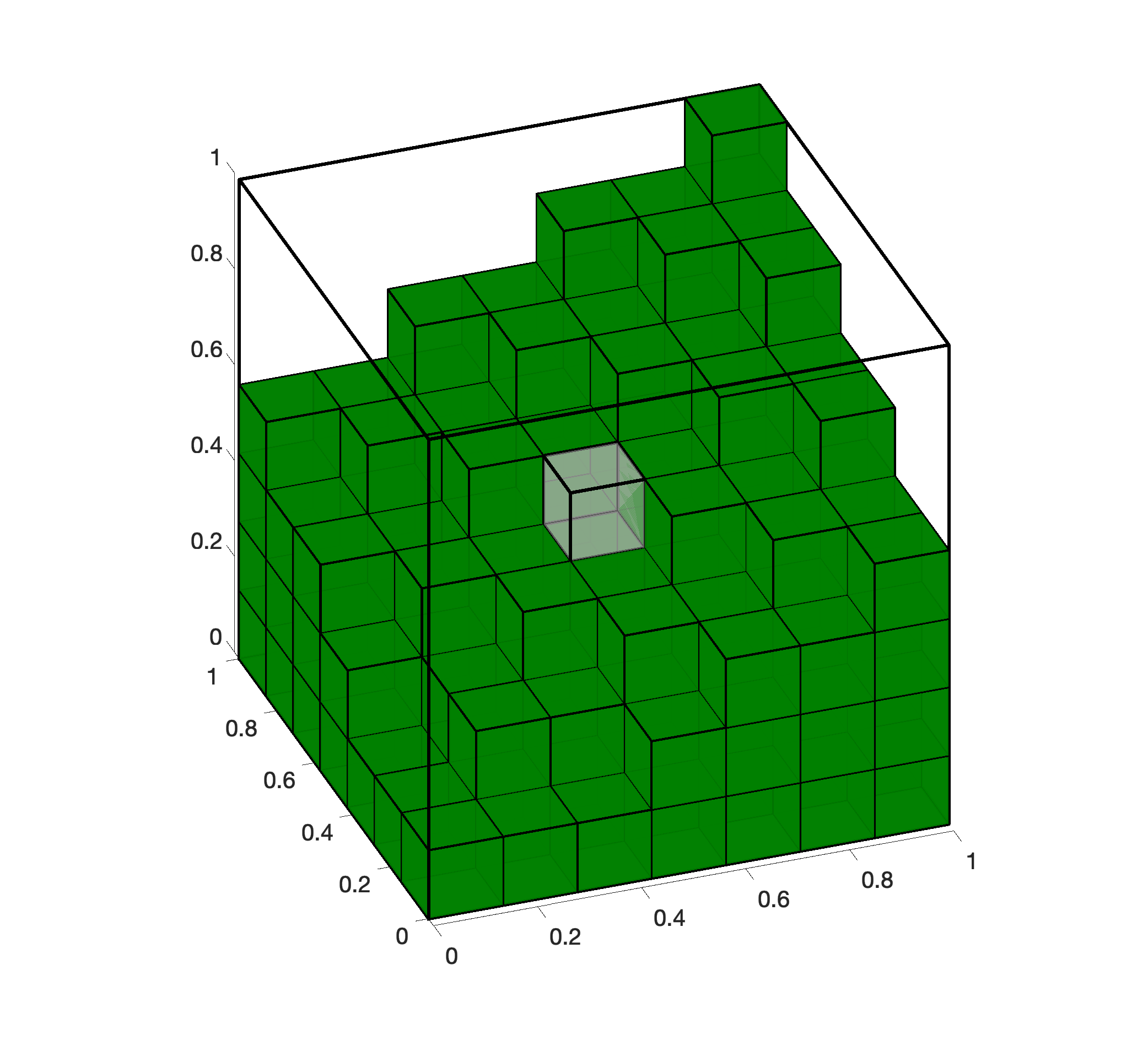}
      \centering
      \caption{\hspace{20pt}Refinement 1}
\end{subfigure}
\begin{subfigure}[b]{0.24\textwidth}
      \centering
      \includegraphics[trim=220 0 220 0,clip,width=1.0\linewidth]{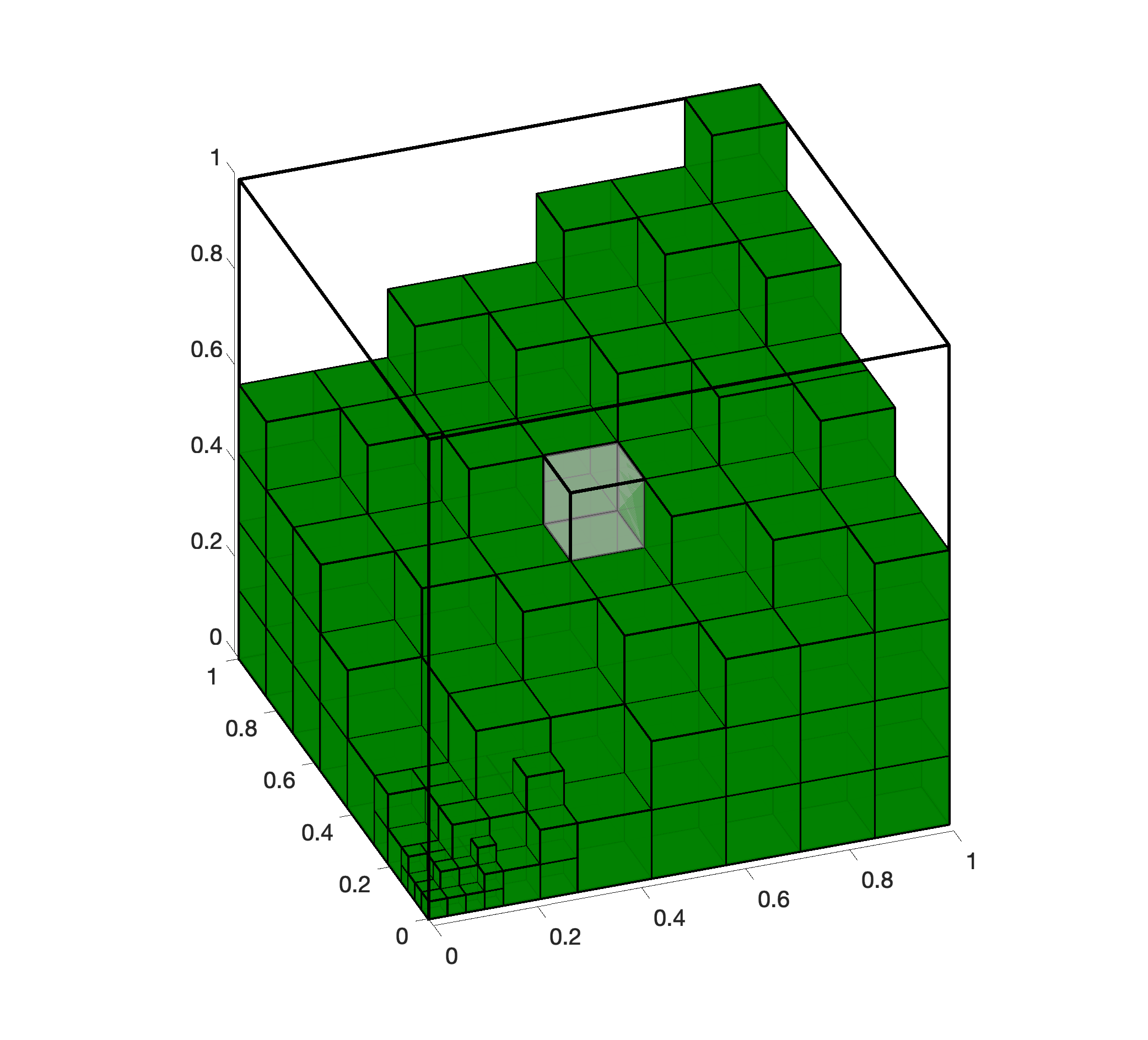}
      \caption{\hspace{20pt}Refinement 3}
\end{subfigure}
\begin{subfigure}[b]{0.24\textwidth}
      \centering
      \includegraphics[trim=220 0 220 0,clip,width=1.0\linewidth]{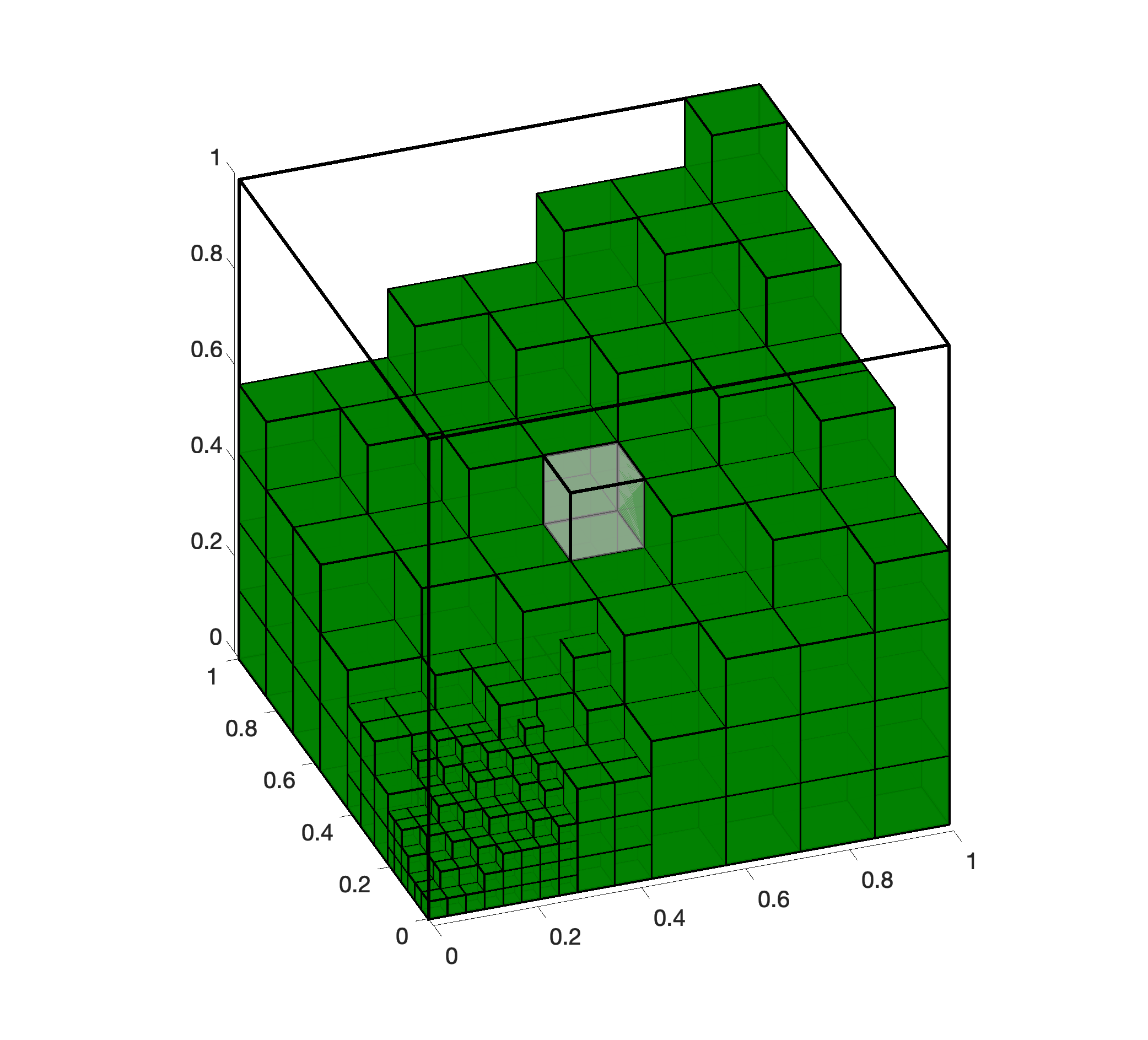}
      \caption{\hspace{20pt}Refinement 5}
\end{subfigure}
\begin{subfigure}[b]{0.24\textwidth}
      \centering
      \includegraphics[trim=220 0 220 0,clip,width=1.0\linewidth]{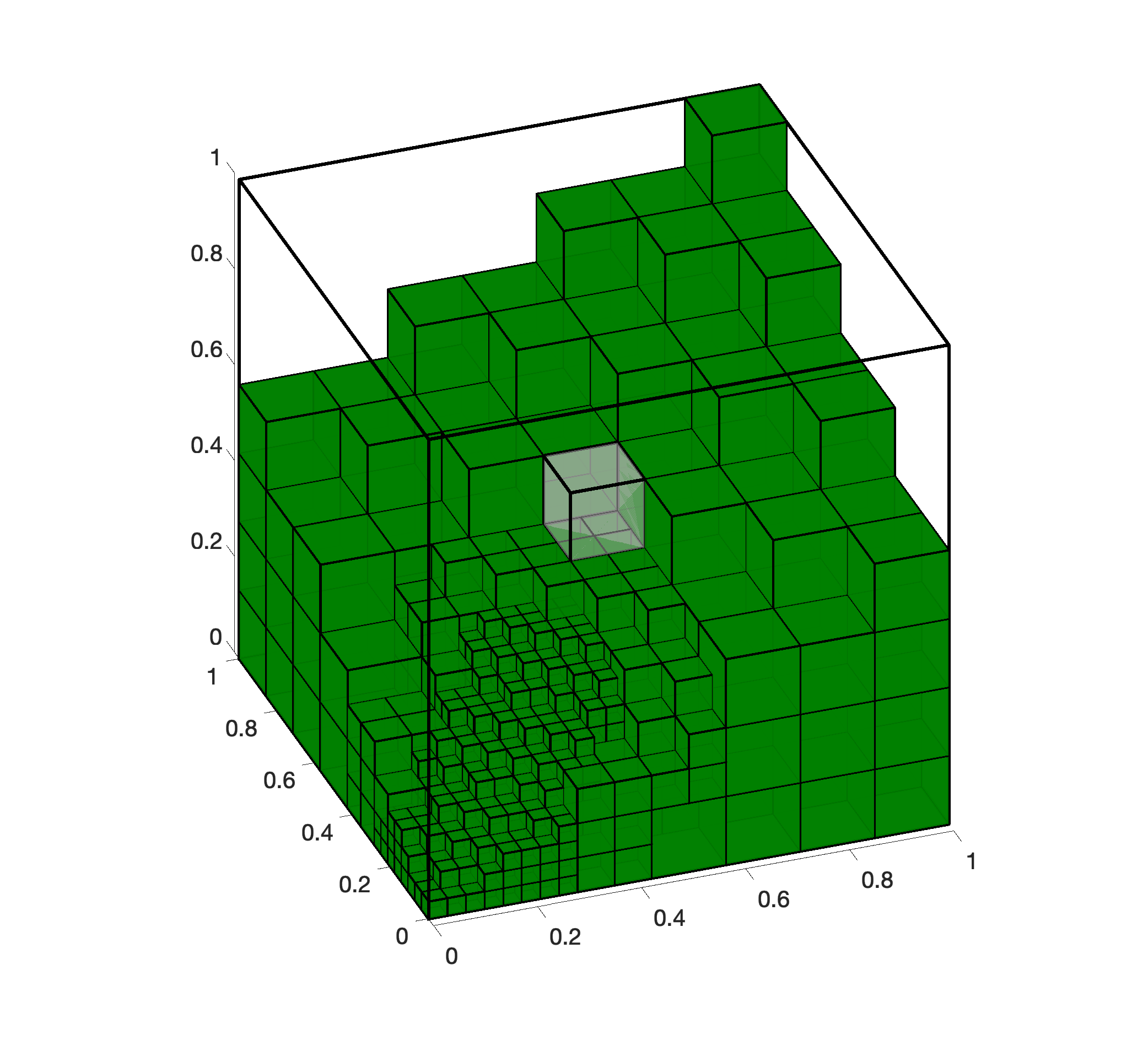}
      \caption{\hspace{20pt}Refinement 7}
\end{subfigure}
%

\begin{subfigure}[b]{0.24\textwidth}
      \centering
      \includegraphics[trim=220 0 220 0,clip,width=1.0\linewidth]{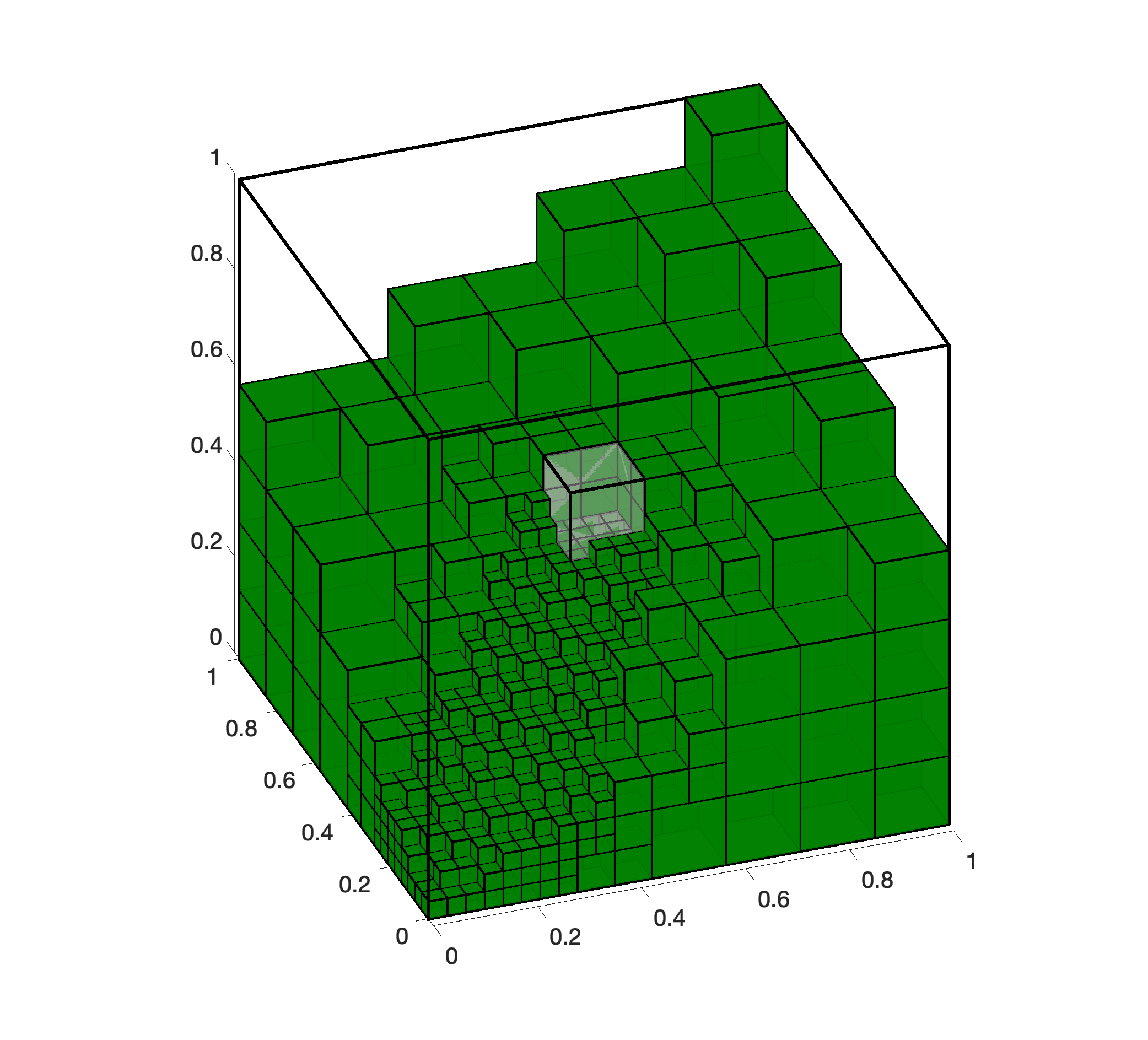}
      \caption{\hspace{20pt}Refinement 9}
\end{subfigure}
\begin{subfigure}[b]{0.24\textwidth}
      \centering
      \includegraphics[trim=220 0 220 0,clip,width=1.0\linewidth]{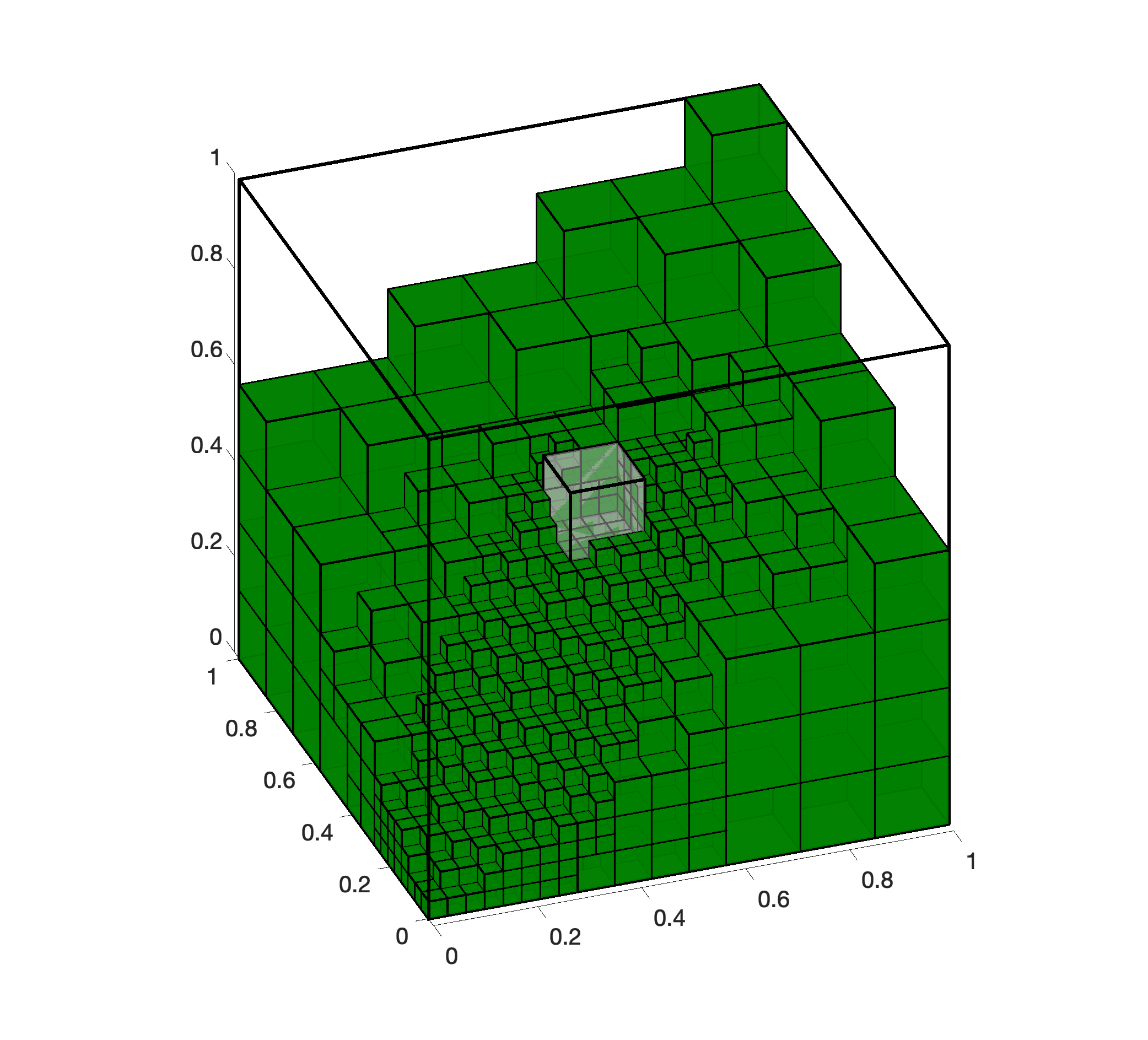}
      \caption{\hspace{20pt}Refinement 11}
\end{subfigure}
\begin{subfigure}[b]{0.24\textwidth}
      \centering
      \includegraphics[trim=220 0 220 0,clip,width=1.0\linewidth]{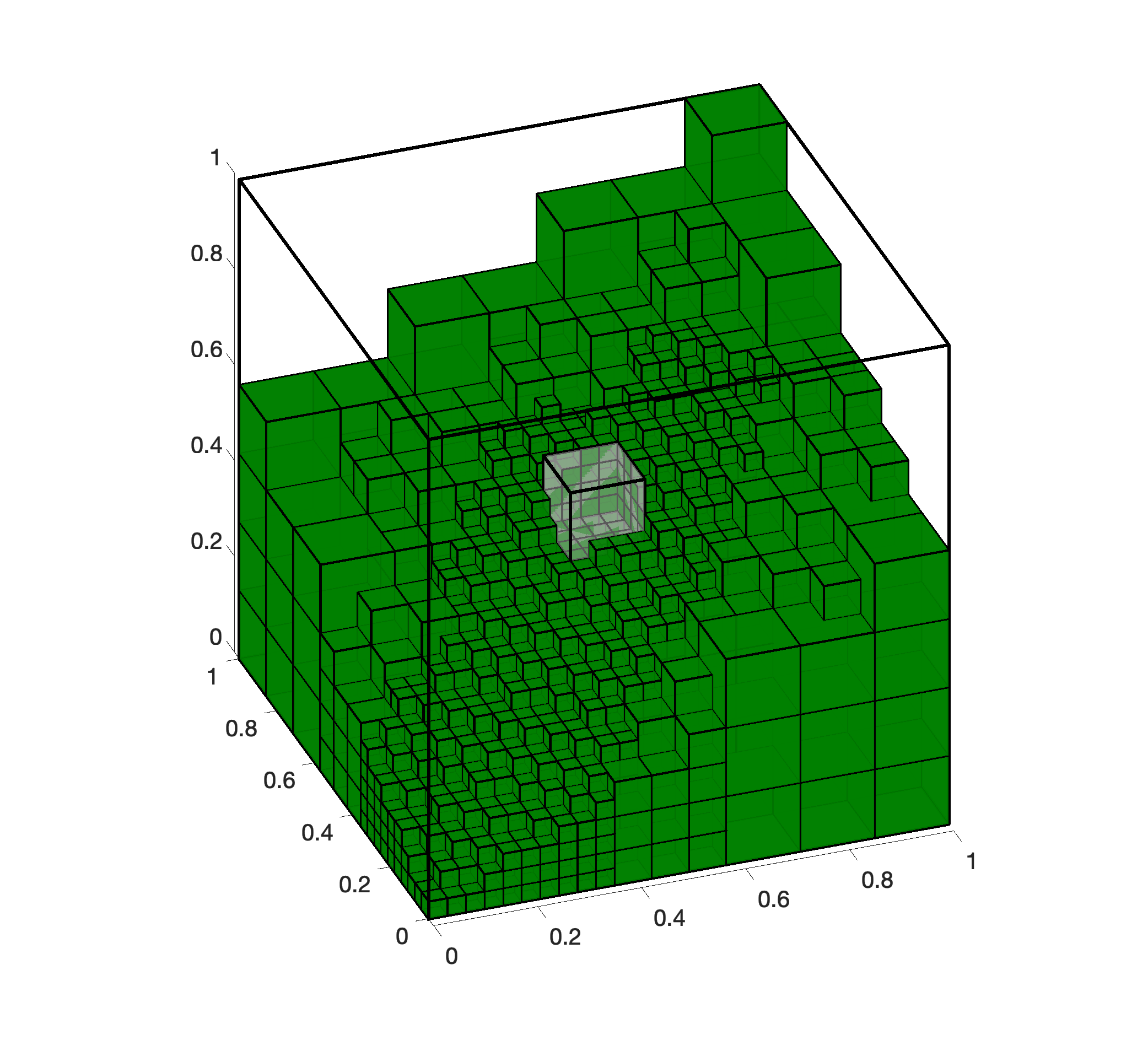}
      \caption{\hspace{20pt}Refinement 13}
\end{subfigure}
\begin{subfigure}[b]{0.24\textwidth}
      \centering
      \includegraphics[trim=220 0 220 0,clip,width=1.0\linewidth]{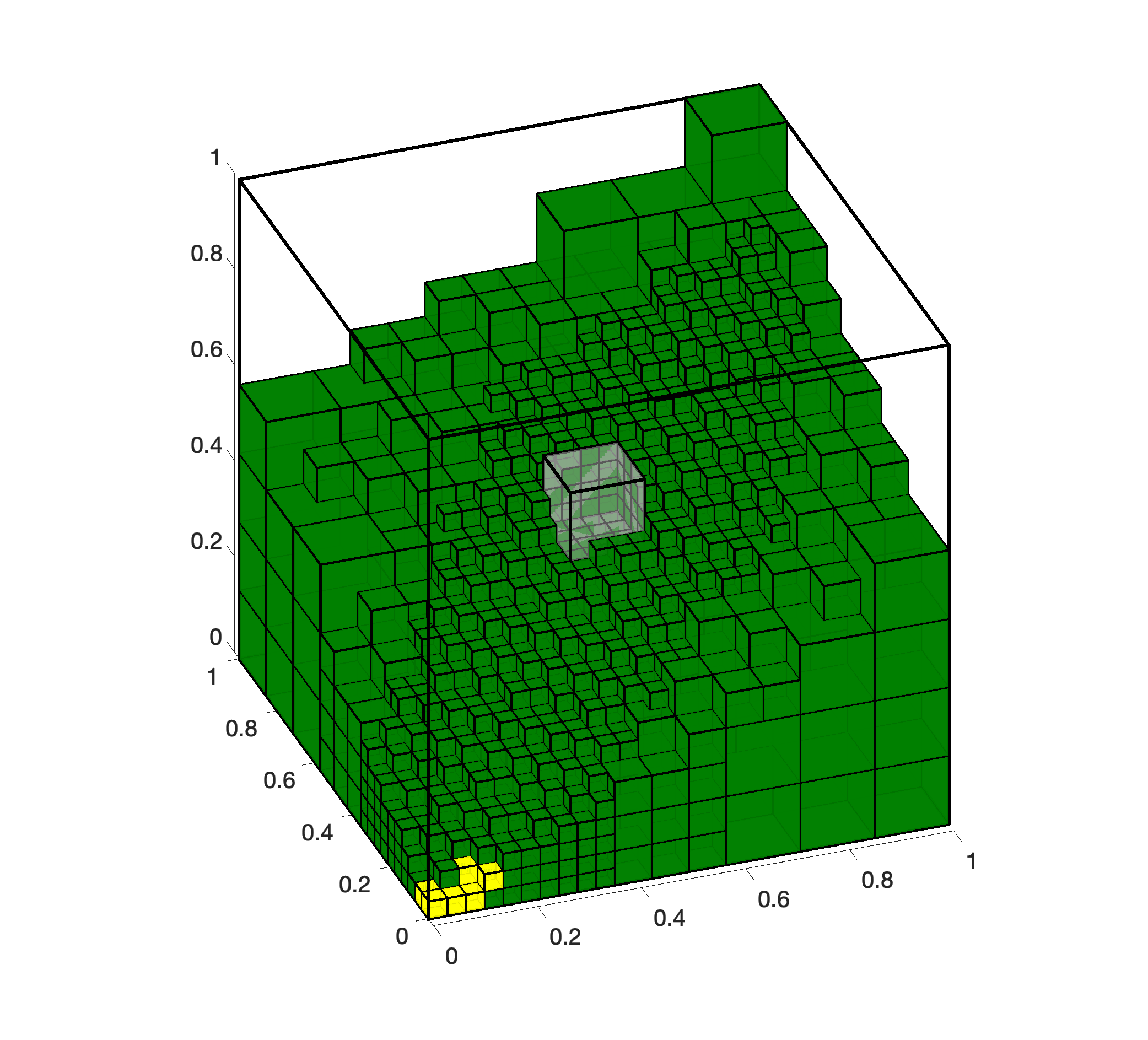}
      \caption{\hspace{20pt}Refinement 15}
\end{subfigure}

\begin{subfigure}[b]{0.24\textwidth}
      \centering
      \includegraphics[trim=220 0 220 0,clip,width=1.0\linewidth]{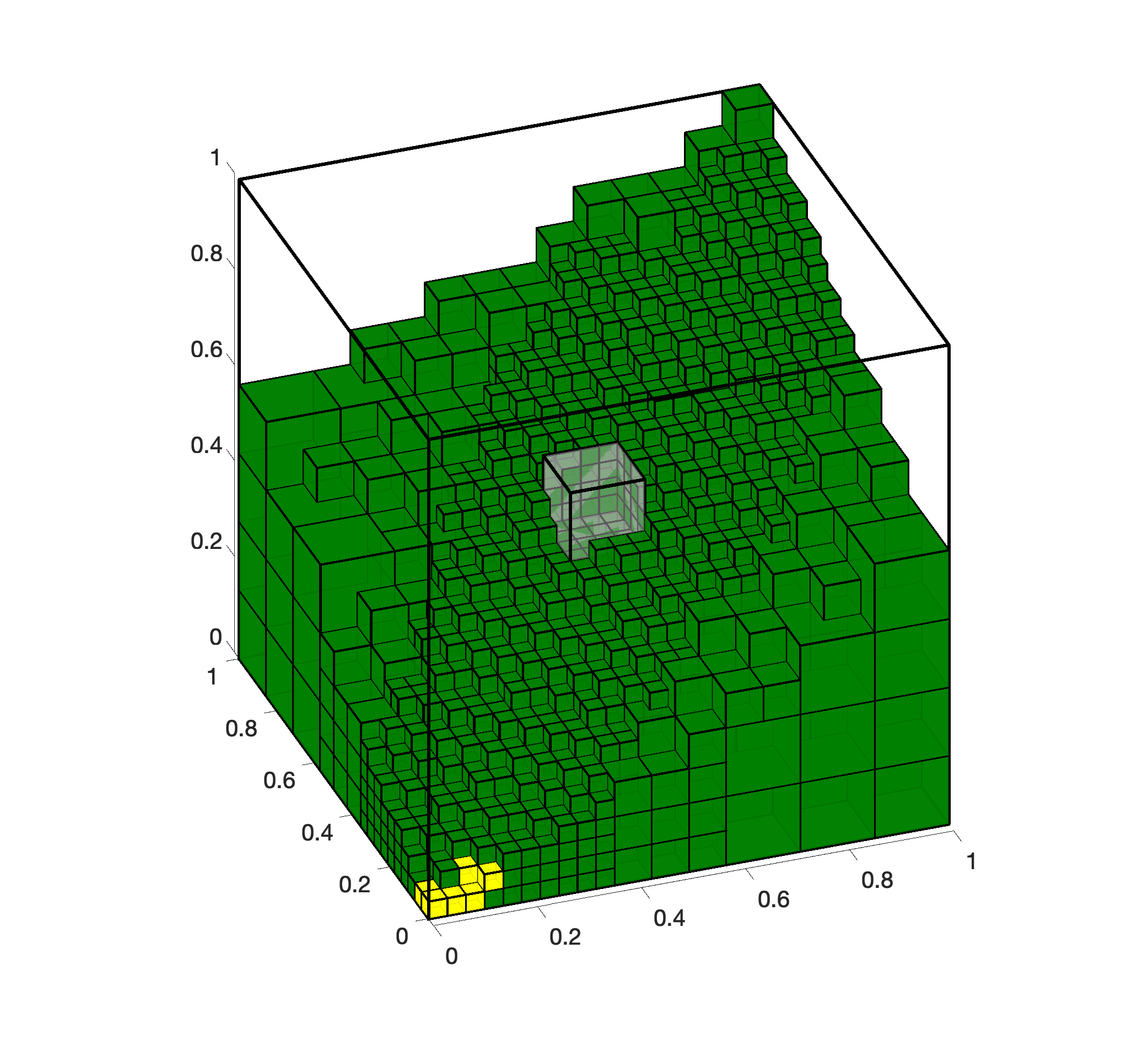}
      \caption{\hspace{20pt}Refinement 17}
\end{subfigure}
\begin{subfigure}[b]{0.24\textwidth}
      \centering
      \includegraphics[trim=220 0 220 0,clip,width=1.0\linewidth]{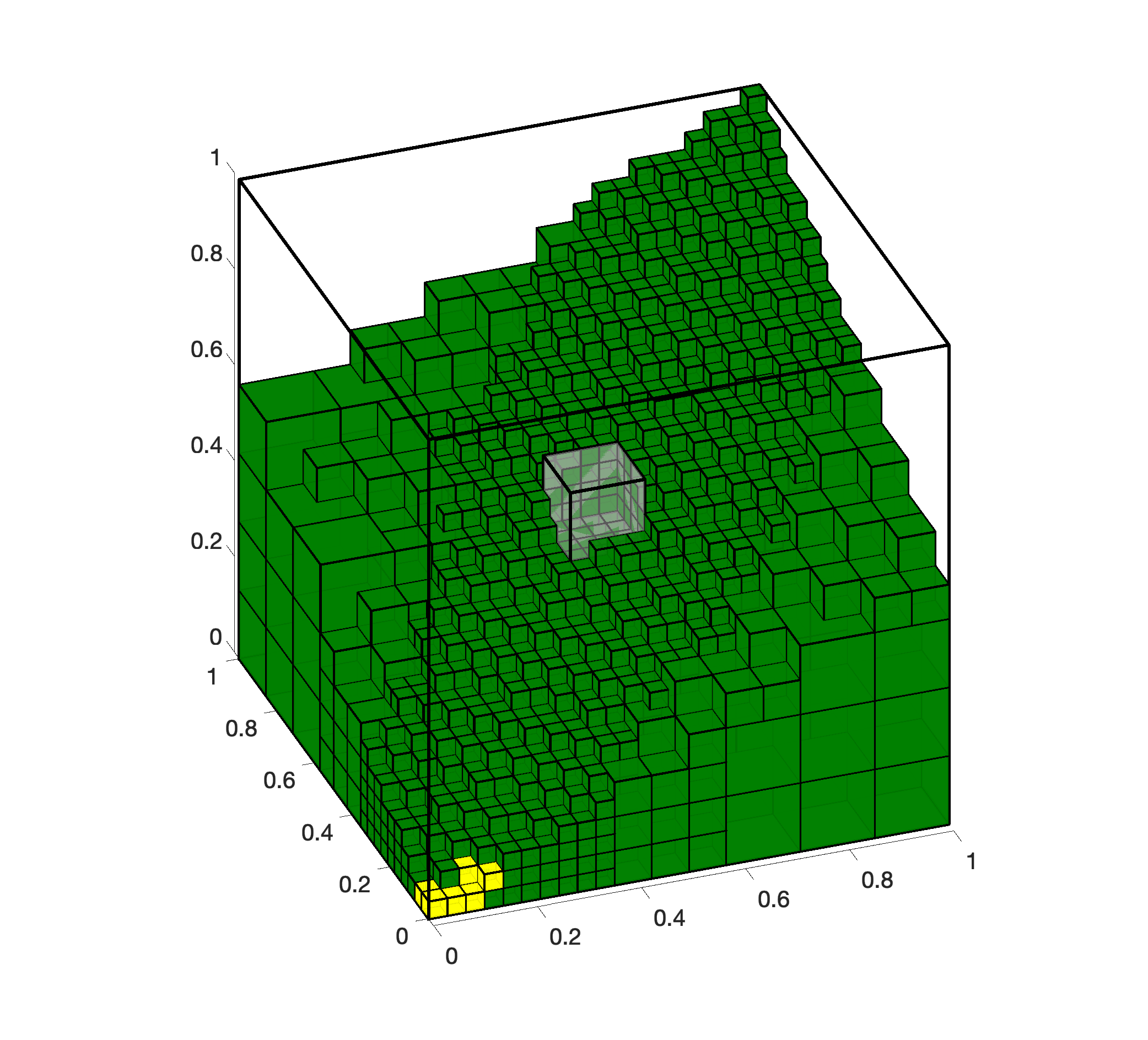}
      \caption{\hspace{20pt}Refinement 19}
\end{subfigure}
\begin{subfigure}[b]{0.24\textwidth}
      \centering
      \includegraphics[trim=220 0 220 0,clip,width=1.0\linewidth]{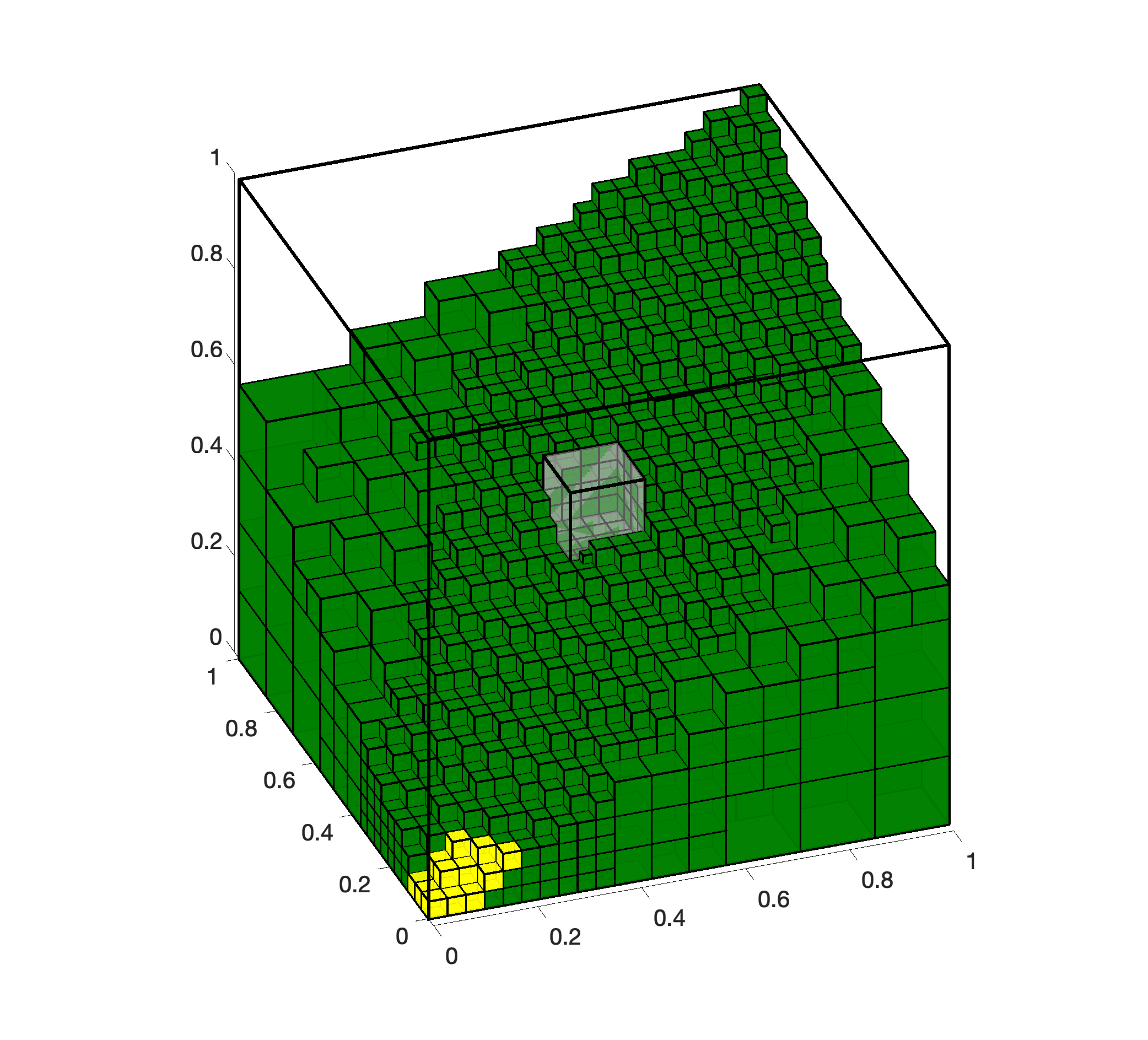}
      \caption{\hspace{20pt}Refinement 21}
\end{subfigure}
\begin{subfigure}[b]{0.24\textwidth}
      \centering
      \includegraphics[trim=220 0 220 0,clip,width=1.0\linewidth]{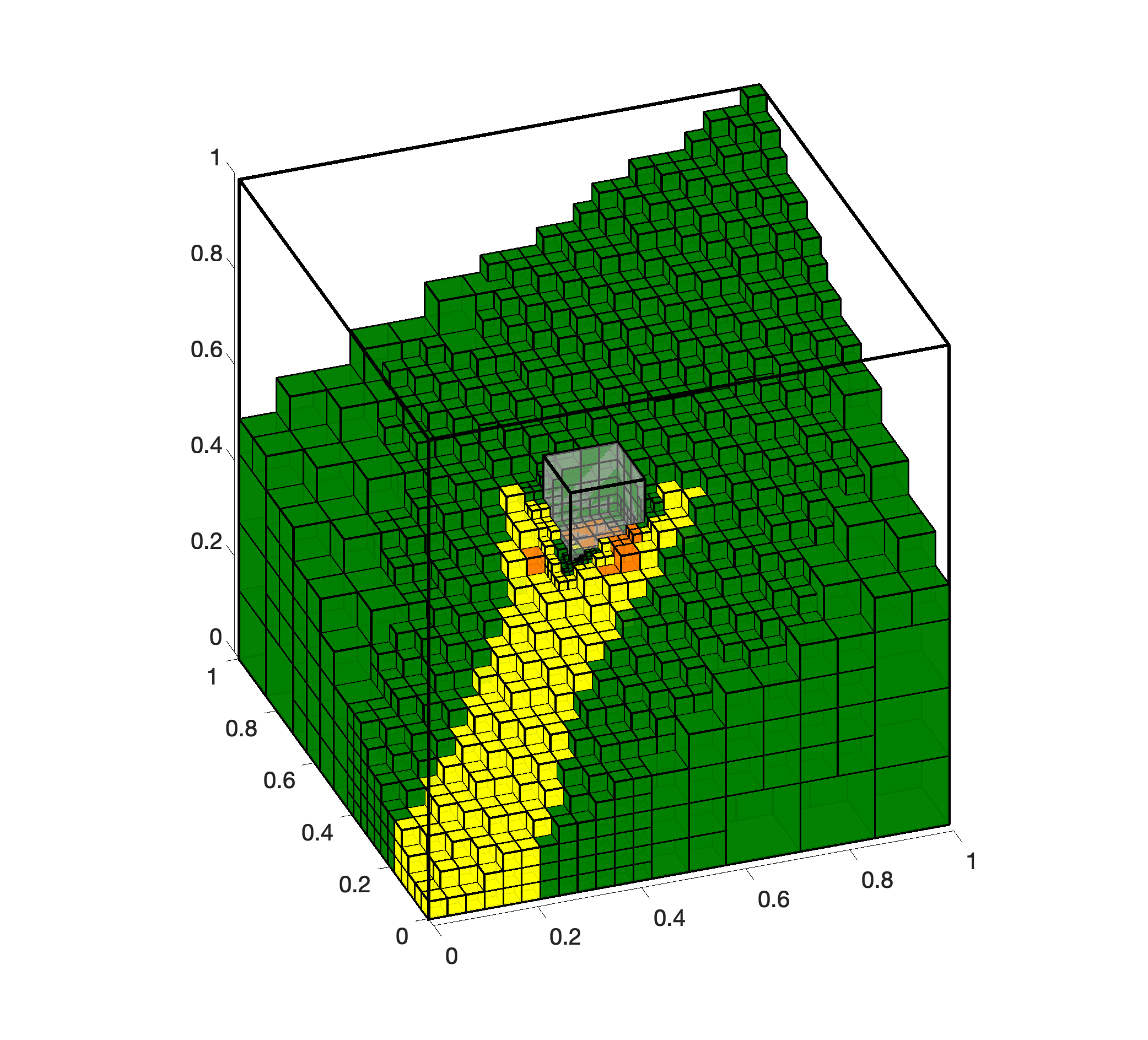}
      \caption{\hspace{20pt}Refinement 23}
\end{subfigure}
\end{figure}

\setcounter{figure}{3}
\begin{figure}[H]
\centering
\begin{tikzpicture}
\node[draw,align=center,fill=green] at (0,0) {$p=3$};
\node[draw,align=center,fill=yellow] at (1.5,0) {$p=4$};
\node[draw,align=center,fill=orange] at (3,0) {$p=5$};
\end{tikzpicture}
      \caption{Evolution of the hp-adaptive meshes.}\label{fig:max_cube2}
\end{figure}
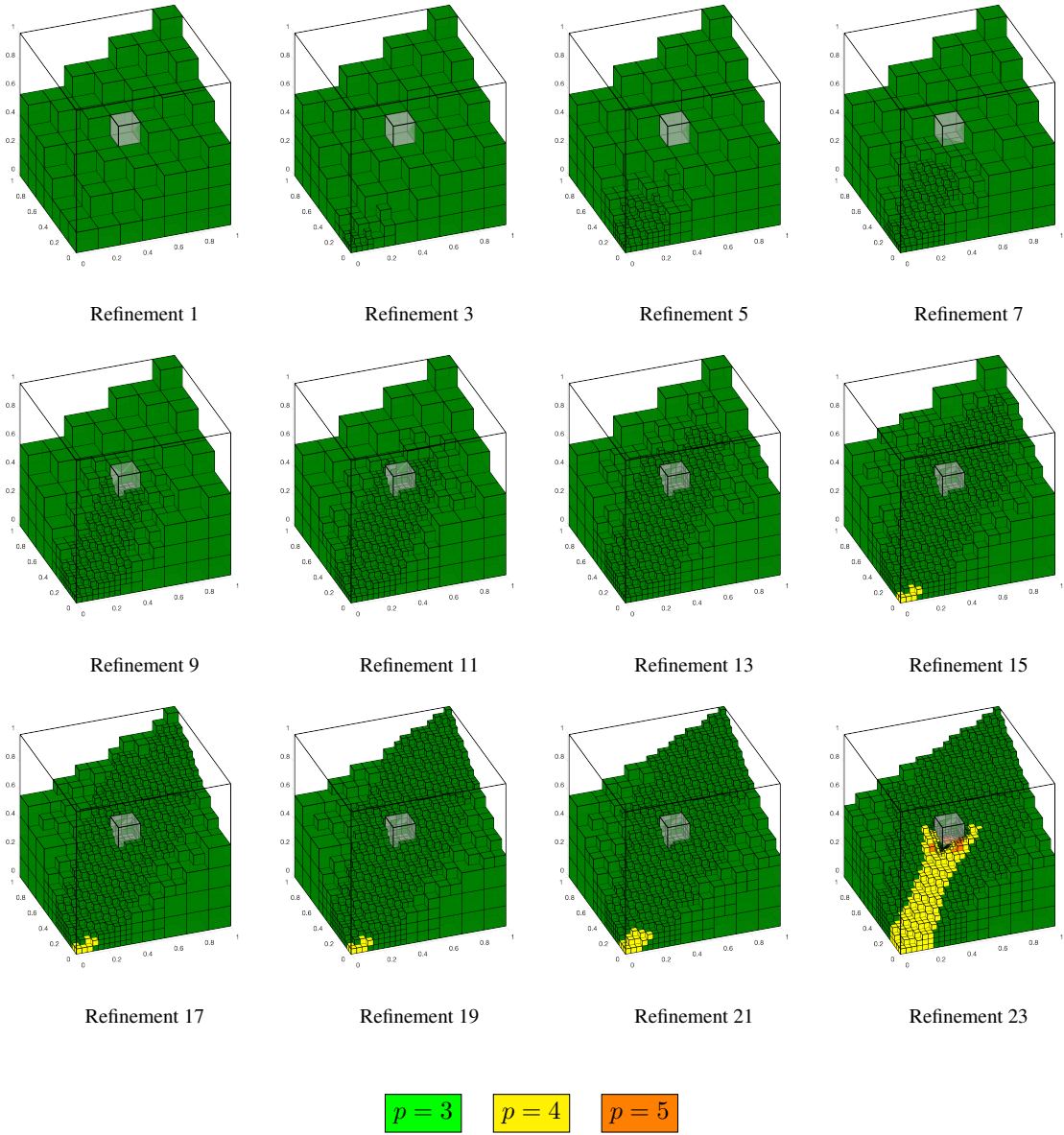
\begin{figure}[H]
\captionsetup[subfigure]{labelformat=empty}
      \centering
\begin{subfigure}[b]{0.24\textwidth}
      \centering
      \includegraphics[trim=80 0 80 0,clip,width=1.0\linewidth]{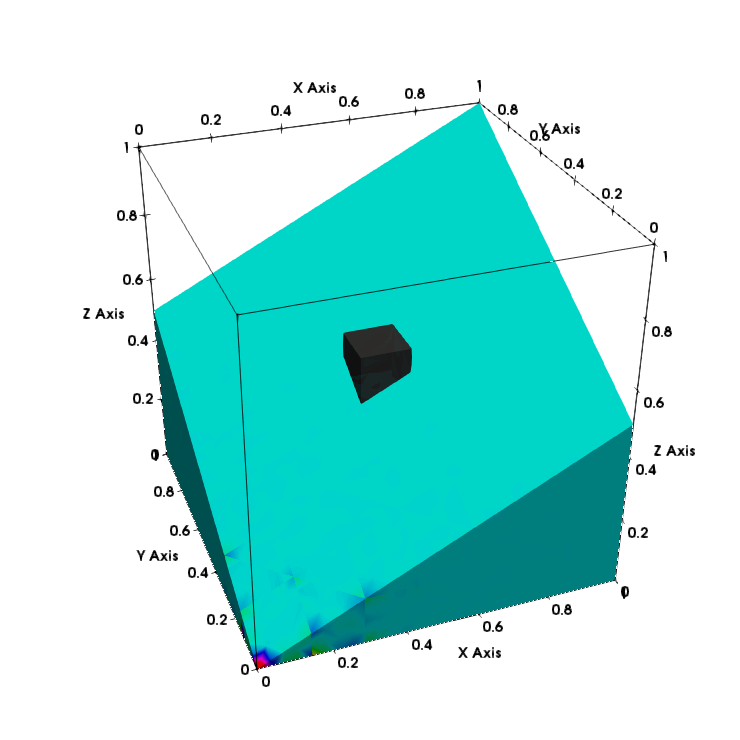}
      \centering
      \caption{\hspace{20pt}Refinement 1}
\end{subfigure}
\begin{subfigure}[b]{0.24\textwidth}
      \centering
      \includegraphics[trim=80 0 80 0,clip,width=1.0\linewidth]{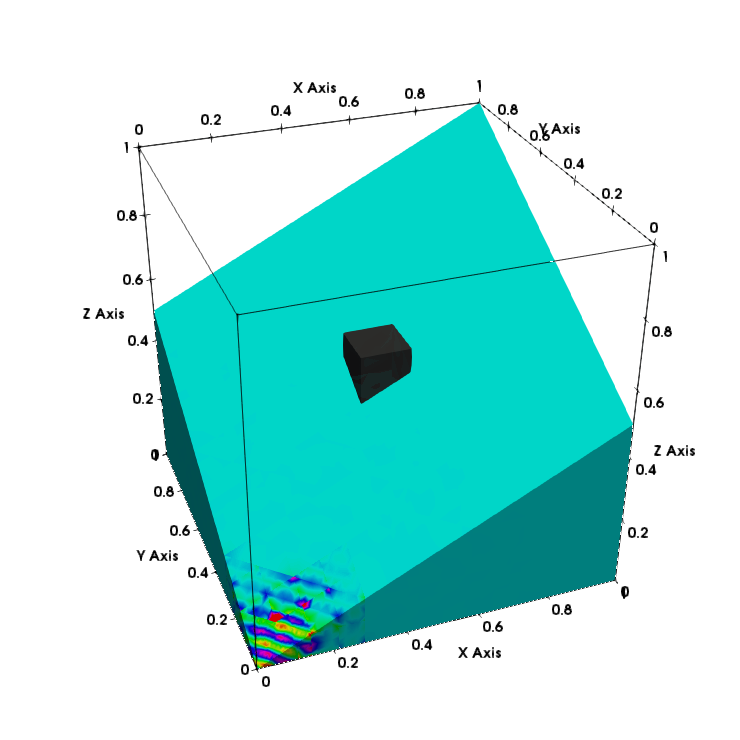}
      \caption{\hspace{20pt}Refinement 3}
\end{subfigure}
\begin{subfigure}[b]{0.24\textwidth}
      \centering
      \includegraphics[trim=80 0 80 0,clip,width=1.0\linewidth]{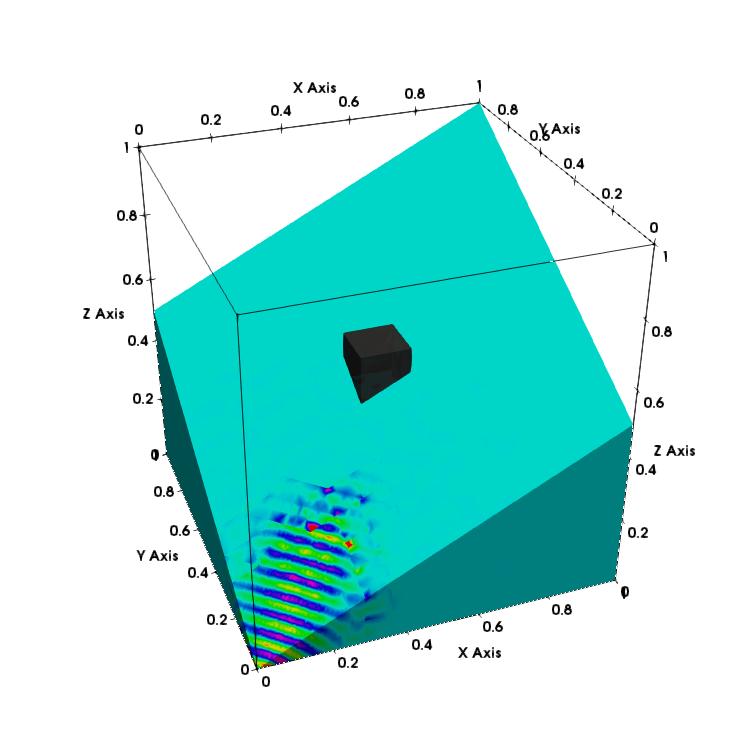}
      \caption{\hspace{20pt}Refinement 5}
\end{subfigure}
\begin{subfigure}[b]{0.24\textwidth}
      \centering
      \includegraphics[trim=80 0 80 0,clip,width=1.0\linewidth]{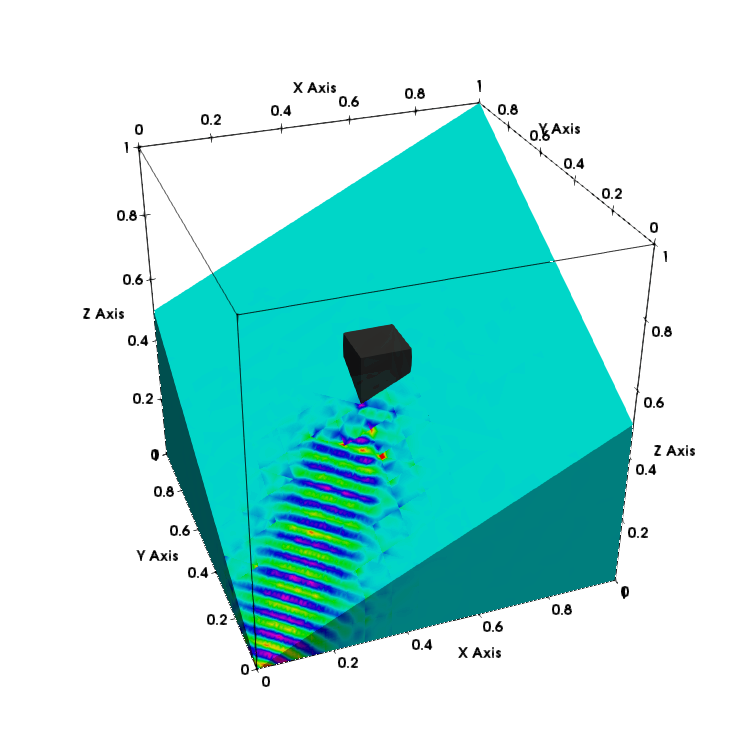}
      \caption{\hspace{20pt}Refinement 7}
\end{subfigure}
%
\begin{subfigure}[b]{0.24\textwidth}
      \centering
      \includegraphics[trim=80 0 80 0,clip,width=1.0\linewidth]{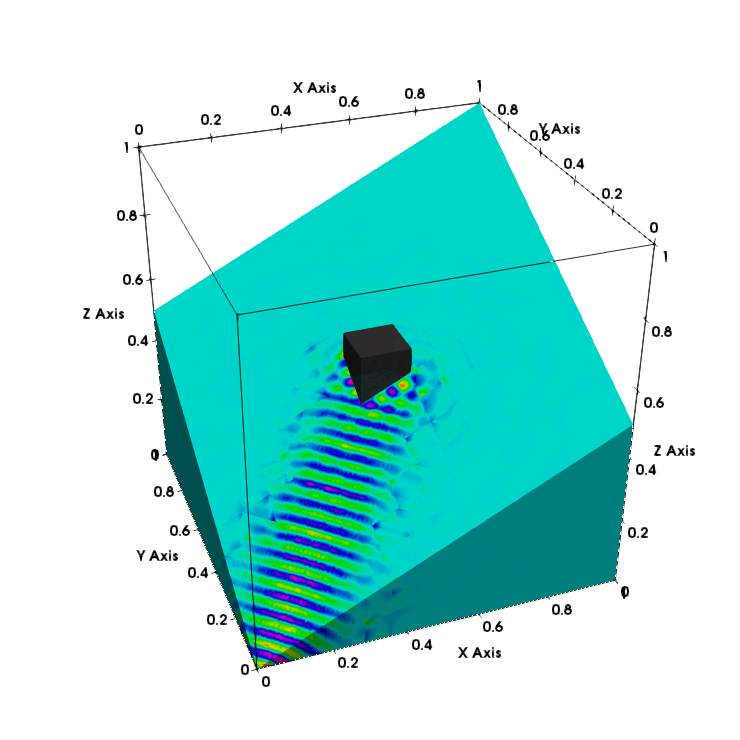}
      \caption{\hspace{20pt}Refinement 9}
\end{subfigure}
\begin{subfigure}[b]{0.24\textwidth}
      \centering
      \includegraphics[trim=80 0 80 0,clip,width=1.0\linewidth]{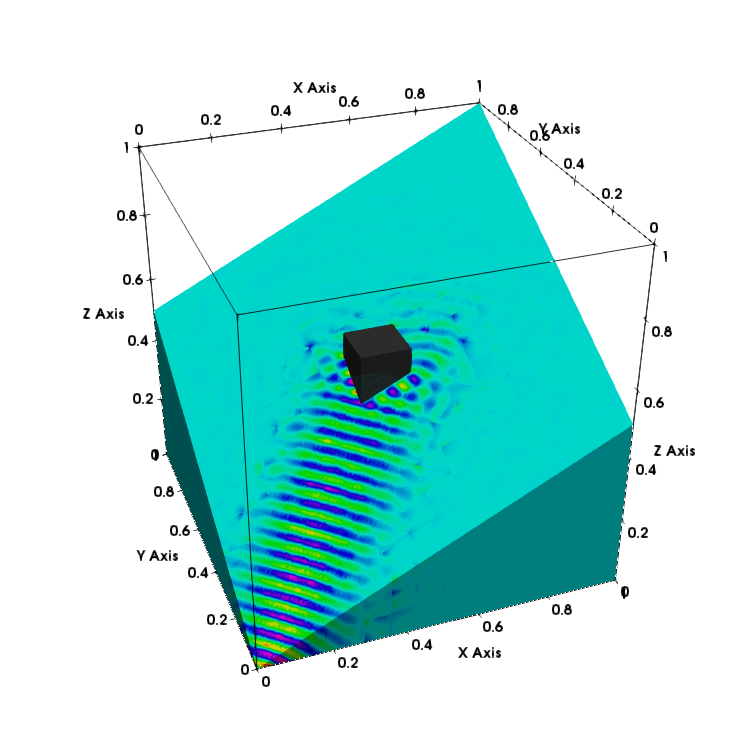}
      \caption{\hspace{20pt}Refinement 11}
\end{subfigure}
\begin{subfigure}[b]{0.24\textwidth}
      \centering
      \includegraphics[trim=80 0 80 0,clip,width=1.0\linewidth]{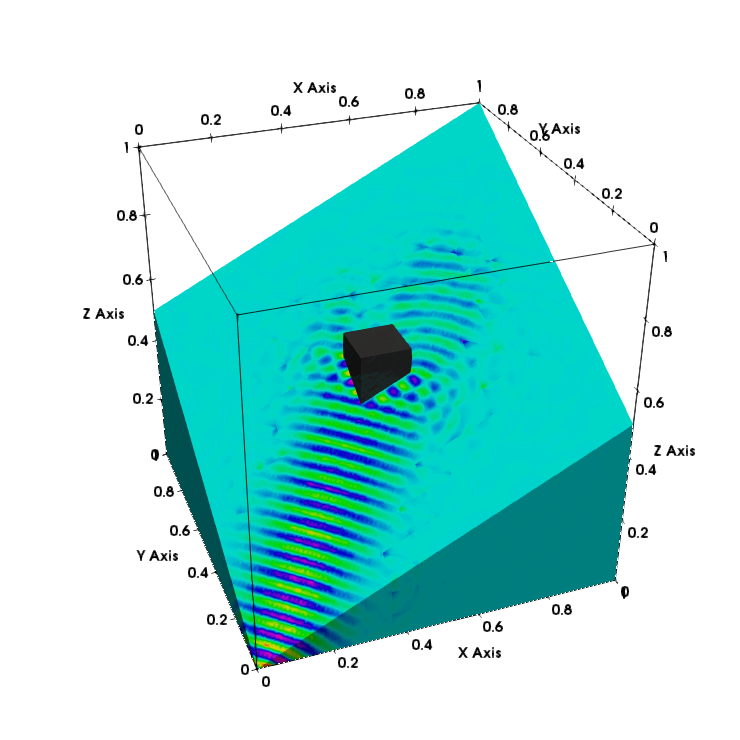}
      \caption{\hspace{20pt}Refinement 13}
\end{subfigure}
\begin{subfigure}[b]{0.24\textwidth}
      \centering
      \includegraphics[trim=80 0 80 0,clip,width=1.0\linewidth]{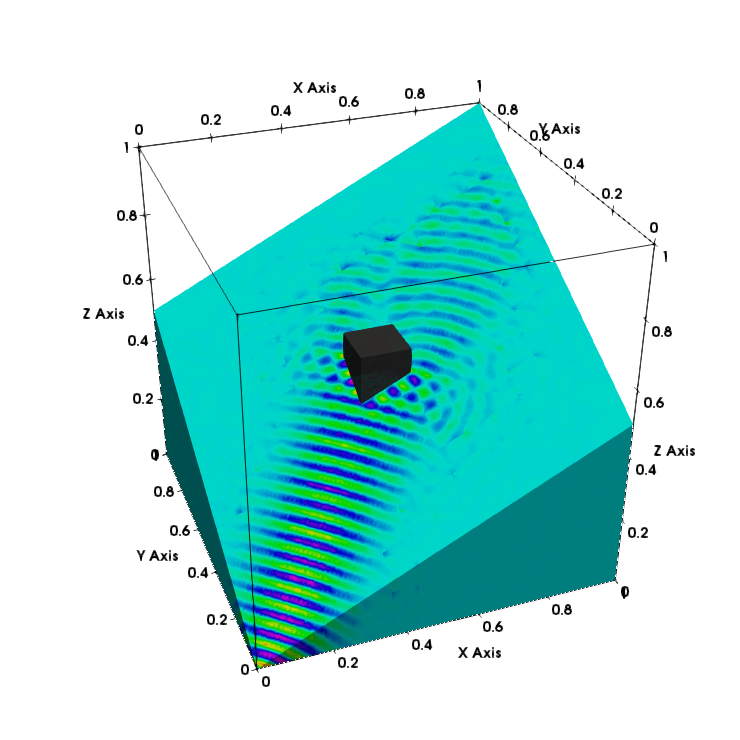}
      \caption{\hspace{20pt}Refinement 15}
\end{subfigure}
%

\begin{subfigure}[b]{0.24\textwidth}
      \centering
      \includegraphics[trim=80 0 80 0,clip,width=1.0\linewidth]{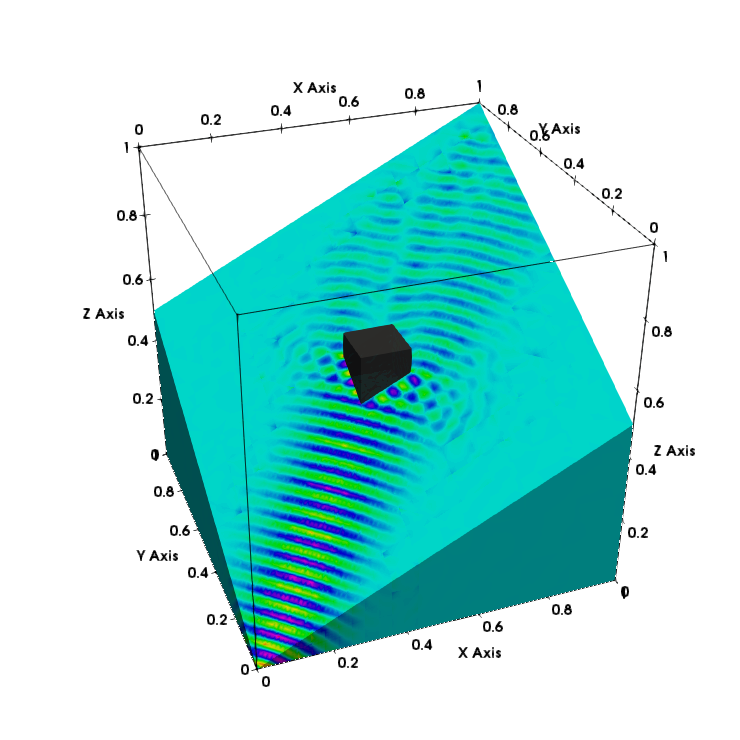}
      \caption{\hspace{20pt}Refinement 17}
\end{subfigure}
\begin{subfigure}[b]{0.24\textwidth}
      \centering
      \includegraphics[trim=80 0 80 0,clip,width=1.0\linewidth]{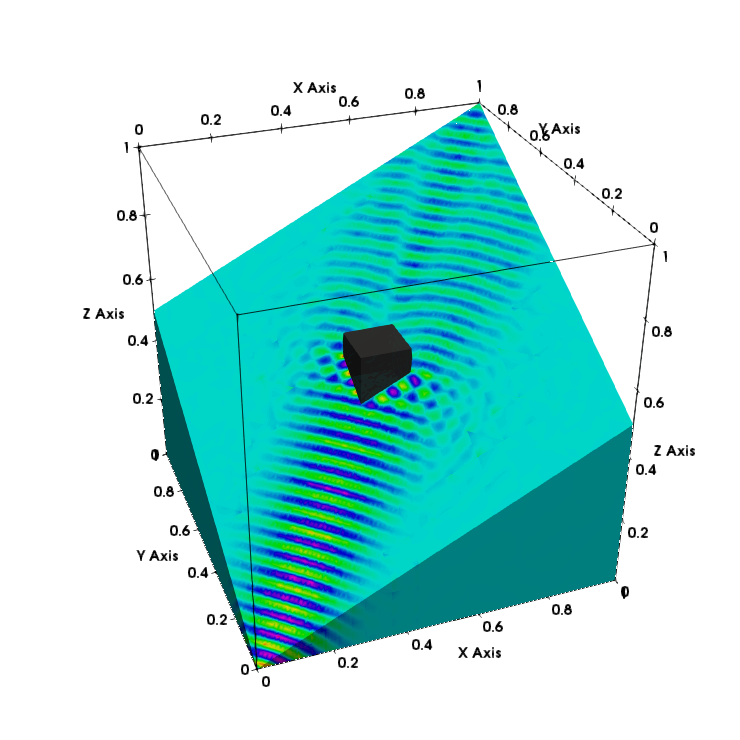}
      \caption{\hspace{20pt}Refinement 19}
\end{subfigure}
\begin{subfigure}[b]{0.24\textwidth}
      \centering
      \includegraphics[trim=80 0 80 0,clip,width=1.0\linewidth]{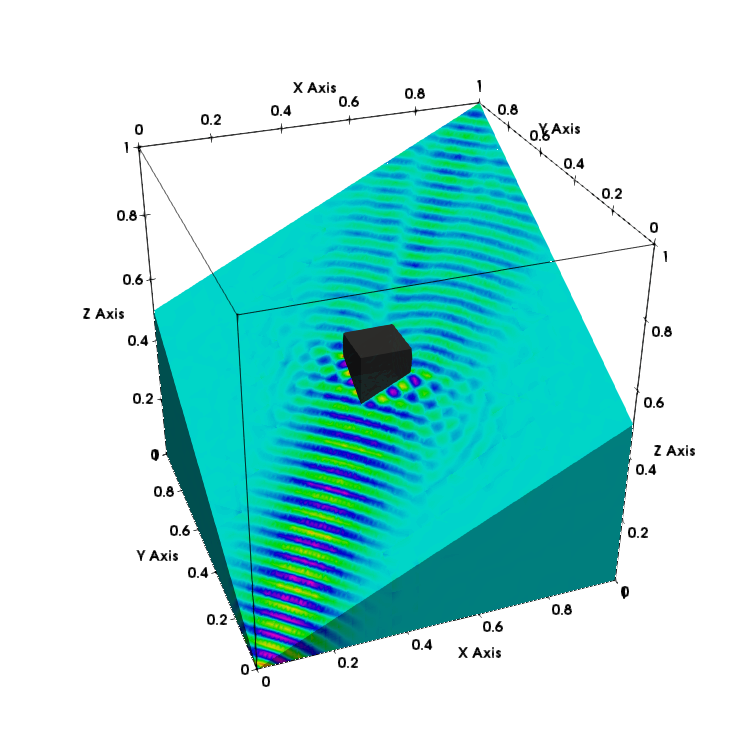}
      \caption{\hspace{20pt}Refinement 21}
\end{subfigure}
\begin{subfigure}[b]{0.24\textwidth}
      \centering
      \includegraphics[trim=80 0 80 0,clip,width=1.0\linewidth]{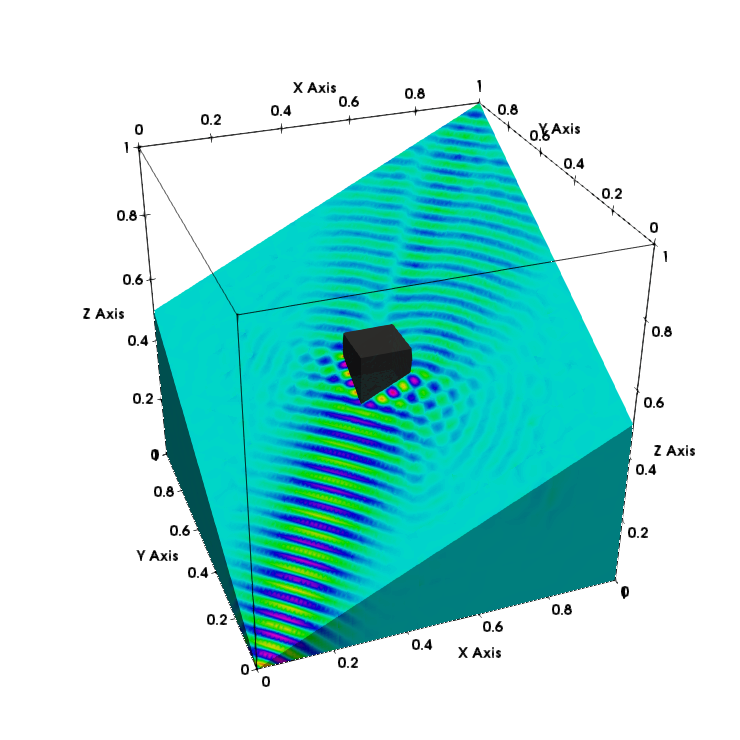}
      \caption{\hspace{20pt}Refinement 23}
\end{subfigure}
\end{figure}

\setcounter{figure}{4}

\begin{figure}[H]
      \centering
      \includegraphics[width=0.5\linewidth]{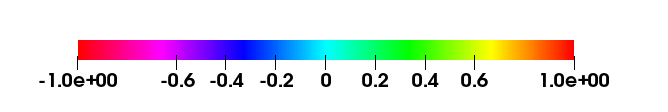}
         \caption{Real part of the numerical solution of the x-component of the electric field. Notice how the DPG adaptive technology refines only in regions of the domain where there is wave activity.}\label{fig:max_efield}
\end{figure}

%% file: results_acoustics.tex
%
%
\subsection{Linear acoustics equations}\label{sec:lin_acoustics_pml}
We continue with three dimensional examples for the linear acoustics problem. Recall the time harmonic form of the acoustics equations, 
\begin{equation*}
   \left\{
      \begin{aligned}
         i\omega p + \text{div}\,u & = f, \\
         i\omega u + \nabla p & = 0,
      \end{aligned}
   \right.
\end{equation*}
where $\omega$ is the angular frequency, $p$ is the pressure and $u$ is the velocity.

\medskip
\subsubsection{Scattering of a plane wave from a sphere}
\begin{figure}[H]
\begin{subfigure}[b]{0.35\textwidth}
\includegraphics[width=1.0\linewidth]{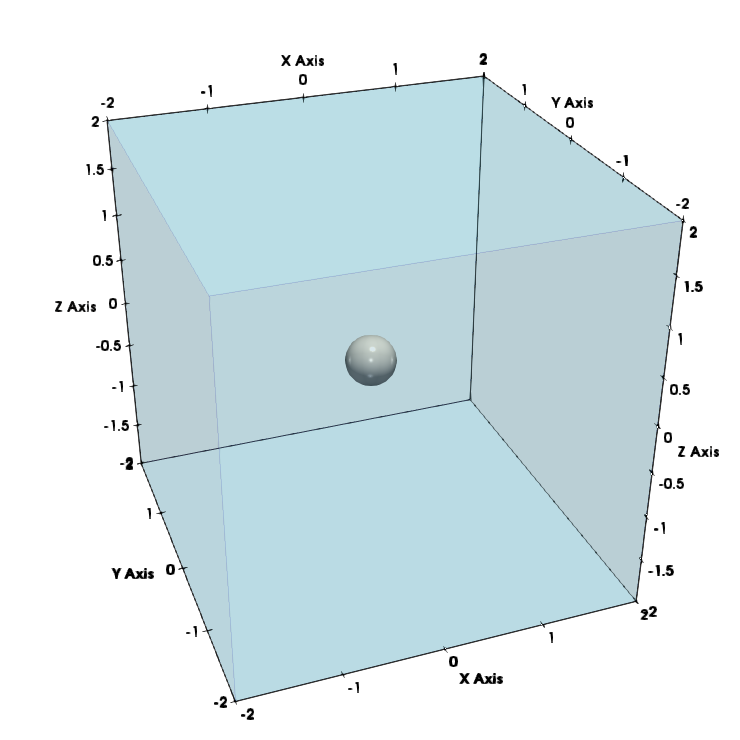}
\caption{Domain}\label{fig:sphere_domain}
\end{subfigure}
\begin{subfigure}[b]{0.35\textwidth}
\includegraphics[trim=50 0 50 0,clip,width=1.0\linewidth]{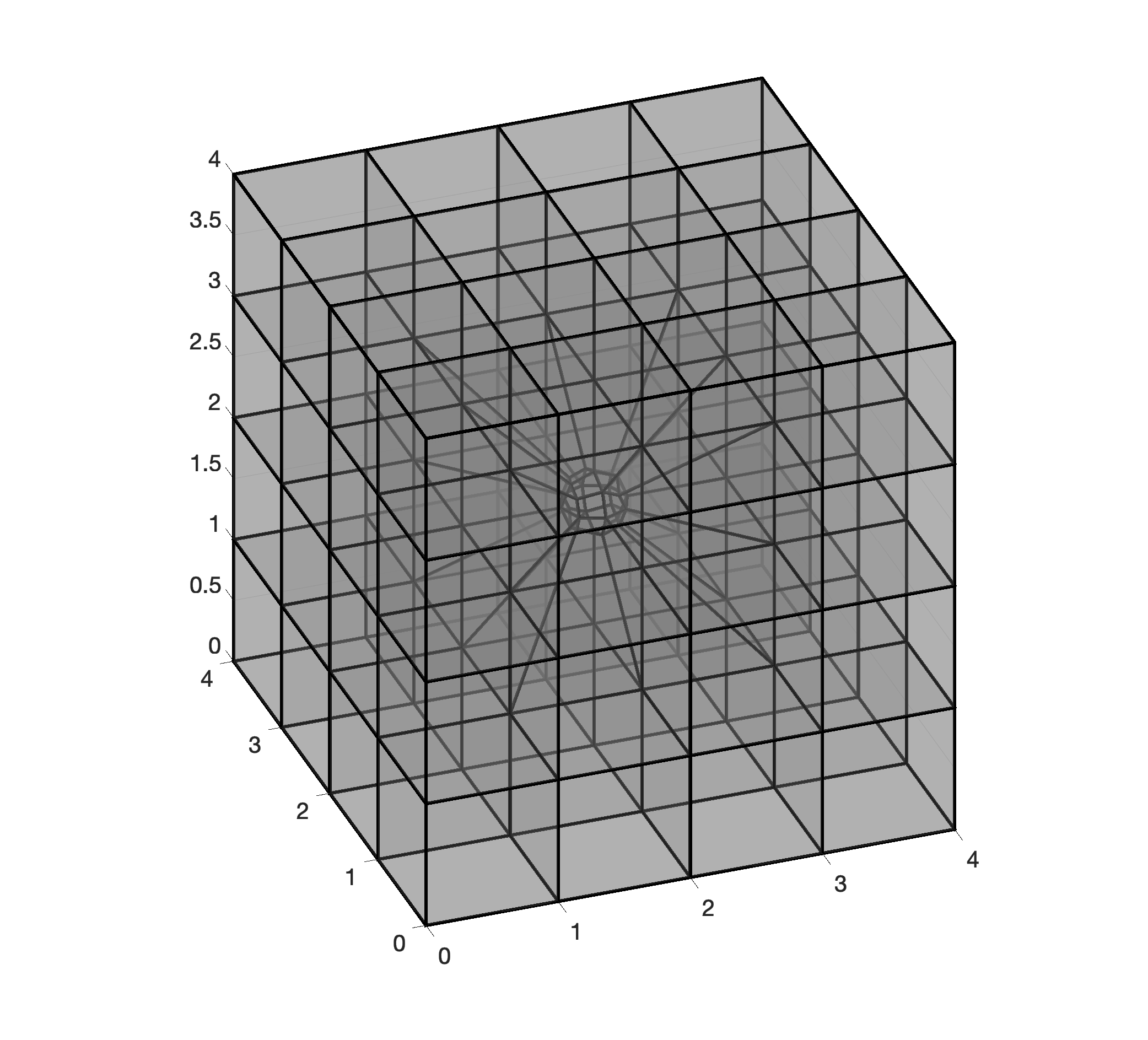}
\caption{Hexahedral mesh}\label{fig:sphere_mesh}
\end{subfigure}
\caption{Computational domain with a spherical scatterer and initial mesh}
\vspace{30pt}
\end{figure}

For our first example we consider the simulation of the scattered wave when a plane wave hits a rigid sphere (see \cref{fig:sphere_domain}). The domain is meshed\footnote{The authors would like to thank Brendan Keith for all his help on constructing the geometry files for this simulation.} by hexahedra as shown in \cref{fig:sphere_mesh} and the spherical surfaces are approximated using transfinite interpolation \cite{Demkbook2}. The incident plane wave is traveling in the direction $(1,1,1)$ with frequency $\omega = 35 \pi$. This corresponds to approximately 30 wavelengths inside the domain. 
The simulation is driven by hard boundary conditions on the spherical boundary ($\Gamma_{\text{sph}}$),
\begin{equation*}
u \cdot n = g = -u_{\text{inc}} \,\, \text{ on }\, \Gamma_{\text{sph}}
\end{equation*}
where $u_{\text{inc}}$ is the incident plane wave. The computational domain is truncated by a homogeneous impedance condition on the outer boundary of the cube $\Gamma_{\text{cube}}$.
\begin{equation*}
p = u \cdot n \,\, \text{ on }\, \Gamma_{\text{cube}}
\end{equation*}
We solve the problem using our multigrid technology in the uniform refinement setting. The fine grid consists of 40960 hexahedra of quartic polynomial order. This results in a linear system of approximately 14 million degrees of freedom. The coarse grid is constructed by two $h$--coarsening steps of the fine grid. It consists of 640 quartic hexahedra and the linear system has size of approximately 200 thousands degrees of freedom. The Conjugate Gradient algorithm starts with zero initial guess and terminates when the residual is less than $10^{-5}$. We use total of 10 smoothing iterations at each level and the smoother relaxation parameter is chosen to be $\theta = 0.2$. In this setting the preconditioned Conjugate Gradient algorithm converges in 12 iterations. We emphasize that we attempted to solve the problem with a sparse multi--frontal solver but it was not possible because of high memory requirements ($\approx 350Gb$). The solution is shown in \cref{fig:sphere_pressure} below. Note that we show the solution on the part of the domain below the plane defined by the point $(0.5,0.5,0.5)$ and the normal vector $(-1,-1,2)$.
\begin{figure}[H]
      \centering
      \includegraphics[width=0.5\linewidth]{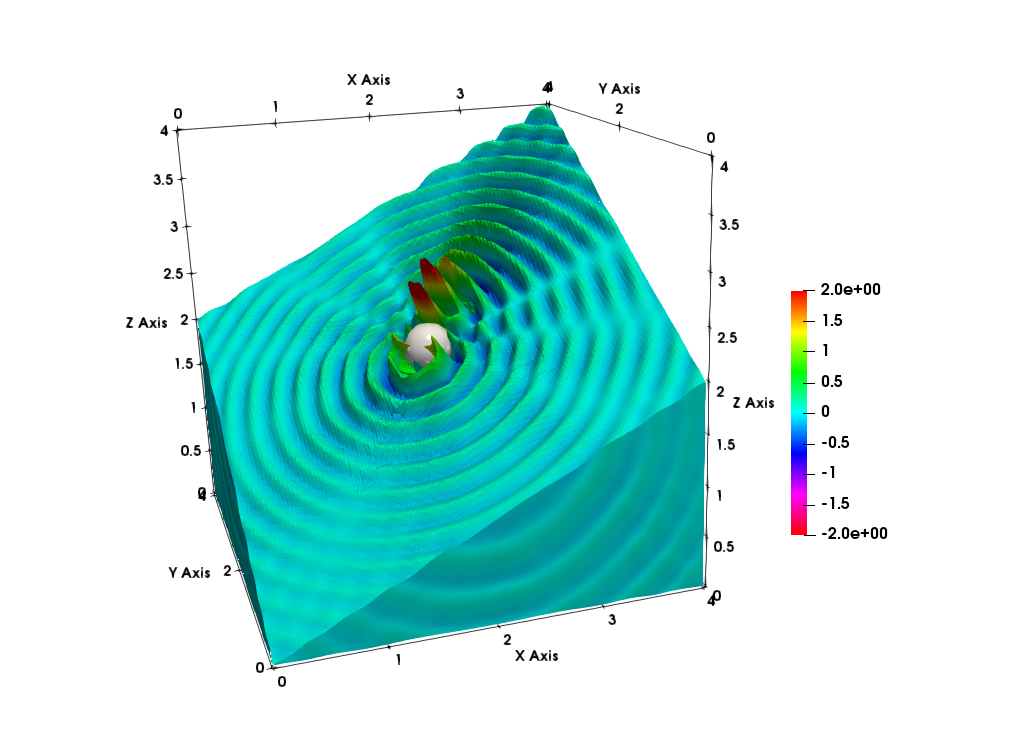}
      \caption{Scattered wave: real part of pressure}\label{fig:sphere_pressure}
\end{figure}

\subsubsection{Scattering of a plane wave from a cube}\label{sec: plane_scat_sec}
Our second acoustics example involves scattering of a plane wave by a cube. The computational domain is $\Omega = (\frac{1}{7},\frac{6}{7})^3\backslash ((\frac{3}{7},\frac{4}{7})^3$. The incident wave has low to medium frequency of $\omega = 16\pi$, and has direction of propagation $(1,1,0)$. The domain is truncated by a \emph{Perfectly Matched Layer} (PML) region of length $L=\frac{1}{7}$ in each direction. The construction of the PML is based on the work described in \cite{astaneh2018perfectly}. Additional details on how the DPG formulation is modified inside the PML region are given in \cite[Appendix D]{petrides2019adaptive}. 

We start the simulation with a mesh consisting of 342 cubes of side size $h=\frac{1}{7}$ and polynomial order $p=3$. Note that the initial mesh is fine enough to control the best approximation and the pollution effect. However, the singularities on the scatterer and the exponential decay of the wave in the PML region are not resolved and therefore the quality of the solution is not good. Starting automatic $h$--adaptivity, we anticipate that the DPG method will perform extensive refinements in order to resolve the singularities at the corners and edges of the cubic scatterer. Additionally, some refinements are expected to occur at the transition from the computational domain into the PML region.

\begin{figure}[H]
\begin{subfigure}[b]{0.49\textwidth}
\includegraphics[width=1.0\linewidth]{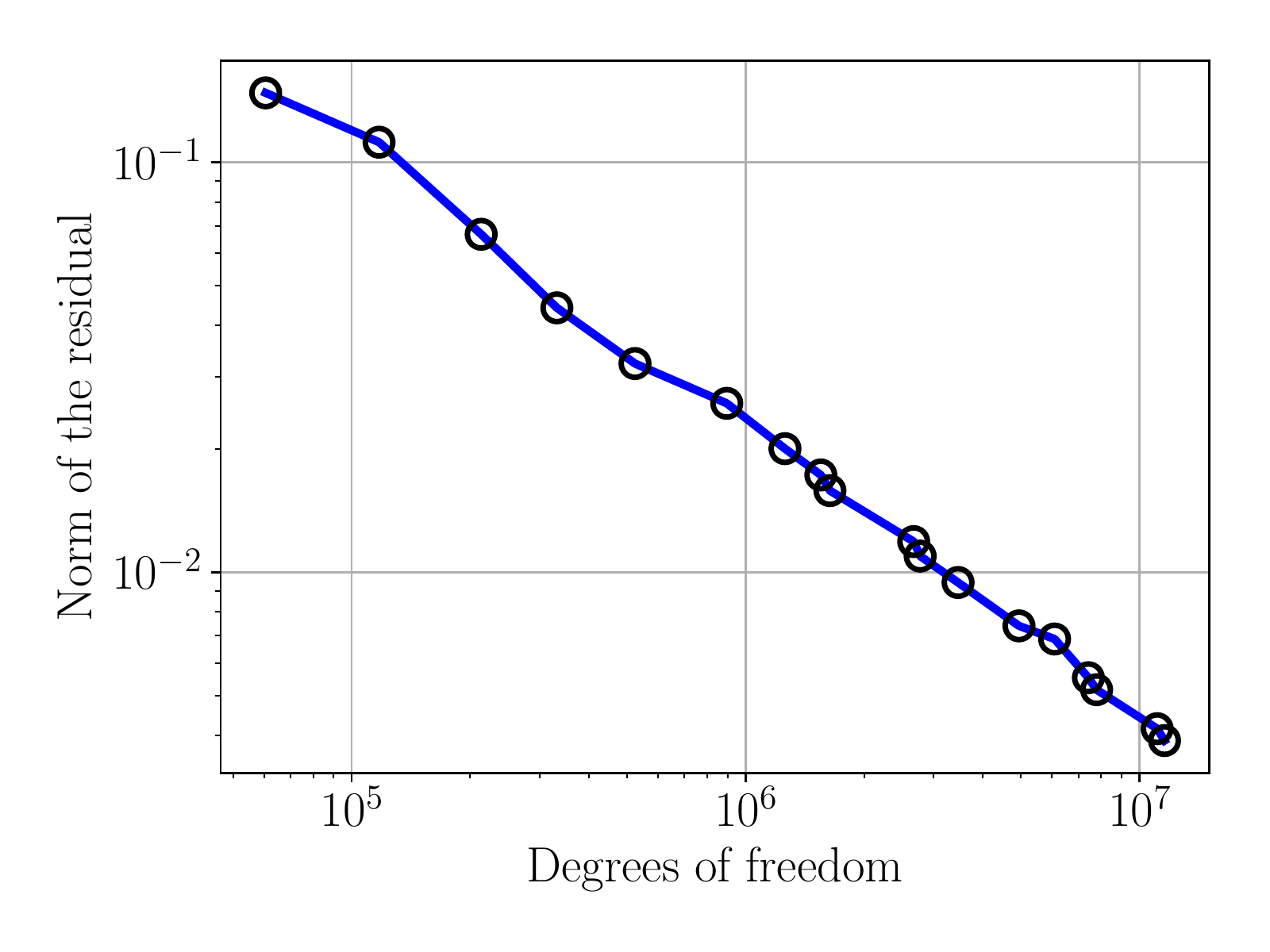}
\caption{Residual convergence}\label{fig:plane_res2}
\end{subfigure}
\begin{subfigure}[b]{0.49\textwidth}
\includegraphics[width=1.0\linewidth]{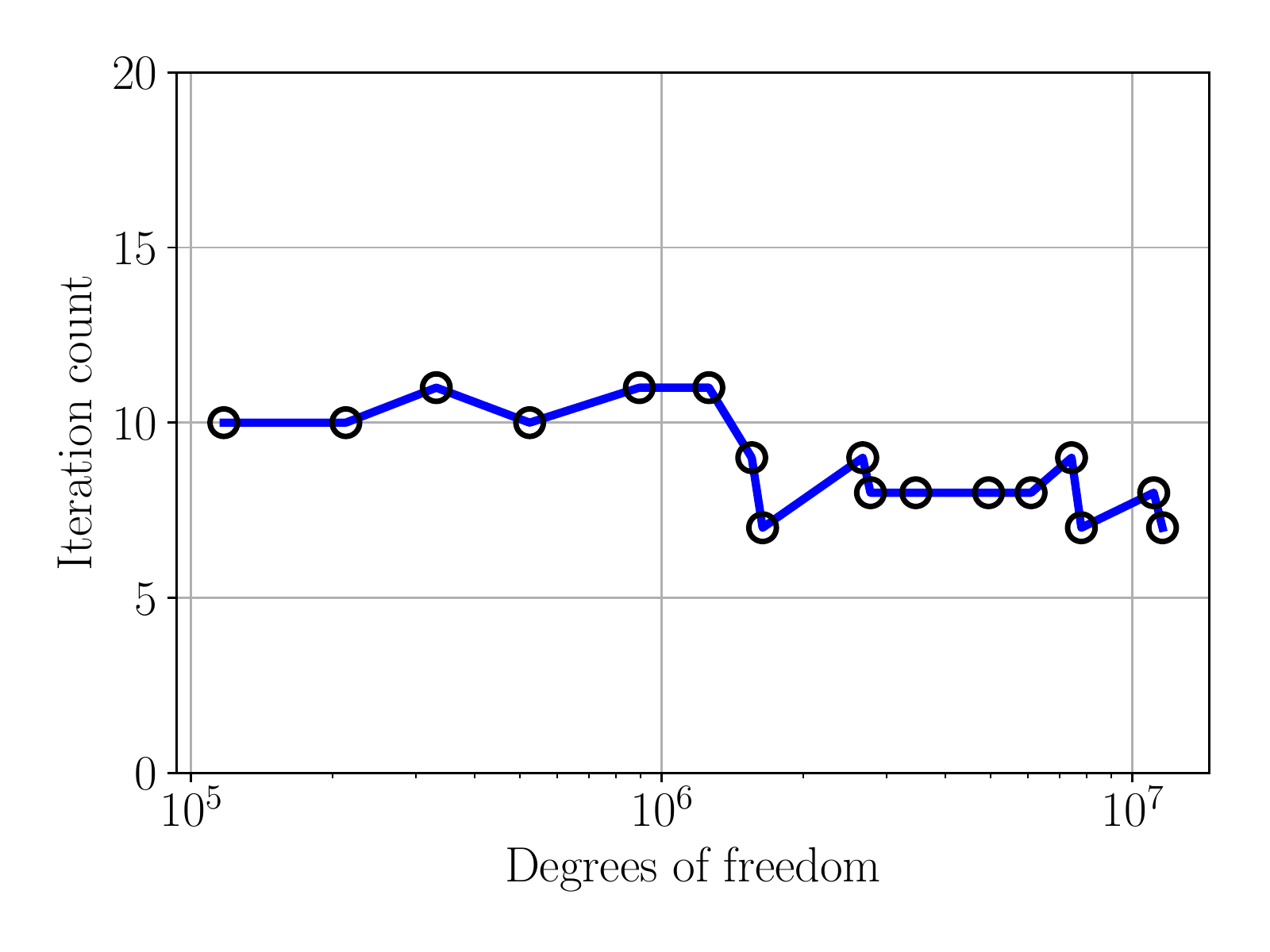}
\caption{PCG iteration count}\label{fig:pcg_plane2}
\end{subfigure}
\caption{Residual and preconditioned CG convergence. Plane wave scattering from a cube. Direction of propagation: (1,1,0).}\label{fig:conv2}
\end{figure}

\begin{figure}[H]
\captionsetup[subfigure]{labelformat=empty}
\begin{subfigure}[b]{0.22\textwidth}
\includegraphics[trim=150 0 150 0,clip,width=1.0\linewidth]{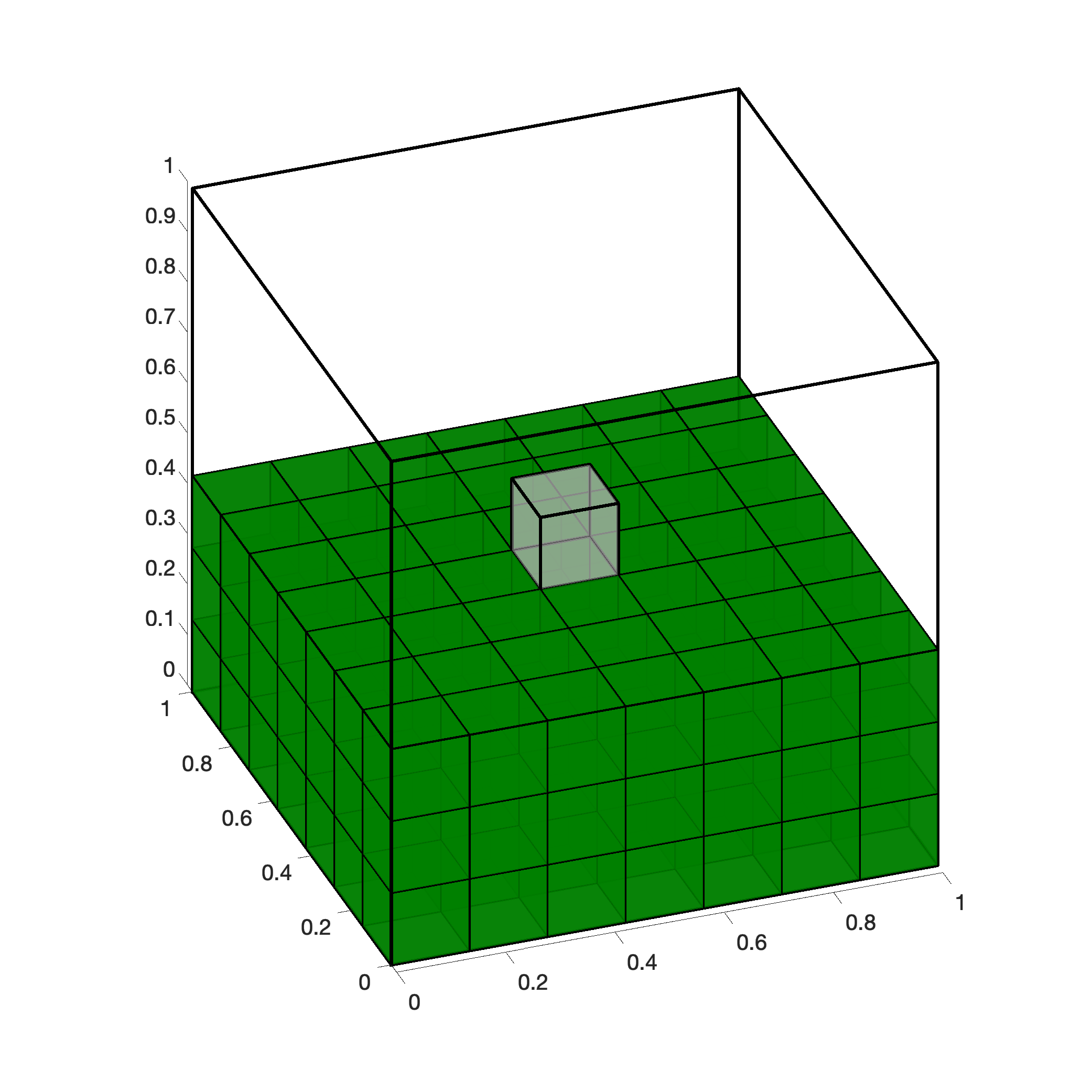}
\caption{Initial mesh}
\end{subfigure}
\begin{subfigure}[b]{0.22\textwidth}
\includegraphics[trim=150 0 150 0,clip,width=1.0\linewidth]{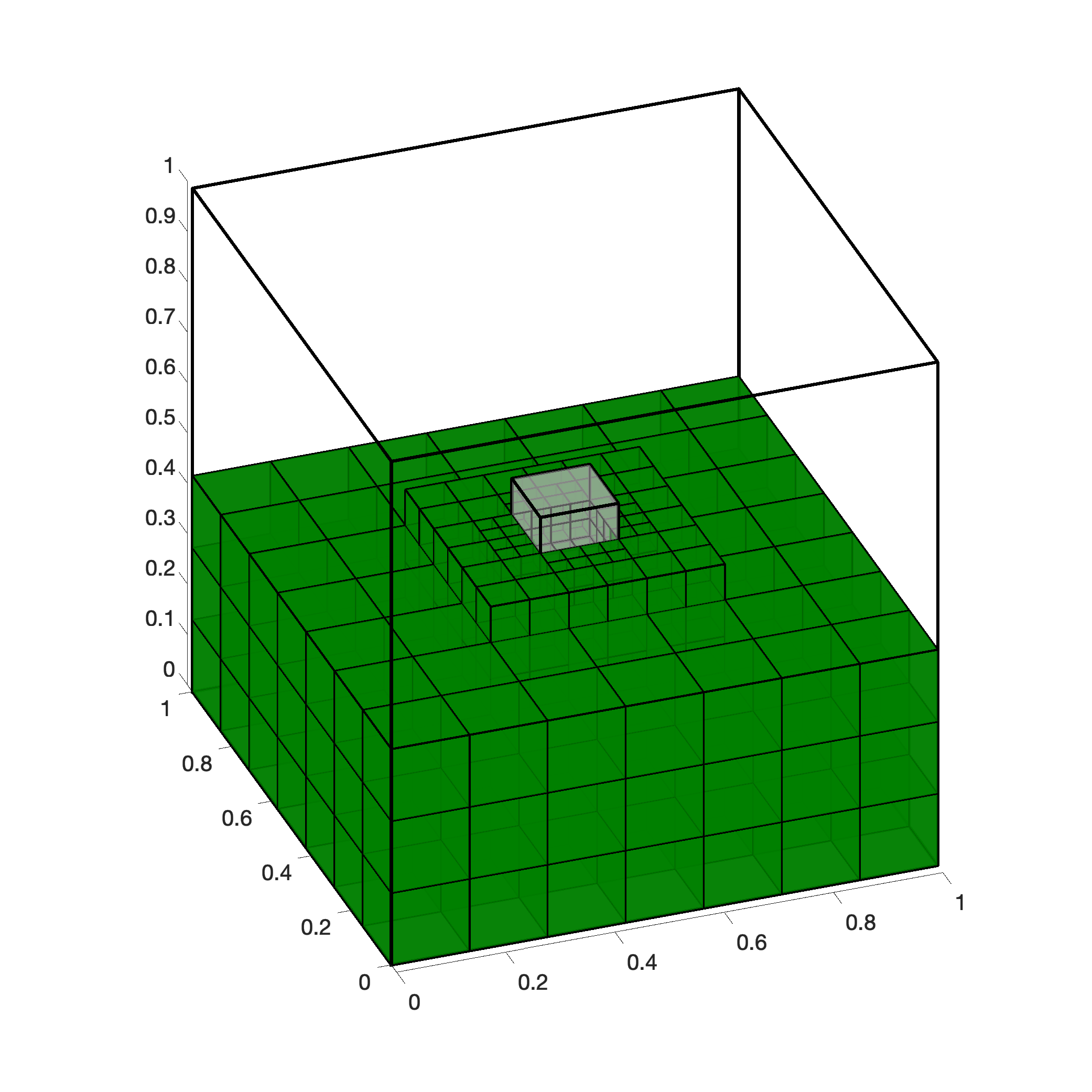}
\caption{Mesh 3}
\end{subfigure}
\begin{subfigure}[b]{0.22\textwidth}
\includegraphics[trim=150 0 150 0,clip,width=1.0\linewidth]{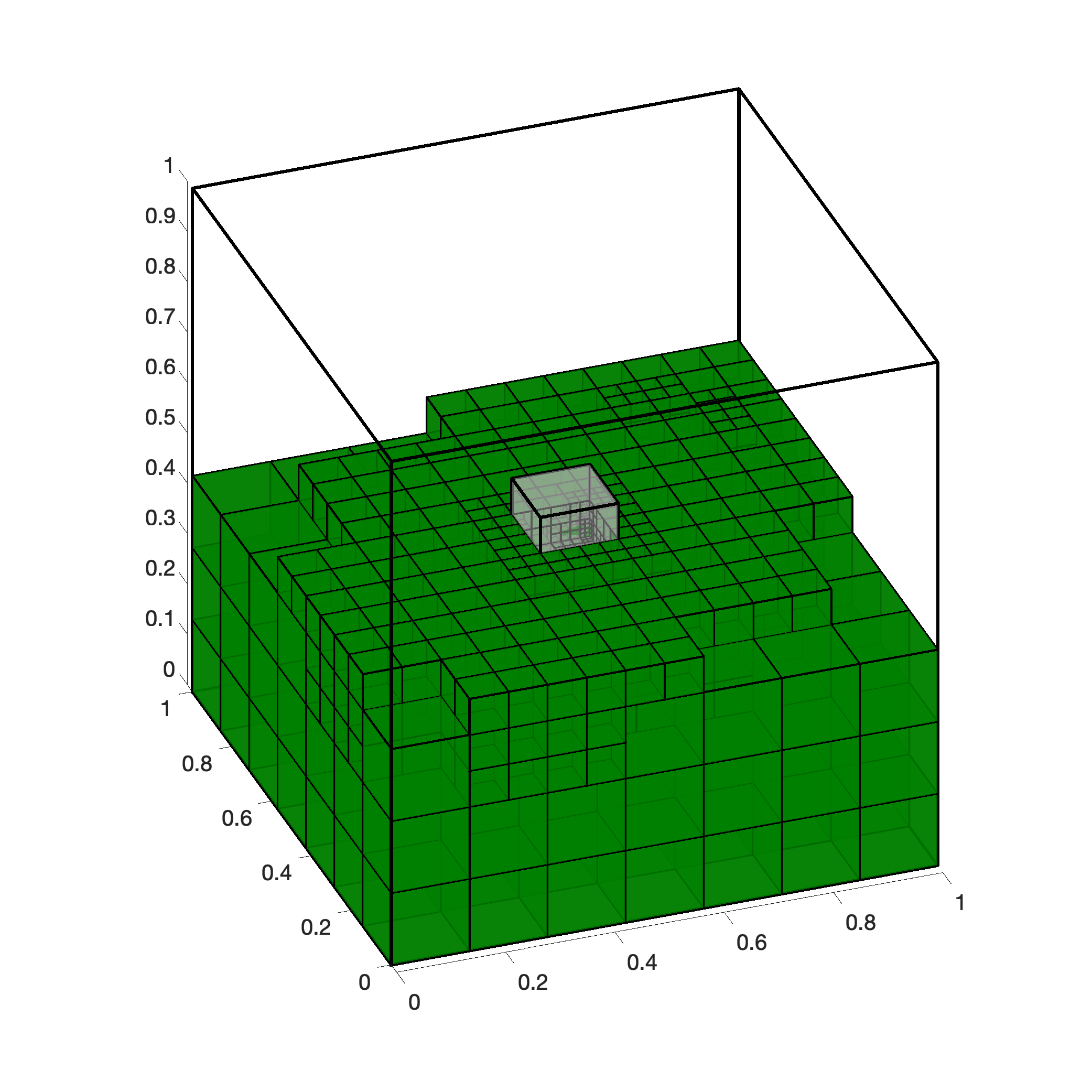}
\caption{Mesh 5}
\end{subfigure}
\begin{subfigure}[b]{0.22\textwidth}
\includegraphics[trim=150 0 150 0,clip,width=1.0\linewidth]{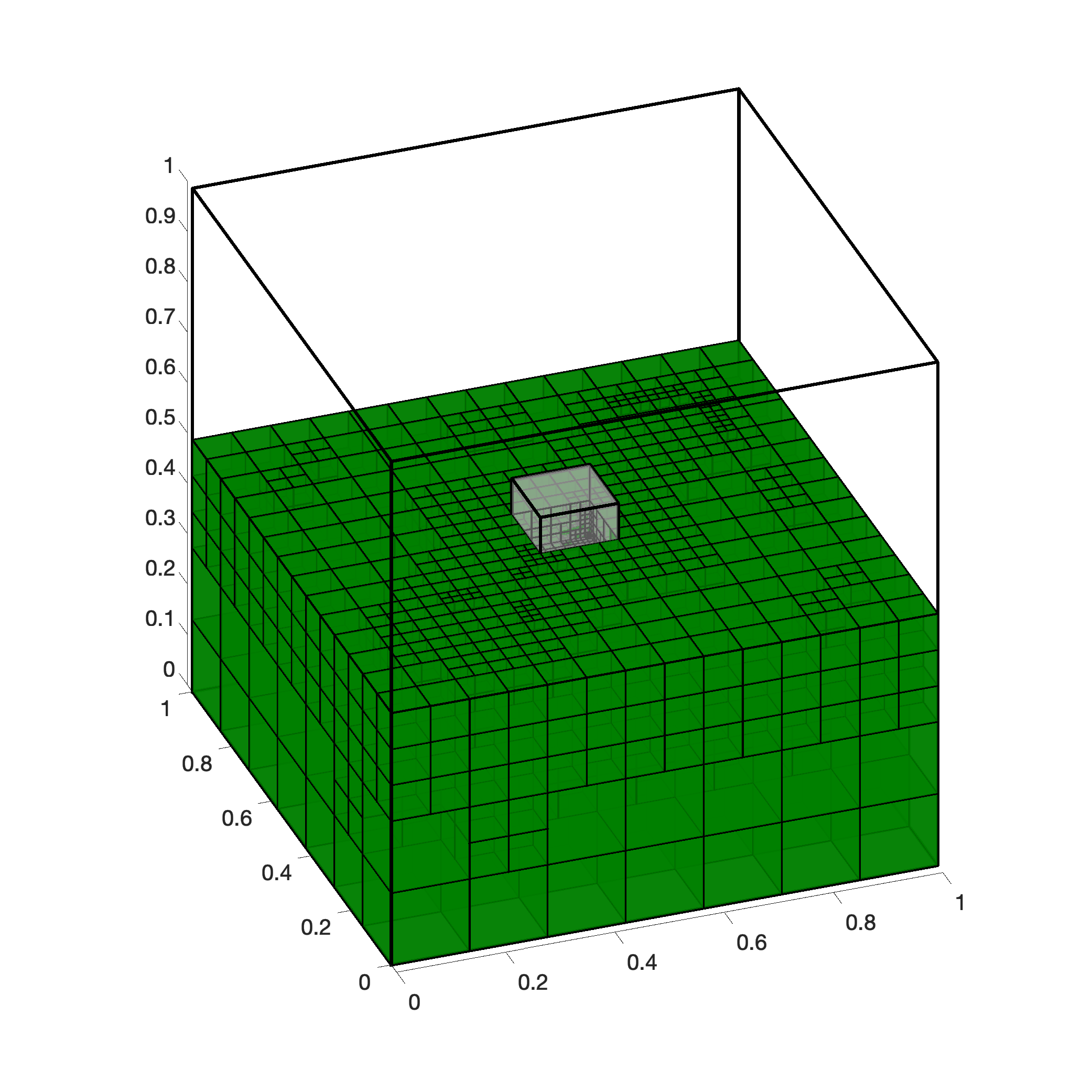}
\caption{Mesh 7}
\end{subfigure}
%
\begin{subfigure}[b]{0.22\textwidth}
\includegraphics[trim=150 0 150 0,clip,width=1.0\linewidth]{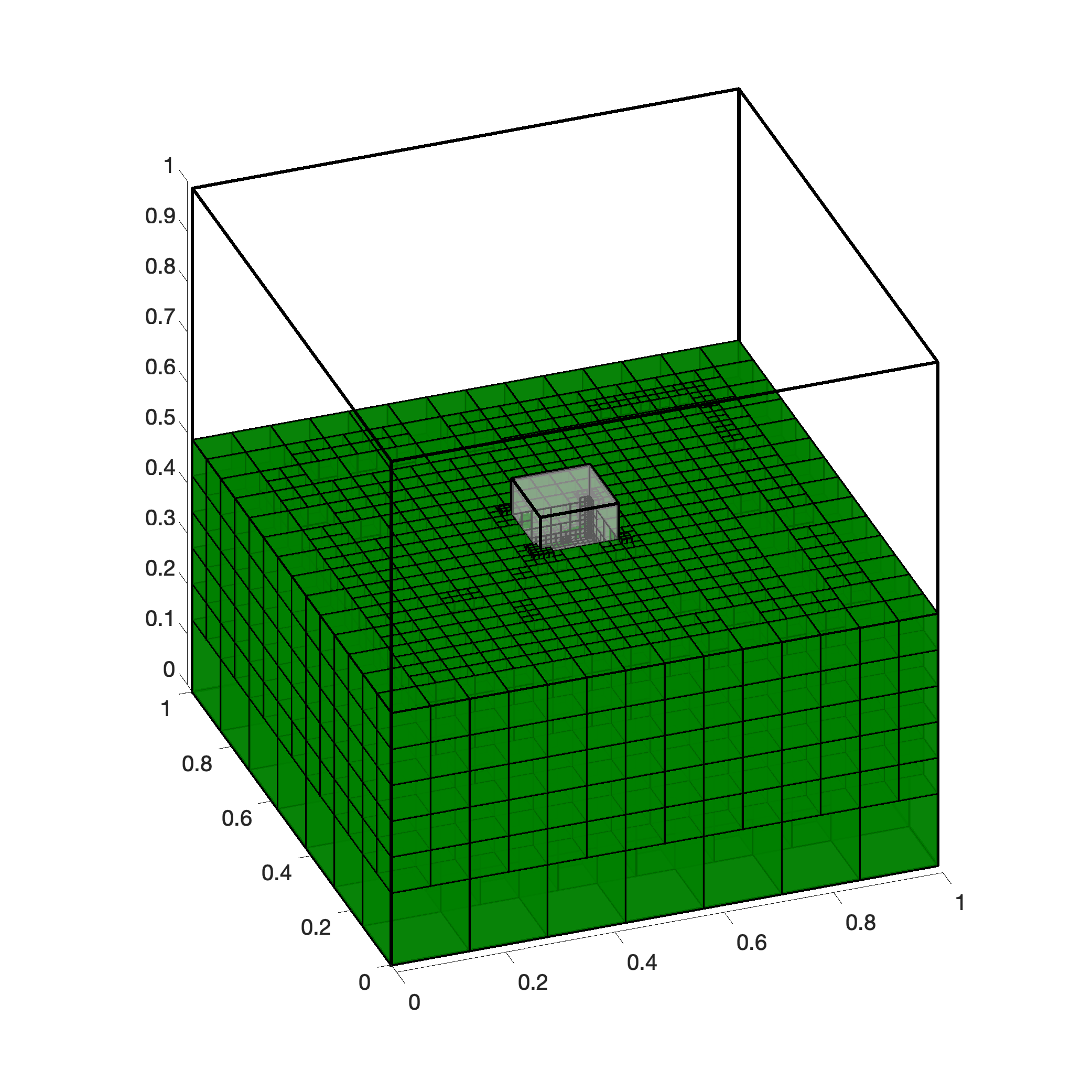}
\caption{Mesh 10}
\end{subfigure}
\begin{subfigure}[b]{0.22\textwidth}
\includegraphics[trim=150 0 150 0,clip,width=1.0\linewidth]{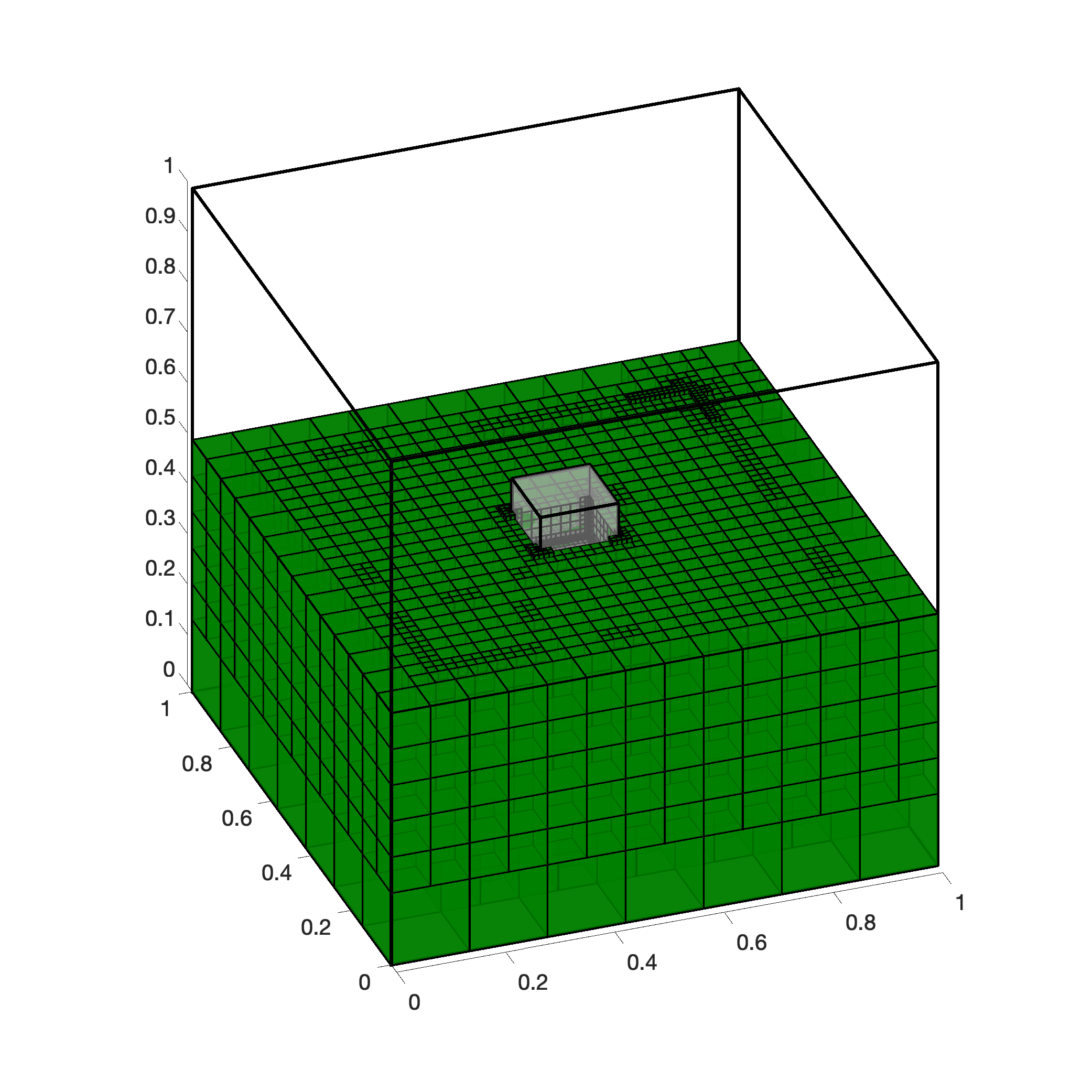}
\caption{Mesh 13}
\end{subfigure}
\begin{subfigure}[b]{0.22\textwidth}
\includegraphics[trim=150 0 150 0,clip,width=1.0\linewidth]{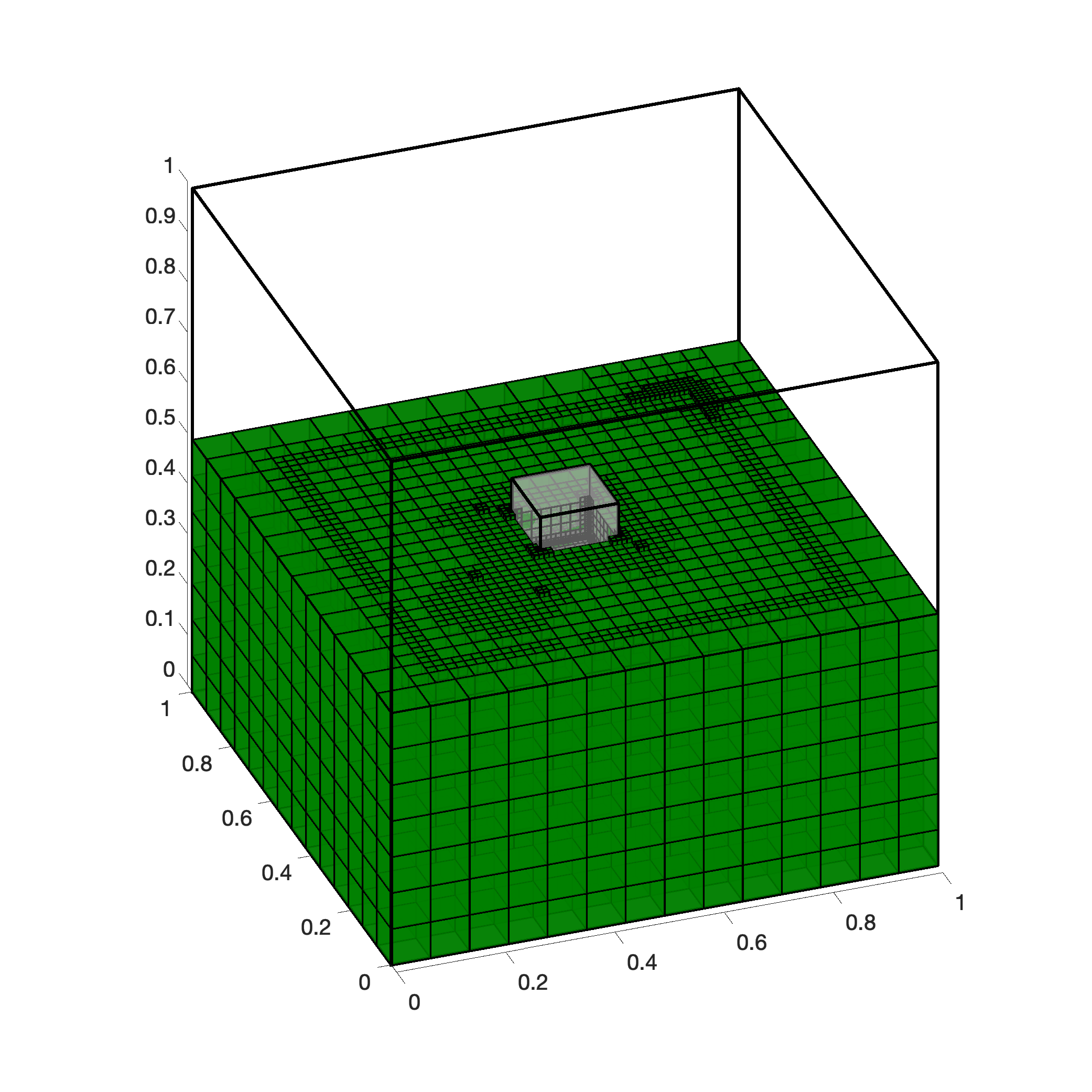}
\caption{Mesh 16}
\end{subfigure}
\begin{subfigure}[b]{0.22\textwidth}
\includegraphics[trim=150 0 150 0,clip,width=1.0\linewidth]{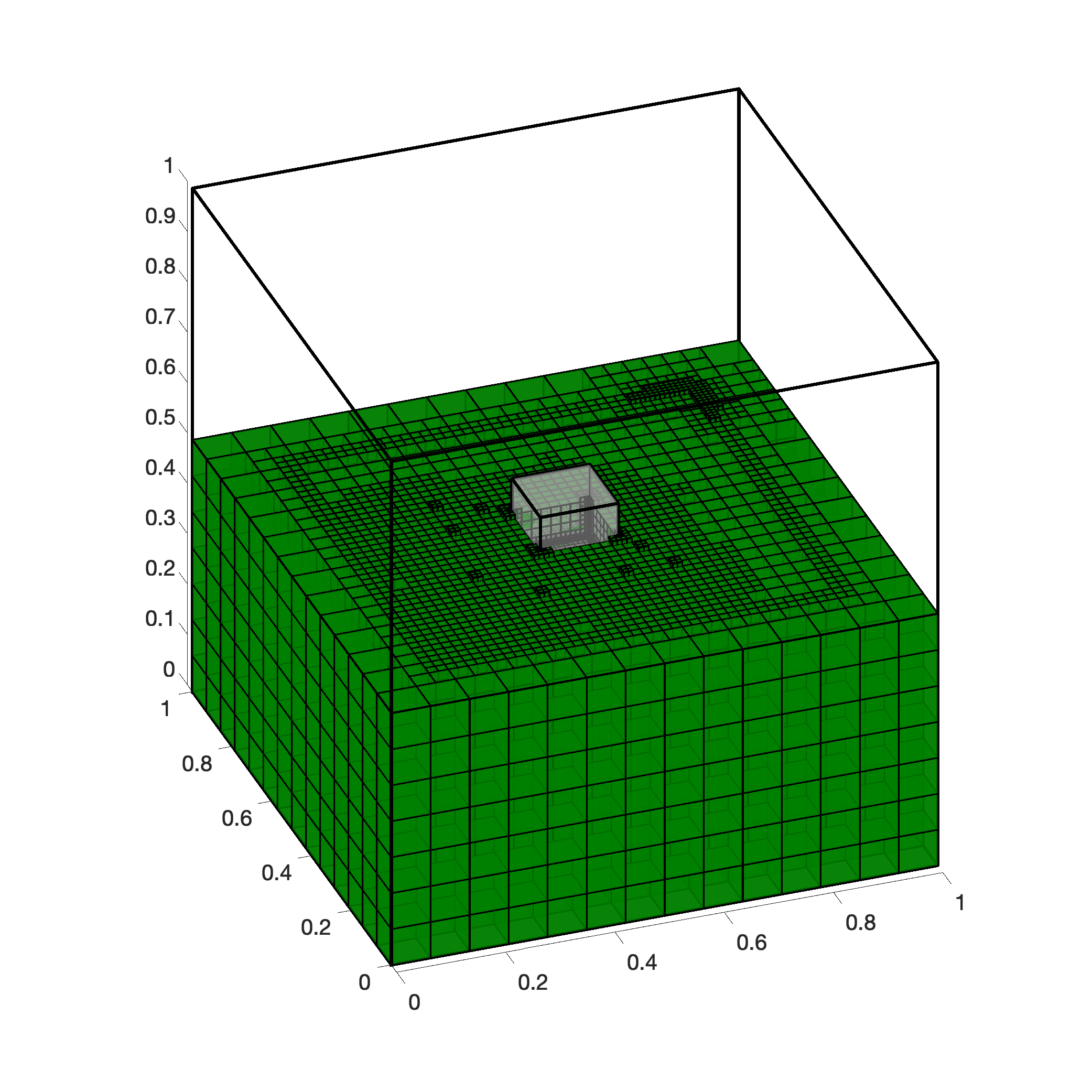}
\caption{Mesh 19}
\end{subfigure}
\caption{Wave propagating in the direction (1,1,0). Evolution of the $h$--adaptive mesh. As expected, a lot of refinements occur in the region close to the scatterer because the singularities have to be resolved. Additional refinements occur in the PML region.}
\label{fig:meshes}
\end{figure}
%
%
%
\begin{figure}[H]
\captionsetup[subfigure]{labelformat=empty}
\begin{subfigure}[b]{0.22\textwidth}
\includegraphics[trim=50 0 50 0,clip,width=1.0\linewidth]{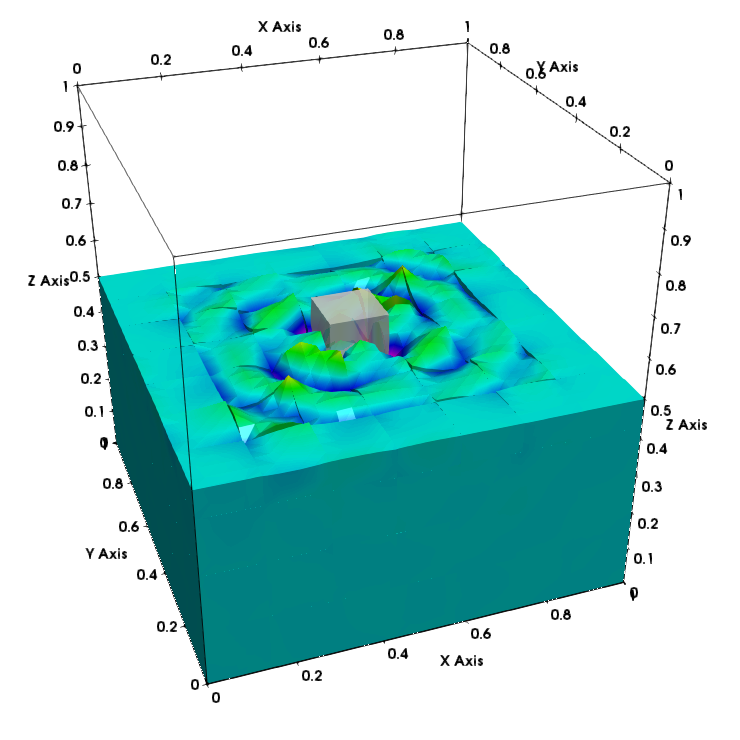}
\caption{Initial mesh}
\end{subfigure}
\begin{subfigure}[b]{0.22\textwidth}
\includegraphics[trim=50 0 50 0,clip,width=1.0\linewidth]{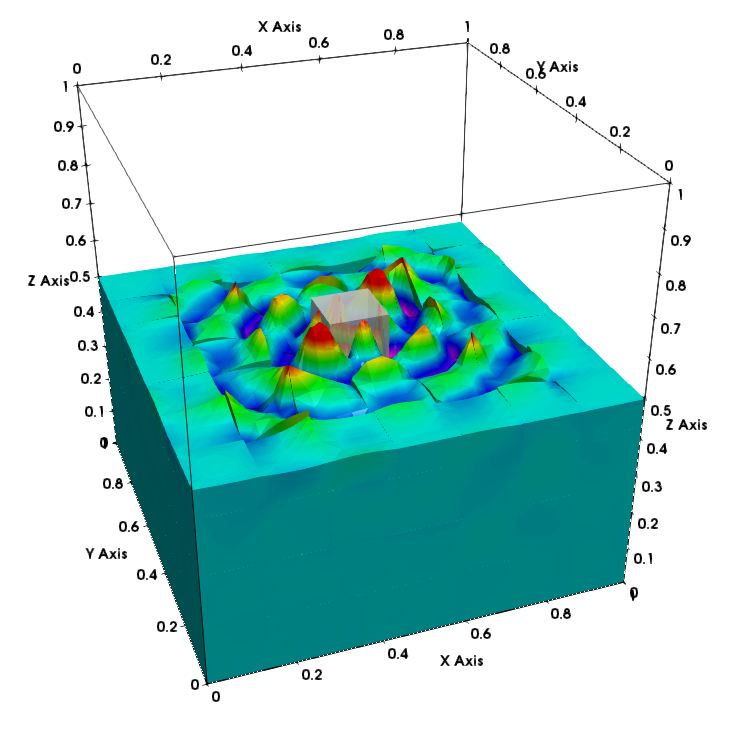}
\caption{Mesh 3}
\end{subfigure}
\begin{subfigure}[b]{0.22\textwidth}
\includegraphics[trim=50 0 50 0,clip,width=1.0\linewidth]{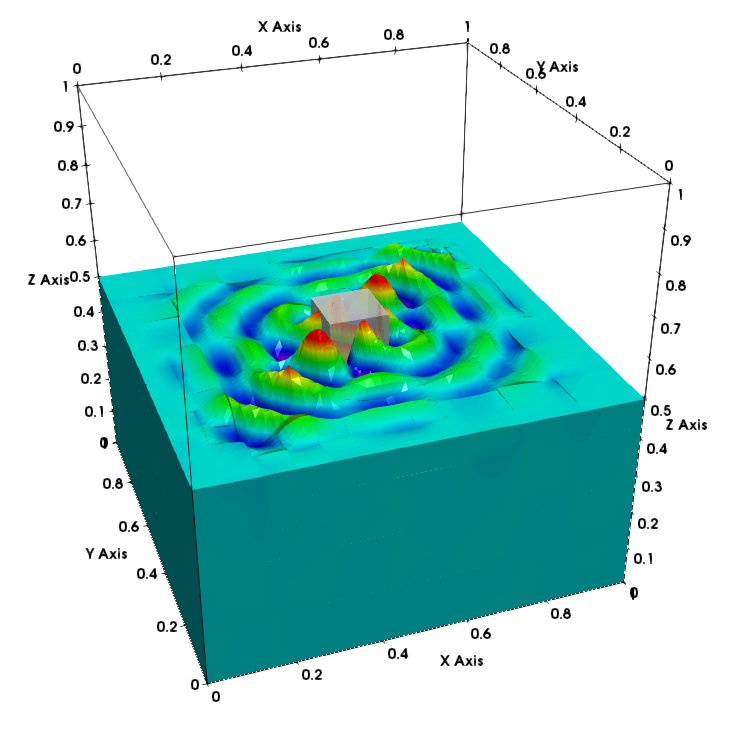}
\caption{Mesh 5}
\end{subfigure}
\begin{subfigure}[b]{0.22\textwidth}
\includegraphics[trim=50 0 50 0,clip,width=1.0\linewidth]{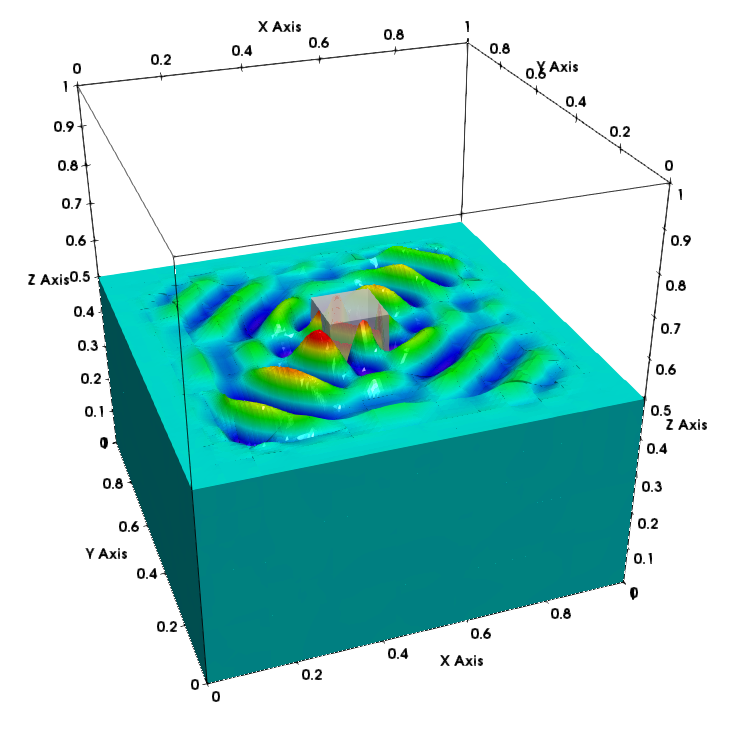}
\caption{Mesh 7}
\end{subfigure}
\begin{subfigure}[b]{0.22\textwidth}
\includegraphics[trim=50 0 50 0,clip,width=1.0\linewidth]{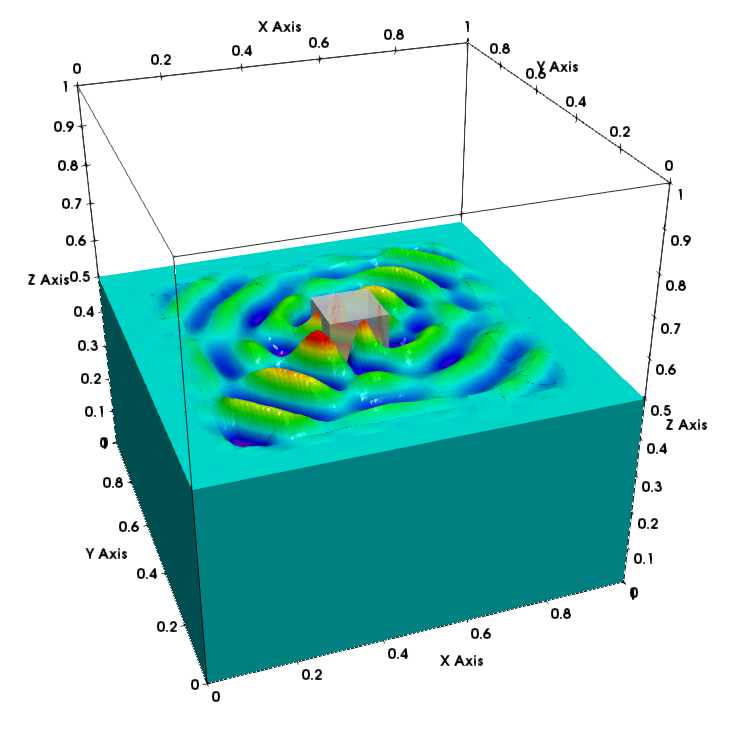}
\caption{Mesh 10}
\end{subfigure}
\begin{subfigure}[b]{0.22\textwidth}
\includegraphics[trim=50 0 50 0,clip,width=1.0\linewidth]{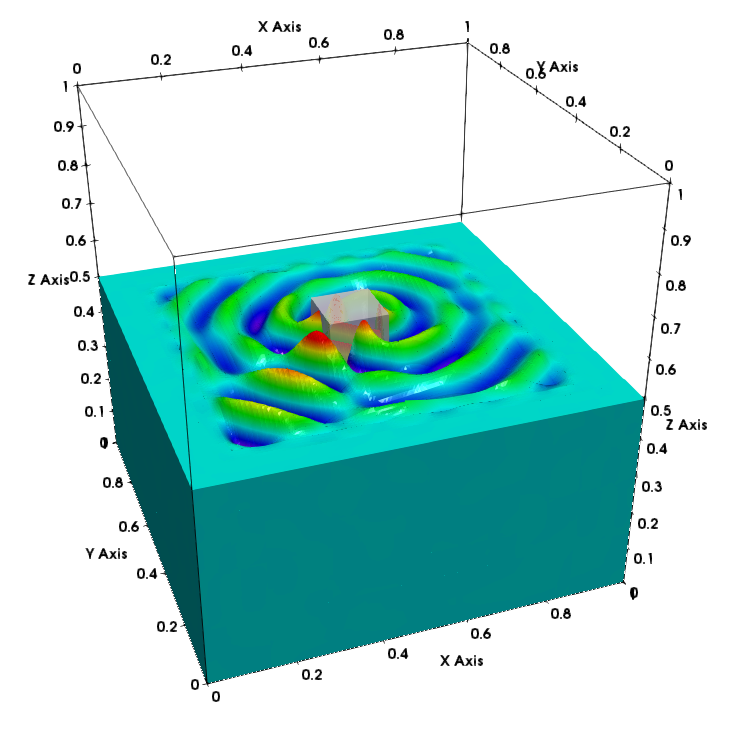}
\caption{Mesh 13}
\end{subfigure}
\begin{subfigure}[b]{0.22\textwidth}
\includegraphics[trim=50 0 50 0,clip,width=1.0\linewidth]{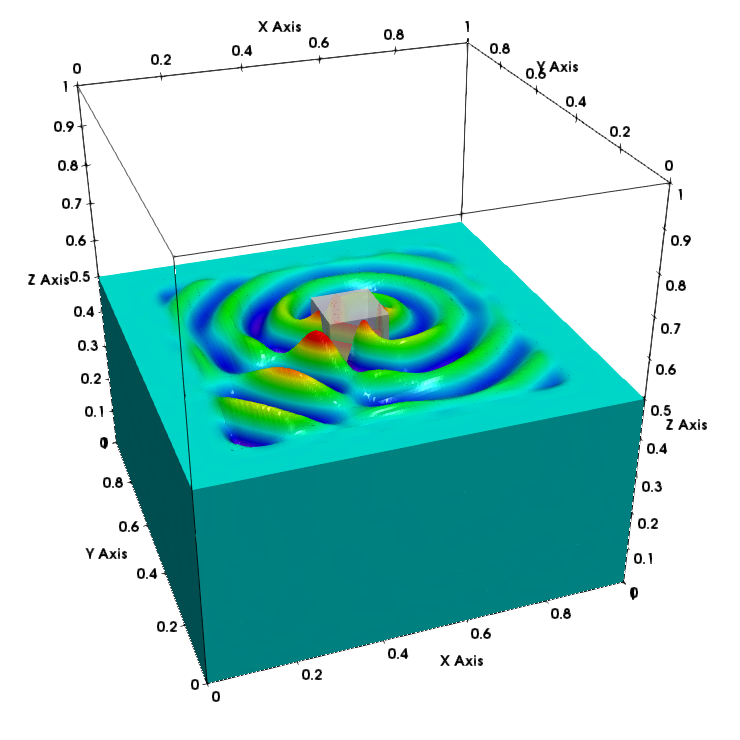}
\caption{Mesh 16}
\end{subfigure}
\begin{subfigure}[b]{0.22\textwidth}
\includegraphics[trim=50 0 50 0,clip,width=1.0\linewidth]{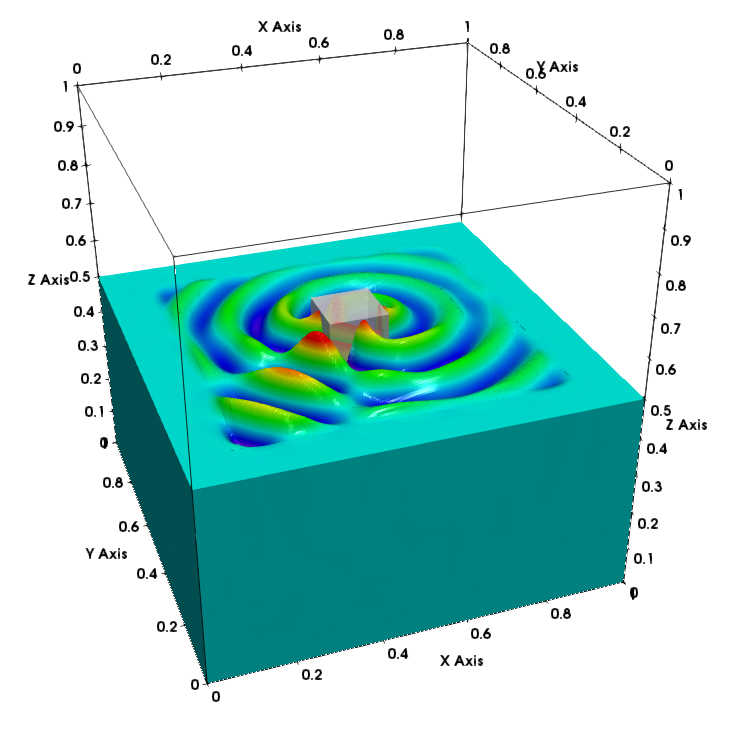}
\caption{Mesh 19}
\end{subfigure}
\label{fig:pressure}
\end{figure}

\setcounter{figure}{10}
\begin{figure}[H]
      \centering
      \includegraphics[width=0.5\linewidth]{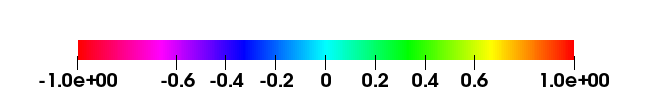}
      \caption{Wave propagating in the direction (1,1,0). Evolution of the solution. The solution rabidly decays in the PML region. Notice that the quality of the solution is affected by the resolution of the singularities}\label{fig:pressure2}
\end{figure}

The multigrid preconditioner setup is the same as in the previous examples. Here,we use at most 6 multigrid levels and the stopping criterion of the CG solver is set to $10^{-5}$. The sequence of meshes along with the solutions are shown in \cref{fig:meshes,fig:pressure}. Notice that the meshes and the solutions are shown only in the part of the domain where $z \le 0.5$. As expected, refinements occur close to the singularities and in the PML region. Observe how the solution changes as the singularities are resolved. Convergence results for the DPG residual and iteration counts for the preconditioned CG solver are given in \cref{fig:conv2}. As in the previous numerical examples the solver shows robust convergence throughout the adaptive refinement process. 

\medskip
\subsection{Scattering of a Gaussian beam from a cube}
Our last experiment involves the simulation of a high frequency Gaussian beam scattering from a cube. The computational domain and the PML region are defined as in \cref{sec: plane_scat_sec}. We run the simulation with angular frequency $\omega = 140\pi$. This corresponds to about 85 wavelengths inside the computational domain. We begin the simulation with a uniform mesh consisting of 342 cubes of polynomial order $p=3$. We then perform adaptive $hp$--refinements with the following strategy. An element inside the computational domain is $h$--refined up to the point where its maximum side size becomes smaller than half a wavelength. An element with size smaller than that is $p$--refined. An exception is made for the elements adjacent to the corners and edges of the cubic scatterer and the elements inside the PML region, where only $h$--refinements are allowed. The simulation is terminated when the DPG residual reduces by an order of magnitude.

The iterative solver setup is as follows. The initial iterate consists of the degrees of freedom corresponding to the prolongation of the solution on the previous mesh to the current mesh. We perform 10 smoothing steps at each multigrid level and we choose the relaxation parameter to be $\theta = 0.2$. Additional computational time is saved by omitting to smooth in areas of the domain where there is no wave activity (the local patch residual is close to zero). Finally, in this simulation we use at most eight multigrid levels and the CG iterations are terminated when the $l^2$--norm of the residual becomes less than $10^{-3}$.

In \cref{fig:beam_meshes,fig:beam_pressure} we show the sequence of meshes and the corresponding numerical solutions for the real part of the pressure respectively. For demonstration purposes a partial region of the mesh is shown, constructed by the union of the regions below three planes. These planes are defined by the point $(0.5,0.5,0.5)$ and the normal vectors $(-0.5,-0.5,1)$, $(0.5,-0.5,1)$ and $(-0.5,0.5,1)$. Notice that the singularities at the corners and edges of the cubic scatterer have to be resolved before the wave can propagate. When the wave reaches the PML region, additional refinements occur in order to capture the steep decay. Convergence results are presented in \cref{fig:beam_convergence}. Observe in \cref{fig:beam_res_convergence} that the norm of DPG residual, which drives the adaptive refinements, tends to zero only after the decay in the PML region is resolved. Lastly, iteration counts for the CG solver preconditioned with our multigrid technology are given in \cref{fig:beam_pcg_convergence}. Notice that the number of iterations remains under control throughout the adaptive process.  
\begin{figure}[H]
\captionsetup[subfigure]{labelformat=empty}
\begin{subfigure}[b]{0.24\textwidth}
\includegraphics[trim=150 0 150 0,clip,width=1.0\linewidth]{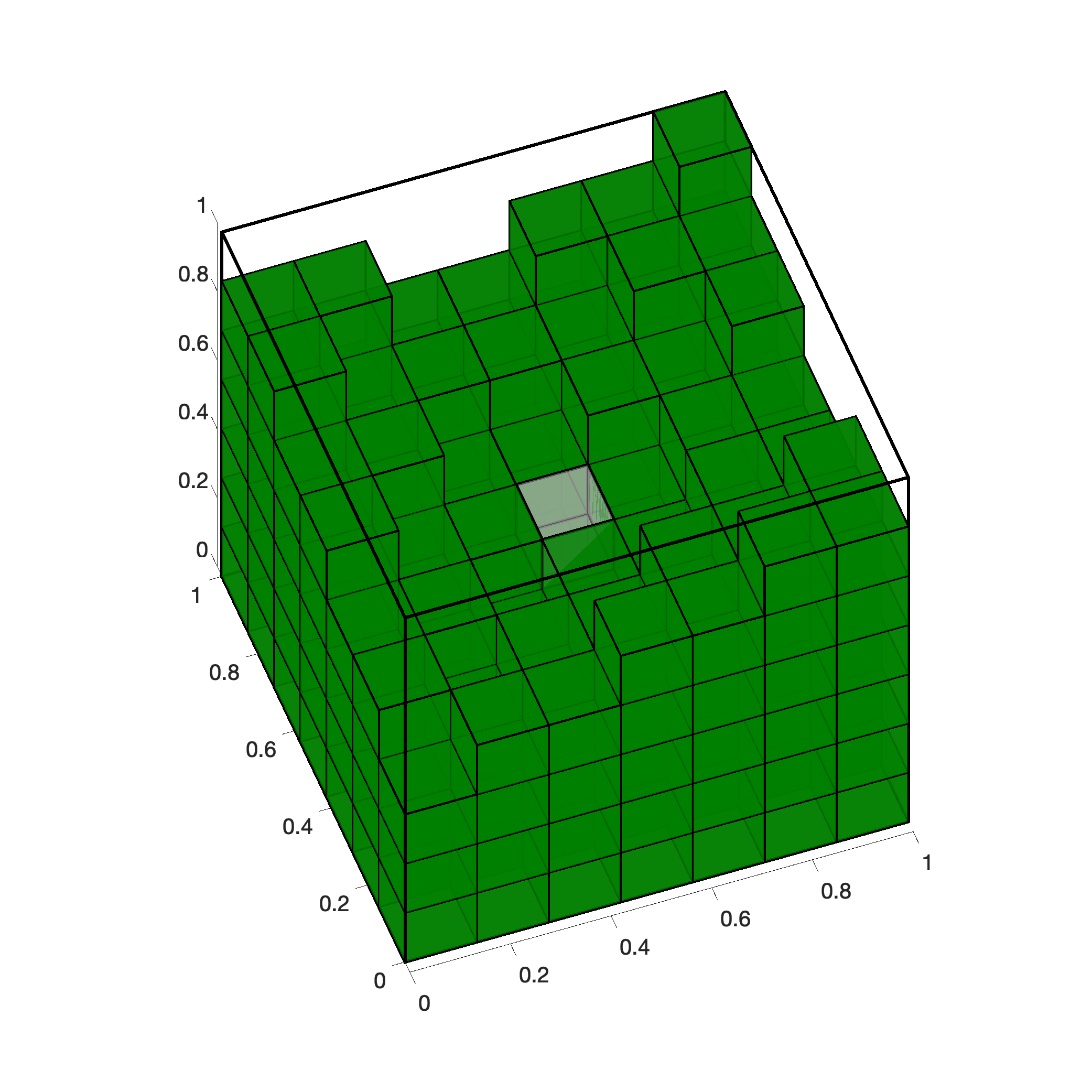}
\caption{Initial mesh}
\end{subfigure}
\begin{subfigure}[b]{0.24\textwidth}
\includegraphics[trim=150 0 150 0,clip,width=1.0\linewidth]{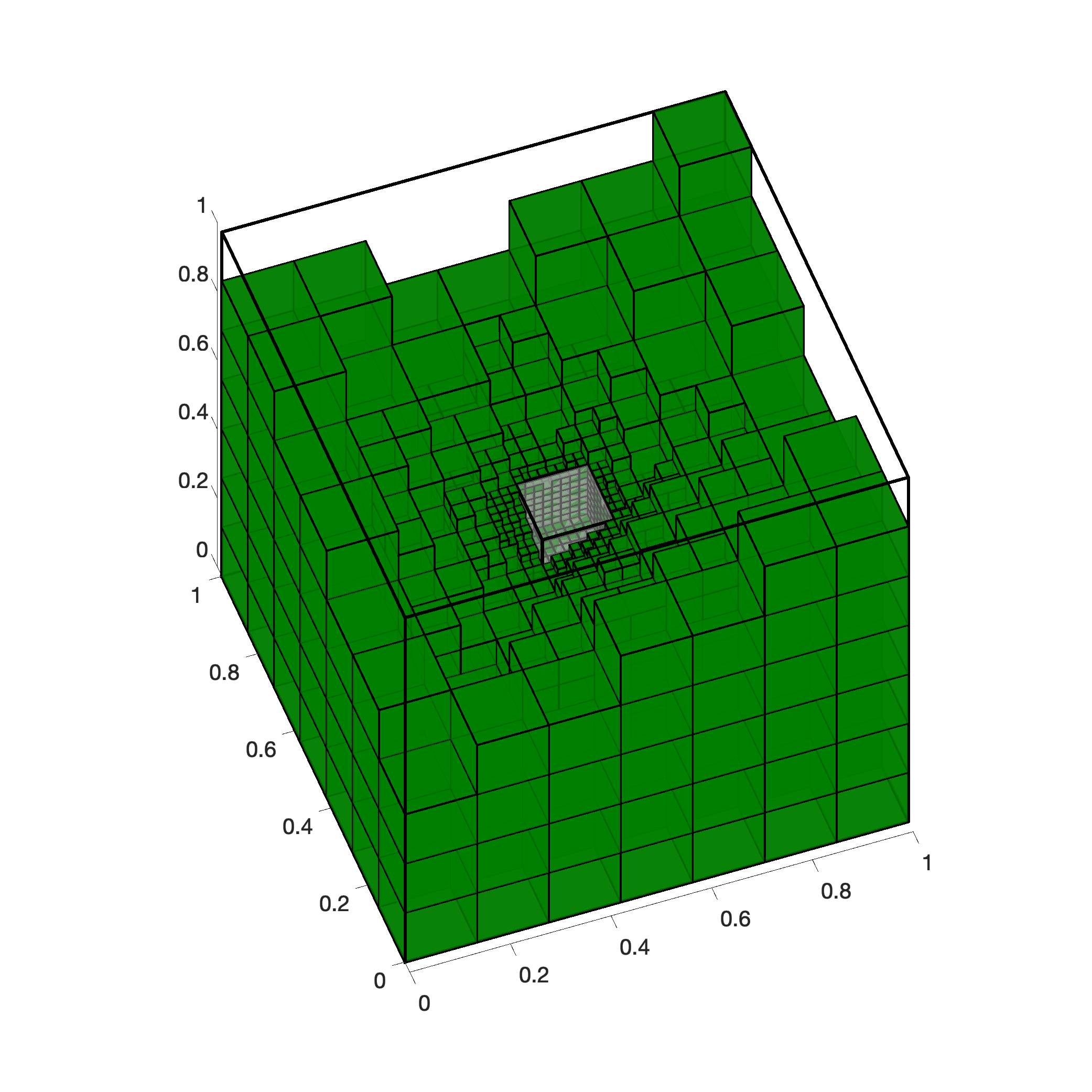}
\caption{Mesh 4}
\end{subfigure}
\begin{subfigure}[b]{0.24\textwidth}
\includegraphics[trim=150 0 150 0,clip,width=1.0\linewidth]{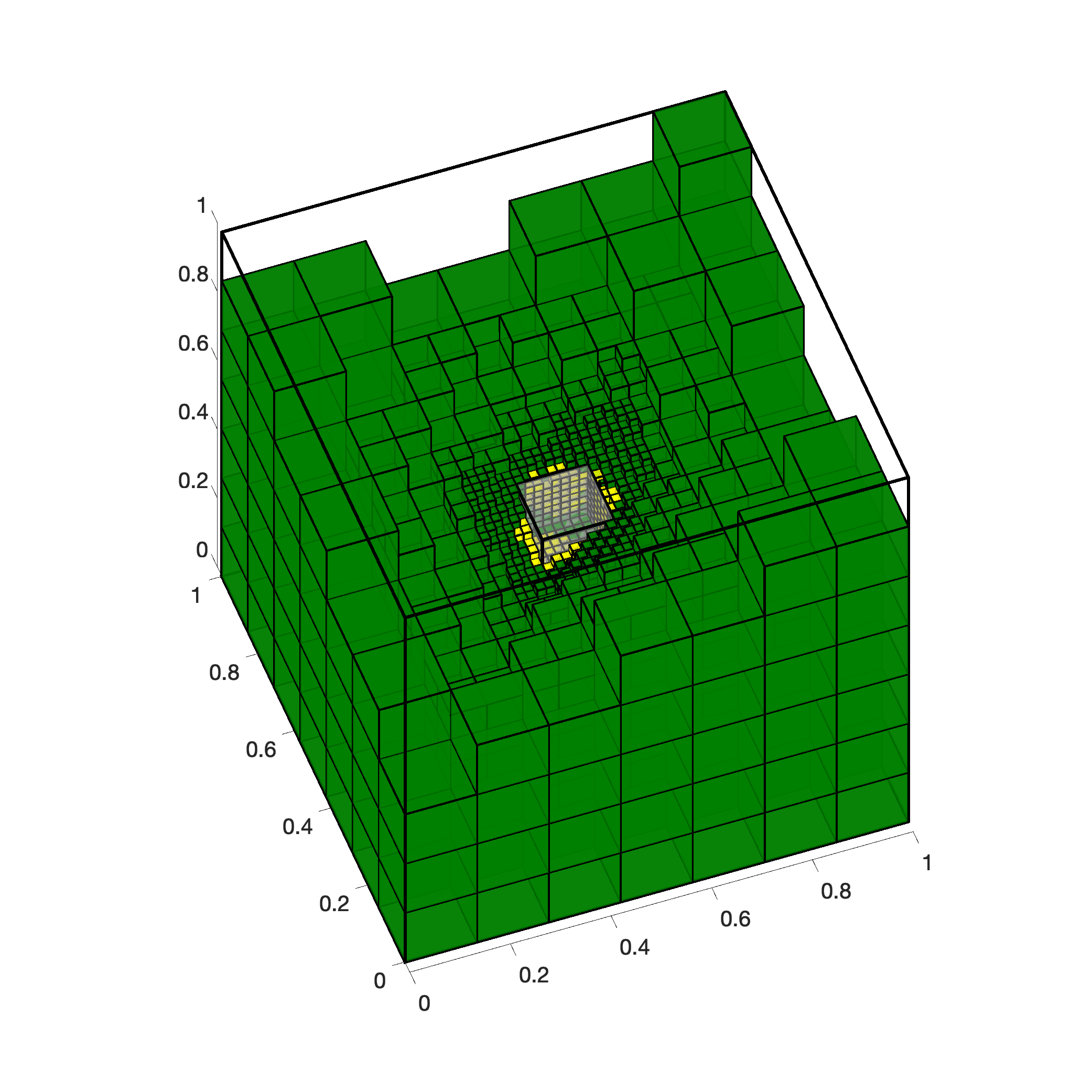}
\caption{Mesh 7}
\end{subfigure}
\begin{subfigure}[b]{0.24\textwidth}
\includegraphics[trim=150 0 150 0,clip,width=1.0\linewidth]{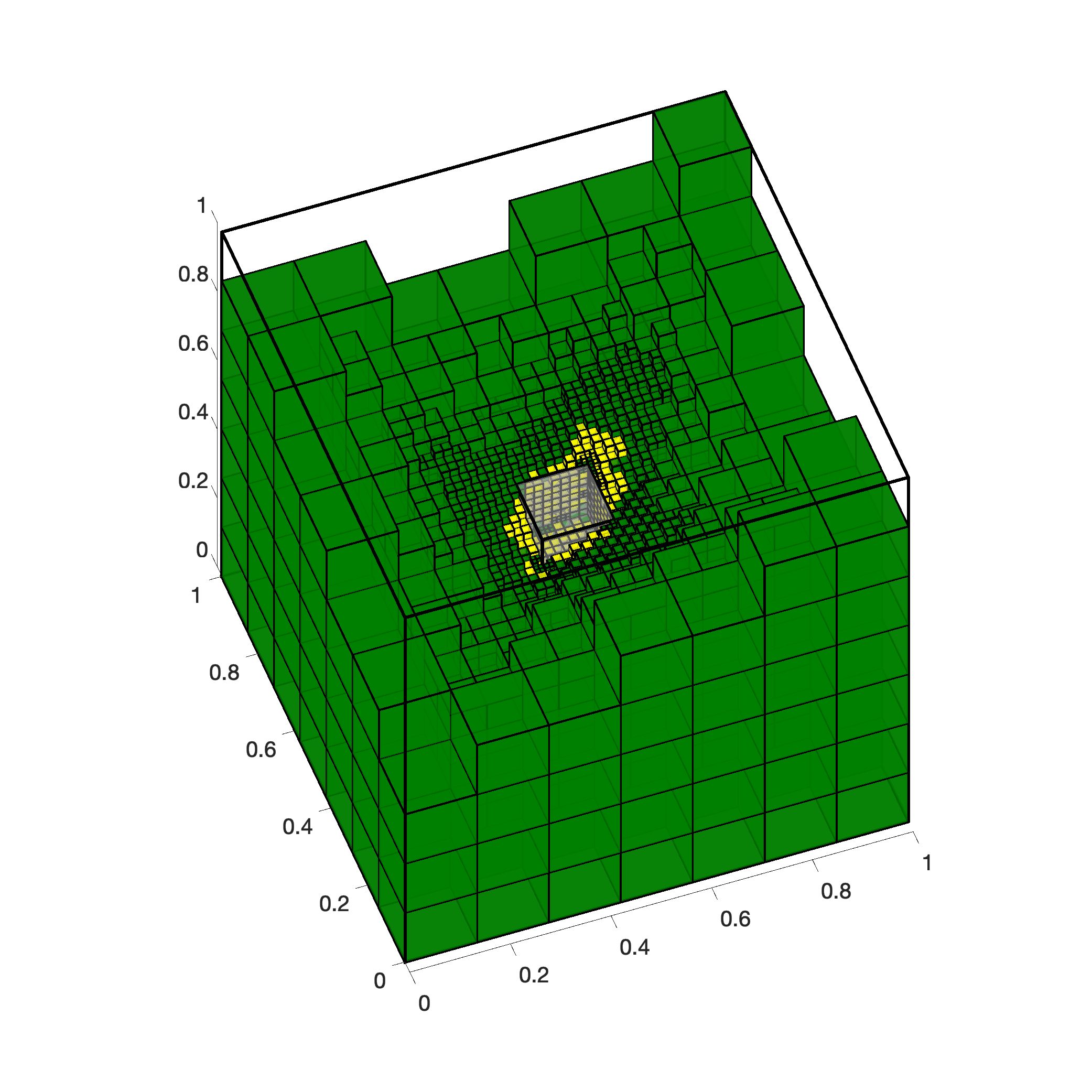}
\caption{Mesh 10}
\end{subfigure}
\begin{subfigure}[b]{0.24\textwidth}
\includegraphics[trim=150 0 150 0,clip,width=1.0\linewidth]{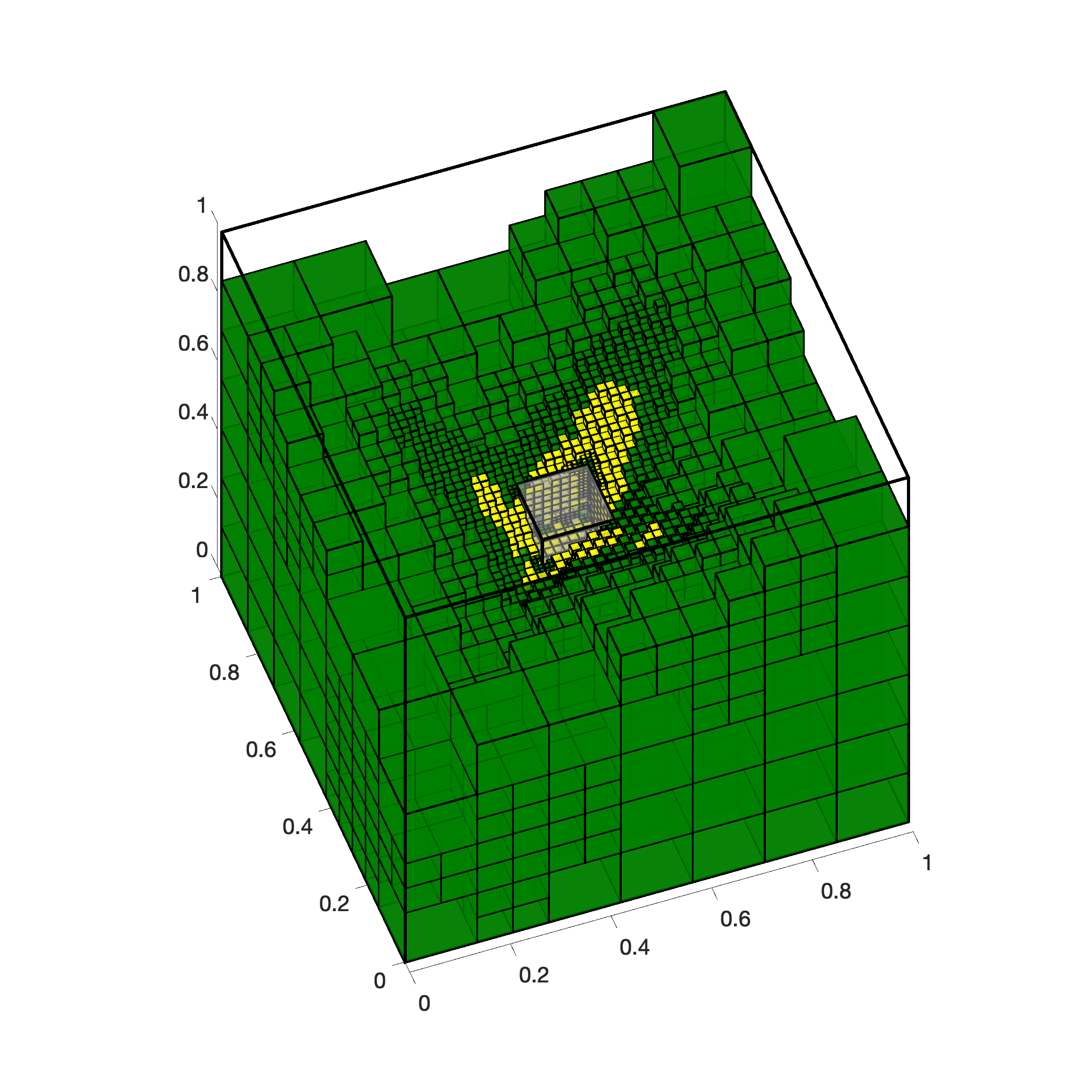}
\caption{Mesh 13}
\end{subfigure}
\begin{subfigure}[b]{0.24\textwidth}
\includegraphics[trim=150 0 150 0,clip,width=1.0\linewidth]{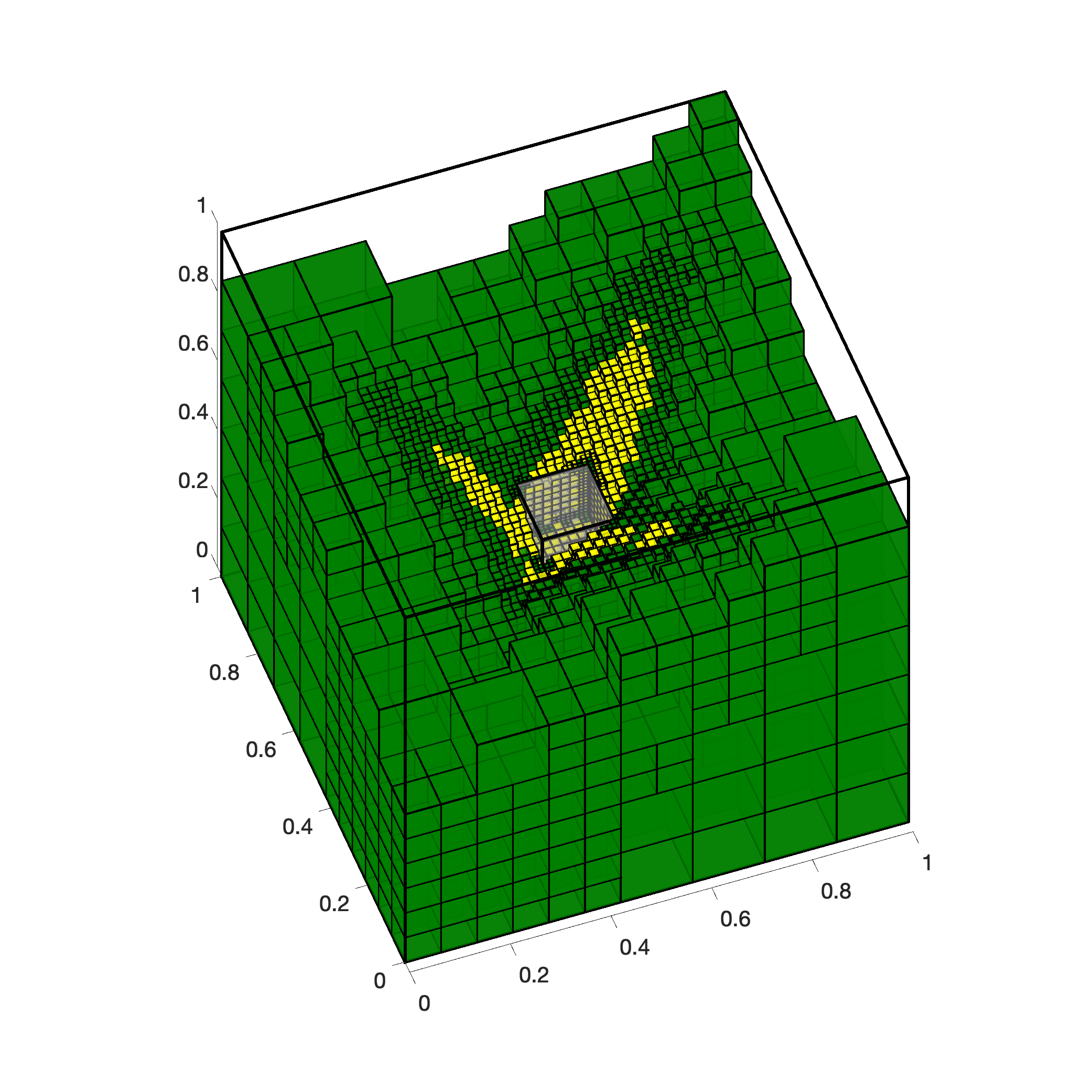}
\caption{Mesh 16}
\end{subfigure}
\begin{subfigure}[b]{0.24\textwidth}
\includegraphics[trim=150 0 150 0,clip,width=1.0\linewidth]{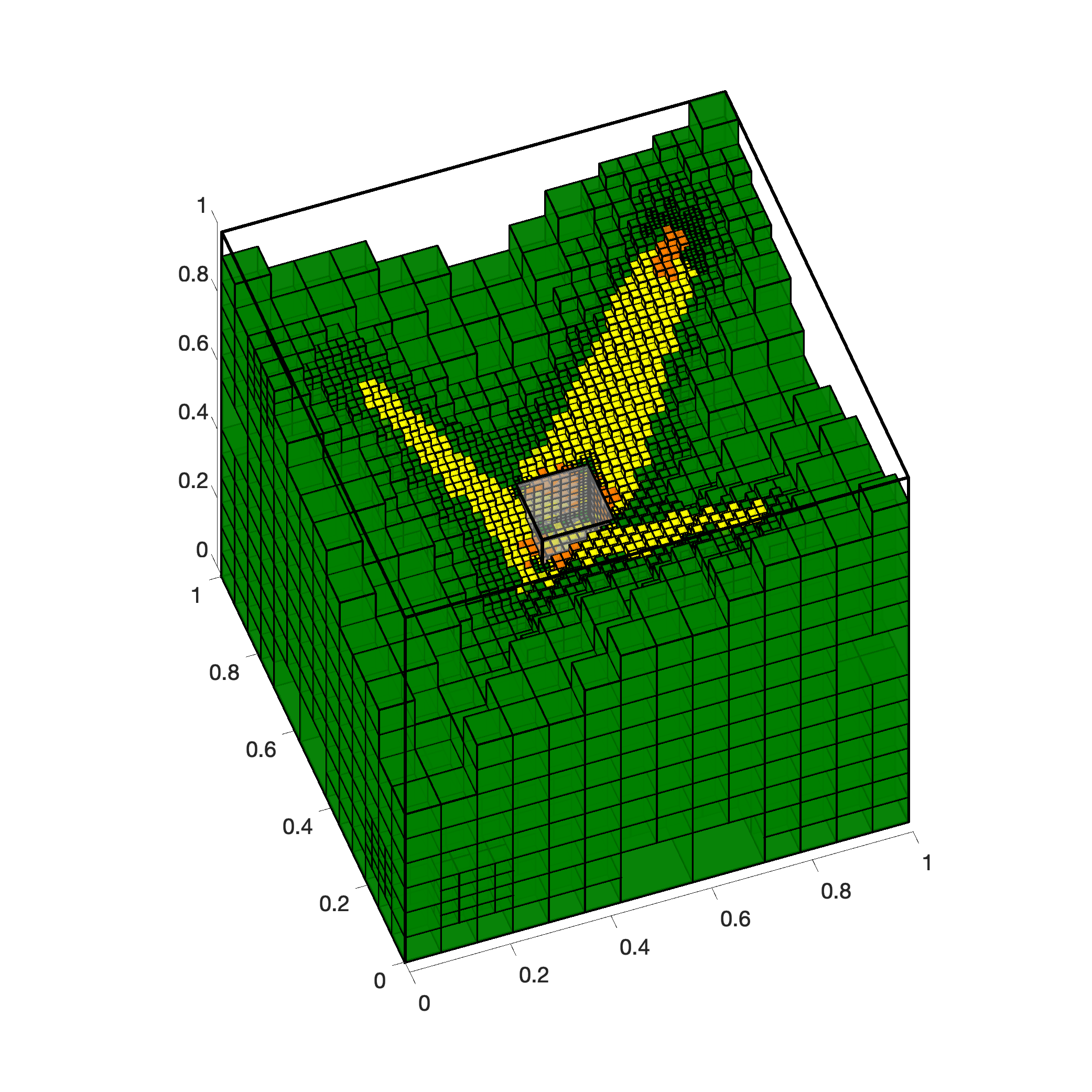}
\caption{Mesh 19}
\end{subfigure}
\begin{subfigure}[b]{0.24\textwidth}
\includegraphics[trim=150 0 150 0,clip,width=1.0\linewidth]{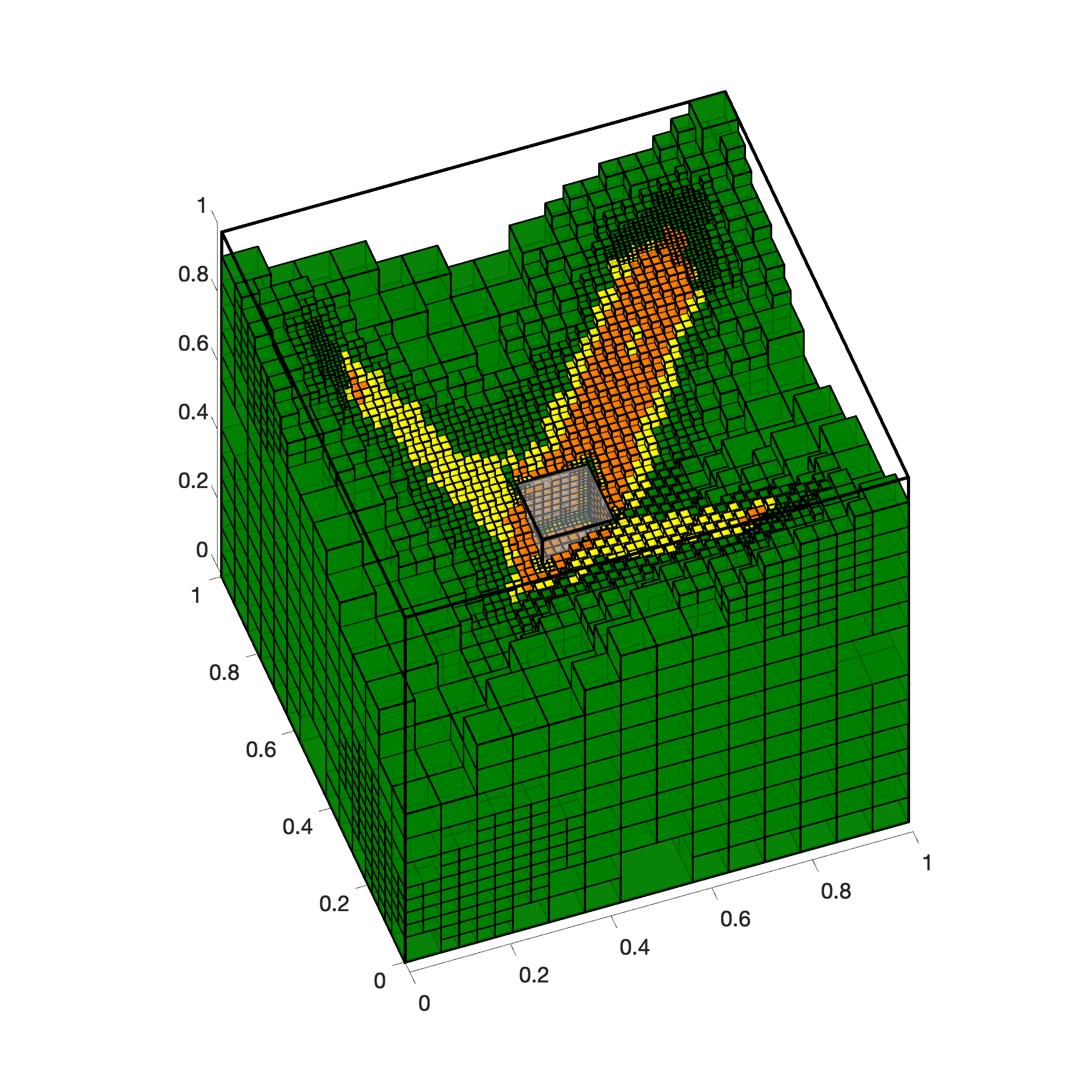}
\caption{Mesh 21}
\end{subfigure}
\end{figure}
\setcounter{figure}{11}
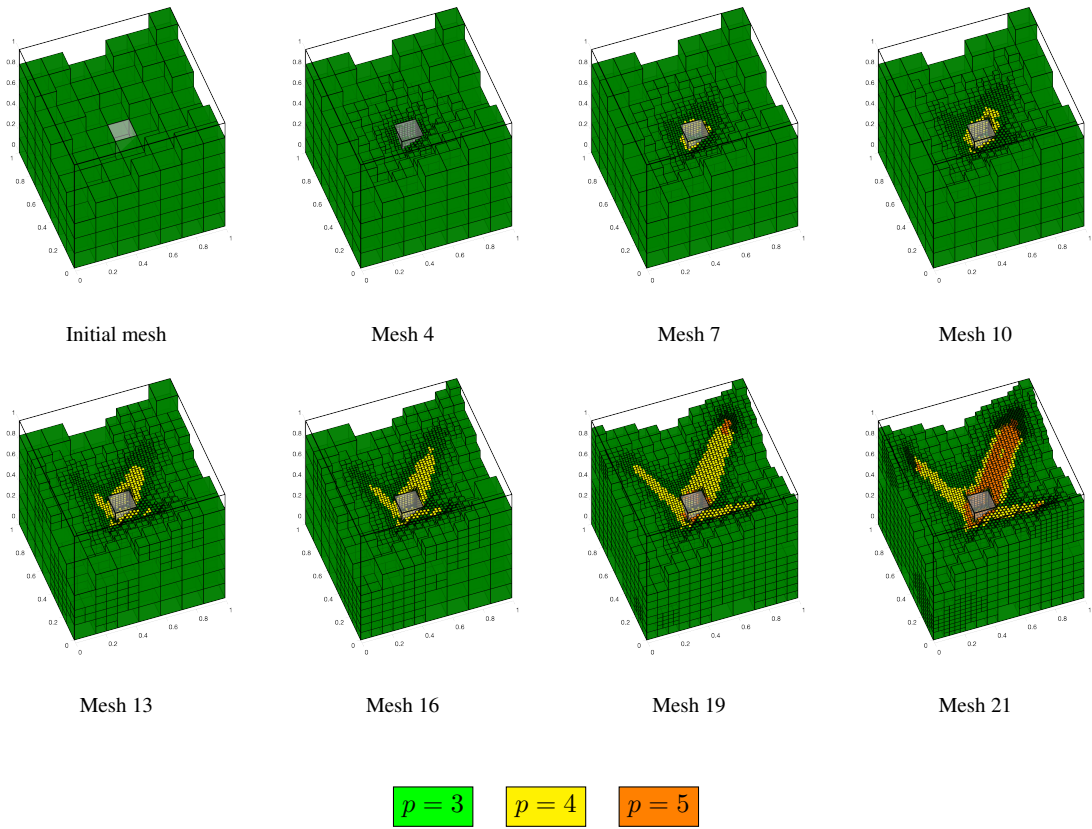
\begin{figure}[H]
\centering
\begin{tikzpicture}
\node[draw,align=center,fill=green] at (0,0) {$p=3$};
\node[draw,align=center,fill=yellow] at (1.5,0) {$p=4$};
\node[draw,align=center,fill=orange] at (3,0) {$p=5$};
\end{tikzpicture}
\caption{Evolution of the $hp$--adaptive mesh. Notice that the singularities at the scatterer have to be resolved before the wave can propagate.}\label{fig:beam_meshes}
\end{figure}
\begin{figure}[H]
\captionsetup[subfigure]{labelformat=empty}
\begin{subfigure}[b]{0.24\textwidth}
\includegraphics[trim=50 0 50 0,clip,width=1.0\linewidth]{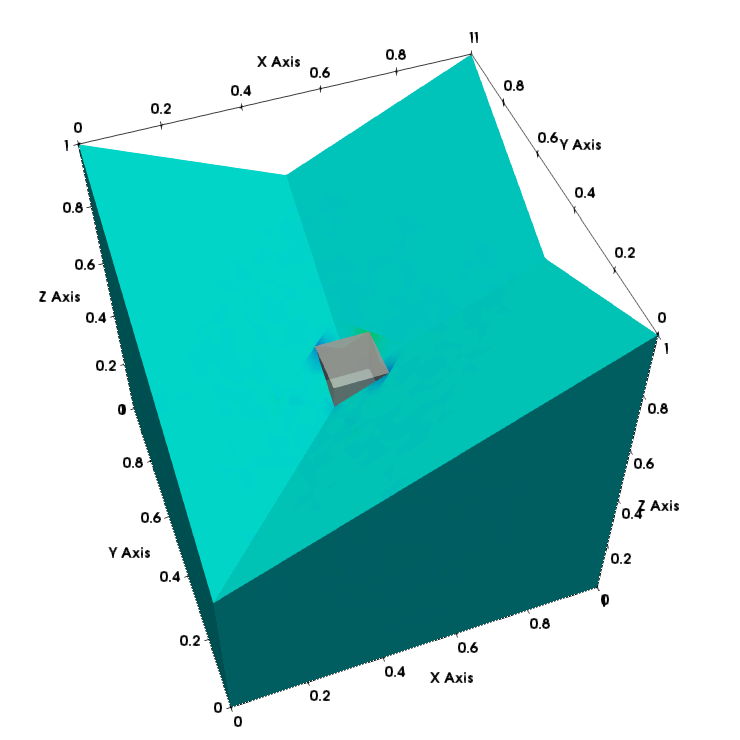}
\caption{Initial mesh}
\end{subfigure}
\begin{subfigure}[b]{0.24\textwidth}
\includegraphics[trim=50 0 50 0,clip,width=1.0\linewidth]{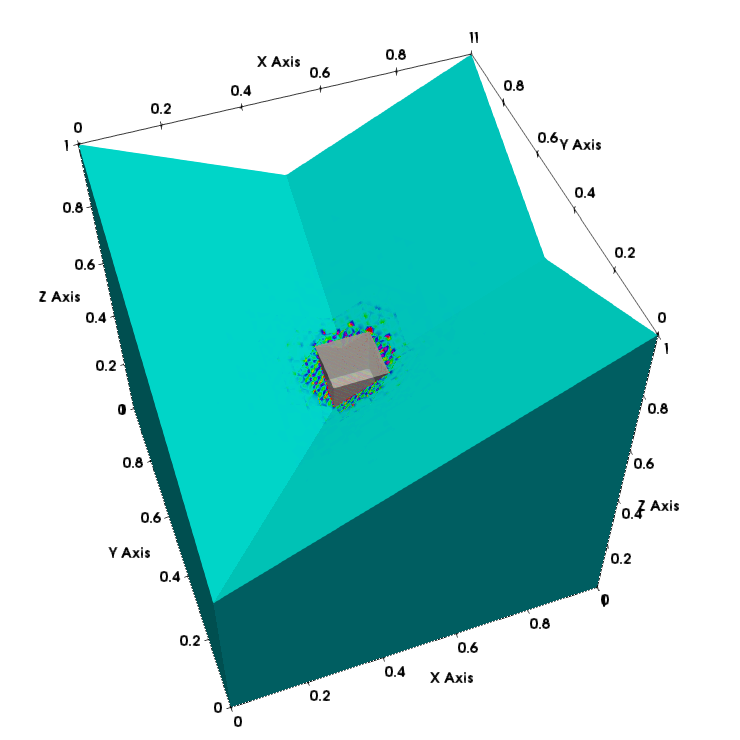}
\caption{Mesh 4}
\end{subfigure}
\begin{subfigure}[b]{0.24\textwidth}
\includegraphics[trim=50 0 50 0,clip,width=1.0\linewidth]{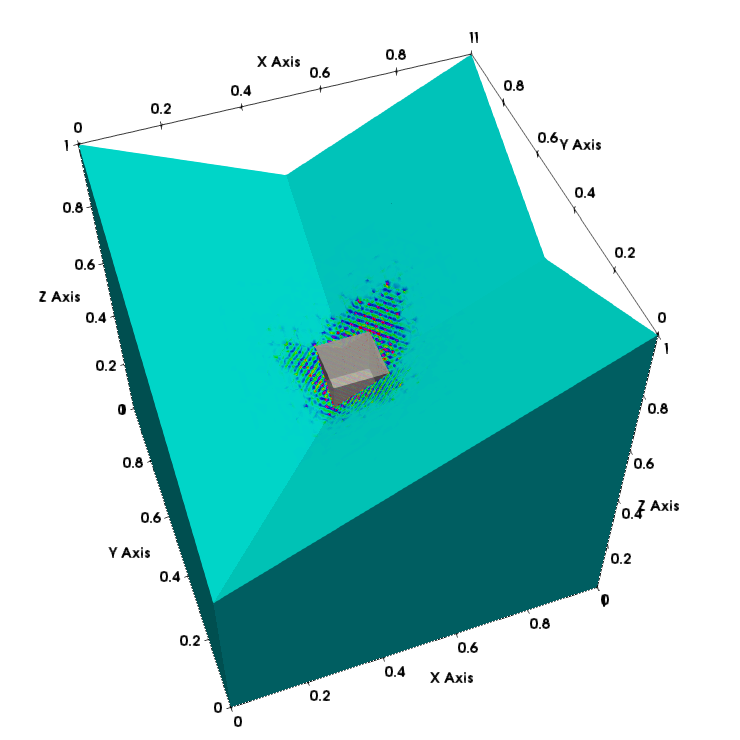}
\caption{Mesh 7}
\end{subfigure}
\begin{subfigure}[b]{0.24\textwidth}
\includegraphics[trim=50 0 50 0,clip,width=1.0\linewidth]{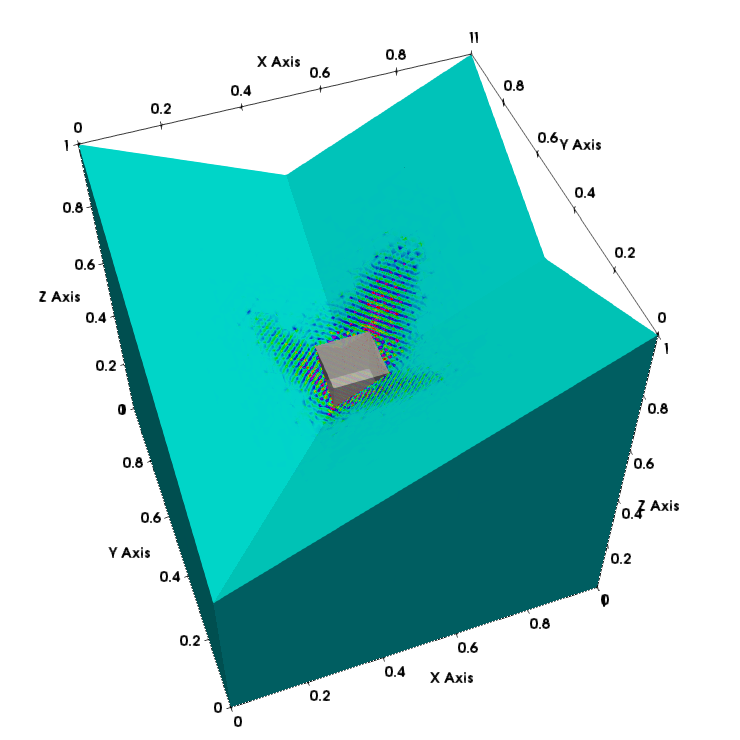}
\caption{Mesh 10}
\end{subfigure}

\begin{subfigure}[b]{0.24\textwidth}
\includegraphics[trim=50 0 50 0,clip,width=1.0\linewidth]{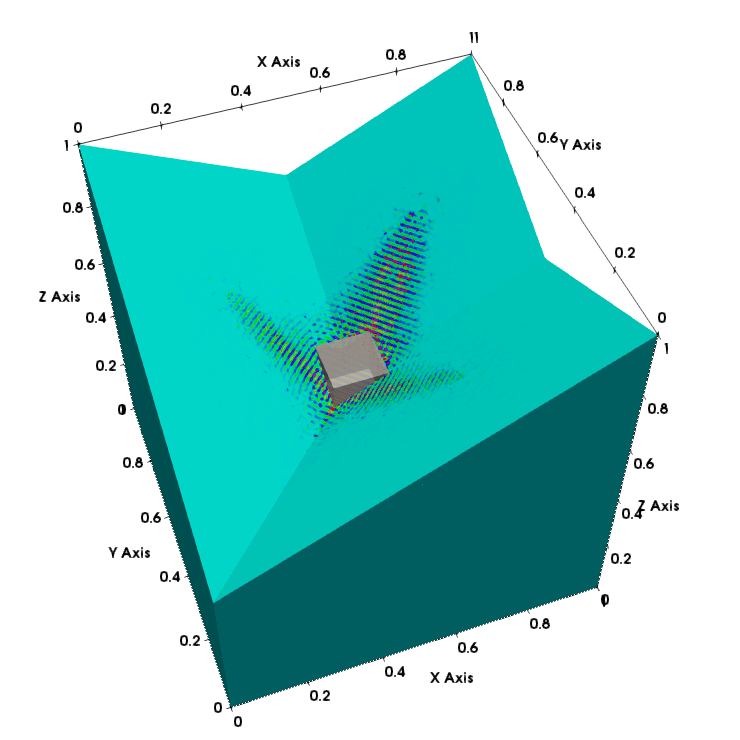}
\caption{Mesh 13}
\end{subfigure}
\begin{subfigure}[b]{0.24\textwidth}
\includegraphics[trim=50 0 50 0,clip,width=1.0\linewidth]{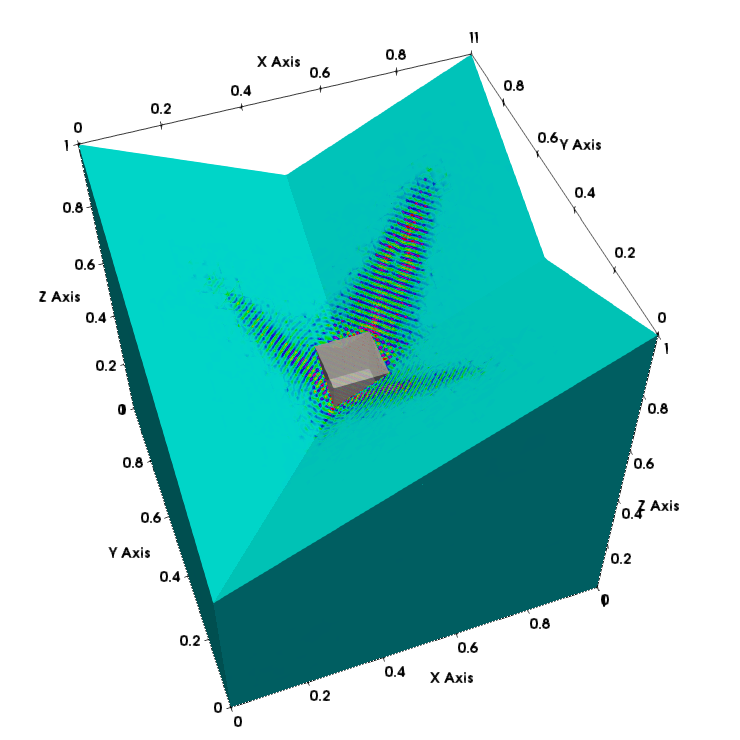}
\caption{Mesh 16}
\end{subfigure}
\begin{subfigure}[b]{0.24\textwidth}
\includegraphics[trim=50 0 50 0,clip,width=1.0\linewidth]{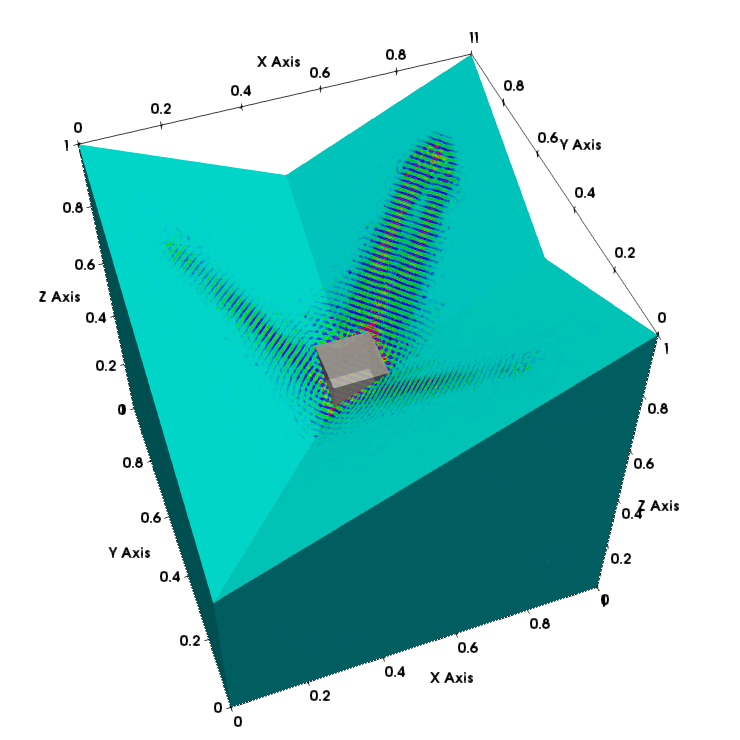}
\caption{Mesh 19}
\end{subfigure}
\begin{subfigure}[b]{0.24\textwidth}
\includegraphics[trim=50 0 50 0,clip,width=1.0\linewidth]{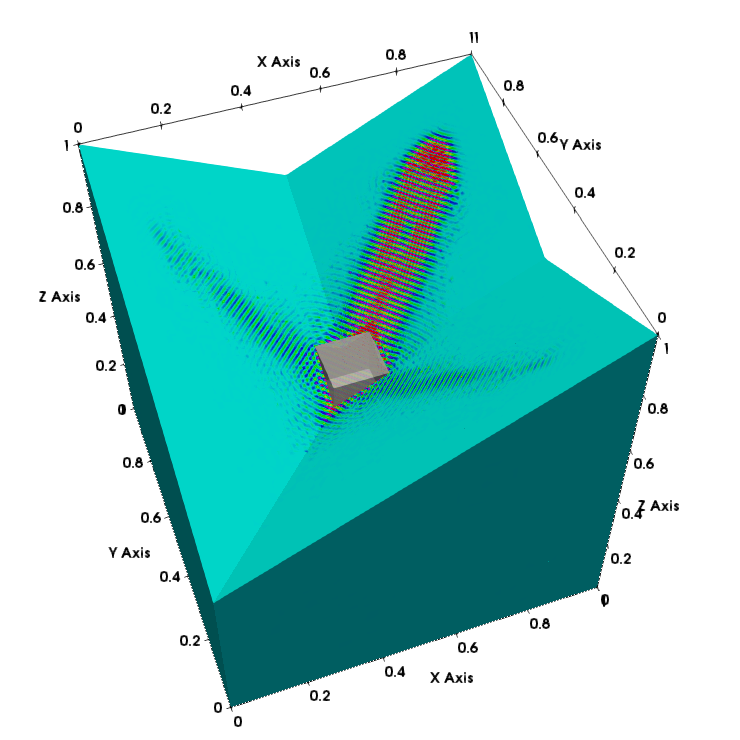}
\caption{Mesh 21}
\end{subfigure}
\end{figure}
\setcounter{figure}{12}
\begin{figure}[H]
      \centering
      \includegraphics[width=0.5\linewidth]{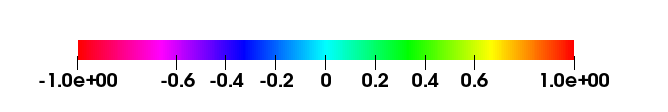}
      \caption{Evolution of the solution. Here the real part of the acoustic pressure is displayed.}\label{fig:beam_pressure}
\end{figure}
\begin{figure}[H]
\begin{subfigure}[b]{0.49\textwidth}
      \centering
\includegraphics[width=1.0\linewidth]{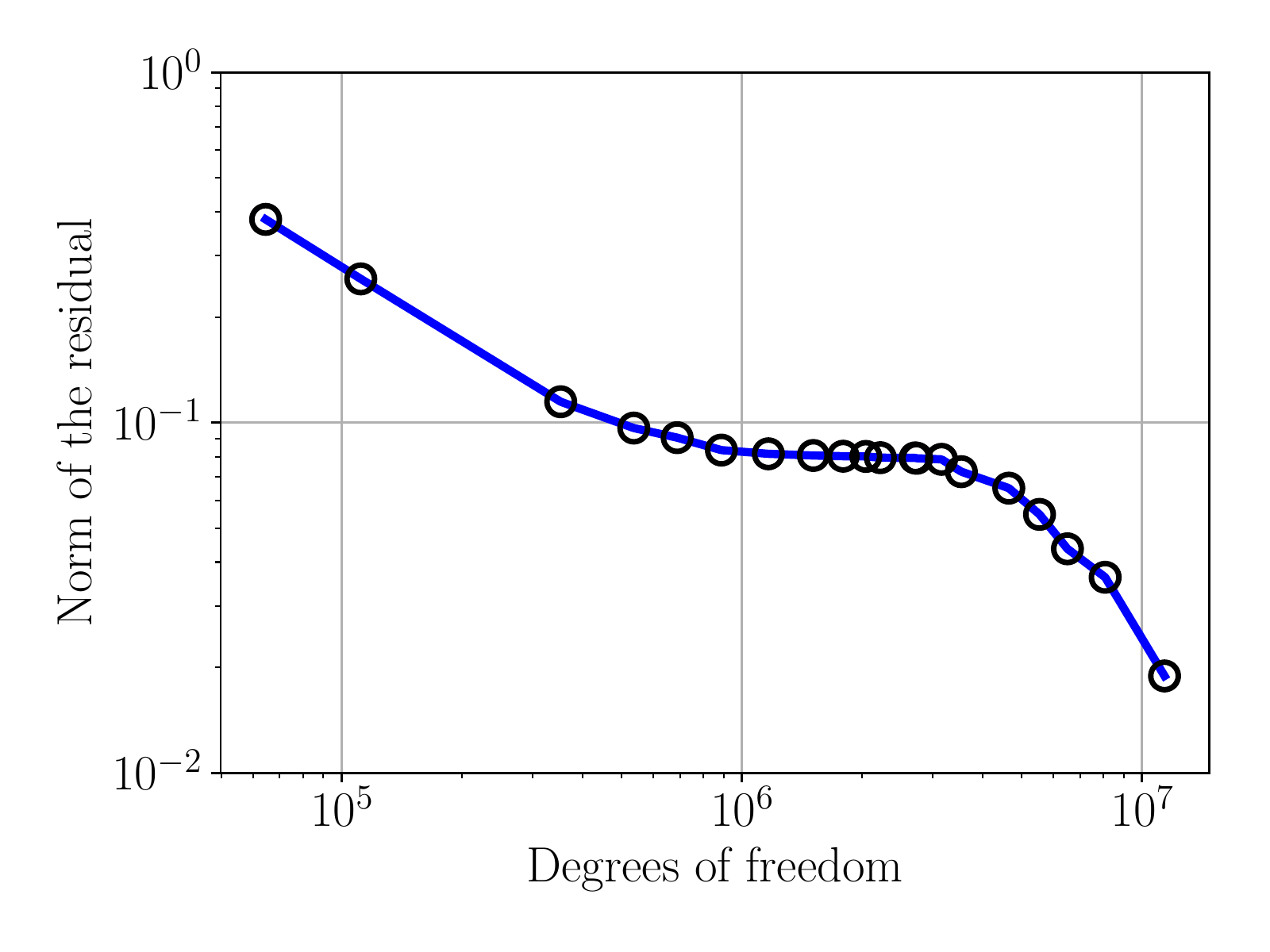}
\caption{Residual convergence}\label{fig:beam_res_convergence}
\end{subfigure}
\begin{subfigure}[b]{0.49\textwidth}
      \centering
\includegraphics[width=1.0\linewidth]{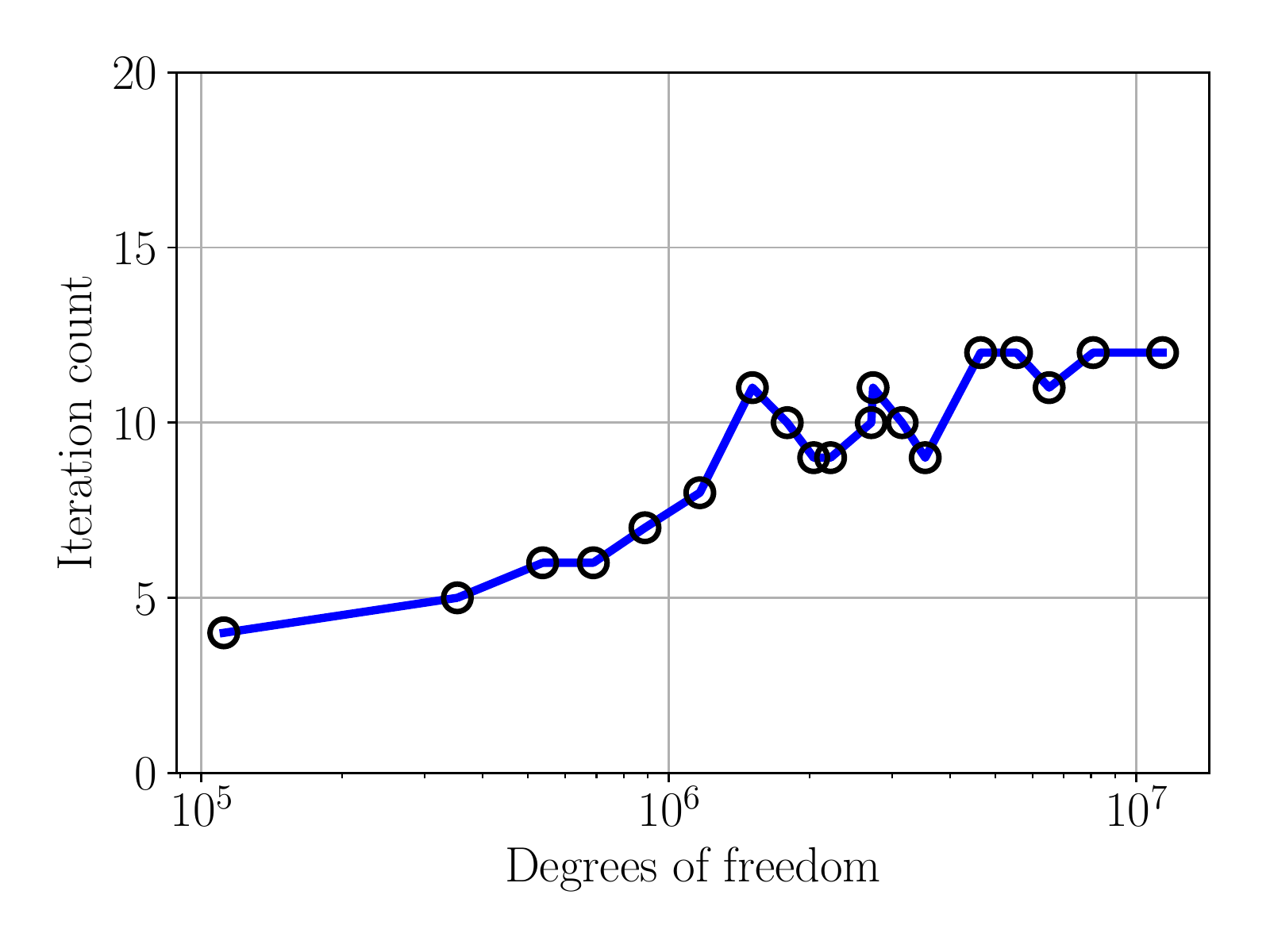}
\caption{Iteration count}\label{fig:beam_pcg_convergence}
\end{subfigure}
\caption{Convergence of the DPG residual and the preconditioned CG solver. Note that the number of iterations of the iterative solver is controlled throughout the adaptive process.}\label{fig:beam_convergence}
\end{figure}

%% file: conclusion.tex

\section{Conclusion}
\label{sec:sec6}

The main accomplishment of this work is the design and implementation of an efficient and robust preconditioner for linear systems arising from DPG discretizations. The construction is heavily based on the attractive properties of the DPG method, but also on well established theory of Schwarz Domain Decomposition and multigrid methods. Special emphasis was given in the simulation of wave problems, specifically acoustics and electromagnetic simulations. As it was showcased from numerous numerical results the method is stable and reliable and it is suitable for adaptive $hp$--meshes. Under certain circumstances, where the coarse grid is ``good'' enough, uniform convergence with respect to the frequency, polynomial order and discretization size can be achieved. Integrating the iterative solver with the adaptive refinement procedure allowed us to solve efficiently, problems in the high--frequency regime, problems which under a uniform--mesh setting would be intractable. Additionally, we presented some theoretical results for the overlapping Schwarz smoother, and demonstrated its dependence on the frequency for the acoustics problem in two--space dimensions.

Our ongoing work focuses on a new implementation of the solver for distributed memory computer architectures. Our aim is to be able to solve high--frequency electromagnetic problems, arising from the simulation of optical laser amplifiers \cite{henneking2019dpg,nagaraj20193d}, problems which require the sufficient resolution of waves of hundred of thousands of wavelengths. Since direct solvers fail to scale efficiently in modern many--core architectures, an efficient iterative solution scheme seems to be the only promising alternative choice. Finally, our theoretical analysis on the one--level additive preconditioner for the acoustics
problem can serve as a stepping stone on a further study in the multilevel setting. While
several other theoretical works exist on preconditioning DPG systems, analytical results on
preconditioning wave operators are yet to be explored.

%% file: interpolant.tex
%
\section{Computing the interpolation norm}
\label{appendix:interpolant}
The continuity constant $C_{\hat{\Pi}} = \|\hat{\Pi}\|_{\hat{U}}$ can be computed using the definition of a norm of an operator, i.e., by solving the maximization problem: 
\begin{equation*}
\|\hat{\Pi}\|_{\hat{U}} =\max_{\hat{v}\in\hat{U}, \hat{v}\ne0} \frac{\|\hat{\Pi} \hat{v}\|_{\hat{U}}}{\|\hat{v}\|_{\hat{U}}}
\end{equation*}
This leads to the generalized eigenvalue problem 
\begin{equation*}
(\hat{\Pi}\hat{v},\hat{\Pi}\widehat{\delta v})_{\hat{U}} = \lambda^2 (\hat{v},\widehat{\delta v})_{\hat{U}}, \qquad  \widehat{\delta v} \in  \hat{U}.
\end{equation*}
Consider now a discrete basis $\{\hat{v}_i\}_{i=1}^n$ for the polynomial subspace $\hat{U}_h \subset \hat{U}$. Then, we need to solve 
\begin{equation*}
\label{eq:eig1}
\sfP^*\sfG\sfP \sfv = \lambda^2 \sfG\sfv 
\end{equation*}
where $\sfv \in \C^n$,  $\sfP$ is the matrix representation of $\hat{\Pi}$ and $\sfG = (\hat{v}_i,\hat{v}_j)_{\hat{U}}$ is the Gram matrix corresponding to the inner product $(\cdot,\cdot)_{\hat{U}}$. The continuity constant is therefore given by the square root of the maximum generalized eigenvalue of ~\eqref{eq:eig1}. 

\bigskip
\paragraph{\bf Computation of matrix $G$}
Given a basis $\{\hat{v}_i\}_{i=1}^n$ of the polynomial subspace $\hat{U}_h \subset \hat{U}$ we can compute the entries of the Gram matrix by using the $(\cdot,\cdot)_{H_\Aop}$ inner product (see \cref{eq:ha}).
Indeed, using the polarization formula \cite[Ch.~6]{oden2010}, it is easy to see that 
\begin{equation}\label{eq:inner_product}
(\hat{v}_i,\hat{v}_j)_{\hat{U}} = (v_i,v_j)_{H_A}
\end{equation}
where $v_i, v_j \in H_\Aop(K)$ are the minimum energy extensions of $\hat{v}_i,\hat{v}_j$ respectively. However, in practice we can only approximate the minimum energy extensions using polynomial extensions. We note that for all the numerical results presented in \cref{sec:sec3}, polynomials of sufficiently large order \footnote{No notable change in the numerical results could be observed when increasing further the polynomial order of approximation} were used. 
\bigskip
\paragraph{\bf Approximating the minimum energy extension.} Given $\hat{v} \in \hat{U}_h$ we can derive the extension $v \in H_A(K)$ that realizes the minimum energy by:
\begin{equation*}
\|\hat{v}\|_{\hat{U}} = \inf_{tr v = \hat{v}} \|v\|_{H_\Aop}
\end{equation*}
Consider the polynomial subspace $H_{\Aop,h}(K) \subset H_{\Aop}(K)$. Then the above minimization problem leads to the following Dirichlet problem.
\begin{equation*}
\left\{
\begin{alignedat}{2}
&\text{Find } v_h\in H_{\Aop,h}(K) \\
&(v_h,\delta v_h)_{H_\Aop} = 0, \quad \delta v_h &&\text{ in } H_{\Aop,h}(K) \\
&v_h = \hat{v},&& \text{ on } \partial K
\end{alignedat}
\right.
\end{equation*}
or equivalently
\begin{equation*}
\left\{
\begin{aligned}
&\text{Find } v_b\in H^b_{\Aop,h}(K) \,\, (\text{bubble functions}) \\
&(\Aop v^b,\Aop\delta v^b) + (v^b,\delta v^b)  = -(\hat{v},\delta v^b)_{H_\Aop},&& \forall \delta v^b\text{ in } H^b_{\Aop,h}(K) \\
\end{aligned}
\right.
\end{equation*}

\bigskip
\paragraph{\bf Definition of the interpolation operator $\Pi$ - computation of matrix $\sfP$.}
Let $\hat{v} \in \hat{U}_h(K)$. We define the projection based interpolant $\hat{v}^p = \hat{\Pi}\hat{v}$ using the following steps: 
\begin{itemize}
\item{Interpolation at vertices:} the interpolant $\hat{v}^p$ matches the function $\hat{v}$ at the vertices
\begin{equation*}
\hat{v}^p(a) = \hat{v}(a), \qquad \forall \text{ vertex } a.   
\end{equation*}
We lift the vertex values using a polynomial extension that lives in the element trace space. This leads to the linear interpolant $\hat{v}_1 \in tr \mathcal{P}^1(K)$. 
\item{Edge projection:} We subtract the linear interpolant $\hat{v}_1$ from the function $\hat{v}$. Now the difference $\hat{v}-\hat{v}_1$ vanishes at the element vertices. Then we project the difference onto the trace space of edge polynomials of order $p_e$ vanishing at the vertices (i.e, the edge bubbles $\mathcal{P}^{p_e}_0$), i.e., 
\begin{equation*}
\left\{
\begin{aligned}
&\hat{v}_{2,e} \in \mathcal{P}^{p_e}_0 \\ 
&\|\hat{v}-\hat{v}_1 - \hat{v}_{2,e}\|_{\hat{U}(e)} \rightarrow \min  \,.  
\end{aligned}
\right.
\end{equation*}
Here, the norm $\|\cdot\|_{\hat{U}(e)}$ is defined by the inner product \eqref{eq:inner_product} on the edge. The edge interpolant is then the sum of the edge projections 
$$\hat{v}_2 = \sum_e \hat{v}_{2,e}.$$
\end{itemize}
The final interpolant is defined by the sum of the vertex and the edge interpolant
$$\hat{v}^p = \hat{v}_1 + \hat{v}_2.$$
We obtain a matrix representation of $\sfP$ by applying it to a basis of the polynomial space. The polynomial spaces are defined as \cite[Ch. 2]{petrides2019adaptive}. In particular, consider

\begin{equation*}
\begin{aligned}
  W^r := \mathcal{Q}^{(r,r)} \subset H^1 (\Omega), \,\, \text{ and } \,\,
  V^r := \mathcal{Q}^{(r,r-1)} \times \mathcal{Q}^{(r-1,r)} \subset H(\text{div},\Omega). 
\end{aligned}
\end{equation*}
Then $\hat{\Pi} : W^{r+1} \times V^{r+1} \rightarrow W^{r} \times V^{r}$, and its matrix representation $\sfP$ is computed as follows. For $\hat{p}\in \,\, tr|_{\partial K} W^{r+1}$, we have vertex shape functions and bubble shape functions. For the vertex shape functions the interpolation operator matches the values at the vertices and so the first 4 columns of matrix P are the orthonormal vectors $\{e_i\}_{i=1}^4$. Note that if $\hat{p}$ is an edge bubble then it vanishes on the vertices. Therefore it is enough to solve the following minimization problem. Let $\hat{p} \in tr|_{\partial K}\mathcal{W}^{r+1}_0$ be an edge bubble vanishing on the boundary. Then we solve
\begin{equation*}
\left\{
\begin{aligned}
&\text{Find } \hat{v}^b \in tr|_{\partial K}\mathcal{W}^{r}_0 \\
&\|\hat{p} - \hat{v}^b\|_{\hat{U}} \rightarrow \min    
\end{aligned}
\right.
\end{equation*}
or equivalently 
\begin{equation*}
\left\{
\begin{aligned}
&\text{Find } \hat{v}^b \in tr|_{\partial K}\mathcal{W}^{r}_0 \\
&(\hat{v}^b,\delta\hat{v}^b)_{\hat{U}} = (\hat{p},\delta\hat{v}^b)_{\hat{U}} \qquad \forall  \delta\hat{v}^b \in tr|_{\partial K}\mathcal{W}^{r}_0
\end{aligned}
\right.
\end{equation*}
or 
\begin{equation}\label{eq:matrixP_comp}
\left\{
\begin{aligned}
&\text{Find } v^b \in \mathcal{W}^{r}_0 \\
&(v^b,\delta v^b)_{H_\Aop} = (p^b,\delta v^b)_{H_\Aop} \qquad \forall \delta v^b \in \mathcal{W}^{r}_0
\end{aligned}
\right .
\end{equation}
where $p^b, v^b, \delta v^b$ are the minimum energy polynomial extensions of $\hat{p}, \hat{v}^b,\delta\hat{v}^b$ respectively. 
Finally the interpolant is defined to be the trace of the solution $v^b$. Notice that both left and right hand sides of \eqref{eq:matrixP_comp} can be computed using the Gram matrix $G$. For the edge bubbles of the variable $\hat{u}_n \in tr|_{\partial K} V^{r+1}$, we follow a similar procedure. 

\bigskip
\paragraph{\bf Computation on the master element.} We can reduce the computation on the master element using appropriate scalings. Assuming that each element $K$ is obtained by a simple scaling of the master element $\bar{K}$, and using the same $H^1$ scaling for both $H^1$ and $H(\text{div})$ functions we have:    

\begin{equation*}
\begin{aligned}
\|\Aop_{\omega}(p,u)\|^2_{L^2(K)} 
:&= \int_K (|i\omega p + \text{div }u|^2 + |i\omega u + \nabla p|^2 ) \,dK \\
& = \int_K h^2 (|i\omega \bar{p} + \frac{1}{h}\bar{\text{div }}\bar{u}|^2 + |i\omega \bar{u} + \frac{1}{h} \bar{\nabla} \bar{p}|^2 ) \,d\bar{K} \\
& = \int_K (|i\omega h \bar{p} + \bar{\text{div }}\bar{u}|^2 + |i\omega h \bar{u} + \bar{\nabla} \bar{p}|^2 ) \,d\bar{K} \\
&=: \|\Aop_{\omega h}(p,u)\|_{L^2(\bar{K})}^2
\end{aligned}
\end{equation*}
and
\begin{equation*}
\|(p,u)\|_{L^2(K)}^2 =  h^2 \|(\bar{p},\bar{u})\|^2_{L^2(\bar{K})} 
\end{equation*}
We can therefore perform all the computations on the master element by using the following norm:
\begin{equation*}
\|(p,u)\|_{H_\Aop(\omega h,h)}^2 =  \|\Aop_{\omega h}(p,u)\|_{L^2(\bar{K})}^2 + h^2 \|(\bar{p},\bar{u})\|^2_{L^2(\bar{K})} 
\end{equation*}

%% file: main.bbl
\begin{thebibliography}{10}

\bibitem{MUMPS2}
{\sc P.~R. Amestoy, I.~S. Duff, J.-Y. L'Excellent, and J.~Koster}, {\em A fully
  asynchronous multifrontal solver using distributed dynamic scheduling}, SIAM
  J. Matrix Anal. Appl., 23 (2001), pp.~15--41 (electronic).

\bibitem{arnold1997preconditioning}
{\sc D.~Arnold, R.~Falk, and R.~Winther}, {\em {Preconditioning in H(div) and
  applications}}, Math. Comp., 66 (1997), pp.~957--984.

\bibitem{babuska}
{\sc I.~Babu{\v{s}}ka}, {\em {Error-bounds for finite element method}}, Numer.
  Math., 16 (1971), pp.~322--333.

\bibitem{Babuska1}
{\sc I.~M. Babu{\v{s}}ka and S.~A. Sauter}, {\em Is the pollution effect of the
  {FEM} avoidable for the {H}elmholtz equation considering high wave numbers?},
  SIAM J. Numer. Anal., 34 (1997), pp.~2392--2423.

\bibitem{Barker2014}
{\sc A.~T. Barker, S.~C. Brenner, E.-H. Park, and L.-Y. Sung}, {\em {A
  one-level additive Schwarz preconditioner for a discontinuous
  Petrov--Galerkin method}}, in Domain Decomposition Methods in Science and
  Engineering XXI, Springer, 2014, pp.~417--425.

\bibitem{barker}
{\sc A.~T. Barker, V.~Dobrev, J.~Gopalakrishnan, and T.~Kolev}, {\em {A
  Scalable Preconditioner for a Primal Discontinuous Petrov--Galerkin Method}},
  SIAM J. Sci. Comput., 40 (2018), pp.~A1187--A1203.

\bibitem{absorption_2}
{\sc M.~Bonazzoli, V.~Dolean, I.~G. Graham, E.~A. Spence, and P.-H. Tournier},
  {\em {Domain decomposition preconditioning for the high-frequency
  time-harmonic Maxwell equations with absorption}}, arXiv preprint
  arXiv:1711.03789,  (2017).

\bibitem{DBLP}
{\sc X.-C. Cai and O.~B. Widlund}, {\em {Domain decomposition algorithms for
  indefinite elliptic problems}}, SIAM J. Sci. Statist. Comput., 13 (1992),
  pp.~243--258.

\bibitem{XavierV2}
{\sc H.~Calandra, S.~Gratton, X.~Pinel, and X.~Vasseur}, {\em An improved
  two-grid preconditioner for the solution of three-dimensional {H}elmholtz
  problems in heterogeneous media}, Numer. Linear Algebra Appl., 20 (2013),
  pp.~663--688.

\bibitem{Carstensen}
{\sc C.~Carstensen, L.~Demkowicz, and J.~Gopalakrishnan}, {\em Breaking spaces
  and forms for the {DPG} method and applications including {M}axwell
  equations}, Comput. Math. Appl., 72 (2016), pp.~494--522.

\bibitem{CH17}
{\sc C.~Carstensen and F.~Hellwig}, {\em {Low-Order Discontinuous
  Petrov--Galerkin Finite Element Methods for Linear Elasticity}}, SIAM J.
  Numer. Anal., 54 (2016), pp.~3388--3410.

\bibitem{chen2013source}
{\sc Z.~Chen and X.~Xiang}, {\em {A source transfer domain decomposition method
  for Helmholtz equations in unbounded domain}}, SIAM J. Numer. Anal., 51
  (2013), pp.~2331--2356.

\bibitem{Demkbook1}
{\sc L.~Demkowicz}, {\em Computing with {$hp$}-adaptive finite elements. {V}ol.
  1}, Chapman \& Hall/CRC Applied Mathematics and Nonlinear Science Series,
  Chapman \& Hall/CRC, Boca Raton, FL, 2007.
\newblock One and two dimensional elliptic and Maxwell problems, With 1 CD-ROM
  (UNIX).

\bibitem{demk_varform}
{\sc L.~Demkowicz}, {\em {Various Variational Formulations and Closed Range
  Theorem}}, ICES Report, 15-03 (2015).

\bibitem{dpg_opt}
{\sc L.~Demkowicz and J.~Gopalakrishnan}, {\em A class of discontinuous
  {P}etrov-{G}alerkin methods. {II}. {O}ptimal test functions}, Numer. Methods
  Partial Differential Equations, 27 (2011), pp.~70--105.

\bibitem{demk2}
{\sc L.~Demkowicz, J.~Gopalakrishnan, I.~Muga, and J.~Zitelli}, {\em Wavenumber
  explicit analysis of a {DPG} method for the multidimensional {H}elmholtz
  equation}, Comput. Methods Appl. Mech. Engrg., 213/216 (2012), pp.~126--138.

\bibitem{Demkbook2}
{\sc L.~Demkowicz, J.~Kurtz, D.~Pardo, M.~Paszy{\'n}ski, W.~Rachowicz, and
  A.~Zdunek}, {\em Computing with {$hp$}-adaptive finite elements. {V}ol. 2},
  Chapman \& Hall/CRC Applied Mathematics and Nonlinear Science Series, Chapman
  \& Hall/CRC, Boca Raton, FL, 2008.
\newblock Frontiers: three dimensional elliptic and Maxwell problems with
  applications.

\bibitem{dryja1994domain}
{\sc M.~Dryja and O.~B. Widlund}, {\em {Domain decomposition algorithms with
  small overlap}}, SIAM J. Sci. Comput., 15 (1994), pp.~604--620.

\bibitem{EngquistY1}
{\sc B.~Engquist and L.~Ying}, {\em Sweeping preconditioner for the {H}elmholtz
  equation: moving perfectly matched layers}, Multiscale Model. Simul., 9
  (2011), pp.~686--710.

\bibitem{ErnstGander}
{\sc O.~G. Ernst and M.~J. Gander}, {\em Why it is difficult to solve
  {H}elmholtz problems with classical iterative methods}, in Numerical analysis
  of multiscale problems, vol.~83 of Lect. Notes Comput. Sci. Eng., Springer,
  Heidelberg, 2012, pp.~325--363.

\bibitem{fuentes2017coupled}
{\sc F.~Fuentes, B.~Keith, L.~Demkowicz, and P.~Le~Tallec}, {\em {Coupled
  variational formulations of linear elasticity and the DPG methodology}}, J.
  Comput. Phys., 348 (2017), pp.~715--731.

\bibitem{keith}
{\sc F.~Fuentes, B.~Keith, L.~Demkowicz, and S.~Nagaraj}, {\em Orientation
  embedded high order shape functions for the exact sequence elements of all
  shapes}, Comput. Math. Appl., 70 (2015).

\bibitem{absorption_1}
{\sc I.~G.~Graham, E.~A.~Spence, and E.~Vainikko}, {\em {Domain decomposition
  preconditioning for High-Frequency Helmholtz problems with absorption}},
  Math. Comp., 86 (2015).

\bibitem{Gander2}
{\sc M.~J. Gander, I.~G. Graham, and E.~A. Spence}, {\em Applying {GMRES} to
  the {H}elmholtz equation with shifted {L}aplacian preconditioning: what is
  the largest shift for which wavenumber-independent convergence is
  guaranteed?}, Numer. Math., 131 (2015), pp.~567--614.

\bibitem{gander2019class}
{\sc M.~J. Gander and H.~Zhang}, {\em {A class of iterative solvers for the
  Helmholtz equation: Factorizations, sweeping preconditioners, source
  transfer, single layer potentials, polarized traces, and optimized Schwarz
  methods}}, Siam Review, 61 (2019), pp.~3--76.

\bibitem{dpg_e}
{\sc J.~Gopalakrishnan, I.~Muga, and N.~Olivares}, {\em {Dispersive and
  Dissipative Errors in the DPG Method with Scaled Norms for Helmholtz
  Equation}}, SIAM J. Sci. Comput., 36 (2014), pp.~A20--A39.

\bibitem{JaySchwarz}
{\sc J.~Gopalakrishnan and J.~E. Pasciak}, {\em Overlapping {S}chwarz
  preconditioners for indefinite time harmonic {M}axwell equations}, Math.
  Comp., 72 (2003), pp.~1--15 (electronic).

\bibitem{Jaymultigrid}
{\sc J.~Gopalakrishnan, J.~E. Pasciak, and L.~F. Demkowicz}, {\em Analysis of a
  multigrid algorithm for time harmonic {M}axwell equations}, SIAM J. Numer.
  Anal., 42 (2004), pp.~90--108 (electronic).

\bibitem{Gopalakrishnan1}
{\sc J.~Gopalakrishnan and W.~Qiu}, {\em {An analysis of the practical DPG
  method}}, Math. Comp., 83 (2014), pp.~537--552.

\bibitem{Gopalakrishnan2015}
{\sc J.~Gopalakrishnan and J.~Sch{\"o}berl}, {\em {Degree and wavenumber [in]
  dependence of Schwarz preconditioner for the DPG method}}, in Spectral and
  High Order Methods for Partial Differential Equations ICOSAHOM 2014,
  Springer, 2015, pp.~257--265.

\bibitem{henneking2020numerical}
{\sc S.~Henneking and L.~Demkowicz}, {\em {A numerical study of the pollution
  error and DPG adaptivity for long waveguide simulations}}, CAMWA,  (2020).

\bibitem{henneking2019dpg}
{\sc S.~Henneking, J.~Grosek, and L.~Demkowicz}, {\em A dpg maxwell approach
  for studying nonlinear thermal effects in active gain fiber amplifiers},
  arXiv preprint arXiv:1912.01185,  (2019).

\bibitem{hiptmair1999multilevel}
{\sc R.~Hiptmair and R.~H. Hoppe}, {\em {Multilevel methods for mixed finite
  elements in three dimensions}}, Numer. Math., 82 (1999), pp.~253--279.

\bibitem{MR3543022}
{\sc B.~Keith, F.~Fuentes, and L.~Demkowicz}, {\em The {DPG} methodology
  applied to different variational formulations of linear elasticity}, Comput.
  Methods Appl. Mech. Engrg., 309 (2016), pp.~579--609.

\bibitem{sp2}
{\sc B.~Keith, S.~Petrides, F.~Fuentes, and L.~Demkowicz}, {\em {Discrete
  least-squares finite element methods}}, Comput. Methods Appl. Mech. Engrg.,
  327 (2017), pp.~226--255.

\bibitem{KimSeungil}
{\sc S.~Kim and H.~Zhang}, {\em Optimized {S}chwarz method with complete
  radiation transmission conditions for the {H}elmholtz equation in
  waveguides}, SIAM J. Numer. Anal., 53 (2015), pp.~1537--1558.

\bibitem{leng2020diagonal}
{\sc W.~Leng and L.~Ju}, {\em {A diagonal sweeping domain decomposition method
  with source transfer for the Helmholtz equation}}, arXiv preprint
  arXiv:2002.05327,  (2020).

\bibitem{LiXu}
{\sc X.~Li and X.~Xu}, {\em {Domain decomposition preconditioners for the
  discontinuous Petrov-Galerkin method}}, ESAIM Math. Model. Numer. Anal., 51
  (2017), pp.~1021--1044.

\bibitem{mumps1}
{\sc J.~W.~H. Liu}, {\em The multifrontal method for sparse matrix solution:
  Theory and practice}, SIAM Rev., 34 (1992), pp.~82--109.

\bibitem{nagaraj20193d}
{\sc S.~Nagaraj, J.~Grosek, S.~Petrides, L.~F. Demkowicz, and J.~Mora}, {\em {A
  3D DPG Maxwell approach to nonlinear Raman gain in fiber laser amplifiers}},
  Journal of Computational Physics: X, 2 (2019), p.~100002.

\bibitem{srir_soc}
{\sc S.~Nagaraj, S.~Petrides, and L.~F. Demkowicz}, {\em {Construction of DPG
  Fortin operators for second order problems}}, Comput. Math. Appl., 74 (2017),
  pp.~1964--1980.

\bibitem{oden2010}
{\sc J.~Oden and L.~Demkowicz}, {\em {Applied Functional Analysis, Third
  Edition}}, Textbooks in Mathematics, CRC Press, Taylor \& Francis Group,
  2018.

\bibitem{petrides2019adaptive}
{\sc S.~Petrides}, {\em {Adaptive multilevel solvers for the discontinuous
  Petrov--Galerkin method with an emphasis on high-frequency wave propagation
  problems}}, PhD thesis, The University of Texas at Austin, 2019.

\bibitem{sp}
{\sc S.~Petrides and L.~F. Demkowicz}, {\em An adaptive {DPG} method for high
  frequency time-harmonic wave propagation problems}, Comput. Math. Appl., 74
  (2017), pp.~1999--2017.

\bibitem{ROBERTS20172018}
{\sc N.~V. Roberts and J.~Chan}, {\em {A geometric multigrid preconditioning
  strategy for DPG system matrices}}, Comput. Math. Appl., 74 (2017), pp.~2018
  -- 2043.

\bibitem{LahayeD}
{\sc A.~H. Sheikh, D.~Lahaye, and C.~Vuik}, {\em On the convergence of shifted
  {L}aplace preconditioner combined with multilevel deflation}, Numer. Linear
  Algebra Appl., 20 (2013), pp.~645--662.

\bibitem{Stolk}
{\sc C.~C. Stolk}, {\em A rapidly converging domain decomposition method for
  the {H}elmholtz equation}, J. Comput. Phys., 241 (2013), pp.~240--252.

\bibitem{Stolk14}
{\sc C.~C. Stolk, M.~Ahmed, and S.~K. Bhowmik}, {\em A multigrid method for the
  {H}elmholtz equation with optimized coarse grid corrections}, SIAM J. Sci.
  Comput., 36 (2014), pp.~A2819--A2841.

\bibitem{taus2020sweeps}
{\sc M.~Taus, L.~Zepeda-N{\'u}{\~n}ez, R.~J. Hewett, and L.~Demanet}, {\em
  {L-Sweeps: A scalable, parallel preconditioner for the high-frequency
  Helmholtz equation}}, J. Comput. Phys., 420 (2020), p.~109706.

\bibitem{astaneh2018perfectly}
{\sc A.~Vaziri~Astaneh, B.~Keith, and L.~Demkowicz}, {\em On perfectly matched
  layers for discontinuous {Petrov-Galerkin} methods}, Comput. Mech.,  (2018).

\bibitem{VionA}
{\sc A.~Vion and C.~Geuzaine}, {\em Double sweep preconditioner for optimized
  {S}chwarz methods applied to the {H}elmholtz problem}, J. Comput. Phys., 266
  (2014), pp.~171--190.

\bibitem{xu1}
{\sc J.~Xu}, {\em Iterative methods by space decomposition and subspace
  correction}, SIAM Rev., 34 (1992), pp.~581--613.

\bibitem{Xu2}
\leavevmode\vrule height 2pt depth -1.6pt width 23pt, {\em An introduction to
  multilevel methods}, in Wavelets, multilevel methods and elliptic {PDE}s
  ({L}eicester, 1996), Numer. Math. Sci. Comput., Oxford Univ. Press, New York,
  1997, pp.~213--302.

\bibitem{Xu}
\leavevmode\vrule height 2pt depth -1.6pt width 23pt, {\em The method of
  subspace corrections}, J. Comput. Appl. Math., 128 (2001), pp.~335--362.
\newblock Numerical analysis 2000, Vol. VII, Partial differential equations.

\bibitem{Zepeda}
{\sc L.~Zepeda-N{\'u}{\~n}ez and L.~Demanet}, {\em The method of polarized
  traces for the 2{D} {H}elmholtz equation}, J. Comput. Phys., 308 (2016),
  pp.~347--388.

\bibitem{Zitelli}
{\sc J.~Zitelli, I.~Muga, L.~Demkowicz, J.~Gopalakrishnan, D.~Pardo, and V.~M.
  Calo}, {\em A class of discontinuous {P}etrov-{G}alerkin methods. {P}art
  {IV}: the optimal test norm and time-harmonic wave propagation in 1{D}}, J.
  Comput. Phys., 230 (2011), pp.~2406--2432.

\end{thebibliography}
